\documentclass[11pt]{article}
\usepackage{geometry, amsmath, amsthm, latexsym, amssymb, graphicx, epstopdf} 
\usepackage{amsthm,amssymb,mathrsfs,setspace,amsmath, bbm}
\usepackage{authblk}
\usepackage{subfiles}
\usepackage[usenames,dvipsnames]{color}
\usepackage{graphicx}
\usepackage{subcaption}
\usepackage[export]{adjustbox}
\usepackage{enumitem}
\numberwithin{equation}{subsection}
\usepackage{comment}
\usepackage{caption}
\usepackage{wrapfig}
\usepackage{sectsty}
\sectionfont{\centering} % Centers all \section headings

\usepackage[dvipsnames]{xcolor}
\definecolor{newred}{rgb}{0.32,0.04,0.04}   %darkred
\definecolor{newblue}{rgb}{0.0,0.06,0.26} %darkblue
\definecolor{newgreen}{rgb}{0.03,0.24,0.15} %darkgreen
\definecolor{newviolet}{rgb}{0.22,0.16,0.32} %darkviolet
\definecolor{newyellow}{rgb}{0.45,0.36,0.08} %darkyellow
\usepackage{xcolor}
\usepackage{hyperref}
    \hypersetup{
    colorlinks=true,
    linkcolor= newblue,
    filecolor= magenta,      
    urlcolor= cyan,
    citecolor= newviolet
    }
\newtheorem{theorem}{Theorem}[section]
\newtheorem{lemma}{Lemma}[section]
\newtheorem{proposition}{Proposition}[section]

\newtheorem{definition}{Definition}[section]
\newtheoremstyle{myremark}
  {3pt}   % Space above
  {3pt}   % Space below
  {\normalfont} % Body font
  {}      % Indent
  {\bfseries} % Head font (bold)
  {}     % Punctuation after head
  {\newline}  % Space after head
  {}
\theoremstyle{myremark}

\newtheoremstyle{myexample}
  {3pt}   % Space above
  {3pt}   % Space below
  {\normalfont} % Body font
  {}      % Indent
  {\bfseries} % Head font (bold)
  {. }     % Punctuation after head
  { }  % Space after head
  {}
\theoremstyle{myexample}

\DeclareMathOperator{\comp}{comp}
\providecommand{\keywords}[1]
{
  \small	
  \textbf{Keywords:} #1
}
\title{\bfseries Synchronisation under quasi-periodic forcing\\ of circle diffeomorphisms with two attracting and two repelling fixed points}
\author{{Kirthana Rajasekar \footnote{acknowledges the support of the Swedish Research Council (grant 2020-03989)\\
\textit{Date:} 23 November 2026}

}}
\date{}

\begin{document}

\maketitle

%%%%%%%%%% ABSTRACT %%%%%%%%%%%%%%%%
\begin{abstract} 
We study the dynamics of circle diffeomorphisms, which have exactly two attracting and two repelling fixed points, under quasiperiodic forcing. For sufficiently small forcing, such a system  exhibit uniformly hyperbolic behaviour. The results focus on systems with larger forcing, where every orbit intersects the basins of attraction of both attracting fixed points infinitely often. For such quasi-periodically forced circle diffeomorphisms, we identify a set of parameters of positive measure for which the corresponding skew-product systems admit a unique attracting and a unique repelling invariant graph, which capture the dynamics of Lebesgue almost every point on $\mathbb{T}^2$. Further, we show that both continuous invariant graphs and invariant graphs that are discontinuous on a dense set can occur within the same class of quasi-periodically forced systems, for different quasi-periodic rotations.
\end{abstract}
%%%%%%%%%%%%%%%%%%%%%%%%%%%%%%%%%%%%%%

%%%%%%%%%% INTRODUCTION %%%%%%%%%%%%%%%%%%%%%%%
\section{Introduction}

In the 1960s--80s, numerous experiments and theoretical studies established the occurrence of attractors with wild fractal geometry in certain dynamical systems, such as the Lorenz and H\'{e}non attractors. These developments motivated the notion of ``strange attractors'', which admits different formulations depending on the context in which it is studied. We consider an attractor strange when it is neither finite nor piecewise differentiable \cite{AlsedaCosta2009,grebogi:1984}. An intriguing question in this context was whether there exists a strange attractor that supports non-chaotic dynamics.

Meanwhile, the study of quasiperiodically forced (qpf) systems was prevalent and motivated by their close connection to non-linear oscillators with incommensurate frequencies \cite{bondeson, Lin2010}, along with numerous other physical systems. The qpf Schr\"{o}dinger operators is one such class of systems, which has grown into a distinctive area of active research, and we refer to the survey \cite{Marx_Jitomirskaya_2017} for an overview of some of the influential results. When a general qpf system has invariant graphs, they are among the simplest invariant structures which can provide crucial insights on the dynamics of the system. We refer to \cite{Stark_1997,Stark_1999} for results on the regularity and stability of invariant graphs of certain forced systems and a survey of some of the previous results in this area.  

In 1984, Grebogi et al. \cite{grebogi:1984} showed in their landmark result the existence of strange non-chaotic attractors (SNA) in quasiperiodically forced systems given by $F: \mathbb{T}^2 \rightarrow \mathbb{T}^2$, where $\mathbb{T} := \mathbb{R}/\mathbb{Z}$ and 
\begin{align*}
    F(\theta, x) = (\theta+\omega,\ 2 \lambda \tanh(x) \cos(\theta)) ,
\end{align*}
where $\omega$ is the golden mean and $\lambda > 0$ is a parameter. The breakthroughs that followed further solidified SNA as a phenomenon in certain non-uniformly hyperbolic systems, which are generally rare but are surprisingly abundant in quasi-periodically forced systems. Following this discovery, Keller \cite{Keller1996} proved the existence of and described SNAs for a broader class of maps of the form 
\begin{align}\label{equation_introduction_Kellermaps}
    F(\theta, x) = (\theta+\omega, f(x)g(\theta)),
\end{align}
where $\omega$ is irrational and specific conditions are imposed on $f$ and $g$, which include the family of maps explored in \cite{grebogi:1984}. In the following years, mechanisms of SNA formation have been studied and further generalised for a broader class of qpf maps called the pinched skew-products, for example in \cite{Haro2006, Glendinning2002, Jager2007}, which encompasses Keller's maps in equation \eqref{equation_introduction_Kellermaps}. An underlying mechanism, called non-smooth saddle-node bifurcation, which creates SNAs also in the case of pinched skew product maps, was described rigorously in \cite{JAGER_2006}. 

The above chain of results, from a specific example to a generalised mechanism, exemplifies a type of development that is common in dynamical systems. We have had numerous results for other special classes such as linear $SL_2(\mathbb{R})$ cocycles \cite{Young:2002}, the logistic map \cite{Bjerklov2009}, and Harper maps \cite{HaroAlex2006, BJERKLOV_2005}.  In many cases, the techniques have obstructions that limit them to specific classes of qpf systems. Thus, there continues to be a demand for more generalised results that capture the underlying mechanisms that generate SNA in qpf systems. The results in this article  contribute to overcoming this gap. We refer to \cite{Jorba2007} for a survey of some of the previous results and their generalisations. 

The mechanisms such as non-smooth saddle-node bifurcations \cite{JAGER_2006,Fuhrmann2018} and fractalisation of invariant curves \cite{AngelJorba2008} were established with a similar incentive. In \cite{JAGER_2006}, for a generalised class of qpf circle maps, a non-smooth collision of the stable and unstable manifolds was proved to create SNA. A result of similar flavour is the fractalisation of invariant curves that was established in \cite{AngelJorba2008}, for a one-parameter family of linear skew-products. Their result demonstrates the progressive destruction of a smooth invariant curve into a fractal structure as its Lyapunov exponent approaches zero. This also serves as a warning for highly wrinkled invariant curves that could be mistaken for a fractal in numerical experiments. The results in \cite{Krikorian2018, Jager_Stark_2006} are part of the effort to classify qpf systems, analogous to Poincaré's classification of circle homeomorphisms, and they are restricted to a certain class of qpf circle homeomorphisms that are homotopic to the identity.

Our techniques are inspired by \cite{BJERKLOV_2005}, where the projective qpf Schr\"{o}dinger cocycles, with specific conditions on the potential, have positive Lyapunov exponents and are minimal for a positive measure of frequencies. Our set up and method is also motivated by \cite{Jager_2009}, where the existence of SNA and minimality of the dynamics was proved for a class of qpf circle maps. These techniques were in turn influenced by \cite{lyapunov_young_1997}, where the existence of positive Lyapunov exponents was proved for an open class of qpf linear $SL_2(\mathbb{R})$ cocycles. This class of techniques, of using critical sets with non-fast return conditions, builds on the techniques in the seminal paper by Benedicks and Carleson on H\'{e}non maps \cite{Benedicks_Carleson_1991}. 

For the systems considered in \cite{BJERKLOV_2005, Jager_2009, lyapunov_young_1997}, the corresponding forced circle maps have exactly one attracting and one repelling regions along the fibre direction. The primary distinguishing characteristic of the result in Theorem \ref{Theorem SNA} is the consideration of Morse–Smale diffeomorphisms with two attracting and two repelling fixed points, which significantly complicates the question of convergence of the orbits along the fibre for the corresponding forced system. It is also worth noting that, unlike many of the previous results, the results in \cite{BJERKLOV_2005}, \cite{Jager_2009} and Theorem \ref{Theorem SNA} do not require pinching or monotone conditions on the circle diffeomorphisms on the fibre.  

%----------------------------- FIGURE circle diffeo and g(theta) ------------------
% trim = left bottom right top
\begin{figure}[htbp] \centering \captionsetup[subfigure]{font=small, skip=2pt} 
\begin{subfigure}[t]{0.5\textwidth} \centering \includegraphics[width=\linewidth, trim=7.5cm 0.5cm 5.5cm 1.5cm,
clip]
    {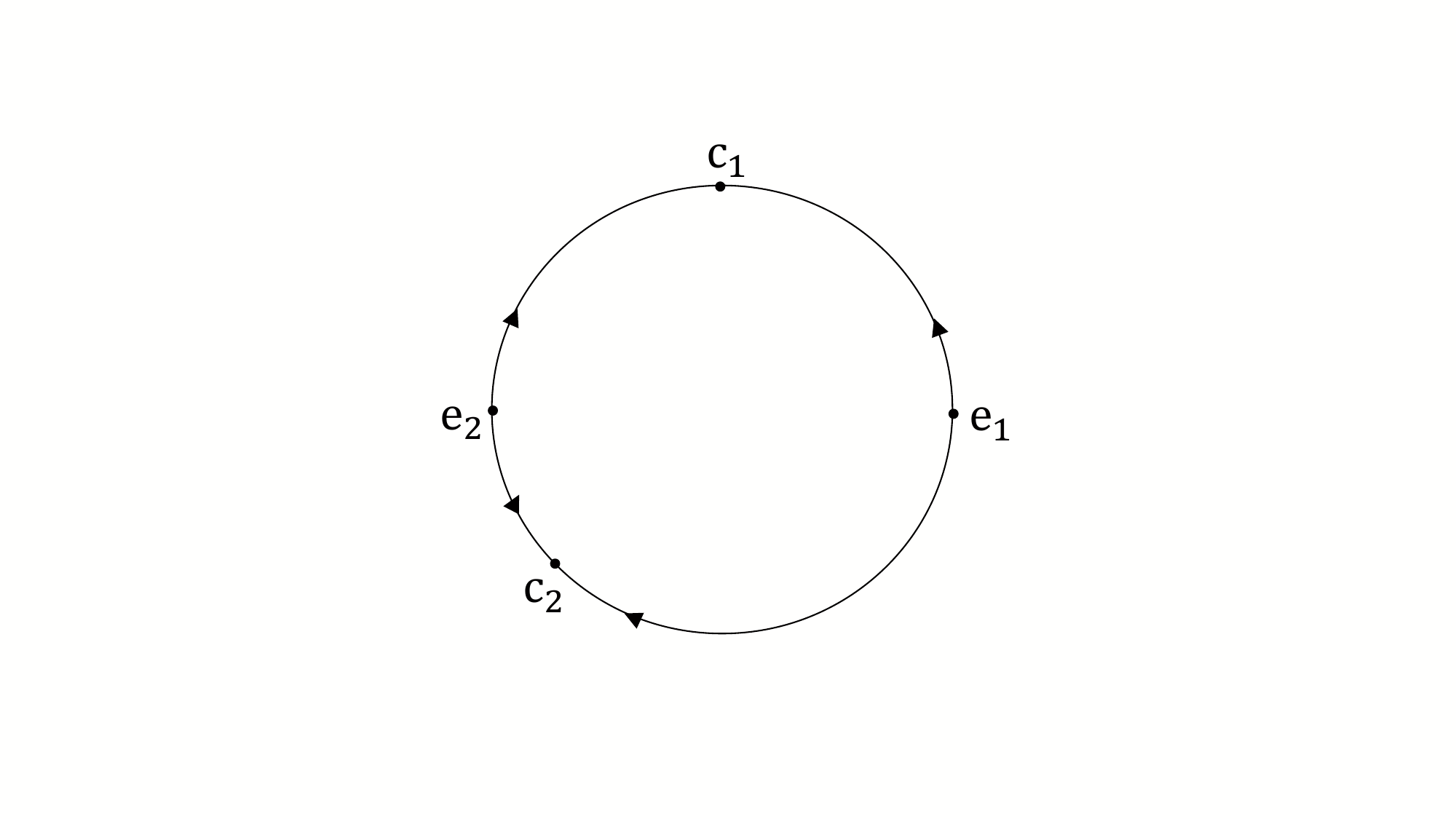} 
\label{fig:circle_diffeomorphism} 
\end{subfigure} \hspace{-0.08\textwidth} 
\begin{subfigure}[t]{0.55\textwidth} \centering \includegraphics[width=\linewidth, trim=1cm 0.45cm 1cm 0.3cm,
clip]{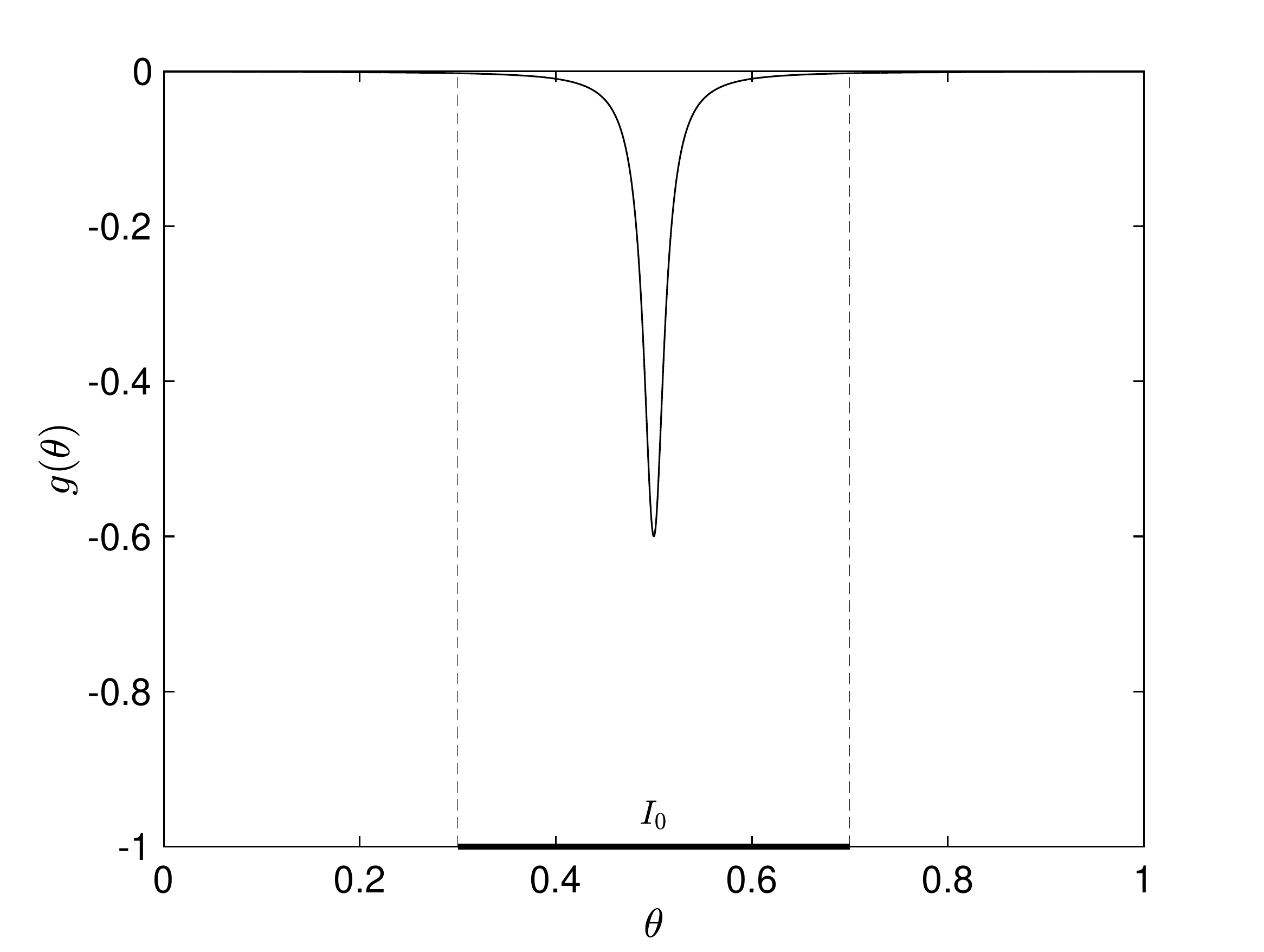} 
\label{fig:critical_set_I0} \end{subfigure} 
\caption{The diagram on the left illustrates a circle diffeomorphism $f$ with two attracting fixed points $\{c_1, c_2\}$ and two repelling fixed points $\{e_1, e_2\}$, for which $d(c_2,e_2)<d(c_1,e_1)$, where $d$ denotes the metric on the circle. The figure on the right contains the plot of the curve $\displaystyle g(\theta) = \frac{-a}{1+\lambda \cos^2(\pi\theta)}$, for the values $a=0.6$ and $\lambda = 625$. Its corresponding critical set $I_0 = [1/2-\lambda^{-1/4}, 1/2+\lambda^{-1/4}]$ is marked on the base. The plot illustrates that the points $\theta \in \mathbb{T}$ for which the forcing $g(\theta)$ on the fibre is the strongest are contained in the critical set $I_0$.} \label{Statementofthm_figure1} 
\end{figure}
%---------------------------------

\section{The Setting}

 Let $f: \mathbb{T} \rightarrow \mathbb{T}$ denote a $C^1$ circle diffeomorphism with exactly two attracting fixed points $\{c_1, c_2\}$ and two repelling fixed points $\{e_1, e_2\}$. Since $\mathbb{T} = [0,1)/\sim$, without loss of generality, let us assume that the fixed points satisfy the relation
 $$0<e_2<c_2<e_1<c_1<1,$$
 when projected onto $[0,1)$. We assume that $\displaystyle {400} (c_2-e_2)< (c_1-e_1)$. Such a circle diffeomorphism belongs to the class of Morse–Smale diffeomorphisms. It is worth mentioning that the constant $400$ in the bound $\displaystyle {400} (c_2-e_2)< (c_1-e_1)$ is chosen for convenience and is not an optimal condition. 

The forcing on the fibre is defined by a function $g:\mathbb{T} \rightarrow \mathbb{R}$ given by 
\begin{align}\label{equation_setting_g func}
    g_{\lambda,a}(\theta) = g(\theta) := \frac{-a}{1+\lambda \cos^2(\pi\theta)},
\end{align}
 where $\lambda >1$ (which is assumed to be large) and $a>0$. The figure \ref{Statementofthm_figure1} contains the graph $x = g(\theta)$ for a fixed $\lambda >1$, and the critical set $I_0$ defined as 
\begin{align*}
    I_0 := [1/2-\lambda^{-1/4}, 1/2+\lambda^{-1/4}],
\end{align*}
which is the region where the forcing $g(\theta)$ is the strongest. Moreover, for any $\theta \notin I_0$ the forcing  $g(\theta)$ is relatively weak. The other essential characteristic of the forcing is that $g(\theta)$ has a unique critical point in $I_0$, where its global minimum is attained. Here, the parameters $a$ and $\lambda$ regulate the strength and localisation of the forcing, respectively. It should be noted that the specific function $g(\theta)$ in equation \eqref{equation_setting_g func} is chosen for convenience. Furthermore, our result also holds true for a broader class of functions which satisfy the property of strong localised forcing in an interval $I_0$ and possess a unique critical point in $I_0$.

For $\omega \in \mathbb{T} \setminus \mathbb{Q}$, we study  quasi-periodically forced circle diffeomorphisms defined by skew-product maps $F_{\omega}: \mathbb{T}^2 \rightarrow \mathbb{T}^2$ of the form 
\begin{align}\label{statementofthm_defn_mapF}
    F_{\omega}(\theta, x) & := \left( \theta+\omega, f(x)+g(\theta) \right) \pmod{1} \\
    & = \left( \theta+\omega, f(x)-\frac{a}{1+\lambda \cos^2(\pi \theta)} \right) \pmod{1},
\end{align}
where $\lambda >1$ (which is assumed to be large), $a>0$ and $\displaystyle {400} (c_2-e_2)< (c_1-e_1)$. Since the map is on $\mathbb{T}^2=(\mathbb{R}/\mathbb{Z})^2$, each coordinate is considered modulo $1$. Note that an orbit $(\theta_n, x_n)_{n\geq 1} := (F^{n}(\theta, x))_{n\geq 1}$ experiences weak forcing when the iterate on the base $\theta_n$ is outside $I_0$ and a stronger forcing when the iterate $\theta_n$ lies in $I_0$. Thus, for larger values of $a>0$ and $\lambda>0$ the orbits of all the points experience a stronger forcing but at fewer iterations.

For a quasi-periodically forced map $F_{\omega}$, the asymptotic behaviour of the orbit of a point $(\theta,x) \in \mathbb{T}^2$ is characterised by the pointwise forward Lyapunov exponents along the fibre, given by
\begin{align}\label{equation_intro_forward Lyapunov}
\overline{L}^{+}(\theta, x) & := \limsup_{n \rightarrow \infty} \frac{1}{n} \log |\partial_x \left(\pi_2\circ F^{ n}(\theta, x)\right)|, \\
    \underline{L}^{+}(\theta, x) & := \liminf_{n \rightarrow \infty} \frac{1}{n} \log |\partial_x \left( \pi_2\circ F^{ n}(\theta, x)\right)|,  \label{equation_intro_lower_forward Lyapunov}
\end{align}
where $\pi_2$ denotes the projection on the second co-ordinate, and the pointwise backward Lyapunov exponents along the fibre, given by
\begin{align}\label{equation_intro_backward Lyapunov}
    \overline{L}^{-}(\theta, x) & := \limsup_{n \rightarrow \infty} \frac{1}{n} \log |\partial_x  \left(\pi_2\circ F^{-n}(\theta, x)\right)|, \\
    \underline{L}^{-}(\theta, x) & := \liminf_{n \rightarrow \infty} \frac{1}{n} \log |\partial_x  \left(\pi_2\circ F^{-n}(\theta, x)\right)|. \label{equation_intro_lower_backward Lyapunov}
\end{align}
A map $\Phi: \mathbb{T} \rightarrow \mathbb{T}$ is an invariant graph of $F_{\omega}$ if it satisfies the condition 
\begin{align*}
    F_{\omega}(\theta, \Phi(\theta)) = (\theta+\omega, \Phi(\theta+\omega)),\ \  \text{   for almost every $\theta \in \mathbb{T}$.}
\end{align*} 
When they exist, Lyapunov exponents capture the local dynamics of the quasi-periodically forced system, while an invariant graph can provide crucial insights on the global dynamics. Further, we call an invariant graph $\Phi$ attracting if $\overline{L}^+(\theta, \Phi(\theta)) < 0$, for Lebesgue almost every $\theta \in \mathbb{T}$, and repelling if $\overline{L}^-(\theta, \Phi(\theta)) < 0$, for Lebesgue almost every $\theta \in \mathbb{T}$.  

\subsection{Discussion of the Setting}

Before we state the main result in section \ref{section_statementoftheThm}, we discuss the various dynamical behaviours exhibited by quasi-periodically forced maps defined in equation \eqref{statementofthm_defn_mapF}, for different strengths of forcing. Consequently, we highlight the difficulties that lie in analysing synchronisation in these quasi-periodically forced circle diffeomorphisms. 

%----------------------------- FIGURE SNA Invariant graphs
% trim = left bottom right top
\begin{figure}[t]
    \centering
    \includegraphics[
        width=0.99\linewidth,
        trim=1cm 0.1cm 0.8cm 0.3cm,
        clip
    ]{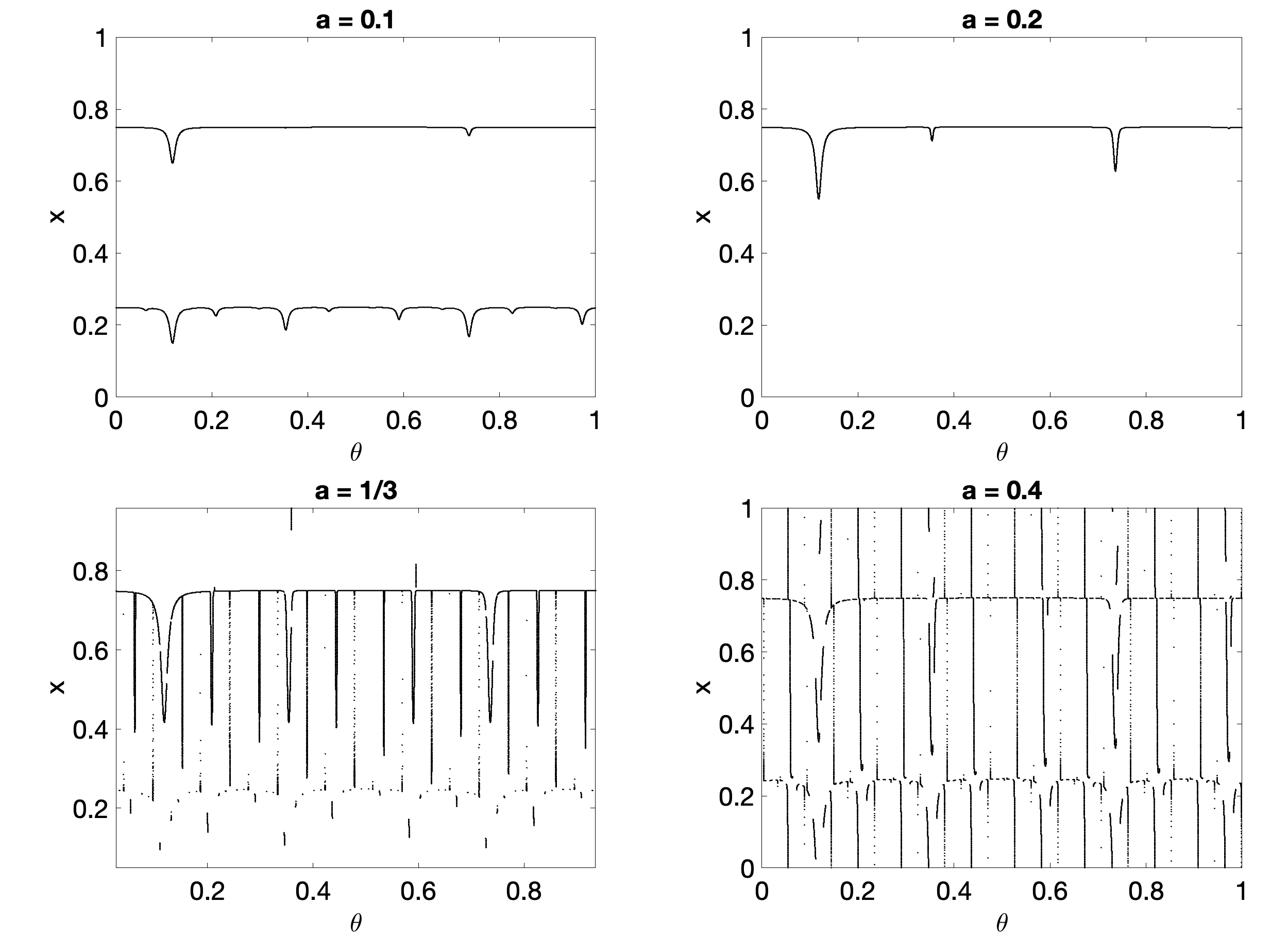}
    \caption{The figure demonstrates the dynamics of the skew-product map $\tilde{F}_{\omega}$ in Example 1 (equation \eqref{equation_examples_qpf map}), for $\omega = (\sqrt{5}-1)/2$ and different values of $a$. In each panel, the value of $a$ is fixed and plot contains $10^6$ iterates of $20$ equidistributed points on a fixed fibre with the first $10^3$ iterates removed. The plots thus illustrate the bifurcation that occurs for increasing values of $a$, where the attracting invariant graphs go from continuous curves to a graph that resembles a SNA.}
    \label{examples_SNA}
\end{figure}
%----------------------------------

Note that an unforced map $(\theta,x) \rightarrow (\theta+\omega, f(x))$ consists of two continuous attracting graphs $\{\mathbb{T} \times \{c_i\}\}_{i=1}^{2}$ and two continuous repelling graphs $\{\mathbb{T} \times \{e_i\}\}_{i=1}^{2}$. Since $(c_1-e_1) \neq (c_2-e_2)$, the circle diffeomorphism $f$ cannot satisfy the symmetric condition $f(x+1/2) = f(x)+1/2$. This ensures that the quasi-periodically forced system is not trivially decomposable into two identical dynamical subsystems, with possibly two attracting invariant graphs. We now demonstrate the different dynamics that can be observed for varying strengths of forcing and for different frequencies $\omega$, in the following subsubsection.

\subsubsection{Examples}\label{subsubsection_examples}

Let us consider quasiperiodically forced circle diffeomorphisms of the form 
\begin{align}\label{equation_examples_qpf map}
    \tilde{F}_{\omega}(\theta, x) = \left(\theta+\omega, \tilde{f}(x) - \frac{a}{1+2000 \cos^2(\pi \theta)}\right),
\end{align}
where $\tilde{f}(x) = x + 0.05\sin(4\pi x) - 0.05\cos(2\pi x)$ and $a>0$. Note that the circle diffeomorphism $\tilde{f}$ has two attracting fixed points $\{c_1,c_2\}$ and two repelling fixed points $\{e_1,e_2\}$ for which $0<e_2<c_2<e_1<c_1<1$ and $0.1< (c_2-e_2) < (c_1-e_1) < 0.4$, where $(c_1-e_1) = 1/3$. Even though the map $\tilde{F}_{\omega}$ is different from the skew-product maps considered in equation \eqref{statementofthm_defn_mapF}, we use these maps to motivate the occurrence of non-chaotic behaviour along the fibre and the existence of attracting invariant graphs of varying regularities for different frequencies $\omega \in \mathbb{T}$. 

\vspace{0.2cm}
\noindent \textbf{Example 1:} Let us fix $\omega = (\sqrt{5}-1)/2$ and consider the corresponding map $\tilde{F}_{\omega}$. For a sufficiently small forcing $a>0$, for which the perturbation can be overcome by the hyperbolicity of $\tilde{f}$, the system consists of invariant curves that are small perturbations of the curves $\{\mathbb{T} \times \{c_i\}\}_{i=1}^{2}$ and $\{\mathbb{T} \times \{e_i\}\}_{i=1}^{2}$. This is illustrated in figure \ref{examples_SNA}, in the plot corresponding to the value $a=0.1$. Here, the plot indicates the possible existence of two continuous attracting invariant graphs. 

Consider an intermediate value of $a$ for which $(c_2-e_2)< a < (c_1-e_1)$. The dynamics of one such system is illustrated in figure \ref{examples_SNA}, in the plot corresponding to the value $a=0.2$. Here, the plot indicates the possible existence of a unique continuous attracting invariant graph. This is likely due to the asymmetry $(c_2-e_2) < (c_1-e_1)$ which leads to a single invariant region to which orbits are eventually forced into by the forcing.

 In the case where $(c_2-e_2)< a = (c_1-e_1)$, unlike the previous case, we now have a tangential intersection of $\{(\theta, c_1 + g_{2000,a}(\theta)): \theta \in \mathbb{T} \}$ and $\mathbb{T} \times \{e_1\}$. The intricacies of the dynamics due to the tangential intersection can be compared to the saddle-node bifurcations studied in \cite{JAGER_2006}. The non-triviality of this case is illustrated in figure \ref{examples_SNA}, in the plot corresponding to $a=1/3$.

For a larger value of forcing where $a>(c_1-e_1)$, orbits can experience strong forcing that kicks the orbit on the fibre from one attracting region to the other. This complicates the dynamics significantly and makes it non-trivial to investigate the existence and regularity of  invariant graphs. In figure \ref{examples_SNA}, the plot corresponding to $a=0.4$ suggests synchronisation of orbits along the fibre and thus the possibility of a unique attracting invariant graph. More specifically, the orbits trace a highly irregular geometric structure. This motivates the coexistence of a geometrically complicated invariant structure and non-chaotic behaviour along the fibre, for quasi-periodically forced circle diffeomorphisms defined in equation \eqref{statementofthm_defn_mapF}.

%-------- FIGURE - Invariant curves ---------------
\begin{figure}[htbp]
    \centering
    \includegraphics[
        width=0.9\linewidth,
        trim=0cm 5cm 0cm 5cm,
        clip
    ]{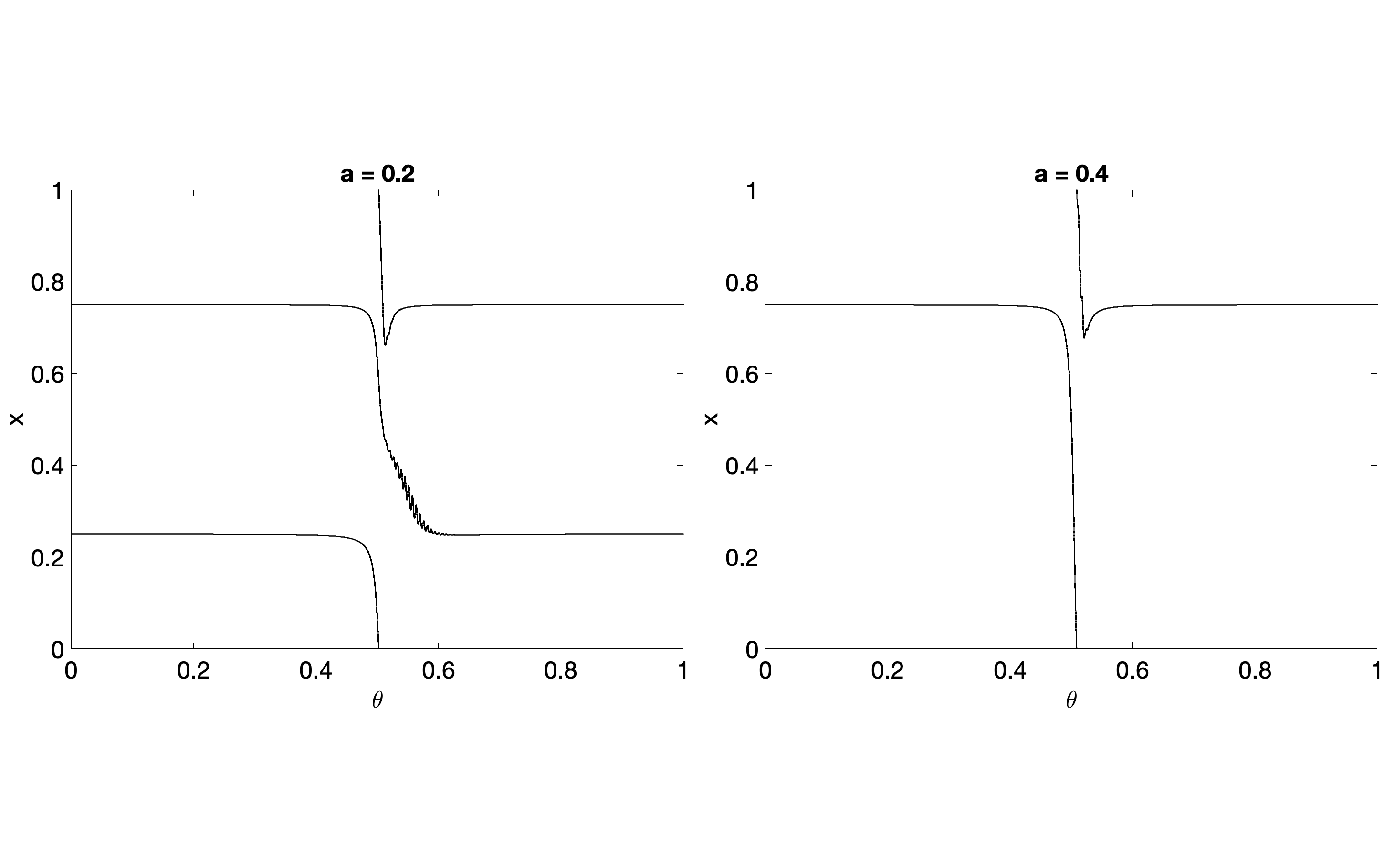}
    \caption{The figure demonstrates the dynamics of the skew-product map $\tilde{F}_{\omega}$ in equation \eqref{equation_examples_qpf map}, for $\omega = \pi/2^9$ and different values of $a$. In each panel, the value of $a$ is fixed and plot contains $10^6$ iterates of $20$ equidistributed points on a fixed fibre with the first $10^3$ iterates removed. The plot illustrates the transition from two attracting invariant curves to a single attracting invariant curve, for increasing values of $a$.}
    \label{examples_curves}
\end{figure}
%----------------------------------------

\vspace{0.2cm}
\noindent \textbf{Example 2:} Let us consider the map $\tilde{F}_{\omega}$ corresponding to the parameters $\displaystyle \omega = {\pi}/{2^9}$ and $a= 0.4$. Similar to Example 1, simulations suggest a synchronisation of orbits along the fibre. In figure \ref{examples_curves}, the plot corresponding to $a=0.4$ is a single synchronised orbit, which suggests a single continuous invariant curve. This supports that the regularity of the invariant graphs depends sensitively on the parameter $\omega$.
%Note that unlike the golden mean $(\sqrt{5}-1)/2$, the irrational $\pi/2^9$ is close to the rational $0$.

\section{Statement of the theorem}\label{section_statementoftheThm}

In this section, we state the main result of this paper, where we establish the non-chaotic behaviour along the fibre in the quasi-periodically forced circle diffeomorphisms defined in equation \eqref{statementofthm_defn_mapF}, under certain open conditions on $\lambda$ and $a$. For any fixed $(\theta,x) \in \mathbb{T}^2$ and $n\in \mathbb{Z}$, we use the notation $(\theta_n, x_n):= F_{\omega}^n(\theta, x)$, to denote the points in the orbit of $(\theta,x)$.
%More specifically, under open conditions on $\lambda$ and $a$, we prove the existence of highly discontinuous invariant graphs and continuous invariant curves, for different parameters $\omega \in \mathbb{T}$.

\begin{theorem}\label{Theorem SNA}
    Let $a>0$ be such that $(c_1-e_1) < a < (c_1-e_1) + (c_2-e_2) + \min \left\{ (e_1-c_2), (1-c_1+e_2) \right\}$. There exists $\lambda_0 > 0$ such that for all $\lambda > \lambda_0$ there exist two disjoint sets $\Omega_1(\lambda), \Omega_2(\lambda) \subset \mathbb{T}\setminus \mathbb{Q}$ of positive measure such that the corresponding maps $F_{\omega}$, as defined in equation \eqref{statementofthm_defn_mapF}, satisfy the following: 
\begin{enumerate}[label=(\alph*)]
    \item\label{Theorem SNA_a} For each $\omega \in \Omega_1$, there exist two single-valued invariant graphs $\Phi, \Psi: \mathbb{T} \rightarrow \mathbb{T}$ such that  
    \begin{align}\label{statementofthem_thm_synceqn1a}
        \lim_{n \rightarrow \infty} d(x_n,  \Phi(\theta_n)) = 0,
    \end{align}
    for Lebesgue a.e. $\theta \in \mathbb{T}$ and $x \neq \Psi(\theta)$, and 
        \begin{align}\label{statementofthem_thm_synceqn1b}
        \lim_{n \rightarrow \infty} d(x_{-n}, \Psi(\theta_{-n})) = 0,
    \end{align}
    for Lebesgue a.e. $\theta \in \mathbb{T}$ and $x \neq \Phi(\theta)$. 
    
    Furthermore, the map has a sink-source orbit, that is, there exists a point $(\theta,x)\in \mathbb{T}^2$ for which we have $\underline{L}^+(\theta,x) > 0$ and $\underline{L}^-(\theta, x) > 0$.
    
    \item\label{Theorem SNA b} For each $\omega \in \Omega_2$, there exist two invariant $C^1$ curves $\Phi, \Psi: \mathbb{T} \rightarrow \mathbb{T}$ for which 
    \begin{align}\label{statementofthem_thm_synceqn2a}
        \lim_{n \rightarrow \infty} d(x_n, \Phi(\theta_n)) = 0,
    \end{align}
    for all $\theta \in \mathbb{T}$ and $x \neq \Psi(\theta)$, and 
    \begin{align}\label{statementofthem_thm_synceqn2b}
        \lim_{n \rightarrow \infty} d(x_{-n}, \Psi(\theta_{-n})) = 0,
    \end{align}
    for all $\theta \in \mathbb{T}$ and $x \neq \Phi(\theta)$.    
\end{enumerate}
    
\end{theorem}

The existence of a sink-source orbit in Theorem \ref{Theorem SNA}\ref{Theorem SNA_a} implies that the graphs of $\Phi$ and $\Psi$ are discontinuous on a dense set of points in $\mathbb{T}$. Thus, Theorem \ref{Theorem SNA} establishes the co-existence of both discontinuous and continuous invariant graphs within the class of qpf circle maps defined in equation \eqref{statementofthm_defn_mapF}, for different parameters $\omega$. 

\subsection{Remarks}\label{statementofthem_remarks} 
\begin{enumerate}[label=(\alph*)]
\item For each $\omega \in \Omega_1\cup \Omega_2$, the result shows that the invariant graphs $\Phi$ and $\Psi$ each support an ergodic measure which is the push-forward of the one-dimensional Lebesgue measure. In particular, these ergodic measures describe the statistical behaviour of Lebesgue almost every point in $\mathbb{T}^2$, under the forward or backward iteration with respect to $F_{\omega}$.
\item In Theorem \ref{Theorem SNA}\ref{Theorem SNA_a}, for each $\omega \in \Omega_1$, the dynamics of the map $F_{\omega}$ is expected to be minimal on $\mathbb{T}^2$. In this case, the graphs $\Phi$ and $\Psi$ will in fact be dense on the torus $\mathbb{T}^2$.
\item The measure of the  parameter set $\Omega_1$ converges to a positive quantity as $\lambda$ tends to infinity, and the set $(\Omega_2 \cup \mathbb{Q})$ is an open set. However, it is worth noting that the sets $\Omega_1$ and $\Omega_2$ derived here are not optimal. Their construction  are specifically for the purpose of Theorem \ref{Theorem SNA}.
\end{enumerate}
\section{Overview of the Paper}

Let us briefly discuss the roadmap to the final proof of Theorem \ref{Theorem SNA} which is in section \ref{section_proof of thm}. Section \ref{section_preliminary notations} contains the preliminary notations, constructions and estimates that will be used in the rest of the paper. The proof of theorem \ref{Theorem SNA} requires three major components, namely the mechanism of slow-returns in section \ref{section_Mechanism of slow-returns}, the description of critical sets in section \ref{section_structure of critical sets} and the construction of parameter sets $\Omega_1$ and $\Omega_2$ in section \ref{section_parameter space}. 

 In each of sections \ref{section_Mechanism of slow-returns} and \ref{section_structure of critical sets}, the propositions form an inductive framework. In section \ref{section_Mechanism of slow-returns}, we introduce a general notion of critical sets and then formulate the corresponding concept of slow-returns, which induces a sequence of arithmetic conditions on parameter $\omega$. In section \ref{section_structure of critical sets}, we fix a specific sequence of critical sets for each admissible $\omega$, and determine its geometric structure using the estimates in section \ref{section_Mechanism of slow-returns}.

The construction in section \ref{section_parameter space} is a sequence of parameter exclusions on $\mathbb{T}$, where we implement conditions of slow-returns imposed in sections \ref{section_Mechanism of slow-returns} and \ref{section_structure of critical sets}. Finally, the three components are combined in section \ref{section_proof of thm} for the construction of invariant graphs $\Phi$ and $\Psi$, and a sink-source orbit. 

\section{Preliminaries}\label{section_preliminary notations}

\subsection{Notations}

This subsection contains the essential notations used in the rest of the article. For any set $B$, the cardinality of the set is denoted by $\# B$. For the set $\{1,2\}$, if $i\in \{1,2\}$ then let $(i+1)$ denote the constant $(i+1) \pmod{2}$ whenever convenient. The one-dimensional torus is given by $\mathbb{T} := \mathbb{R}/\mathbb{Z}$. Thus, there is a natural correspondence between the one-dimensional torus and the interval $[0,1)$. For any two points $x,y \in [0,1)$ such that $0<x<y<1$, note that we have two distinct intervals in $[0,1)$ given by $(x,y) \subset [0,1)$ and 
\begin{align*}
    (y,x) := [0,x) \cup (y,1).
\end{align*}
Then, we use the notation $(x,y)$ and $(y,x)$ to also denote the corresponding intervals in the torus $\mathbb{T}$. Let $d$ denote the natural metric endowed on it. When convenient, we also use the notation
\begin{align*}
    |x| := d(0,x), \text{ for each } x\in \mathbb{T}.
\end{align*}

Let $\mathbb{T}^2$ denote the two-dimensional torus for which $\pi_1:\mathbb{T}^2 \rightarrow \mathbb{T}$ represents the projection on the first coordinate given by $\pi_1(\theta,x) = \theta$, for all $(\theta,x) \in \mathbb{T}^2$. Similarly, the projection on the second coordinate $\pi_2:\mathbb{T}^2 \rightarrow \mathbb{T}$ is given by $\pi_2(\theta, x) = x$, for all $(\theta,x) \in \mathbb{T}^2$. For any Lebesgue measurable subset $A \subset \mathbb{T}^2$, its two-dimensional Lebesgue measure is denoted by $|A|$. Similarly, when convenient, for any Lebesgue measurable set $I \subset \mathbb{T}$, let $|I|$ denote its one-dimensional Lebesgue measure. For each interval $I = (a,b)$ and constant $\varepsilon > 0$, the $\varepsilon$-neighbourhood of the interval $I$ is given by
$$B_{\varepsilon}(I) := (a-\varepsilon, b+\varepsilon).$$
For each $n\geq 1$, we also use the notation $nI$ for the set $B_{n|I|/2}(I)$, which is a neighbourhood of $I$, and which satisfies $|nI|=n|I|$. For an interval $I(\omega) \subset \mathbb{T}$ which is of the form
\begin{align*}
    I(\omega) = [a(\omega), b(\omega)], 
\end{align*}
where $a$ and $b$ are continuously differentiable functions on subsets of $\mathbb{T}$, we define the derivative of $I(\omega)$ with respect to $\omega$ by the following:
\begin{align*}
    \left|\partial_{\omega} I(\omega) \right|:= \max\{|a'(\omega)|, |b'(\omega)|\}.
\end{align*}
Consider a set $J(\omega) = \bigcup_{i=1}^{\mathcal{N}}I^i(\omega)$, where $\{I^{i}(\omega):1\leq i \leq \mathcal{N}\}$ is a finite collection of intervals for which the endpoints are continuously differentiable in $\omega$. Then the definitions discussed so far can be naturally extended to $J$. That is, for any $\varepsilon >0$, $n\geq 1$, and an admissible $\omega \in \mathbb{T}$, 
\begin{align*}
    B_{\varepsilon}(J) & := \bigcup_{i=1}^{\mathcal{N}} B_{\varepsilon}(I^i), \\
   nJ & := \bigcup_{i=1}^{\mathcal{N}} nI^i,\\
    \left|\partial_{\omega} J(\omega) \right| & := \max\{\left|\partial_{\omega} I^i(\omega) \right|: i\in \{1,\dots , \mathcal{N}\}\}.
\end{align*}
For any interval $I \subset \mathbb{T}$ and continuous function $\phi:I \rightarrow \mathbb{T}$, we use the following notations:
\begin{align*}
    |\phi|_{\sup} & := \sup \left\{ |\phi(\theta)|: \theta \in I\right\} \\
    |\phi|_{\inf} & := \inf \left\{ |\phi(\theta)|: \theta \in I\right\}.
\end{align*}

Consider a fixed $\omega \in \mathbb{T}$ and a quasi-periodically forced circle diffeomorphism $F_{\omega}$ as defined in equation \eqref{statementofthm_defn_mapF}. Then, we use the notation $F=F_{\omega}$ for convenience. For each $\theta \in \mathbb{T}$, let $f_{\theta}(x)$ denote the circle diffeomorphism $\pi_2 \circ F(\theta, x)$. For every $(\theta,x) \in \mathbb{T}^2$ and $n\geq 0$, the projection of the $n$th iterate $F^n(\theta, x)$ onto the second coordinate is denoted by $f^{n}_{\theta}(x)$, that is,
\begin{align*}
    F^n(\theta, x) = (\theta+n\omega, f^{n}_{\theta}(x)).
\end{align*}
Similarly, for every $\theta \in \mathbb{T}$ and $x\in \mathbb{T}$, we use the notation $f_{\theta}^{-1}(x) := f^{-1}(x-g(\theta))$, and thus 
\begin{align*}
    F^{-1}(\theta,x) & = (\theta-\omega, f^{-1}_{\theta-\omega}(x)) = (\theta-\omega, f^{-1}(x-g(\theta-\omega))).
\end{align*}
For every $(\theta,x) \in \mathbb{T}^2$ and $n\geq 0$, the projection of $F^{-n}(\theta, x)$ onto the second coordinate is denoted by $f^{-n}_{\theta}(x)$, that is,
$F^{-n}(\theta, x) = (\theta-n\omega, f^{-n}_{\theta}(x))$.

%%%%%%%%%%%%%%%%%%%%%%%%
%       CONSTRUCTIONS
%%%%%%%%%%%%%%%%%%%%%%%%%

\subsection{Constructions}

In this subsection, we construct specific subsets in $\mathbb{T}^2$ which will be essential for our analysis of the dynamics in $(\mathbb{T}^2, F_{\omega})$. By the hypotheses of Theorem \ref{Theorem SNA}, we assume, without loss of generality, the existence of four distinct fixed points that satisfy $$0<e_2<c_2<e_1<c_1<1,$$ where $\{c_i\}_{i=1}^{2}$ are the attracting fixed points, and $\{e_i\}_{i=1}^{2}$ are the repelling fixed points. For each $i \in \{1,2 \}$, we describe the distribution of the fixed points on the torus in terms of the following constants:
\begin{align*}
  d_i & :=  (c_i-e_i) \\
\delta_1 & := (e_1 - c_2)\\
\delta_2 & := (1 - c_1 +e_2),
\end{align*}
which then satisfy the conditions ${400} d_2 < {d_1}$ and $ d_2 < a < d_1 + d_2 + \min \left\{ \delta_1, \delta_2 \right\}$. The constant $a>0$ then satisfies the following inequalities:
\begin{align*}
    d_1 & < |a| < d_1 + \delta_1 + d_2 \\
    d_2 & < |a| < d_2 + \delta_2 + d_1.
\end{align*}
The above conditions allow the construction of regions in $\mathbb{T}^2$ where orbits can experience either a uniform expansion or uniform contraction along the fibre.

%-------------------------- FIGURE Ci Ei Regions 
% trim = left bottom right top
\begin{figure}[t]
    \centering
    \includegraphics[
        width=1\linewidth,
        trim=0cm 0.1cm 0cm 0cm,
        clip
    ]{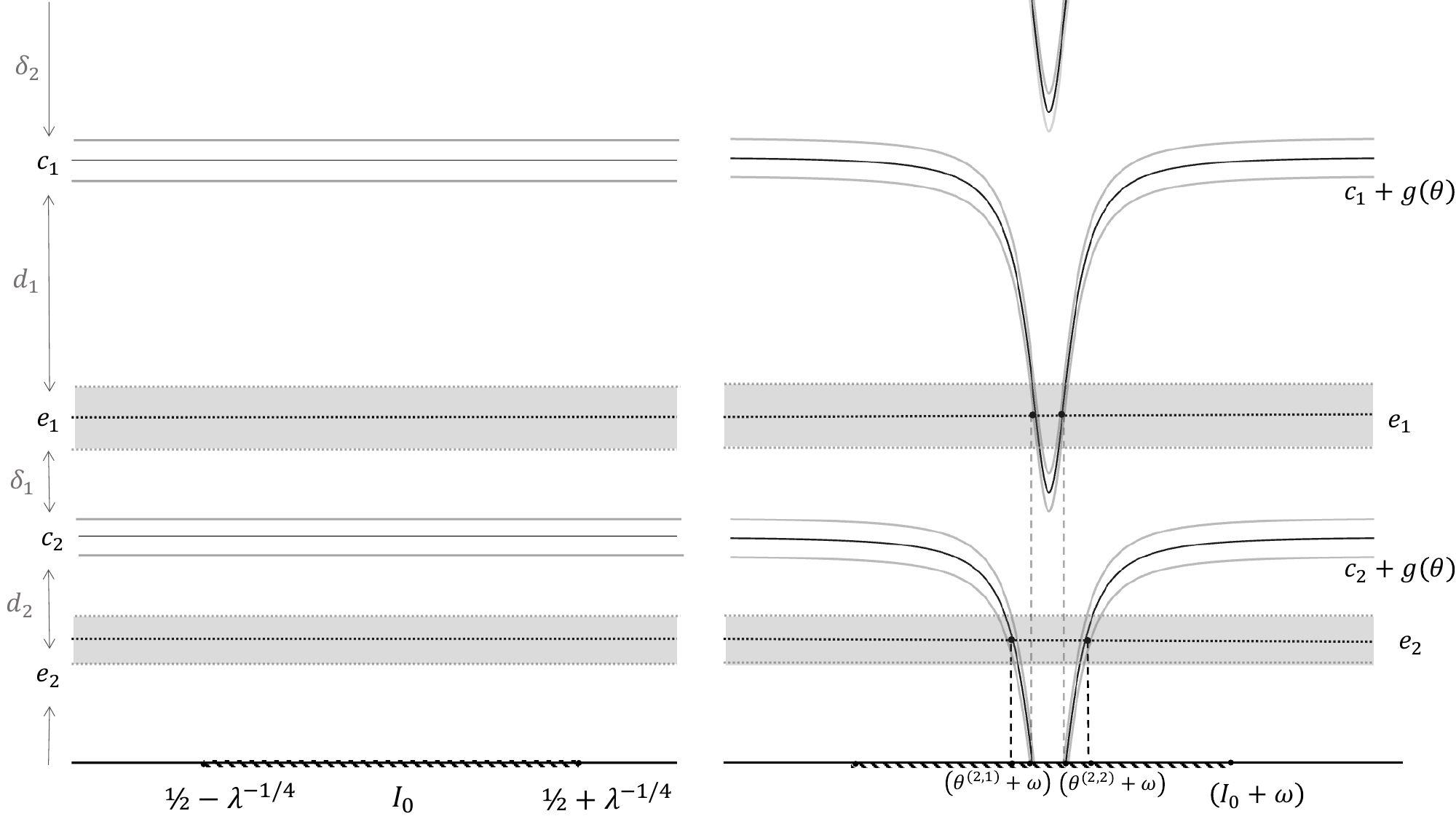}
    \caption{The above figure demonstrates the implications of the conditions in definition \ref{definition_preliminaries_expanding contracting regions}, for the contracting and expanding regions $\{C_i,E_i\}_{i=1}^{2}$. The figure on the right illustrates the sets $\{(\theta+\omega,x): \theta \in I_0,\ x\in [c_i^-+g(\theta), c_i^{+}+g(\theta)]\}$, for $i=1,2$, and their intersection with $\{(I_0+\omega) \times E_i\}_{i=1}^{2}$.}
    \label{figure_expanding contracting regions}
\end{figure}
%-----------------------------

\begin{definition}\label{definition_preliminaries_expanding contracting regions}
    
For the fixed circle diffeomorphism $f$, let $\alpha_u>0$ be such that $\alpha_u^{-1} < |f'(x)| < \alpha_u$, for all $x\in \mathbb{T}$. For the fixed constant $a>0$, we consider a constant $\alpha_c < 1$ and a collection of disjoint compact subsets of $\mathbb{T}$ denoted by $\{E_i, C_i:i=1,2\}$, which satisfy the following conditions:
\begin{enumerate}[label=\alph*)]
    \item For each $i,j\in \{1,2\}$ and $i\neq j$, let $E_i \subset \mathbb{T}$ be such that 
    \begin{enumerate}[label=(\arabic*)]
    \item $E_i = [e_i^{-}, e_i^{+}] := [e_i-\varepsilon_{e_i}, e_i+\varepsilon_{e_i}]$, where $\displaystyle \varepsilon_{e_i} < \min \left\{ \frac{d_i+\delta_{j}+d_{j} -|a|}{2}, \frac{|a|-d_i}{2} \right\}$. Note that $e_i \in E_i$.
        \item For every $x\in E_i$, we have $|f'(x)| > \alpha_c^{-1}$ and equivalently, $|(f^{-1})'(f(x))| < \alpha_c$.
    \end{enumerate}
    \item For each $i,j\in \{1,2\}$ and $i\neq j$, let $C_i\subset \mathbb{T}$ be chosen with respect to the fixed $E_i$ such that 
    \begin{enumerate}[label=(\arabic*)]
    \item $C_i = [c_i^{-}, c_i^{+}] := [c_i-\varepsilon_{c_i}, c_i + \varepsilon_{c_i}]$, where $\displaystyle \varepsilon_{c_i} < \min \left\{\frac{d_i+\delta_i+d_{j} -|a|}{2}, \frac{|a|-d_i-\varepsilon_{e_i}}{2} \right\}$. Note that $c_i \in C_i$.
        \item For every $x\in C_i$, we have $|f'(x)|< \alpha_c$.
    \end{enumerate}
\end{enumerate}
%For ease, we choose $\displaystyle \varepsilon_{c_i} := \frac{\varepsilon_{e_i}}{4}$.
 Since the fixed constants $\alpha_c$ and $\alpha_u$ satisfy $1<\alpha^{-1}_c<\alpha_u$, we fix an integer $b\geq 2$ to be such that
\begin{align}\label{setup_constant b}
  \alpha_c^{b} \alpha_u < 1.
\end{align}
\end{definition}
The above defined contracting and expanding regions will be used extensively to analyse the dynamics on the fibres. Note that by definition \ref{definition_preliminaries_expanding contracting regions}, we have 
\begin{align}\label{preliminaries_f contraction and expansion}
    f([c_i^{-}, c_i^{+}]) & \subset (c_i^{-}, c_i^{+}), \\
f^{-1}([e_i^{-}, e_i^{+}]) & \subset (e_i^{-}, e_i^{+}),
\label{preliminaries_f contraction and expansion_1}
\end{align}
for each $i\in \{1,2\}$. Similarly, there exists an integer $S_0 > 0$ such that 
\begin{align}
f^{S_0}([e_1^{+}, e_2^{-}]) &\subset (c_1^{-}, c_1^{+}) \label{preliminaries_f contraction and expansion S01}\\ f^{S_0}([e_2^{+}, e_1^{-}]) &\subset (c_2^{-}, c_2^{+}) \label{preliminaries_f contraction and expansion S02}\\ f^{-S_0}([c_2^+, c_1^-]) &\subset (e_1^-, e_1^+) \label{preliminaries_f contraction and expansion S03}\\ f^{-S_0}([c_1^+, c_2^-]) &\subset (e_2^-, e_2^+).
\label{preliminaries_f contraction and expansion S04} 
\end{align}

%%%%%%%%%%%%%%%%%%%%%%%%%%%%%%%%%
%     ESTIMATES 
%%%%%%%%%%%%%%%%%%%%%%%%%%%%%%%%%

\subsection{Estimates}

For each $\lambda > 2$, recall that the critical set $I_0$ is given by 
$$ I_0 = \left[ 1/2 - \lambda^{-1/4}, 1/2 + \lambda^{-1/4} \right].$$
In this section, we establish some preliminary and essential estimates for $g(\theta)$ and $g'(\theta)$. Furthermore, we describe specific regions in $I_0$ which are characterised by distinctive dynamical behaviour along their fibre. 
For a sufficiently large $\lambda$, by the definition  $\displaystyle g(\theta)= \frac{-a}{1+\lambda \cos^2(\pi\theta)}$, we have 
\begin{align}\label{equation_preliminaries_g outside I0}
    |g(\theta)| < \lambda^{-1/2}, \text{ for all } \theta \notin I_0.
\end{align}
Then, by equations (\ref{preliminaries_f contraction and expansion} - \ref{preliminaries_f contraction and expansion_1}), we can choose $\lambda$ sufficiently large such that for every $\theta \notin I_0$, 
\begin{align}\label{setup_I_0 contraction charac}
     f_{\theta}([c_i^-, c_i^+]) \subset (c_i^-,c_i^+), \text{ for all } i\in \{1,2\}, 
\end{align}
and for every $\theta \notin (I_0 + \omega)$  
\begin{align}\label{setup_I_0 expansion charac} 
     f^{-1}_{\theta}([e_i^{-}, e_i^{+}]) \subset (e_i^{-}, e_i^{+}), \text{ for all } i\in \{1,2\}. 
\end{align}
Similarly, by equations (\ref{preliminaries_f contraction and expansion S01} - \ref{preliminaries_f contraction and expansion S04}), for sufficiently large $\lambda$ and for every $\omega \in \mathbb{T}\setminus \mathbb{Q}$, if $\{\theta+k\omega\}_{k=0}^{S_0-1} \cap I_0 = \emptyset$ then  
\begin{align}
   f^{S_0}_{\theta}([e_1^{+}, e_2^{-}]) \subset (c_1^{-}, c_1^{+}) \label{equation_preliminaries_f contraction S0 1} \\
    f^{S_0}_{\theta}([e_2^{+}, e_1^{-}]) \subset (c_2^{-}, c_2^{+}), \label{equation_preliminaries_f contraction S0 2}
\end{align}
and if $\{\theta -k\omega\}_{k=-S_0}^{0} \cap (I_0+\omega) = \emptyset$ then  
\begin{align}
     f^{-S_0}_{\theta}([c_2^{+}, c_1^{-}]) \subset (e_1^{-}, e_1^{+})\label{equation_preliminaries_f expansion S0 1} \\
    f^{-S_0}_{\theta}([c_1^{+}, c_2^{-}]) \subset (e_2^{-}, e_2^{+}).
    \label{equation_preliminaries_f expansion S0 2}
\end{align}

%----------

\begin{lemma}\label{lemma_preliminaries_g' estimate}
    Let $c,d > 0$ be such that $0 < c < d < a$. For a sufficiently large $\lambda$, we have the following estimates:
\begin{enumerate}[label=(\alph*)]
    \item For all $\theta \notin I_0$,     \begin{align}\label{equation_preliminaries_g' outside I0}
    |g'(\theta)| < \lambda^{-1/4+1/128}.
    \end{align}
    \item For each $\theta \in I_0$ such that $c  < |g(\theta)| < d < |a|$,  
\begin{align}\label{equation_preliminaries_g' steep in I0}
    \lambda^{1-1/128} & <|g'(\theta)|< \lambda^{1+1/128}.
\end{align}
\end{enumerate}
\end{lemma}

\begin{proof}
By the definition of $g(\theta)$, 
\begin{align*}
|g'(\theta)| & = \left| \frac{\lambda a \pi \sin(2\pi \theta)}{(1+\lambda \cos^{2}(\pi\theta))^2} \right| = \left|\frac{\lambda g(\theta)^{2} \pi \sin(2\pi \theta)}{a} \right|.
\end{align*}
Then, for a sufficiently large $\lambda$, we obtain $ |g'(\theta)| < \lambda^{-1/4+1/128}$ for all $\theta \notin I_0$. For each $\theta \in I_0$ such that $0< c  < |g(\theta)| < d < |a|$, this gives us
\begin{align}
\frac{\lambda \pi |\sin(2\pi\theta)|}{a} c^2 & < |g'(\theta)| < \frac{\lambda \pi |\sin(2\pi\theta)|}{a} d^2.
\end{align}
For a sufficiently large $\lambda$, this proves that $\lambda^{1-1/128} <|g'(\theta)|< \lambda^{1+1/128}$. 
\end{proof}

%---------

For the fixed regions $\{E_i, C_i\}_{i=1}^{2}$, by the conditions in definition \ref{definition_preliminaries_expanding contracting regions}, there exist constants $0< c < d <|a|$ such that 
\begin{align*}
    \left\{ \theta \in \mathbb{T}: x + g(\theta) \in E_1 \cup E_2, \text{ for some } x\in C_1\cup C_2 \right\} \subset \left\{ \theta \in \mathbb{T}: c< |g(\theta)|< d \right\} \subset I_0.
\end{align*}
Thus, for a sufficiently large $\lambda$ with respect to the sets in $\{E_i, C_i\}_{i=1}^{2}$, we have the estimates in equation \eqref{equation_preliminaries_g' steep in I0} (in lemma \ref{lemma_preliminaries_g' estimate}). This in turn refers to the transverse intersections of the contracting and expanding regions under the forcing. It follows from the definition of  $g(\theta)$ that, for each $r \in \mathbb{R}$ which satisfies $-a<r<0$, the function has exactly two preimages given by
\begin{align}
    g^{-1}(r) = \frac{1}{2} \pm \frac{1}{\pi} \sin^{-1} \sqrt{\frac{1}{\lambda}\left(\frac{a}{r}-1\right)}.
\end{align}  
This enables the following definition.

\begin{definition}\label{definition_preliminaries_theta(i,j)}
    For the function $\displaystyle g(\theta) = \frac{-a}{1+\lambda \cos^2( \pi \theta )}$, by the conditions on the constant $a>0$, for each $i\in \{1,2\}$ there exists distinct points $\theta^{(i,1)}$ and $\theta^{(i,2)}$ such that $\theta^{(i,1)} < 1/2 < \theta^{(i,2)}$ and
\begin{align*}
     g(\theta^{(i,1)}) + c_i = g(\theta^{(i,2)}) + c_i = e_i.
\end{align*}
\end{definition}
Further, since $d_2< d_1/400$, the points in $\{\theta^{(i,j)}\}_{i,j=1}^{2}$ satisfy the relation
$$\theta^{(2,1)} < \theta^{(1,1)} < 1/2 < \theta^{(1,2)} < \theta^{(2,2)}.$$ 
By the definition of $g(\theta)$, for each $j\in\{1,2\}$ we have $g'(\theta^{(j,1)}) < 0$ and $g'(\theta^{(j,2)}) > 0$. 

\begin{definition}\label{definition_preliminaries F(1,1)_G and S sets}
    For $i,j \in \{1,2\}$, let the points $\theta^{(i,j)}$ be as defined in definition \ref{definition_preliminaries_theta(i,j)}. Then, we define the following subsets of $I_0$:
\begin{align*}
    G^{(1)} & = (\theta^{(2,1)}, \theta^{(1,1)}) \\
    G^{(2)} & = (\theta^{(1,2)}, \theta^{(2,2)}) \\
     G & = G^{(1)} \cup G^{(2)} \\
    S  & = (\theta^{(1,1)}, \theta^{(1,2)}).
\end{align*}
\end{definition}

Note that if an orbit intersects the set $S$ on the base, then it experiences a strong forcing. In particular, if such an orbit also intersects the contracting region $C_1$ (region $C_2$) along the fibre, under forward iteration it is necessarily forced to the basin of attraction of the point $c_2$ (point $c_1$). Similarly,  the set $G$ is also characterised by the dynamics on the fibre, which is established in subsection \ref{subsection_characteristics G and S}.

%------------------- FIGURE - G S sets
% trim = left bottom right top
\begin{figure}[t]
    \centering
    \includegraphics[
        width=0.9\linewidth,
        trim=1cm 0.1cm 0.8cm 15cm,
        clip
    ]{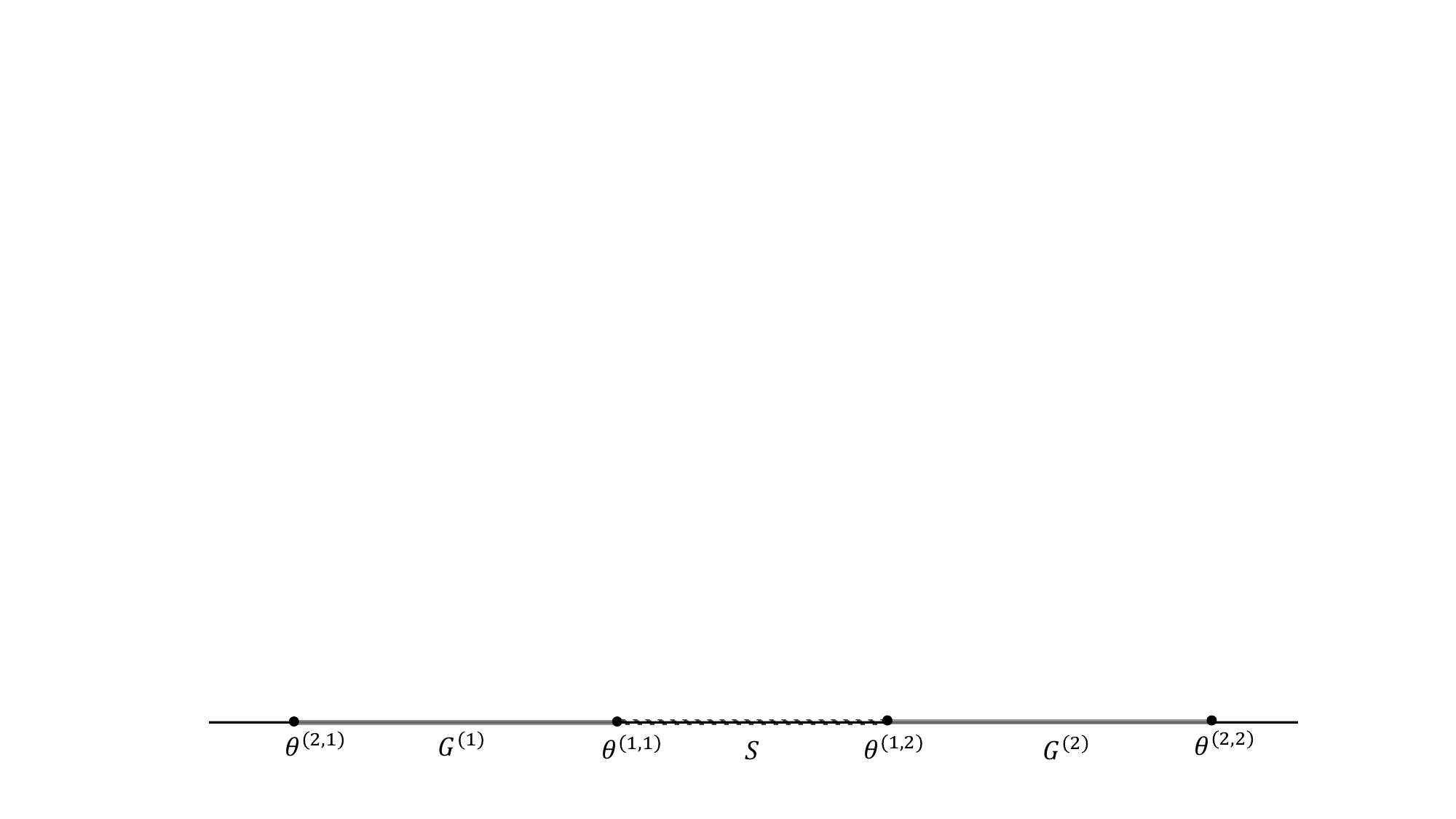}
    \caption{The above figure illustrates the distribution of the points defined in definition \ref{definition_preliminaries_theta(i,j)} and the subsets of $I_0$ defined in definition \ref{definition_preliminaries F(1,1)_G and S sets}.}
    \label{figure_subsets}
\end{figure}
%------------------------

\begin{lemma}\label{lemma_preliminary estimates_|G| and |S|}
    For the fixed $a>0$ as in the hypotheses of theorem \ref{Theorem SNA} and for sufficiently large $\lambda$,
    \begin{align*}
        |G^{(1)}| = |G^{(2)}| > 9|S| > 9\lambda^{-1/2-1/32}.
    \end{align*}
\end{lemma}

\begin{proof}
    By symmetry, we have $|G^{(1)}| = |G^{(2)}|$. From the definition of the points $\theta^{(i,j)}$, we have 
\begin{align}\label{equation_preliminaries_ratio btw distances}
    \frac{d(\theta^{(1,1)}, \theta^{(1,2)})}{d(\theta^{(2,1)}, \theta^{(2,2)})} = \frac{\sin^{-1}\left(\sqrt{\frac{1}{\lambda}(\frac{a}{d_1}-1)}\right)}{\sin^{-1}\left(\sqrt{\frac{1}{\lambda}(\frac{a}{d_2}-1)}\right)}.
\end{align}
By the conditions $d_1< a$, $d_2/d_1 < 1/400$, and the definition of the function $g(\theta)$, 
\begin{align}\label{equation_preliminaries_bound}
    \frac{\sqrt{\frac{1}{\lambda}(\frac{a}{d_1}-1)}}{\sqrt{\frac{1}{\lambda}(\frac{a}{d_2}-1)}} = \sqrt{\frac{a-d_1}{a-d_2}} \sqrt{\frac{d_2}{d_1}} < \frac{1}{20}.
\end{align}
For the sets $(\theta^{(1,1)}, \theta^{(1,2)}) = S$,  $(\theta^{(2,1)}, \theta^{(1,1)}) = G^{(1)}$ and $(\theta^{(1,2)}, \theta^{(2,2)}) = G^{(2)}$, by equations \eqref{equation_preliminaries_ratio btw distances} and \eqref{equation_preliminaries_bound}, we obtain 
\begin{align*}
    \frac{|S|}{|S| + 2|G^{(i)}|} < \frac{1}{20} \implies |G^{(i)}| & > 9|S|.
\end{align*}
Since $g(\theta^{(1,1)}) = d_1$ and $a>d_1$, 
\begin{align*}
    |S| = \frac{2}{\pi} \sin^{-1}\left( \sqrt{\frac{1}{\lambda} \left( \frac{a}{d_1} - 1 \right)} \right).
\end{align*}
Then, for a sufficiently large $\lambda$,
\begin{align}\label{equation_preliminaries_lower estimate for |G| and |S|}
    |S| = \frac{2}{\pi} \sin^{-1}\left( \sqrt{\frac{1}{\lambda} \left( \frac{a}{d_1} - 1 \right)} \right) > \lambda^{-1/2 - 1/32}.
\end{align}

\end{proof}

%\implies |G| > 18  \lambda^{-1/2 - 1/32}
%%%%%%%%%%%%%%%%%%%%%%%%%%%%%%%%%%%%%%%%%%%%%%%%%

%%%%%%%%%%%%%%%%%%%%% BLOCK Mechanism %%%%%%%%%%%%%%%%%%%%
\section{Mechanism of slow-returns}\label{section_Mechanism of slow-returns}
%------------------- BLOCK -1 ---------------------------------------

This section contains the first of the three key components in the proof of theorem \ref{Theorem SNA}. We develop an inductive framework for estimating the rates of contraction along increasingly long segments of orbits. By equations (\ref{setup_I_0 contraction charac}) and (\ref{setup_I_0 expansion charac}), we can estimate the rate of contraction along the fibre until an orbit intersects the set $I_0$ in the base. Since an irrational circle rotation is a minimal, every orbit $(\theta+n\omega)_{n\geq 1}$ intersects $I_0$ infinitely often. This motivates our analysis of parameters and orbits for which intersections with $I_0$ occur only at widely spaced iterations. In subsection \ref{subsection_block 0}, we discuss the mechanism of slow-returns for the critical set $I_0$. In subsection \ref{subsection_block n}, we develop an inductive scheme with respect to a generalised notion of critical sets.

For each $(\theta_0, x_0) \in \mathbb{T}^2$ and $k\in \mathbb{Z}$, we use the notation $(\theta_k, x_k) := F^k(\theta, x)$ for the points in its orbit. For completeness, we include the rates of contraction along parts of an orbit that do not intersect the critical set $I_0$ in the following lemmas.

\begin{lemma}\label{block -1 lemma}\label{lemma_block -1 contraction}
     Let $\omega \in \mathbb{T}\setminus \mathbb{Q}$, $\theta \in \mathbb{T}\setminus I_0$ and $x \in C_1 \cup C_2$. Let $N\in \mathbb{N}$ be the smallest number such that $\theta_N \in I_{0}$. Then, for any $0 \leq k \leq N$,
         \begin{description}
        \item[$(C1)_{-1}$]  \ \ \ \ \ \ \ \ \  
         $\left| f'(x_{N})\dots f'(x_k) \right| \leq \alpha_c ^{N-k+1}$
    \end{description}
\end{lemma}

\begin{proof}
    Since $\{\theta+m\omega\}_{m=0}^{N-1} \cap I_0 = \emptyset$, by equation \eqref{setup_I_0 contraction charac}, for each $k \in \{0,\dots ,N\}$ we have $x_k \in C_1 \cup C_2$ and thus, $|f'(x_k)| < \alpha_c$. This proves that $|f'(x_N)\dots f'(x_k)| \leq \alpha_c^{N-k+1}$, for all $k \in \{0,\dots ,N\}$.
\end{proof}

The following lemma can be proved similarly.

\begin{lemma}\label{lemma_block -1 expansion}
        Let $\omega \in \mathbb{T}\setminus \mathbb{Q}$, $\theta \in \mathbb{T}\setminus (I_0+\omega)$ and $x \in E_1 \cup E_2$. Let $N\in \mathbb{N}$ be the smallest number such that $\theta_{-N} \in I_0$. Then, for any $0 \leq k \leq N-1$ we have
    \begin{description}
        \item[$(E1)_{-1}$]  \ \ \ \ \ \ \ \ \  
       $\left| (f^{-1})'(f(x_{-N+1}))\dots (f^{-1})'(f(x_{-k})) \right| \leq \alpha_c ^{N-k}$
    \end{description}
\end{lemma}

%------------------ BLOCK 0 -------------------------------------
%%%%%%%%%%%%%%%%%%%%%%%%%%%%%%%%%%%%%%%%%%%%%%%%%%%%%%%%%%%%%%%%%%
%------------------- BLOCK 0 ---------------------------------

\subsection{Slow-returns to the critical set $I_0$}\label{subsection_block 0}

In this subsection, we estimate the rate of contraction along the fibre for orbits which return to the critical set $I_0$ slowly. The mechanism requires three main constructions, which are introduced in definitions \ref{definition_block 0_disjointness condn}. \ref{definition_Block 0_Big Theta} and \ref{definition_Block 0 I1 set}. In definition \ref{definition_block 0_disjointness condn}, we introduce slow-return conditions for which every orbit with respect to a feasible parameter $\omega$ intersects $I_0$ at widely spaced iterations. Under such conditions, we specifically focus on orbits of points introduced in definition \ref{definition_Block 0_Big Theta}.  In this step of the inductive framework, our estimates will be limited to segments of orbit that do not intersect a  subset of $I_0$, as in definition \ref{definition_Block 0 I1 set}.

\begin{definition}[Slow return conditions]\label{definition_block 0_disjointness condn}
    Let us fix the constants $\tilde{M_0} = M_0 = \lambda^{1/16}$ and $K_0 = (b+1)2^8$, where $\lambda>2$ and $b>0$ is such that $\alpha_c^b\alpha_u < 1$, as in equation \eqref{setup_constant b}. Then, we define $\mathcal{F}_0 = \mathcal{F}_0(M_0, f, \lambda)$ as the set of all parameters in $\mathbb{T}$ for which the following condition is satisfied:
\begin{description} 
\item[$(\mathcal{F}1)_0$] \label{desc:F1_0}
\begin{align*}
    d\left(3I_0, \bigcup_{m=1}^{2K_0M_0}3(I_0+m\omega)\right) > 0.
\end{align*}
\end{description}
As a convention, let $\mathcal{F}_{-1} := \mathbb{T}$
\end{definition}
For $\lambda$ big enough such that $2S_0<M_0$, for every $\omega \in \mathcal{F}_0$ we also have 
\begin{align*}
    d \left( I_0\pm (M_0+S_0)\omega, \bigcup_{|m|\leq M_0} (I_0+m\omega) \right) > 0,
\end{align*}
which will prove to be a useful consequence of the condition $(\mathcal{F}1)_0$.

\begin{definition}\label{definition_Block 0_Big Theta}
For the fixed $M_0$, let
\begin{align*}
\Theta_{0} := \mathbb{T}\setminus \bigcup_{|m|\leq M_0 - 1}(I_0 + m\omega) 
\end{align*}
and as a convention, let $\Theta_{-1} := \mathbb{T}$.   
\end{definition}

We observe that for every $\omega \in \mathcal{F}_0(M_0)$ and $\theta \in \Theta_0$, any two successive visits of the orbit $(\theta+m\omega)_{m\in \mathbb{Z}}$ to the critical set $I_0$ are spaced by at least $2K_0M_0$ iterations. This formalises the idea of slow-return to critical set $I_0$. 

\begin{definition}[Critical set $I_1$]\label{definition_Block 0 I1 set}
    For any fixed $\omega \in \mathcal{F}_0$, a critical set $I_1(\omega)$ is a finite union of open intervals, which satisfies the following conditions:
\begin{enumerate}
    \item $I_1(\omega) \subset I_0$.
    \item The set $I_1(\omega)$ consists of at most four intervals.
    \item $I_1(\omega) \supset \pi_1 \left( F^{M_0}(A_0) \bigcap F^{-M_0}(B_0) \right)$, where
\begin{align*}
    A_0 & := (I_0 - M_0 \omega) \times (C_1\cup C_2) \\
    B_0 & := (I_0 + M_0\omega) \times (E_1 \cup E_2).
\end{align*}
\end{enumerate}
\end{definition}

The above definition thus captures a critical behaviour within $I_0$ and broadens the notion of a critical set. The critical sets $I_1$ also  mark the limitation of the estimates developed in this subsection, which is established in the propositions in the following sub-subsection.

%%%%%%%%%%%%%%%%%%%%%%%%%%%%%%%%%%%%%%%%%%%%%%%%%%%%%%%%%%%%%%%%%%
%------------------- MAIN PROPOSITIONS OF I0 ---------------------

\subsubsection{Main propositions corresponding to  $I_0$}

We now establish the rates of contraction corresponding to the points in $\displaystyle \Theta_0$ and parameters that satisfy the disjointness condition $(\mathcal{F}1)_0$. Let us fix $$\displaystyle \beta_0 := \frac{(M_0+S_0)}{(K_0M_0)}.$$ 
Since $2S_0 < M_0$ and $K_0 = (b+1)2^8$, note that $\displaystyle 0< \beta_0 < \frac{1}{4(b+1)}$, where $b>0$ satisfies the condition $\alpha_u < \alpha_c^{-b}$.

%------------------- PROP BEGIN -------------------

\begin{proposition}\label{proposition_block 0 contraction}
        Let $\omega \in \mathcal{F}_{0}$ and  $I_1(\omega)$ be a non-empty critical set as in definition \ref{definition_Block 0 I1 set}. Let $\theta_0 \in \Theta_{0}$, $x_0 \in C_1 \cup C_2$, and $N\geq 0$ be the smallest integer such that $\theta_N \in I_{1}(\omega)$. Then, we have the following:
    \begin{description}
        \item[$(C1)_0$] For all $k \in \{0,\dots,N\}$, 
        \begin{align}\label{equation_block 0 contraction_C1_0}
            x_{k} \notin C_1 \cup C_2 \implies \theta_k \in \bigcup_{m=1}^{M_0+S_0-1}(I_0+m\omega).
        \end{align}
        \item[$(C2)_0$] For all $k \in \{0,\dots,N\}$ and $\displaystyle \beta_0 = {(M_0+S_0)}/{K_0M_0}$, 
        \begin{align*}
             \left| f'(x_{N})\dots f'(x_k) \right| \leq \left( \alpha_c^{1-(1+b)\beta_0} \right)^{N-k+1}.
        \end{align*}
         \item[$(C3)_0$] $\left|\partial_x (\pi_2 \circ F^{N+1}(\theta_0, x_0))\right| \leq \left( \alpha_c^{1-(1+b)\beta_0} \right)^{N+1}$, where  $\displaystyle \beta_0 = {(M_0+S_0)}/{K_0M_0}$.
    \end{description}
\end{proposition}

\begin{proof}
    For fixed $\omega \in \mathcal{F}_0$ and $\theta_0 \in \Theta_0$, there exists a sequence of all the integers $0 \leq k\leq N$ for which $\theta_{k}\in I_0$. Let us denote it by $0<N_1 < N_2 <\dots < N_J=N$. We fix a point $x_0 \in (C_1\cup C_2)$.
    
\vspace{0.2cm}
\noindent \textit{Proof of $(C1)_0$:} Since $N_1$ is the smallest positive integer for which $\theta_{N_1} \in I_0$, by $(C1)_{-1}$ (in lemma \ref{lemma_block -1 contraction}) we have $\pi_2\circ F^k(\theta_0,x_0) = x_k \in (C_1 \cup C_2)$, for all $0\leq k \leq N_1$. If $J=1$, then this is proves $(C1)_0$.  Let us consider the case $J>1$, where it is thus sufficient to prove the implication for the segment $(\theta_k,x_k)_{k=N_1}^{N_{J}}$. For each $1 \leq j \leq J-1$, we now define a new property $(C1)_0^{N_{j+1}}$ for the stretch $\{(\theta_m, x_m)\}_{m=N_j}^{N_{j+1}}$ given by the following: 
\begin{description}
    \item[$(C1)_0^{N_{j+1}}$] For any $N_j \leq k \leq N_{j+1}$,
    \begin{align}\label{equation_block 0 contraction_C1_0[N_j-1]}
        x_k \notin C_1 \cup C_2 \implies \theta_k \in \bigcup_{m=1}^{M_0+S_0-1}(I_0+m\omega).
    \end{align}
\end{description}
Note that to prove $(C1)_0$, it is thus sufficient to prove $(C1)_0^{N_{j+1}}$ for all $j=1,\dots, J-1$. Since $N_1\neq N_J$, we infer that $\theta_{N_1} \in I_0\setminus I_1$. Then, from the definition of set $I_1$ and since  $x_{N_1-M_0} \in C_1 \cup C_2$, we can conclude that $x_{N_1+M_0} \notin E_1 \cup E_2$. By the disjointness condition $(\mathcal{F}1)_0$ (in definition \ref{definition_block 0_disjointness condn}), it follows that $$(\theta+m\omega)_{m=N_1+M_0}^{N_1+M_0+S_0} \bigcap I_0 = \emptyset,$$ and thus, by equations (\ref{equation_preliminaries_f contraction S0 1} - \ref{equation_preliminaries_f contraction S0 2}), we obtain $x_{N_1+M_0+S_0} \in C_1 \cup C_2$. Meanwhile on the base, by the definition of $\Theta_0$ (in definition \ref{definition_Block 0_Big Theta}) and the disjointness condition $(\mathcal{F}1)_0$ (in definition \ref{definition_block 0_disjointness condn}), we have $\theta_{N_1+M_0+S_0} \in \Theta_0$. By the choice of constants $N_1$ and $N_2$ we have $$(\theta +m\omega)_{m=N_1+M_0+S_0}^{N_2-1} \cap I_0 = \emptyset,$$ and thus it follows that $x_k \in C_1\cup C_2$ for all $0 \leq k \leq N_2$. This proves the property $(C1)_0^{N_{j+1}}$, as defined in equation \eqref{equation_block 0 contraction_C1_0[N_j-1]}, for the case $j=1$. 
    
    Let us consider the case when $J>2$. Say there exists $1\leq {j'} \leq (J-2)$ such that the property $(C1)_0^{N_{ {j'}+1}}$ holds true.  Since $\theta_{N_{j'}-M_0} \in (I_0-M_0\omega)$, we infer that $x_{N_{j'}-M_0} \in C_1 \cup C_2$. Since $\theta_{N_{j'}} \in I_0 \setminus I_1$, it follows from the definition of $I_1$ that $x_{N_{j'}+M_0} \notin E_1 \cup E_2$. By the disjointness condition $(\mathcal{F}1)_0$ and equations (\ref{equation_preliminaries_f contraction S0 1} - \ref{equation_preliminaries_f contraction S0 2}), we can conclude that  $\theta_{N_{j'}+M_0+S_0} \in \Theta_0$ and $x_{N_{j'}+M_0+S_0} \in C_1\cup C_2$. Then, by $(C1)_{-1}$, we have $x_k \in C_1\cup C_2$ for all $(N_{j'}+M_0+S_0) \leq k \leq N_{{j'}+1}$. This proves $(C1)_0^{N_{{j'}+1}}$, that is for all $0 \leq k \leq N_{{j'}+1}$ 
    \begin{align*}
        x_k \notin C_1 \cup C_2 \implies \theta_k \in \bigcup_{m=1}^{M_0+S_0-1}(I_0+m\omega).   
    \end{align*}
This ends the proof of induction for the property $(C1)_0^{N_{j+1}}$, for $1\leq j \leq J-1$. 

\vspace{0.2cm}
\noindent \textit{Proof of $(C2)_0$:} By the disjointness condition $(\mathcal{F}1)_0$ (in definition \ref{definition_block 0_disjointness condn}), it follows that $\theta_k \notin I_0,$ for all $0 \leq k < N_1$, and this gives us $x_m \in C_1 \cup C_2$ for all $0 \leq m \leq N_1$. Then, by definition of the regions $\{C_i\}_{i=1}^{2}$ (in definition \ref{definition_preliminaries_expanding contracting regions}), we obtain the bound
    \begin{align*}
          |f'(x_{N_1})\dots f'(x_k)| \leq \alpha_c^{N_1-k+1} \leq \left( \alpha_c^{1-\beta_0} \alpha_u^{\beta_0}\right)^{N_1-k+1}.
    \end{align*}
 In the case $J=1$, this is proves $(C2-3)_0$. Let us consider when $J>1$. For each $1 \leq j \leq J-1$, we now define a new property for the stretch $\{(\theta_m, x_m)\}_{m=N_j}^{N_{j+1}}$, given by the following: 
\begin{description}
    \item[$(C2)_0^{N_{j+1}}$] For any $N_j \leq k \leq N_{j+1}$,
    \begin{align}\label{equation_block 0 contraction_C2_0[N_j-1]}
        |f'(x_{N_{j+1}})\dots f'(x_k)| \leq \left( \alpha_c^{1-\beta_0} \alpha_u^{\beta_0} \right)^{N_{j+1}-k+1}.
    \end{align}
\end{description}
To prove $(C2)_0$, it is thus sufficient to prove $(C2)_0^{N_{j+1}}$ for all $j \in \{1,\dots, (J-1)\}$. Let us fix some $1\leq j \leq J-1$. By disjointness condition $(\mathcal{F}1)_0$, it follows that $|N_{j+1}-N_j|\geq 4K_0M_0$. Since $\theta_{N_j} \in I_0\setminus I_1$, by $(C1)_0$ and disjointness condition $(\mathcal{F}1)_0$, we deduce that $x_k\in C_1\cup C_2$, for all $N_j+M_0+S_0 \leq k \leq N_{j+1}$.
Thus, for any $N_j+M_0+S_0 \leq k \leq N_{j+1} $, 
\begin{align}\label{equation_block 0 contraction_C2_0_1}
    |f'(x_{N_{j+1}})\dots f'(x_{k})| \leq \alpha_c^{N_{j+1}-k+1} \leq \left( \alpha_c^{1-\beta_0} \alpha_u^{\beta_0} \right)^{N_{j+1}-k+1} < \left( \alpha_c^{1-(1+b)\beta_0} \right)^{N_{j+1}-k+1}.
\end{align}
Since $\theta_{N_j} \in I_0\setminus I_1$, for all $N_j\leq k \leq N_j+M_0+S_0-1$,
\begin{align*}
    \theta_{k} \in \bigcup_{m=0}^{M_0+S_0-1}(I_{0}+m\omega).
\end{align*}
For each $N_j\leq k \leq N_j+M_0+S_0-1$, by the disjointness condition $(\mathcal{F}1)_0$, we can deduce that $N_j+M_0+S_0-k \leq M_0+S_0 < K_0M_0$, $2K_0M_0 \leq N_{j+1}-N_j$ and $K_0M_0 \leq N_{j+1}-k+1$. This gives us the following bounds:
\begin{align*}
\frac{N_j+M_0+S_0-k}{N_{j+1}-k+1} & \leq \beta_0 = \frac{M_0+S_0}{K_0M_0}.\\
    \frac{N_{j+1}-N_j-M_0-S_0+1}{N_{j+1}-k+1} & \geq 1-\beta_0. 
\end{align*}
We now derive the minimum rate of contraction along the fibre, by considering the upper bound $|f'(x_k)|< \alpha_u$, for all $N_j\leq k \leq N_j+M_0+S_0-1$. This gives us 
\begin{align*}
    |f'(x_{N_{j+1}})\dots f'(x_{N_j+M_0+S_0})f'(x_{N_j+M_0+S_0-1})\dots & f'(x_k)| 
     \leq \alpha_c^{N_{j+1}-N_j-M_0-S_0+1} \alpha_u^{N_j+M_0+S_0-k}\\
    & \leq \left( \alpha_c^{\frac{N_{j+1}-N_j-M_0-S_0+1}{N_{j+1}-k+1}} \alpha_u^{\frac{N_j+M_0+S_0-k}{N_{j+1}-k+1}} \right)^{N_{j+1}-k+1},
\end{align*}
for all $N_j\leq k \leq N_j+M_0+S_0-1$. 
The estimates above prove that for all $k=N_j,\dots ,N_{j}+M_0+S_0-1$,
\begin{align}\label{equation_block 0 contraction_C2_0_2}
    |f'(x_{N_{j+1}})\dots f'(x_k)| \leq \left( \alpha_c^{1-\beta_0} \alpha_u^{\beta_0} \right)^{N_{j+1}-k+1} <  \left( \alpha_c^{1-(1+b)\beta_0} \right)^{N_{j+1}-k+1}.
\end{align}
The equations \eqref{equation_block 0 contraction_C2_0_1} and \eqref{equation_block 0 contraction_C2_0_2} conclude the proof of $(C2)_0^{N_{j+1}}$, for all $j\in \{1,\dots ,J-1\}$. 

\vspace{0.2cm}
\noindent\textit{Proof of $(C3)_0$:} For $(\theta_0,x_0)$ and any $n\geq 1$, by an application of the chain rule we have
\begin{align*}
    \partial_x \left( \pi_2 \circ F^{n} (\theta_0, x_0) \right) =  f'(x_{n-1}) \dots f'(x_0).
\end{align*}
By $(C1)_0$, we have $|f'(x_N)|<\alpha_c$. Then, by $(C2)_0$, it follows that
\begin{align*}
    \left|\partial_x \left( \pi_2 \circ F^{N+1} (\theta_0, x_0) \right)\right| =  |f'(x_{N}) \dots f'(x_0)| < \left( \alpha_c^{1-(1+b)\beta_0} \right)^{N+1}.
\end{align*}

\end{proof}

%------------------- PROP END -------------------

For every $\theta \in \Theta_0$ and $x,y \in C_1\cup C_2$, the above proposition describes the local behaviour for certain finite segments of the forward orbits of $(\theta,x)$ and $(\theta,y)$. The following proposition clarifies the global behaviour of orbits along the fibre, which will be essential for results on synchronisation.

%------------------- PROP BEGIN -------------------

\begin{proposition}\label{proposition_Block 0 sync contraction}
Let $\omega \in \mathcal{F}_{0}$ and $I_1(\omega)$ be a non-empty critical set, as defined in definition \ref{definition_Block 0 I1 set}. Let $\theta_0 \in \Theta_{0}$, $x_0,y_0 \in C_i$ for some $i\in \{1,2\}$, and $N\in \mathbb{N}$ be the smallest integer such that $\theta_N \in I_{1}(\omega)$. Then, we have the following:
    \begin{description}
        \item[$(Sc1)_0$]  For all $k \in \{0,1,\dots,N\}$,
        \begin{align}\label{equation_block 0 contraction_C2_0[N_j-1]}
            \theta_k \notin \bigcup_{m=1}^{M_0+S_0-1}(I_0+m\omega) \ \ \Longrightarrow \ \  x_k,y_k \in C_i, \ \ \text{ for some } i\in \{1,2\}.
        \end{align}
    \end{description}
\end{proposition}

\begin{proof}
    Let $0<N_1<\cdots <N_J = N$ be the sequence of all integers for which $\theta_{N_j}\in I_0$. Let us fix $i \in \{1,2\}$ and some $x \in C_i$. Since $(\theta + m\omega)_{m=0}^{N_1} \cap I_0 = \emptyset$, we also have $x_k \in C_i$ for all $k=0,\dots ,N_1$. Since $(\{\theta\} \times C_i) \subset \mathbb{T}^2$ is a connected set and $F$ is continuous, this gives us
    \begin{align}\label{equation_block 0 contraction_S0_1}
       \left( F^{k}\left(\{\theta\} \times C_i\right)\right) \subset C_i,
    \end{align}
    for all $k \in \{0,\dots ,N_1\}$. In the case when $J=1$, this is proves $(Sc1)_0$. 
    
    Let us now consider the case when $J>1$ and assume that $x_{N_{j'}-M_0}, y_{N_{j'}-M_0} \in C_i$, for some $i\in\{1,2\}$, for some fix ${j'} \in \{1,\dots ,(J-1)\}$. Without loss of generality, we consider $x_{N_{j'}-M_0} < y_{N_{j'}-M_0}$, which gives us
\begin{align}\label{equation_block 0 contraction_S0_2}
    \{\theta_{N_{j'}-M_0}\} \times [x_{N_{j'}-M_0}, y_{N_{j'}-M_0}] \subset \{\theta_{N_{j'}-M_0}\} \times C_i.
\end{align}
    Since $\theta_{N_{j'}} \in I_0 \setminus I_1$, by disjointness condition $(\mathcal{F}1)_0$ (in definition \ref{definition_block 0_disjointness condn}), we have $\theta_{N_{j'}+M_0+S_0} \notin \cup_{m=1}^{M_0+S_0-1}(I_0+m\omega)$ and by $(C1)_0$ (in lemma \ref{proposition_block 0 contraction}), this gives us 
    $$x_{N_{j'}+M_0+S_0},y_{N_{j'}+M_0+S_0} \in (C_1\cup C_2).$$ 
    Say there exists $k\in [N_{j'}+M_0+S_0, N_{{j'}+1}]$ such that $x_k\in C_2$ and $y_k\in C_1$, then we also have $x_{N_{j'}+M_0+S_0} \in C_2$ and $y_{N_{j'}+M_0+S_0} \in C_1$. Since the sets in equation (\ref{equation_block 0 contraction_S0_2}) are connected, and $F$ is continuous and orientation preserving, 
\begin{align}
    F^{2M_0+S_0}\left(\theta_{N_{j'}-M_0}, [x_{N_{j'}-M_0}, y_{N_{j'}-M_0}]\right) & = \{\theta_{N_{j'}+M_0+S_0}\} \times [x_{N_{j'}+M_0+S_0}, y_{N_{j'}+M_0+S_0}].
\end{align}
By the distribution of the regions $\{C_i,E_i\}_{i=1}^{2}$, this gives us
\begin{align*}
  \left( \{\theta_{N_{j'}+M_0+S_0}\} \times E_1 \right) \subset \left( \{\theta_{N_{j'}+M_0+S_0}\} \times [x_{N_{j'}+M_0+S_0},y_{N_{j'}+M_0+S_0}]\right).
\end{align*}
By the condition $(\mathcal{F}1)_0$, we have $F^{-S_0}((I_0+(M_0+S_0)\omega) \times E_1) \subset ((I_0+M_0\omega) \times E_1)$, which then implies
\begin{align*}
    F^{M_0}\left(\theta_{N_{j'}-M_0} \times C_i \right) \cap F^{-M_0} \left(\theta_{N_{j'}+M_0} \times E_1\right) & \neq \emptyset,
\end{align*}
due to equation \eqref{equation_block 0 contraction_S0_2}. However, this leads to a contradiction as $\theta_{N_{j'}}\notin I_1$. A similar reasoning follows for when $x_k\in C_1$ and $y_k\in C_2$, for some $k\in [N_{j'}+M_0+S_0, N_{{j'}+1}]$. This proves that, for each $k\in[N_{j'}+M_0+S_0, N_{{j'}+1}]$ there exists some $l\in \{1,2\}$ such that $x_{k}, y_{k} \in C_l$. This, along with equation \eqref{equation_block 0 contraction_S0_1}, concludes our inductive proof in $j \in \{0,\dots ,(J-1)\}$ and thus the proof of the lemma.
\end{proof}

%------------------- PROP END -------------------

Similar to propositions \ref{proposition_block 0 contraction} and \ref{proposition_block 0 sync' contraction}, the following propositions \ref{proposition_block 0 expansion} and \ref{proposition_block 0 sync' expansion} can be proved similarly for the backward orbits of points in $\Theta_0$, and with respect to the expanding regions $\{E_i\}_{i=1}^{2}$.

\begin{proposition}\label{proposition_block 0 expansion}
       Let $\omega \in \mathcal{F}_{0}$ and $I_1(\omega)$ be a non-empty critical set, as in definition \ref{definition_Block 0 I1 set}. Let $\theta_0 \in \Theta_{0}$, $x_0 \in E_1 \cup E_2$ and $N\in \mathbb{N}$ be the smallest positive integer such that $\theta_{-N} \in I_{1}(\omega)$. Then, we have the following:  
    \begin{description}
        \item[$(E1)_0$] For all $k\in \{0,\dots,N\}$,  
        \begin{align}\label{equation_block 0 expansion_E1_0}
        x_{-k} \notin E_1 \cup E_2 \implies \theta_{-k} \in \bigcup_{m=0}^{M_0+S_0-1}(I_0-m\omega).
        \end{align}
        \item[$(E2)_0$] For all $k\in \{0,\dots,N-1\}$, and $\displaystyle \beta_0 = {(M_0+S_0)}/{K_0M_0}$, 
        \begin{align*}
        \left| (f^{-1})'(f(x_{-N+1}))\dots (f^{-1})'(f(x_{-k})) \right| \leq \left( \alpha_c^{1-(1+b)\beta_0} \right)^{N-k}.
        \end{align*}
        \item[$(E3)_0$] $\left|\partial_x (\pi_2 \circ F^{-N}(\theta_0, x_0))\right| \leq \left( \alpha_c^{1-(1+b)\beta_0} \right)^{N}$, where $\displaystyle \beta_0 = {(M_0+S_0)}/{K_0M_0}$.
    \end{description}
\end{proposition}

\begin{proposition}\label{proposition_Block 0 sync expansion}
    Let $\omega \in \mathcal{F}_{0}$ and $I_1(\omega)$ be a non-empty critical set as in definition \ref{definition_Block 0 I1 set}. Let $\theta_0 \in \Theta_{0}$, $x_0,y_0 \in E_i$ for some $i\in \{1,2\}$, and let $N\in \mathbb{N}$ be the smallest positive integer such that $\theta_{-N} \in I_{1}(\omega)$. Then,
        \begin{description}
        \item[$(Se0)_0$]  For any $k \in \{0,1,\dots,N-1\}$, 
        \begin{align*}
            \theta_{-k} \notin \bigcup_{m=0}^{M_0+S_0-1}(I_0-m\omega) \ \ \implies \ \  x_{-k},y_{-k} \in E_i, \ \ \text{ for some } i\in \{1,2\}.
        \end{align*}
        \end{description}
\end{proposition}

%------------------------ COMMENT ---------------------
\begin{comment}
To summarise, if we consider certain points iterated forward until it intersects $I_1$, denoted by $\{(\theta_i,x_i)\}_{i=0}^{N}$. The constant $\beta_0$ satisfies the following conditions, for any $0 \leq k \leq N$:
\begin{align*}
    \#\left\{m\in [k,N]: x_m\notin C_1\cup C_2 \right\} & \leq \beta_0 (N-k+1) \\
    \#\left\{m\in [k,N]: x_m\in C_1\cup C_2 \right\} & \geq (1- \beta_0) (N-k+1) 
\end{align*}
Similarly, if we consider points iterated backwards until it hits $I_1+\omega$, that is of the form $\{(\theta_{-i},x_{-i})\}_{i=0}^{N}$.
\begin{align*}
    \#\left\{m\in [k,N]: y_{-m}\notin E_1\cup E_2 \right\} & \leq \beta_0 (N-k+1) \\
    \#\left\{m\in [k,N]: y_{-m}\in E_1\cup E_2 \right\} & \geq (1- \beta_0) (N-k+1)
\end{align*}    
\end{comment}

%-------------------------- Block n ----------------------
%%%%%%%%%%%%%%%%%%%%%%%%%%%%%%%%%%%%%%%%%%%%%%%%%%%%%%%%%%%%%%%%%%
%------------------- BLOCK n ---------------------------------

\subsection{Slow-returns to the critical set $I_n$}\label{subsection_block n}

In this subsection, we establish an inductive framework for estimating the rate of convergence along the fibre, for increasingly long stretches of orbits. Analogous to subsection \ref{subsection_block 0}, we impose slow-return conditions on parameter $\omega$ in definition \ref{definition_block n_disjointness condn}. In definition \ref{definition_block n_Big theta}, we introduce the set of points for which the orbits can be studied under the slow-return conditions. In definition \ref{definition_block n_sync condn parameter}, we introduce conditions for the parameters $\omega$ which would specifically enable the analysis of synchronisation of orbits along the fibre direction.

We now generalise the notion of critical set and construct a sequence of nested sets in $I_0$, which are characterised by a strong critical behaviour. When the critical sets are non-empty, they can be the source of non-uniform dynamical behaviour in the system. 

\begin{definition}[Critical Sets]
\label{definition_block n_Critical sets In}
    For $n\geq 1$ and $\omega \in \mathbb{T}$, let $0< M_{-1} < M_0 < M_1 < \cdots < M_{n-1} $ be positive integers, where $M_i \in [\tilde{M_i}, 2\tilde{M_i}]$ and $\tilde{M}_i := \alpha_c^{-\tilde{M}_{i-1}/8}$ for all $1\leq i\leq n-1$. Then, we can define a nested sequence of sets $I_0 \supset I_1 \supset \dots \supset I_n$ such that, for each $1\leq i \leq n$ the set $I_i = I_i(\omega)$ is a closed set that satisfies the following conditions:
    \begin{description}
    \item[$(\mathcal{I}1)_i$]\label{(F1)_i} $I_i \subset I_{i-1}$.
    \item[$(\mathcal{I}2)_i$]\label{(F2)_i} The set $I_i$ consists of at most four intervals.
    \item[$(\mathcal{I}3)_i$]\label{(F3)_i} $I_i \supset \pi_1 \left( F^{M_{i-1}}(A_{i-1}) \bigcap F^{-M_{i-1}}(B_{i-1}) \right)$, where
\begin{align*}
    A_{i-1} & = (I_{i-1} - M_{i-1} \omega) \times C_1\cup C_2 \\
    B_{i-1} & = (I_{i-1} + M_{i-1}\omega) \times E_1 \cup E_2.
\end{align*}
\end{description}
Each set $I_i$ in such a sequence is called a critical set.
\end{definition}

Note that the trivial sequence $(I_n)_{n\geq 1} = (I_0)_{n\geq 1}$ also satisfies the conditions in definition \ref{definition_block n_Critical sets In}. Thus, the mechanism of slow-returns in section \ref{section_Mechanism of slow-returns} is applicable for broad class of critical sets. However, we will primarily focus on certain non-trivial critical sets. More specifically, for the proof of theorem \ref{Theorem SNA}, we will be interested in a unique sequence of critical sets corresponding to each fixed $\omega$, which will be discussed later in definition \ref{definition_structure In_critical set In}.

\begin{definition}[Slow return conditions]\label{definition_block n_disjointness condn}
    Let $n \geq 1$ be an integer. Consider the  increasing sequences of positive integers $\displaystyle (K_n)_{n\geq 1}$ and $\displaystyle  (M_n)_{n\geq 1}$, where $\displaystyle K_n := (b+1)2^{8+n}$, $M_n \in [\tilde{M}_n, 2\tilde{M}_n]$ where $\displaystyle  \tilde{M}_n := \alpha_c^{-\tilde{M}_{n-1}/8}$ for all $n\geq 1$, and $\displaystyle \tilde{M}_0 = M_0 = \lambda^{1/16}$. Then, we define $\displaystyle \tilde{\mathcal{F}}_n(M_n,\dots,M_0,f,\lambda)$ as the set of all parameters in $\mathbb{T}$ for which the sequence of critical sets $(I_j)_{j=0}^{n}$ (as in definition \ref{definition_block n_Critical sets In}) satisfy the following conditions:         
\begin{description} 
\item[$(\mathcal{F}1)_n$]\label{(F1)_n} For each $j\in \{0, 1,\dots ,n\}$,
\begin{align*}
    d\left(3I_j, \bigcup_{m=1}^{2K_jM_j}3(I_j+m\omega)\right) > 0.
\end{align*}

\item[$(\bar{\mathcal{F}}1)_n$] For each $j \in \{1,\dots,n\}$,
         \begin{align*}
    d\left(I_j, \bigcup_{m=1}^{2K_jM_j}(I_j+m\omega)\right) > \varepsilon_j,
\end{align*}
where $\varepsilon_j := 3\alpha_c^{\tilde{M}_{j-1}/3}$ and $\varepsilon_1 := 3\alpha_c^{M_0}$.

\item[$(\mathcal{F}2)_n$] For each $j \in \{1,\dots,n\}$,
         \begin{align*}
             d\left(I_j \pm (M_j+S_0)\omega, \bigcup_{l=0}^{j-1} \bigcup_{|m|\leq 2M_l} (I_l + m\omega)\right) > 0.
         \end{align*}
         
\end{description}
Further, let us define the set
\begin{align*}
    \tilde{\mathcal{F}}_n(f,\lambda) = \left\{ \omega\in \mathbb{T}: \omega \in  \tilde{\mathcal{F}}_n(M_n,\dots,M_0,f,\lambda), \text{ for some } M_i \in [\tilde{M}_i, 2\tilde{M_i}] \text{ for } i=1,\dots ,n\right\}.
\end{align*}
\end{definition}
Observe that for every $n\geq 1$ we have $\tilde{\mathcal{F}}_{n+1} \subset \tilde{\mathcal{F}}_n \subset {\mathcal{F}}_0$. For each $\omega \in \tilde{\mathcal{F}}_n$, since $S_0<2M_0<\dots <M_n$, the following disjointness properties are satisfied: 
\begin{align*}
    d \left(I_j\pm M_j \omega, \bigcup_{|m|\leq M_0}(I_0+m\omega) \bigcup_{l=1}^{j-1}\bigcup_{|m|\leq M_l+M_1}(I_l+m\omega) \right) & > 0, \text{ for all } 1\leq j \leq n.
\end{align*}
It is worth noting that the disjointness condition $(\bar{\mathcal{F}}1)_n$ can in fact be proved to be stronger than $({\mathcal{F}}1)_n$, as a consequence of the geometric structure of the critical sets. This will discussed further in section \ref{section_parameter space}.

\begin{definition}\label{definition_block n_Big theta}
    For each $n\geq 1$ and $\omega \in \tilde{\mathcal{F}}_n$, let us define the following subsets of $\mathbb{T}$:
\begin{align}\label{equation_block n_Big theta tilde n}
    \Theta_{n} & := \mathbb{T}\setminus \bigcup_{l=0}^{n}\bigcup_{|m|\leq M_l -1}(I_l + m\omega), \\
    \tilde{\Theta}_{n} & := \Theta_0 \setminus \bigcup_{l=1}^{n}\bigcup_{|m|\leq 2M_l}(I_l + m\omega).
    \stepcounter{equation}\tag{\theequation} \label{equation_block n_Big theta n}
\end{align}
\end{definition}

Then, for each $n\geq 1$ we have $\Theta_n \subset \Theta_{n-1}$ and $\tilde{\Theta}_{n+1} \subset \tilde{\Theta}_n$. For any $\omega \in \mathcal{F}_n(M_n,\dots,M_0,f,\lambda)$ and $\theta \in \mathbb{T}$, any two successive visits of the orbit $(\theta+m\omega)_{m\in \mathbb{Z}}$ to the critical set $I_n$ are spaced by at least $2K_nM_n$ iterations. This formalises the idea of slow-return to critical set $I_n$. 

 \begin{definition}[Synchronisation condition]\label{definition_block n_sync condn parameter}
    Let $\bar{G}$ be a non-empty open subset of $I_0 \setminus 3I_1(\omega)$, for each $\omega \in \mathcal{F}_0$. We define a parameter set $\mathcal{F}^s = \mathcal{F}^s(\bar{G}) \subset \mathbb{T}$ to be the set of all $\omega \in \mathbb{T}$ that satisfies the following:
    \begin{description}
    \item[$(\mathcal{G})$]\label{(F3)_1} For all $\theta \in \mathbb{T}$, 
\begin{align*}
    \left\{\theta+k\omega \right\}_{k=0}^{M_0^{56}} \cap \bar{G} \neq \emptyset.
\end{align*}
\end{description}
For each $n\geq 1$, let us define the parameter set $\mathcal{F}_n(M_n,\dots,M_0,f,\lambda,{\bar{G}})$  as   
\begin{align}
    \mathcal{F}_n(M_n,\dots,M_0,f,\lambda,{\bar{G}}) := \tilde{\mathcal{F}}_n(M_n,\dots,M_0,f,\lambda) \cap \left\{\omega \in \mathbb{T}: \omega \text{ satisfies } (\mathcal{G}) \text{ for the set } \bar{G}\right\},
\end{align}
and the parameter set $\mathcal{F}_n$ as
\begin{align}
    \mathcal{F}_n = \mathcal{F}_n(f,\lambda,{\bar{G}}) := \tilde{\mathcal{F}}_n(f,\lambda) \cap \left\{\omega \in \mathbb{T}: \omega \text{ satisfies } (\mathcal{G}) \text{ for the set } \bar{G}\right\}. 
\end{align}
As a convention, let $\tilde{\mathcal{F}}_0 := \mathcal{F}_0$.
\end{definition}

The condition $(\mathcal{G})$ in definition \ref{definition_block n_sync condn parameter} thus imposes a condition of partial equidistribution of every orbit, under the irrational circle rotation with respect to $\omega$.

\begin{definition}[Synchronisation property]\label{definition_block n_sync condn set_g1 g2}
    Let $\omega \in \tilde{\mathcal{F}}_1$ and $\bar{G}$ be a non-empty open subset of $I_0 \setminus 3I_1(\omega)$. Then, we define a set $\hat{G}$ to be a non-empty union of open intervals that satisfies the following conditions:
\begin{description}
\item[$(g1)$] $B_{\varepsilon}(\bar{G}) \subset \hat{G} \subset I_0 \setminus I_1(\omega)$, for every $\varepsilon < |I_1(\omega)|$.
    \item[$(g2)$] $\pi_2 \circ F^{2M_0+S_0} \left( (\hat{G}-M_0 \omega) \times \left( C_1 \cup C_2 \right) \right) \subset C_1$.
    \item[$(g3)$]  $\pi_2 \circ F^{-2M_0-S_0} \left( (\hat{G}+M_0 \omega) \times \left( E_1 \cup E_2 \right) \right) \subset E_1$.
\end{description}
Recall that, for each interval $I = (a,b)$, we use the notation $B_{\varepsilon}(I) := (a-\varepsilon, b+\varepsilon)$.
\end{definition}
Thus, the set $\hat{G}$ is characterised by the dynamics on the fibre which forces the orbits from two different contracting (expanding) regions to synchronise and converge into a single contracting (expanding) region, under forward (backward) iterations. It is worth mentioning that the sets $\bar{G}$ and $\hat{G}$ are closely related to the set $G$ as in definition \ref{definition_preliminaries F(1,1)_G and S sets}. The existence of non-empty sets $\bar{G}$ and $\hat{G}$ will be proved later in subsection \ref{subsection_characteristics G and S}.

%%%%%%%%%%%%%%%%%%%%%%%%%%%%%%%%%%%%%%%%%%%%%%%%%%%%%%%%%%%%%%%%%%
%------------------- MAIN PROPOSITIONS OF In ---------------------

\subsubsection{Main propositions corresponding to  $I_n$}

For each $n\geq 1$, and for the parameters in $\tilde{\mathcal{F}}_n$ and ${\mathcal{F}}_n$ we now establish the rate of contraction for orbits that return slowly to the critical sets $(I_j)_{j=1}^{n}$. 

\begin{definition}\label{definition_block n_beta n's}
For each integer $n\geq 1$, let us define the constant
\begin{align*}
    \beta_n & := \left( 1 - \prod_{i=0}^{n}\left(1-\frac{M_i+S_0}{K_iM_i} \right) \right),
\end{align*} 
 and as a convention let $\beta_{-1} := 0$.
\end{definition}
Note that for each $n\geq 1$, by the choice of constants $(K_i)_{i\leq n}$ and $(M_i)_{i\leq n}$, we have  
\begin{align*}
1-\beta_n & = \prod_{i=0}^{n}\left(1-\frac{M_i+S_0}{K_iM_i}\right) \\
& \geq  \prod_{i=0}^{n}\left(1-\frac{1}{2^{7+i}(b+1)}\right) > 1 - \frac{1}{4(b+1)}. 
\end{align*}
The sequence $(\beta_n)_{n\geq 1}$ is a bounded increasing sequence, where for each $n\geq 1$ the following recursive relation is satisfied:
\begin{align}\label{equation_block n_beta n recursive}
    (1-\beta_n) = \left(1-\frac{M_n+S_0}{K_nM_n}\right)(1-\beta_{n-1}).
\end{align}
Further, for each $n\geq 1$ we have $\displaystyle \beta_n < \frac{1}{4(b+1)}$, and thus
\begin{align}\label{equation_block n_beta n bound}
    \alpha_c^{1-(1+b)\beta_n} & < \alpha_c^{1/2}.
\end{align}

%------------------- PROP BEGIN -------------------

\begin{proposition}\label{proposition_block n contraction}
    Let $n\geq 0$, $\omega \in \tilde{\mathcal{F}}_{n}$, $\theta_0 \in \Theta_{n}$, $x_0 \in C_1 \cup C_2$. Then the following conditions are satisfied: 
    \begin{description}
        \item[$(C1)_n$] If $N\geq 1$ is the smallest positive integer such that $\theta_N \in I_{n+1}$, then for all $k \in \{0,\dots,N\}$, 
\begin{align}\label{equation_block n contraction_Cn_0}
     x_{k} \notin C_1 \cup C_2 \implies \theta_k \in \bigcup_{l=0}^{n}\bigcup_{m=1}^{M_l+S_0-1}(I_l+m\omega).
\end{align}
        \item[$(C2)_n$] If $N\geq 1$ is the smallest positive integer such that $\theta_N \in I_{n+1}$, then for all $k\in \{0,\dots,N\}$,  
\begin{align*}
    \left|f'(x_{N})\dots f'(x_k) \right| \leq \left( \alpha_c^{1-(1+b)\beta_n} \right)^{N-k+1}.
\end{align*}
        \item[$(C3)_n$] If $M \geq 1$ is a positive integer such that $\theta_M \in I_{n}\setminus I_{n+1}$, then for all $k\in \{0,\dots,M\}$,
\begin{align}\label{equation_block n contraction_C3 n_1}
x_{k} \notin C_1 \cup C_2 \implies \theta_k \in \bigcup_{l=0}^{n}\bigcup_{m=1}^{M_l+S_0-1}(I_l+m\omega)  \\
   \text{ and } \left|f'(x_{M})\dots f'(x_k) \right| \leq \left( \alpha_c^{1-(1+b)\beta_n} \right)^{M-k+1}.
     \stepcounter{equation}\tag{\theequation} \label{equation_block n contraction_C3 n_2}
\end{align} 
        \item[$(C4)_n$] If $N\geq 1$ is the smallest positive integer such that $\theta_N \in I_{n+1}(\omega)$, then 
        \begin{align}
            \left|\partial_x (\pi_2 \circ F^{N+1}(\theta_0, x_0))\right| \leq \left( \alpha_c^{1-(1+b)\beta_n} \right)^{N+1}.
        \end{align}
    \end{description}
\end{proposition}

\begin{proof}
The proof proceeds by an induction on $n\geq 1$. Recall that for every $n\geq 1$ we have $\mathcal{F}_{n} \subset \mathcal{F}_{n-1}$. For the case $n=0$, the proof follows directly from proposition \ref{proposition_block 0 contraction}, where proof of $(C3)_0$ also follows from the proof of $(C1-2)_0$ (in proposition \ref{proposition_block 0 contraction}). For an arbitrary $n\geq 1$, let us now assume the results $(C1-4)_{n-1}$ to be true, and prove $(C1-4)_n$. Let $\omega \in \mathcal{F}_n$, $\theta_0 \in \Theta_n$ and $x\in C_1 \cup C_2$ be fixed. Let $0<N_1 < N_2< \cdots < N_J = N$ be the sequence of all integers for which $\theta_{N_j}\in I_n$, and $N$ is the smallest positive integer for which $\theta_{N} \in I_{n+1}$.

\vspace{0.2cm}
\noindent\textit{Proof of $(C1)_n$:} By our assumption, for every $\omega \in \tilde{\mathcal{F}}_{n-1}$, $\theta_0 \in \Theta_{n-1}$, $x_0 \in C_1 \cup C_2$ and smallest integer $N'$ such that $\theta_{N'} \in I_n$, 
    \begin{description}
        \item[$(C1)_{n-1}$] For all $k=0,\dots , N'$, we have $x_k \notin C_1 \cup C_2 \implies \theta_k \in \bigcup_{l=0}^{n-1} \bigcup_{m=1}^{M_l+S_0-1}(I_l+m\omega)$.
    \end{description}
For each $1 \leq j \leq J$, let us define the following new condition for the orbit stretch $\{(\theta_m, x_m)\}_{m=0}^{N_{j}}$:
\begin{description}
    \item[$(C1)_n^{N_{j}}$] For any $k = 0,\dots ,N_{j}$,
    \begin{align}\label{equation_block 0 contraction_C1_n[N_j-1]}
        x_k \notin C_1 \cup C_2 \implies \theta_k \in \bigcup_{l=0}^{n} \bigcup_{m=1}^{M_l+S_0-1}(I_l+m\omega).
    \end{align}
\end{description}
Since $\omega \in \tilde{\mathcal{F}}_{n} \subset \tilde{\mathcal{F}}_{n-1} \subset \mathcal{F}_0$, $\theta_0 \in \Theta_n \subset \Theta_{n-1}$ and $N_1$ is the smallest integer for which $\theta_{N_1} \in I_n$, by the assumption $(C1)_{n-1}$ it follows that for all $k\in \{0,\dots ,N_1\}$
\begin{align*}
    x_k \notin C_1 \cup C_2 & \implies \theta_k \in \bigcup_{l=0}^{n-1} \bigcup_{m=1}^{M_l+S_0-1}(I_l+m\omega).
\end{align*}
This proves the condition $(C1)_n^{N_1}$ and also settles the case when $J=1$. 

Let us consider the case when $J>1$ and assume the condition $(C1)_n^{N_j}$ to be true for some $1\leq j \leq (J-1)$. Since $(\theta + (N_j-M_n)\omega) \in (I_n - M_n\omega)$, by the disjointness conditions $(\mathcal{F}2)_n$ (in definition \ref{definition_block n_disjointness condn}) and $(\mathcal{F}1)_0$ (in definition \ref{definition_block 0_disjointness condn}), it follows that
\begin{align}
(\theta + (N_j-M_n)\omega) \notin  \bigcup_{l=0}^{n-1} \bigcup_{m=1}^{M_l+S_0-1}(I_l+m\omega), 
\end{align}
and by $(C1)_n^{N_j}$, this gives us $x_{N_j-M_n} \in C_1\cup C_2$. Since $j\leq J-1$ and thus $\theta_{N_j} \in I_n\setminus I_{n+1}$, by the definition of the critical sets (in definition \ref{definition_block n_Critical sets In}), we can conclude that $x_{N_j+M_n} \notin E_1\cup E_2$. Then, the disjointness condition $(\mathcal{F}2)_n$ and equations (\ref{equation_preliminaries_f contraction S0 1} - \ref{equation_preliminaries_f contraction S0 2}) yield us $x_{N_j+M_n+S_0} \in C_1\cup C_2$. For all $k=N_j,\dots ,N_j+M_n+S_0-1$, we have $\displaystyle \theta_k \in \bigcup_{m=1}^{M_n+S_0-1}(I_n+m\omega)$ and thus 
\begin{align*}
    \theta_k \in \bigcup_{l=0}^{n} \bigcup_{m=1}^{M_l+S_0-1}(I_l+m\omega).
\end{align*}
By disjointness conditions $(\mathcal{F}1-2)_0$ and definition of $\Theta_n$, it follows that $\theta_{N_j+M_n+S_0} \in \Theta_n$. For the forward iterates of $(\theta_{N_j+M_n+S_0}, x_{N_j+M_n+S_0})$, it is clear that $\theta_{N_{j+1}} \in I_n$ and $\theta_{k} \notin I_n$, for all $N_j+M_n+S_0 \leq k < N_{j+1}$. Then, it follows from $(C1)_{n-1}$ that, for all $k = N_j+M_n+S_0, \dots ,N_{j+1}$,
\begin{align}\label{equation_block n_C1 n_1}
    x_{k} \notin C_1 \cup C_2 & \implies \theta_k \in \bigcup_{l=0}^{n-1} \bigcup_{m=1}^{M_l+S_0-1}(I_l+m\omega) \subset \bigcup_{l=0}^{n} \bigcup_{m=1}^{M_l+S_0-1}(I_l+m\omega).
\end{align}
Thus, by $(C1)_n^{N_j}$ and equation \eqref{equation_block n_C1 n_1}, we have proved that  
\begin{align*}
    x_k \notin C_1 \cup C_2 & \implies \theta_k \in \bigcup_{l=0}^{n} \bigcup_{m=1}^{M_l+S_0-1}(I_l+m\omega),
\end{align*}
for all $k = 0,\dots , N_{j+1}$, which is statement of $(C1)_n^{N_{j+1}}$. Note that the proof of $(C1)_n^{N_j}$, as defined in equation \eqref{equation_block 0 contraction_C1_n[N_j-1]}, is in fact independent of the existence of integer $N$ such that $\theta_N \in I_{n+1}$. Hence, this concludes the inductive proof of $(C1)_n$ as well as equation \eqref{equation_block n contraction_C3 n_1} in $(C3)_n$.

%%%%%%%%%%%%%%%%%%%%%%%%%%%%%%%%%%%%%%%%%%%%%%%%%%%%%%%%%

%%%%%%%%%%%%%%%%%%%%%%%%%%%%%%%%%%%%%%%%%%%%%%%%%%%%%%%%%

\vspace{0.2cm}
\noindent\textit{Proof of $(C2)_n$}: For every $\omega \in \tilde{\mathcal{F}}_{n-1}$, $\theta_0 \in \Theta_{n-1}$, $x_0\in C_1 \cup C_2$ and the smallest integer $N'$ such that $\theta_{N'} \in I_n$, we assume that  
    \begin{description}
        \item[$(C2)_{n-1}$] $ |f'(x_{N'})\dots f'(x_k)| \leq ( \alpha_c^{1-\beta_{n-1}} \alpha_u^{\beta_{n-1}})^{N'-k+1}$, where $\displaystyle \beta_{n-1} = 1 - \prod_{i=0}^{n-1}\left(1-\frac{M_i+S_0}{2K_iM_i}\right)$. 
    \end{description}
For each $1 \leq j \leq J$, let us define the following condition for the orbit stretch $\{(\theta_m, x_m)\}_{m=0}^{N_{j}}$:
\begin{description}
    \item[$(C2)_n^{N_{j}}$] For any $0 \leq k \leq N_{j}$,
    \begin{align}\label{equation_block 0 contraction_C2_n[N_j-1]}
        |f'(x_{N_{j}})\dots f'(x_k)| \leq \left( \alpha_c^{1-(1+b)\beta_{n}} \right)^{N_{j}-k+1}.
    \end{align}
\end{description}
Since $\omega \in \tilde{\mathcal{F}}_{n} \subset \tilde{\mathcal{F}}_{n-1} \subset \mathcal{F}_0$, $\theta_0 \in \Theta_n \subset \Theta_{n-1}$ and $N_1$ is the smallest integer for which $\theta_{N_1} \in I_n$, by the inductive assumption $(C2)_{n-1}$ and the inequality $\beta_{n-1} < \beta_{n}$, we deduce that 
\begin{align*}
    |f'(x_{N_1})\dots f'(x_k)| \leq \left( \alpha_c^{1-(1+b)\beta_{n-1}} \right)^{N_1-k+1} < \left( \alpha_c^{1-(1+b)\beta_{n}} \right)^{N_1-k+1}, \text{ for all } k=0,\dots ,N_1. 
\end{align*}
This proves the condition $(C2)_n^{N_1}$ and also proves $(C2)_n$ for the case  $J=1$. 

Let us consider the cases when $J>1$ and assume $(C2)_n^{N_{j}}$ to be true for some $1\leq j \leq J-1$. Since $\theta_{N_j-M_n} \in (I_n-M_n\omega)$, by the disjointness condition $(\mathcal{F}1-2)_{n}$ (in definition \ref{definition_block n_disjointness condn}) and $(\mathcal{F}1)_0$ (in definition \ref{definition_block 0_disjointness condn}), it follows that
\begin{align*}
    \theta_{N_j-M_n} \notin \bigcup_{m=0}^{n}\bigcup_{m=0}^{M_l+S_0-1}(I_l+m\omega).
\end{align*}
and, by $(C1)_n$ this gives us $x_{N_j-M_n} \in C_1\cup C_2$. Since $j \leq (J-1)$ and thus $\theta_{N_j} \in I_n\setminus I_{n+1}$, by the definition of the critical sets (in definition \ref{definition_block n_Critical sets In}) $x_{N_j+M_n} \notin E_1 \cup E_2$. Then, the disjointness condition $(\mathcal{F}2)_n$ and equations (\ref{equation_preliminaries_f contraction S0 1} - \ref{equation_preliminaries_f contraction S0 2}) yield us $x_{N_j+M_n+S_0} \in C_1 \cup C_2$. Note that $\theta_{N_j+M_n+S_0} \in \Theta_n$, $\theta_{N_{j+1}} \in I_n$ and $\theta_{k} \notin I_n$, for all $N_j+M_n+S_0 \leq k < N_{j+1}$. Then, by $(C2)_{n-1}$ we obtain the estimates 
\begin{align}\label{equation_block n_C2n_2}
    |f'(x_{N_{j+1}})\dots f'(x_{k})| < \left( \alpha_c^{1-(b+1)\beta_{n-1}} \right)^{N_{j+1}-k+1},
\end{align}
for all $N_j+M_n+S_0 \leq k \leq N_{j+1}$. By the disjointness condition $(\mathcal{F}1)_n$, we can deduce that $N_{j+1}-N_{j}+1 > 2K_nM_n$ and thus, $N_{j+1}-(N_j+M_n+S_0) > 2K_nM_n-(M_0+S_0)$. Then, by equation \eqref{equation_block n_C2n_2}, we obtain
\begin{multline*}
    |f'(x_{N_{j+1}})\dots f'(x_{N_j+M_n+S_0}) f'(x_{N_j+M_n+S_0-1})\dots f'(x_k)|\\ < \left( \alpha_c^{1-(b+1)\beta_{n-1}} \right)^{N_{j+1}-N_j-M_n-S_0+1} 
    \alpha_u^{N_j+M_n+S_0-k}, 
\end{multline*}
for every $N_{j} \leq k \leq N_{j}+M_0+S_0-1$. Since $\alpha_u < \alpha_c^{-b}$, the above equation can be reduced to 
\begin{align}\label{equation_block n_C2n_3}
    |f'(x_{N_{j+1}})\dots f'(x_k)| < \left( \alpha_c^{P} \right)^{N_{j+1}-k+1}, \text{ for all } N_{j} \leq k \leq N_{j}+M_0+S_0-1,
\end{align}
where $P$ denotes the expression
\begin{multline}\label{equation_block n_C2n_power}
    P := \frac{N_{j+1}-N_j-M_n-S_0+1}{N_{j+1}-k+1} - \frac{\beta_{n-1}(b+1)(N_{j+1}-N_j-M_n-S_0+1)}{N_{j+1}-k+1} \\ 
    - \frac{b(N_j+M_n+S_0-k)}{N_{j+1}-k+1}.
\end{multline}
The disjointness condition $(\mathcal{F}1)_n$, gives us the bounds $N_j+M_n+S_0-k \leq M_n+S_0$,  $N_{j+1}-N_j \geq 2K_nM_n$ and $N_{j+1}-k+1 \geq K_nM_n$, for every $N_j\leq k \leq N_{j+1}$. The expression then satisfies the bound 
\begin{align*}
    P & > 1 - (b+1)\left(\frac{\beta_{n-1}(N_{j+1}-N_j-M_n-S_0+1)}{N_{j+1}-k+1} + \frac{b(N_j+M_n+S_0-k)}{(b+1)(N_{j+1}-k+1)}\right) \\
    & >  1 - (b+1)\left( 1 - (1 - \beta_{n-1}) \left(1 - \frac{b}{(b+1)(1-\beta_{n-1})} \frac{M_n+S_0}{K_nM_n}\right)\right). 
\end{align*}
Since $\displaystyle \beta_{n-1} < \frac{1}{4(b+1)}$, we have $\displaystyle \frac{b}{(b+1)(1-\beta_{n-1})} < 1$. By equation \eqref{equation_block n_beta n recursive}, we have the recursive relation
$$\beta_n = 1 - (1-\beta_{n-1})\left(1-\frac{M_n+S_0}{K_nM_n}\right).$$
This yields us $P > 1-(b+1)\beta_n$. Thus, the equation \eqref{equation_block n_C2n_3} reduces to
\begin{align}\label{equation_block n_C2n_4}
    |f'(x_{N_{j+1}})\dots f'(x_k)| < \left( \alpha_c^{1-(b+1)\beta_n} \right)^{N_{j+1}-k+1}, \text{ for all } N_{j} \leq k \leq N_{j}+M_0+S_0-1.
\end{align}
The estimates in equations \eqref{equation_block n_C2n_2} and \eqref{equation_block n_C2n_4} thus proves the condition $(C2)_n^{N_{j+1}}$, for all $j=1,\dots ,J-1$. Note that the proof is in fact independent of the existence of integer $N$ such that $\theta_N \in I_{n+1}$. Hence, this concludes the inductive proof of $(C2)_n$ as well as equation \eqref{equation_block n contraction_C3 n_2} in $(C3)_n$. The proof of $(C4)_n$ follows from the estimate in $(C2)_n$ for the case $k=N$.

\end{proof}

%------------------- PROP END -------------------

%------------------- PROP BEGIN -------------------

\begin{proposition}\label{proposition_block n sync contraction}
    Let $n\geq 0$, $\omega \in \tilde{\mathcal{F}}_{n}$, $\theta_0 \in \Theta_{n}$, $x_0, y_0 \in C_i$ for some $i\in \{1,2\}$ and $N\geq 1$ be the smallest positive integer such that $\theta_N \in I_{n+1}(\omega)$. Then the following condition is satisfied:
\begin{description}
            \item[$(Sc1)_n$] For all $k\in \{0,1,\dots,N\}$, 
        \begin{align*}
            \theta_k \notin \bigcup_{l=0}^{n}\bigcup_{m=1}^{M_l+S_0-1}(I_l+m\omega) \Longrightarrow  x_k,y_k \in C_j, \text{ for some } j\in \{1,2\}.
        \end{align*}
\end{description}
\end{proposition}

\begin{proof}
The proof proceeds by an induction on $n\geq 0$. For $n=0$, the proof of $(Sc1)_0$ follows from  proposition \ref{proposition_Block 0 sync contraction}. 
For an arbitrary $n\geq 1$, let us now assume the result to be true for $(n-1)$ case. That is, for every $\omega \in \tilde{\mathcal{F}}_{n-1}$, $\theta_0 \in \Theta_{n-1}$, $x_0,y_0 \in C_i$ for some $i\in\{1,2\}$, we assume that:
    \begin{description}
        \item[$(Sc1)_{n-1}$] For all $k=0,\dots , N'$, where $N'>0$ is the smallest integer such that $\theta_{N'} \in I_n$,
        $$\theta_k \notin \bigcup_{l=0}^{n-1} \bigcup_{m=1}^{M_l+S_0-1}(I_l+m\omega) \implies x_k,y_k \notin C_l, \text{ for some } l\in \{1,2\}.$$
    \end{description}
Let $\omega \in \mathcal{F}_n$, $\theta_0 \in \Theta_n$ and $x\in C_1 \cup C_2$ be fixed. Let $0<N_1 < N_2< \cdots < N_J = N$ be the sequence of all integers for which $\theta_{N_j}\in I_n$, where $N$ is the smallest positive integer for which $\theta_{N} \in I_{n+1}$. For each $1 \leq j \leq J$, let us define the following new condition for the orbit stretch $\{(\theta_m, x_m)\}_{m=0}^{N_{j}}$:
\begin{description}
    \item[$(Sc1)_n^{N_{j}}$] For any $k = 0,\dots ,N_{j}$,
    \begin{align}\label{equation_block 0 contraction_Sc1_n[N_j-1]}
        \theta_k \notin \bigcup_{l=0}^{n-1} \bigcup_{m=1}^{M_l+S_0-1}(I_l+m\omega) \implies x_k,y_k \notin C_l, \text{ for some } l\in \{1,2\}.
    \end{align}
\end{description}
Since $\omega \in \tilde{\mathcal{F}}_{n} \subset \tilde{\mathcal{F}}_{n-1} \subset \mathcal{F}_0$, $\theta_0 \in \Theta_n \subset \Theta_{n-1}$ and $N_1$ is the smallest integer for which $\theta_{N_1} \in I_n$, by the inductive hypothesis $(C1)_{n-1}$, we have
\begin{align}\label{equation_block n_Sc1_n_1}
    \theta_k \notin \bigcup_{l=0}^{n-1} \bigcup_{m=1}^{M_l+S_0-1}(I_l+m\omega), \implies x_k,y_k \notin C_l, \text{ for some } l\in \{1,2\},
\end{align}
for all $k=0,\dots ,N_1$. This proves the condition $(Sc1)_n^{N_1}$ and also settles the case when $J=1$. 

Let us consider the case when $J>1$ and assume that the condition $(Sc1)_n^{N_j}$ is satisfied for some $1\leq j \leq J-1$. Note that we have $\theta_{N_j-M_n} \in (I_n-M_n\omega)$ and thus, by disjointness condition $(\mathcal{F}1-2)_n$ (in definition \ref{definition_block n_disjointness condn}) and by $(Sc1)_n^{N_j}$ (in equation \eqref{equation_block 0 contraction_Sc1_n[N_j-1]}), it follows that $\displaystyle x_{N_j-M_n},y_{N_j-M_n} \in C_i$, for some $i\in \{1,2\}$. Without loss of generality, let $x_{N_j-M_n}< y_{N_j-M_n}$ and thus
\begin{align}\label{equation_block n contraction_Sc1_2}
    \{\theta_{N_j-M_n}\} \times [x_{N_j-M_n}, y_{N_j-M_n}] \subset \{\theta_{N_j-M_n}\} \times C_i.
\end{align}
Since $\theta_{N_j} \in I_n\setminus I_{n+1}$, by disjointness condition $(\mathcal{F}1-2)_n$ and by $(C1)_n$ (in proposition \ref{proposition_block n contraction}), we can deduce that $x_{N_j+M_n+S_0},y_{N_j+M_n+S_0} \in C_1\cup C_2$.

Say there exists $k \in [N_j+M_n+S_0, N_{j+1}]$ such that $x_k \in C_{i_1}$ and $y_k \in C_{i_2}$ for two distinct integers $i_1,i_2\in \{1,2\}$. By the disjointness condition $(\mathcal{F}1-2)_n$, we have $\theta_{N_j+M_n+S_0} \in \Theta_{n}$ and  $\theta_{N_{k}} \notin I_n$ for all $k\in [N_j+M_n+S_0, N_{j+1}-1]$. Then, by $(Sc1)_{n-1}$, the points $x_{N_j+M_n+S_0}$ and $y_{N_j+M_n+S_0}$ must also belong to two distinct contracting regions and by the distribution of the regions $\{C_i,E_i\}_{i=1}^{2}$, this gives us 
\begin{align}\label{equation_block n contraction_Sc1_3}
    \left(\{\theta_{N_j+M_n+S_0}\} \times [x_{N_j+M_n+S_0}, y_{N_j+M_n+S_0}] \right) \cap \left(\{\theta_{N_j+M_n+S_0}\} \times (E_1\cup E_2)\right) & \neq \emptyset. 
\end{align}
By $(\mathcal{F}2)_n$, it follows that $\displaystyle \{\theta_{N_j+M_n+k}\}_{k=0}^{S_0} \cap I_0 = \emptyset$ and as a result,
$$F^{-S_0}\left(\{\theta_{N_j+M_n+S_0}\} \times (E_1\cup E_2)\right) \subset \left(\{\theta_{N_j+M_n}\} \times (E_1\cup E_2)\right).$$ Further, since $F$ is continuous and orientation preserving
\begin{align}
    F^{2M_n+S_0}\left(\theta_{N_j-M_n}, [x_{N_j-M_n}, y_{N_j-M_n}]\right) & = \{\theta_{N_j+M_n+S_0}\} \times [x_{N_j+M_n+S_0}, y_{N_j+M_n+S_0}],
\end{align}
it follows from equations \eqref{equation_block n contraction_Sc1_2} and \eqref{equation_block n contraction_Sc1_3} that
\begin{align*}
F^{M_n}\left(\{\theta_{N_j-M_n}\} \times [x_{N_j-M_n}, y_{N_j-M_n}]\right) \cap F^{-M_n}\left(\{\theta_{N_j+M_n}\} \times (E_1\cup E_2)\right) & \neq \emptyset \\
\implies F^{M_n}\left(\{\theta_{N_j-M_n}\} \times (C_1\cup C_2)\right) \cap F^{-M_n}\left(\{\theta_{N_j+M_n}\} \times (E_1\cup E_2)\right) & \neq \emptyset.
\end{align*}
 However, this leads to a contradiction as $\theta_{N_j} \notin I_{n+1}(\omega)$. Thus, for the fixed $x_0,y_0 \in C_i$ and each $k\in[N_j+M_n+S_0, N_{j+1}]$, there exists some $l\in \{1,2\}$ such that $x_{k}, y_{k} \in C_l$. This concludes our proof of $(Sc1)_1^{N_{j+1}}$ and consequently, the proof of induction of the lemma.

\end{proof}

%------------------- PROP END -------------------

%------------------- PROP BEGIN -------------------

\begin{proposition}\label{proposition_block 0 sync' contraction}
    For $\omega \in \mathcal{F}_1(M_1,M_0,f,\lambda,{\bar{G}})$, where $\bar{G}$ is a non-empty set in $I_0 \setminus 3I_1(\omega)$, assume there exists a set $\hat{G}$ as in definition \ref{definition_block n_sync condn set_g1 g2}. Let $\theta_0 \in \Theta_{1}$, $x_0 \in C_1$, $y_0 \in C_2$ and $N$ be the smallest integer such that $\theta_N \in I_{1}$. Then, the following condition is satisfied:
        \begin{description} 
         \item[$(S1)_0$] If there exists a constant $N_s < M_1$ such that $\theta_{N_s-M_0-S_0} \in \hat{G}$, then we have $\theta_{N_s} \in \Theta_0$ and $x_{N_s}, y_{N_s} \in C_1$. Further, for any $N_s \leq  k \leq N$, 
            \begin{align}\label{equation_block 0 sync' contraction}
                \theta_k \notin \bigcup_{m=1}^{M_0+S_0-1}(I_0+m\omega) \Longrightarrow x_k, y_k \in C_i, \text{ for some } i\in \{1,2\}.
            \end{align}
        \end{description}
\end{proposition}

\begin{proof}

By definition \ref{definition_block n_sync condn set_g1 g2}, since $\hat{G} \subset I_0\setminus I_1$, it is clear that $\theta_{N_s-M_0-S_0} \in I_0 \setminus I_1$ and $\theta_{N_s-2M_0-S_0} \in (I_0-M_0\omega)$. Then, from the definition of $\Theta_0$ and $(\mathcal{F}1)_0$ (in definition \ref{definition_block 0_disjointness condn}), we deduce that $\theta_{N_s} \in \Theta_0$. By $(C1)_0$ (in proposition \ref{proposition_block 0 contraction}), the orbits along the fibre satisfy $x_{N_s-2M_0-S_0}, y_{N_s-2M_0-S_0} \in (C_1 \cup C_2)$.
By property $(g2)$ (in definition \ref{definition_block n_sync condn set_g1 g2}), which applies here to the forward iterate of the point 
$$(\theta_{N_s-2M_0-S_0}, x_{N_s-2M_0-S_0}) \in ((\hat{G}-M_0\omega)\times (C_1\cup C_2)),$$ this proves that $x_{N_s}, y_{N_s} \in C_1$. Note that, as a consequence of $N_s<M_1$ and $\theta_0 \in \Theta_0$ we have $N_s < N$. Here, the implication in equation \eqref{equation_block 0 sync' contraction}  follows from $(Sc1)_0$ (in proposition \ref{proposition_Block 0 sync contraction}).
    
\end{proof}

%------------------- PROP END -------------------

%------------------- PROP BEGIN -------------------

\begin{proposition}\label{proposition_block n sync' contraction}
For each $n\geq 1$ and $\omega \in \mathcal{F}_n(M_n,\dots,M_0,f,\lambda,{\bar{G}})$, where $\bar{G}$ is a non-empty set in $I_0 \setminus 3I_1(\omega)$, assume there exists a set $\hat{G}$ as in definition \ref{definition_block n_sync condn set_g1 g2}. Let $\theta_0 \in \Tilde{\Theta}_{n}$, $x_0 \in C_1$, $y_0 \in C_2$ and $N$ be the smallest integer such that $\theta_N \in I_{n+1}$. Then, the following condition is satisfied:
        \begin{description}          
         \item[$(S1)_n$] If there exists $N_s < M_1$ such that $\theta_{N_s-M_0-S_0} \in \hat{G}$, then we have $\theta_{N_s}\in \Theta_n$ and $x_{N_s}, y_{N_s} \in C_1$. Further, for any $N_s \leq  k \leq N$, 
            \begin{align}\label{equation_block n sync' contraction_S1_n}
                \theta_k \notin \bigcup_{l=0}^{n}\bigcup_{m=1}^{M_l+S_0-1}(I_l+m\omega) \Longrightarrow x_k, y_k \in C_i, \text{ for some } i\in \{1,2\}.
            \end{align}
        
        \end{description}
\end{proposition}

\begin{proof}
By definition \ref{definition_block n_sync condn set_g1 g2}, since $\hat{G} \subset I_0\setminus I_1$, it is clear that $\theta_{N_s-M_0-S_0} \in I_0 \setminus I_1$ and $\theta_{N_s-2M_0-S_0} \in (I_0-M_0\omega)$.
From the definition of $\tilde{\Theta}_n$, we can conclude that $\theta_k \notin I_1$, for all $k\in \{0,\dots ,M_1\}$. Then, by $(C1)_0$ (in proposition \ref{proposition_block 0 contraction}), the orbits along the fibre satisfy
$x_{N_s-2M_0-S_0}, y_{N_s-2M_0-S_0} \in (C_1\cup C_2).$ Then, by property $(g2)$ (in definition \ref{definition_block n_sync condn set_g1 g2}), it follows that $x_{N_s}, y_{N_s} \in C_1$. Since $N_s < M_1$ and $\theta_{N_s} \in (I_0+(M_0+S_0)\omega)$, by $(\mathcal{F}1)_0$ and definition of set $\tilde{\Theta}_n$, we deduce that $\theta_{N_s} \in \Theta_n$. Now the proof of equation \eqref{equation_block n sync' contraction_S1_n} follows from $(Sc1)_n$ (in proposition \ref{proposition_block n sync contraction}).
\end{proof}

%------------------- PROP END -------------------

Analogous to the propositions discussed above, estimates for the contraction rates of certain forward iterates, we have the following propositions for the contraction rates of certain backward iterates. The proof of these results also follow similarly.

%------------------- PROPs Expanding region BEGIN -------------------

\begin{proposition}\label{proposition_block n expansion}
Let $n\geq 1$, $\omega \in \tilde{\mathcal{F}}_{n}$, $\theta_0 \in \Theta_{n}$, $x_0 \in E_1 \cup E_2$. Then the following conditions are satisfied: 
    \begin{description}
        \item[$(E1)_n$] If $N\geq 1$ be the smallest such that $\theta_{-N} \in I_{n+1}(\omega)$, then for all $k \in \{0,\dots,N-1\}$, 
        \begin{align}\label{equation_block n expansion_E1_0} 
        x_{-k} \notin E_1 \cup E_2 \Longrightarrow \theta_{-k} \in \bigcup_{l=0}^{n}\bigcup_{m=0}^{M_l+S_0-1}(I_l-m\omega).
        \end{align}
        
        \item[$(E2)_n$] If $N\geq 1$ be the smallest such that $\theta_{-N} \in I_{n+1}(\omega)$, then for all $k \in \{0,\dots,N-1\}$, 
        \begin{align*}
            \left|(f^{-1})'(f(x_{-N+1}))\dots (f^{-1})'(f(x_{-k})) \right| \leq \left( \alpha_c^{1-(1+b)\beta_n} \right)^{N-k}.
        \end{align*}

 \item[$(E3)_n$] If $M$ is a positive integer such that $\theta_{-M} \in I_{n}(\omega)\setminus I_{n+1}(\omega)$, then for all $k\in \{0,\dots,M-1\}$,
\begin{align}\label{equation_block n expansion_E3 n_1}
x_{-k} \notin E_1 \cup E_2 \implies \theta_{-k} \in \bigcup_{l=0}^{n}\bigcup_{m=0}^{M_l+S_0}(I_l-m\omega) \text{ and } \\
    \left|f'(x_{-M+1})\dots f'(x_{-k}) \right| \leq \left( \alpha_c^{1-(1+b)\beta_n} \right)^{M-k}
     \stepcounter{equation}\tag{\theequation} \label{equation_block n expansion_E3 n_2}
\end{align} 
        
        \item[$(E4)_n$] $\left|\partial_x (\pi_2 \circ F^{-N}(\theta_0, x_0))\right| \leq \left( \alpha_c^{1-(1+b)\beta_n} \right)^{N}$. 
    \end{description}
\end{proposition}

\begin{proposition}\label{proposition_block n sync expansion}
    Let $n\geq 1$, $\omega \in \tilde{\mathcal{F}}_{n}$, $\theta_0 \in \Theta_{n}$, $x_0, y_0 \in E_i$ for some $i\in \{1,2\}$, and let $N\geq 1$ be the smallest such that $\theta_{-N} \in I_{n+1}(\omega)$. Then, the following condition is satisfied: 
\begin{description}
            \item[${(Se1)}_n$] For any $k\in \{0,1,\dots,N-1\}$, 
        \begin{align*}
            \theta_k \notin \bigcup_{l=0}^{n}\bigcup_{m=0}^{M_l+S_0-1}(I_l-m\omega) \Longrightarrow x_k,y_k \in E_i, \text{ for some } i\in \{1,2\}.
        \end{align*}
\end{description}
\end{proposition}

\begin{proposition}\label{proposition_block 0 sync' expansion}
    For $\omega \in \mathcal{F}_1(M_1,M_0,f,\lambda,{\bar{G}})$, where $\bar{G}$ is some non-empty set in $I_0 \setminus 3I_1(\omega)$, assume there exists a set $\hat{G}$ as in definition \ref{definition_block n_sync condn set_g1 g2}. Let $\theta_0 \in \Theta_{1}$, $x_0 \in E_1$, $y_0 \in E_2$ and $N$ be the smallest integer such that $\theta_{-N} \in I_{1}$. Then, the following condition is satisfied:
        \begin{description} 
         \item[$(\bar{S}1)_0$] If there exists $N_s<M_1$ such that $\theta_{-N_s+M_0+S_0} \in \hat{G}$, then we have $\theta_{-N_s} \in \Theta_0$ and $x_{-N_s}, y_{-N_s} \in E_1$.
            Further, for any $N_s \leq  k \leq N$, 
            \begin{align*}
                \theta_{-k} \notin \bigcup_{m=0}^{M_0+S_0-1}(I_0-m\omega) \Longrightarrow x_{-k}, y_{-k} \in E_i, \text{ for some } i\in \{1,2\}.
            \end{align*}
        \end{description}
\end{proposition}

\begin{proposition}\label{proposition_block n sync' expansion}
    For each $n\geq 2$ and $\omega \in \mathcal{F}_n(M_n,\dots,M_0,f,\lambda,{\bar{G}})$, where $\bar{G}$ is some non-empty set in $I_0 \setminus 3I_1(\omega)$, assume there exists a set $\hat{G}$ as in definition \ref{definition_block n_sync condn set_g1 g2}. Let $\theta_0 \in \tilde{\Theta}_{n}$,  $x_0 \in E_1$, $y_0 \in E_2$ and $N$ be the smallest integer such that $\theta_{-N} \in I_{n+1}$. Then, the following condition is satisfied:
        \begin{description} 
         \item[${(\bar{S}1)}_n$] If there exists $N_s<M_1$ such that $\theta_{-N_s+M_0+S_0} \in \hat{G}$, then we have $\theta_{-N_s} \in \Theta_n$ and $x_{-N_s}, y_{-N_s} \in E_1$. Further, for any $N_s \leq  k \leq N$ we have 
            \begin{align*}
                \theta_{-k} \notin \bigcup_{l=0}^{n}\bigcup_{m=0}^{M_l+S_0-1}(I_l-m\omega) \Longrightarrow x_{-k}, y_{-k} \in E_i, \text{ for some } i\in \{1,2\}.
            \end{align*}
        \end{description}
\end{proposition}

%------------------- PROPs Expanding region END -------------------

%------------------------------------------------
%%%%%%%%%%%%%%%%%%%%%%%%%%%%%%%%%%%%%%%%%%%%%%%%%

%%%%%%%%%%%%%%%%%%%%% STRUCTURE of critical sets %%%%%%%%%%%%%%%%%%%%%%%%%%%%%%%%%%%%%%
\section{Structure of critical sets}\label{section_structure of critical sets}
In this section, we address the second of the three main components of the proof of Theorem~\ref{Theorem SNA}. For each admissible parameter $\omega$, we introduce and describe a particular sequence of sets $(I_n(\omega))_{n\geq 1}$, where each set satisfies the conditions of a critical set in definition \ref{definition_block n_Critical sets In}. We then determine their geometric structure and analyse their dependence on the parameter $\omega$. 

More specifically, subsection \ref{subsection_structure of critical set I1} contains estimates describing the non-empty critical set $I_1(\omega)$ associated with each $\omega \in \mathcal{F}_0$. Recall that the conditions of synchronisation discussed in subsection \ref{subsection_block n} rely on the existence of well-defined synchronisation sets $\hat{G}$ and $\bar{G}$ (in definition \ref{definition_block n_sync condn set_g1 g2} and \ref{definition_block n_sync condn parameter}). In relation to this, we identify suitable choices for the sets $\hat{G}$ and $\bar{G}$, and describe their structure in subsection \ref{subsection_characteristics G and S}. In subsection \ref{subsection_structure of critical set I2}, we focus on the critical set $I_2(\omega)$ and in particular, in subsubsection \ref{subsubsection_empty vs nonempty I2}, we identify the properties that determine whether the set $I_2(\omega)$ is empty or non-empty. When $I_2$ is non-empty, we establish its structure in subsubsection \ref{subsubsection_estimate nonempty I2}. 
For each $n\geq 3$, analogous results for the critical set $I_n(\omega)$ are established in subsection \ref{subsection_structure of critical set I(n+1)}. 

We also note that the estimates in this section are derived under finitely many lower bound conditions on $\lambda$. That is, the results in this section hold true for all $\lambda > \lambda'$, for some $\lambda'>1$. We now introduce certain subsets of $I_0$, which will be used to describe critical sets $I_n$ for $n\geq 1$, in subsections \ref{subsection_structure of critical set I1}, \ref{subsection_structure of critical set I2} and \ref{subsection_structure of critical set I(n+1)}.

%--------------------- Critical Set I_0 ---------------------------

\begin{definition}[Iterating Blocks] \label{definition_structure I0_iterating blocks}
    For each $\omega \in \mathbb{T}$ and $i \in \{1,2\}$, let us define the following sets in $\mathbb{T}^2$:
\begin{align*}
    A_{-1}^{(i)} & := \left\{(\theta+\omega,x): \theta \in I_0,\ x\in (c_i^{-}+g(\theta), c_i^{+}+g(\theta)) \right\} \\
    B_{-1}^{(i)} & := (I_0+\omega) \times E_i.
\end{align*}
\end{definition}

%------------------- LEMMA BEGIN -------------------

\begin{lemma}\label{lemma_structure I0_match Ai Bi}
    For all sufficiently large $\lambda$ and any $i,j \in \{1,2\}$, we have $A_{-1}^{(i)} \cap B_{-1}^{(j)} \neq \phi$ if and only if  $i=j$.
\end{lemma}

\begin{proof}
By the definition of $g(\theta)$, for a sufficiently large $\lambda > 1$,
$$|g(1/2 -\lambda^{-1/4})| = |g(1/2 +\lambda^{-1/4})| < \lambda^{-1/2} < |c_i^+ - c_i^-|/4 = \varepsilon_{c_i}/2.$$ 
For any $i\in \{1,2\}$, by the conditions on $a>0$ in the hypothesis of Theorem \ref{Theorem SNA}, we can conclude that $|c_i-e_i| < |a| = |g(1/2)|$. Thus, there exists $\theta \in I_0$ such that 
\begin{align*}
    [c_i^{-}+g(\theta), c_i^{+}+g(\theta)] \cap [e_i^{-}, e_i^{+}] \neq \phi.
\end{align*}
This gives us $A_{-1}^{(i)} \cap B_{-1}^{(i)} \neq \phi$, for all $i\in \{1,2\}$.

If $i,j \in \{1,2\}$ such that $i\neq j$, then by the conditions on $a$ in the hypothesis of Theorem \ref{Theorem SNA} and properties of the sets $\{E_i,C_i\}_{i=1}^{2}$ in definition \ref{definition_preliminaries_expanding contracting regions}, we have $e_{j}^{+} < c_i^{-} + a < c_i^{+}+a$ and thus
\begin{align*}
 [c_i^{-}+g(\theta), c_i^{+}+g(\theta)] \cap [e_{j}^{-}, e_{j}^{+}] & = \phi, \text{ for all } \theta \in I_0.
\end{align*}
Thus, for $i\neq j$, we have $B_{-1}^{(i)}\cap A_{-1}^{(j)} = \phi$. 
\end{proof}

%------------------- LEMMA END -------------------

For the purpose of the following proposition, recall the points $\{\theta^{(i,j)}\}_{i,j=1}^{2} \subset I_0$, as  in definition \ref{definition_preliminaries_theta(i,j)}, which satisfy the condition
\begin{align*}
     g(\theta^{(i,1)}) + c_i = g(\theta^{(i,2)}) + c_i = e_i, \text{ for each } i\in \{1,2\},
\end{align*}
and the relation
$$\theta^{(2,1)} < \theta^{(1,1)} < 1/2 < \theta^{(1,2)} < \theta^{(2,2)}.$$ 
By the definition of $g(\theta)$, for each $j\in\{1,2\}$ we have $g'(\theta^{(j,1)}) < 0$ and $g'(\theta^{(j,2)}) > 0$. 

%------------------- PROP BEGIN -------------------

\begin{proposition}\label{proposition_structure I0_I0 subsets}
For all sufficiently large $\lambda$, and for each $i \in \{1,2\}$, there exists non-empty disjoint intervals $\{I_0^{(i,j)}\}_{j=1}^2$ such that
    \begin{align*}
        \pi_1 (A_{-1}^{(i)} \cap B_{-1}^{(i)}) = \bigcup_{j=1}^2 (I_0^{(i,j)}+\omega),
    \end{align*}
    and $\theta^{(i,j)} \in I_{0}^{(i,j)}$, for each $j\in \{1,2\}$. Further, for the intervals in $\{I_0^{(i,j)}\}_{i,j=1}^2$ we have the following estimates:
\begin{enumerate}[label=\alph*)]
    \item \label{lemma item_structure I0_I0 components_b} For any $\theta \in \bigcup_{i=1}^{2} I_{0}^{(i,1)}$,  
    \begin{align*}
         - \lambda^{1+1/128} <g'(\theta)< -\lambda^{1 - 1/128} 
    \end{align*}
    and for any $\theta \in \bigcup_{i=1}^{2} I_{0}^{(i,2)}$, 
    \begin{align*}
         \lambda^{1 - 1/128} < g'(\theta)< \lambda^{1 + 1/128}.
    \end{align*}
    \item \label{lemma item_structure I0_I0 components_c} $|\bigcup_{i,j=1}^{2}I_0^{(i,j)}| < \lambda^{-1+1/64}.$
        \item \label{lemma item_structure I0_I0 components_d} For any $i,i',j,j' \in \{1,2\}$, we have
    \begin{align*}
        (i,j) \neq (i',j') \Longleftrightarrow I_0^{(i,j)}\cap I_0^{(i',j')} = \phi, 
    \end{align*}
    which implies that the intervals in $\{I_0^{(i,j)}\}_{i,j=1}^{2}$ are all disjoint.
\end{enumerate}
\end{proposition}

\begin{proof}  

Let us fix $i\in \{1,2\}$. By lemma \ref{lemma_structure I0_match Ai Bi}, for all sufficiently large $\lambda$ the sets $A_{-1}^{(i)}$ and $B_{-1}^{(i)}$ have a non-empty intersection over the set $(I_0+\omega)$. It follows from the choice of sets $\{C_i,E_i\}_{i=1}^{2}$ in definition \ref{definition_preliminaries_expanding contracting regions} and the conditions on $a$ that $a > |c_i^{+}-e_i^{-}|$. Since $1/2$ is the unique critical point of $g$ and for which $g(1/2) = a$, this gives us $1/2 \notin \pi_1 (A_{-1}^{(i)}\cap B_{-1}^{(i)})$. From the estimate of $g'(\theta)$ in lemma \ref{lemma_preliminaries_g' estimate}, we deduce that the sets $A_{-1}^{(i)}$ and $B_{-1}^{(i)}$ must intersect transversally of the intervals $(1/2-\lambda^{-1/4}, 1/2)$ and $(1/2, 1/2+\lambda^{-1/4})$. This proves that $\pi_1(A_{-1}^{(i)} \cap B_{-1}^{(i)})$ consists of exactly two closed intervals denoted by $(I_0^{(i,j)}+\omega)$, where $\theta^{(i,j)} \in I_0^{(i,j)}$, for each $j \in\{1,2\}$.

%------------------ Part 3 -------------
\vspace{0.2cm}
\noindent\textit{Proof of \ref{lemma item_structure I0_I0 components_b}}: For the points $\{\theta^{(i,j)}\}_{i,j=1}^{2}$ as defined in \ref{definition_preliminaries_theta(i,j)}, we have $g(\theta^{(j,1)}) < 0$ and $g(\theta^{(j,2)}) > 0$, for each $j \in \{1,2\}$. 
Since each interval $I_0^{(i,j)}$ is a closed neighbourhood of $\theta^{(i,j)}$, there exists values $0<c<d<|a|$ such that $c<|g(\theta)|<d< |a|$, for all $\theta \in \bigcup_{i,j=1}^{2} I_0^{(i,j)}$. By lemma \ref{lemma_preliminaries_g' estimate}, for a big enough $\lambda$, this gives us
\begin{align}\label{equation_structure of I0_negative bound for g' on I0 components}
     - \lambda^{1+1/128} <g'(\theta)< -\lambda^{1 - 1/128},  \text{ for all }  \theta \in \bigcup_{i=1}^{2} I_0^{(i,1)},
\end{align}
and 
\begin{align}\label{equation_structure of I0_positive bound for g' on I0 components}
     \lambda^{1 - 1/128} < g'(\theta)< \lambda^{1 + 1/128},  \text{ for all }  \theta \in \bigcup_{i=1}^{2} I_0^{(i,2)}.
\end{align}

%---------------------- Proof of 4 ---------------------

\vspace{0.2cm}
\noindent\textit{Proof of \ref{lemma item_structure I0_I0 components_c}:} By construction, $(I_0^{(i,j)}+\omega)$ is the projection on the first co-ordinate of the transverse intersection of $A_{-1}^{(i)}$ and $B_{-1}^{(i)}$ over a closed interval. Then, for each $i\in \{1,2\}$ the set $I_0^{(i,1)}$ is of the form $(s,t)$, where $s,t \in I_0$ satisfies the equations $c_i^{-}+g(s) = e_i^{+}$ and $c_i^{+}+g(t) = e_i^{-}$. For each $i,j \in \{1,2\}$, by the Mean Valued Theorem, we have 
\begin{align*}
    |I_0^{(i,j)}| \leq \frac{2\max \{|c_i^{+}-c_i^{-}|, |e_i^{+}-e_i^{-}|: i=1,2\}}{|g'(\theta)|_{inf}},
\end{align*}
where $|g'(\theta)|_{inf} = \inf\{|g'(\theta)|: \theta \in \bigcup_{i,j=1}^{2}I_0^{(i,j)}\}$, and by the bound on $|g'(\theta)|$ in equations \eqref{equation_structure of I0_negative bound for g' on I0 components} and \eqref{equation_structure of I0_positive bound for g' on I0 components}, it gives us 
\begin{align*}
    |I_0^{(i,j)}| \leq \frac{2\max \{|c_i^{+}-c_i^{-}|, |e_i^{+}-e_i^{-}|: i=1,2\}}{\lambda^{1-1/128}}.
\end{align*}
 For all sufficiently large $\lambda >0$, this proves that
\begin{align}\label{equation_structure of I0_bound measure of I0 components}
     \left|\bigcup_{i,j=1}^{2} I_{0}^{(i,j)} \right| \leq \lambda^{-1+1/64}.
\end{align}

%-------------------------- PROOF OF 5 ---------------------------------

\vspace{0.2cm}
\noindent\textit{Proof of \ref{lemma item_structure I0_I0 components_d}:} 
It is now sufficient to prove that $I_0^{(2,1)} \cap I_0^{(1,1)} = \emptyset$ and $I_0^{(1,2)} \cap I_0^{(2,2)} = \emptyset$. By lemma \ref{lemma_preliminary estimates_|G| and |S|}, we have the bound
\begin{align}
   |G^{(1)}| = d(\theta^{(2,1)}, \theta^{(1,1)}) > \lambda^{-\frac{1}{2}-\frac{1}{32}}.
\end{align}
Thus, for $\lambda$ big enough such that $2\lambda^{-1+1/64} < \lambda^{-1/2-1/32}$, we obtain
\begin{align*}
 2(|I_0^{(1,1)}|+ |I_0^{(2,1)}|) & < d(\theta^{(1,1)}, \theta^{(2,1)}) \implies I_0^{(1,1)} \cap I_0^{(2,1)} = \emptyset.
\end{align*}
Similarly, from the bound $ d(\theta^{(1,2)}, \theta^{(2,2)}) > \lambda^{-\frac{1}{2}-\frac{1}{32}}$ in lemma \ref{lemma_preliminary estimates_|G| and |S|}, it follows that $I_0^{(2,1)} \cap I_0^{(2,2)} = \emptyset$.
\end{proof}

%------------------- PROP END -------------------

%--------------------- Critical Set I_1 ---------------------------

\subsection{Structure of critical set $I_1$}\label{subsection_structure of critical set I1}

In this subsection, for each $\omega \in \mathcal{F}_0$, we define a specific critical set $I_1(\omega)$ satisfying the conditions in definition \ref{definition_block n_Critical sets In}. For the particular critical set, we establish its geometric structure and analyse its dependence on the parameter $\omega$. We emphasize that, throughout the remainder of the article, the notation $I_1(\omega)$ is reserved exclusively for the sets introduced in the following definition.
\begin{definition}[Critical subsets]\label{definition_structure I1_I1 definition}
    For any $\omega \in \mathcal{F}_{0}(M_0)$, we fix the set $I_1(\omega)$ to be the following: 
\begin{align*}
I_{1}(\omega) & := \pi_2 \left( F^{M_0}(A_0) \bigcap F^{-M_0}(B_0) \right),
\end{align*}
where
\begin{align*}
    A_0 & = (I_0 - M_0 \omega) \times \left( C_1\cup C_2 \right) \\
    B_0 & = (I_0 + M_0\omega) \times \left( E_1 \cup E_2 \right).
\end{align*}
Then, for each $i,j \in \{1,2\}$, we define
    \begin{align*}
        I_1^{(i,j)}(\omega) = I_1(\omega) \cap I_0^{(i,j)}(\omega), 
    \end{align*}
    where the set $I_0^{(i,j)}(\omega)$ is as described in proposition \ref{proposition_structure I0_I0 subsets}.
\end{definition}
For convenience, we use the notation $I_1:= I_1(\omega)$ when appropriate. Let us denote the components of the sets $A_0$ and $B_0$, by the following:
\begin{align*}
    A_0^{(i)} & :=  (I_0 - M_0\omega) \times C_i \\
    B_0^{(i)} & := (I_0 + M_0\omega) \times E_i, 
\end{align*}
where $i\in \{1,2\}$.

%------------------- LEMMA BEGIN -------------------
\begin{lemma}\label{lemma_stucture I1_blocks and curves relation}
For any $\omega \in \mathcal{F}_0$ and $i\in \{1,2\}$, the sets $F^{-M_0+1}(B_0^{(i)})$ and $F^{M_0+1}(A_0^{(i)})$ can be expressed as the following:
\begin{align*}
    F^{-M_0+1}(B_0^{(i)}) & = \left\{(\theta + \omega,x): \theta \in I_0 ,\ x\in \left[\psi_{0,i}^{-}(\theta), \psi_{0,i}^{+}(\theta)\right] \right\}, \\
    F^{M_0+1}(A_0^{(i)}) & = \left\{(\theta+\omega,x): \theta \in I_0,\ x\in \left[\phi_{0,i}^{-}(\theta), \phi_{0,i}^{+}(\theta)\right] \right\},
\end{align*}
where 
\begin{align*}
    \psi_{0,i}^{\pm}(\theta) & = \pi_2 \circ F^{-M_0+1}(\theta+M_0\omega, e_i^{\pm}),\\
    \phi_{0,i}^{\pm}(\theta) & = \pi_2 \circ F^{M_0+1}(\theta - M_0\omega, c_i^{\pm}),
\end{align*}
for all $\theta \in I_0$.
\end{lemma}

\begin{proof}
For each $\theta \in I_0$, since $f_{\theta}(x)$ is a continuous and orientation preserving diffeomorphism on the torus, we obtain
   \begin{align*}
       F^{-M_0+1}\left( \{\theta + M_0\omega \} \times [e_i^{-}, e_i^{+}]\right) & = \{\theta + \omega \} \times [\psi_{0,i}^{-}(\theta), \psi_{0,i}^{+}(\theta)] \\
       F^{M_0+1}\left( \{\theta - M_0\omega \} \times [c_i^{-}, c_i^{+}]\right) & = \{\theta + \omega \} \times [\phi_{0,i}^{-}(\theta), \phi_{0,i}^{+}(\theta)],
   \end{align*}
for each $i \in \{1,2\}$. This proves the lemma.
\end{proof}

%------------------- LEMMA END -------------------

The preceding lemma shows that the geometry of the critical set $I_1(\omega)$ is closely connected to the structure of the curves $\phi_{0,i}^{\pm}$ and $\psi_{0,i}^{\pm}$, where $i \in \{1,2\}$. In the following propositions we establish the necessary estimates for these curves.

%------------------- PROP BEGIN -------------------

\begin{proposition}\label{proposition_structure I1_backward iterate curves}
    For all sufficiently large $\lambda$, $\omega \in \mathcal{F}_0$ and for each $i \in \{1,2\}$, the maps $\psi_{0,i}^{\pm}(\theta)$ satisfy the following conditions:
\begin{enumerate}[label=\alph*)]
\item \label{lemma item_structure I1_backward iterate curves_a} For every $\theta \in I_0$, we have $ [\psi_{0,i}^{-}(\theta), \psi_{0,i}^{+}(\theta)]  \subset E_i.$

\item \label{lemma item_structure I1_backward iterate curves_b} For every $\theta \in I_0$, 
\begin{align*}
    \alpha_c^{b(M_0-1)}|e_i^+ - e_i^-| & \leq |\psi_{0,i}^{+}(\theta) - \psi_{0,i}^{-}(\theta)| \leq \alpha_c^{M_0-1}|e_i^+ - e_i^-|.
\end{align*}
\item  \label{lemma item_structure I1_backward iterate curves_c} For every $\theta \in I_0$,
$$|\partial_{\omega}\psi_{0,i}^{\pm}(\theta)|, |\partial_{\theta}\psi_{0,i}^{\pm}(\theta)| < \lambda^{-1/4+1/64}.$$ 
\end{enumerate}

\end{proposition}

\begin{proof}
%\textit{Proof of \ref{lemma item_structure I1_backward iterate curves_a}:} 
Recall that, by the disjointness condition $(\mathcal{F}1)_0$ (in definition \ref{definition_block 0_disjointness condn}), the intervals in the set $\left\{I_0+m\omega\right\}_{m=0}^{M_0}$ are pairwise disjoint. For each $i \in \{1,2\}$, the expanding region $E_i$ satisfies the invariance condition discussed in equation \eqref{setup_I_0 expansion charac}. Then, for each  $\theta \in (I_0+M_0 \omega)$ and $x\in E_i$, this gives us
\begin{align*}
    \pi_2 \circ F^{-k}(\theta, x) \in E_i, 
\end{align*}
for all $k \in \{0,\dots, M_0-1\}$.
Since the map $f_{\theta}^{-1}$ is orientation preserving for any $\theta \in \mathbb{T}$, this proves that
\begin{align*}
    [\psi_{0,i}^{-}(\theta),\psi_{0,i}^{+}(\theta)] \subset E_i, \text{ for all } i\in \{1,2\}.
\end{align*}
%\textit{Proof of \ref{lemma item_structure I1_backward iterate curves_b}:} 
For any $x \in (E_1\cup E_2)$ and $\theta \in (I_0+M_0\omega)$, we can derive the expression
\begin{align*}
    \partial_x\left(\pi_2 \circ F^{-M_0+1}(\theta, x) \right) = \prod_{j=1}^{M_0-1} \partial_xf^{-1}(f(x_{-j})).
\end{align*}
Then, from $(E1)_{-1}$ (in lemma \ref{lemma_block -1 expansion}) and the bound $\alpha_c^b < \alpha_u^{-1}$, we obtain
\begin{align*}
     \alpha_c^{b(M_0-1)} & < |\partial_x \left( \pi_2 \circ F^{-M_0+1}(\theta, x)\right)| < \alpha_c^{M_0-1},  
\end{align*}
for all $x_0 \in E_1\cup E_2$ and $\theta_0 \in (I_0+M_0\omega)$. This proves that, for each $i \in \{1,2\}$ and $\theta \in I_0$,
\begin{align*}
\alpha_c^{b(M_0-1)} |e_i^{+}-e_i^{-}| < |\psi_{0,i}^{+}(\theta) - \psi_{0,i}^{-}(\theta)| < \alpha_c^{M_0-1}|e_i^{+}-e_i^{-}|.
\end{align*}
For each $\theta \in I_0$, we have the bound  $|g'(\theta+k\omega)|<  \lambda^{-1/4 + 1/128}$ for all $k \in \{1,\dots ,M_0-1\}$, due to the estimate of $g'(\theta)$ outside the critical set $I_0$ in lemma \ref{lemma_preliminaries_g' estimate}. For each $i\in \{1,2\}$ and $\theta \in I_0$, we then derive
\begin{align*}
    | \partial_{\theta} \psi_{0,i}^{\pm}(\theta)| & \leq \sum_{k=1}^{M_0-1} |g'(\theta+k\omega)|\prod_{j=1}^k |(f^{-1})'(f(x_{-M_0+j}))| \\
& \leq \sum_{k=1}^{M_0-1} \lambda^{-1/4 + 1/128} \alpha_c^k < \lambda^{-1/4+1/64}.
\end{align*}
Similarly, we obtain  
\begin{align*}
    | \partial_{\omega} \psi_{0,i}^{\pm}(\theta)| & \leq \sum_{k=1}^{M_0-1} |k g'(\theta+k\omega)|\prod_{j=1}^k |(f^{-1})'(f(x_{-M_0+j}))|  < \lambda^{-1/4+1/64}.
\end{align*}

\end{proof}

%------------------- PROP END -------------------

%------------------- PROP BEGIN -------------------
\begin{proposition}\label{proposition_structure I1_forward iterate curves}
    For all sufficiently large $\lambda$ and for each $i \in \{1,2\}$, the curves $\phi_{0,i}^{\pm}(\theta)$ satisfy the following conditions: 
    
\begin{enumerate}[label=\alph*)]
\item \label{lemma item_structure I1_forward iterate curves_a} For all $\theta \in I_0$, we have $\pi_2 \circ F^{M_0}(A_0^{(i)}) \subset C_i$.

\item \label{lemma item_structure I1_forward iterate curves_b} For all $\theta \in I_0$,
\begin{align*} 
\alpha_c^{b(M_0+1)}|c_i^+ - c_i^-| & \leq |\phi_{0,i}^+(\theta) - \phi_{0,i}^-(\theta)| \leq \alpha_c^{M_0+1}|c_i^{+}-c_i^{-}| 
\end{align*}

\item \label{lemma item_structure I1_forward iterate curves_c} For any $i\in \{1,2\}$, if $\theta \in I^{(i,1)}_0$ then
\begin{align*}
     -\lambda^{1+1/64} & < \partial_{\theta} \phi_{0,i}^{\pm}(\theta) < -\lambda^{1-1/64}, 
\end{align*}
and if $\theta \in I^{(i,2)}_0$ then
\begin{align*}
    \lambda^{1-1/64} & < \partial_{\theta} \phi_{0,i}^{\pm}(\theta) < \lambda^{1+1/64}.
\end{align*}
For any $\theta \in \bigcup_{i,j=1}^2 I_0^{(i,j)}$,  
\begin{align*}
    \lambda^{-1/4-1/64} & < |\partial_{\omega} \phi_{0,i}^{\pm}(\theta)| < \lambda^{-1/4+1/64}.
\end{align*}
\end{enumerate}
\end{proposition}

\begin{proof}
%\textit{Proof of \ref{lemma item_structure I1_forward iterate curves_a}:} 
Recall that, by the disjointness condition $(\mathcal{F}1)_0$ (in definition \ref{definition_block 0_disjointness condn}), the intervals in the set $\left\{I_0+m\omega\right\}_{m=0}^{M_0}$ are pairwise disjoint. For each $i \in \{1,2\}$, the contracting region $C_i$ satisfies the invariance condition discussed in equation \eqref{setup_I_0 contraction charac}. Then, for each  $\theta \in (I_0-M_0 \omega)$ and $x\in C_i$, this gives us
\begin{align}\label{equation_structure I1_forward iterate_1}
    \pi_2 \circ F^{k}(\theta, x) \in C_i,
\end{align}
for all $k\in \{0,\dots, M_0\}$ and this proves that $\pi_2 \circ F^{M_0}(A_0^{(i)}) \subset C_i$.
%\textit{Proof of \ref{lemma item_structure I1_forward iterate curves_b}:} 
From the definition of the curves $\phi_{0,i}^{\pm}(\theta)$, for any $\theta \in I_0$ and $x\in C_i$, we obtain the expression
\begin{align*}
     \partial_{x} \left(\pi_2 \circ F^{M_0+1}(\theta-M_0\omega, x_0) \right) & = \prod_{k=0}^{M_0}f'(x_k).
\end{align*}
Then, from $(C1)_{-1}$ (in lemma \ref{lemma_block -1 contraction}) and the bound $\alpha_c^b < \alpha_u^{-1}$, we obtain
\begin{align*}
\alpha_c^{b(M_0+1)} < |\partial_x F^{M_0+1}(\theta-M_0\omega, x_0) | & < \alpha_c^{M_0+1},
\end{align*}
for all $x\in C_i$. This proves that, for every $\theta \in I_0$,
    \begin{align*}
 \alpha_c^{b(M_0+1)}|c_i^+ - c_i^-| \leq |\phi_{0,i}^+(\theta) - \phi_{0,i}^-(\theta)| \leq \alpha_c^{M_0+1}|c_i^{+}-c_i^{-}|. 
\end{align*}
%\textit{Proof of \ref{lemma item_structure I1_forward iterate curves_c}:} 
For each $i \in \{1,2\}$ and $\theta \in I_0$, we can derive the expression
\begin{align}\label{equation_structure I1_forward iterate curves_phi}
    \partial_{\theta} \phi_{0,i}^{\pm}(\theta) = g'(\theta) + \sum_{k=0}^{M_0-1} g'(\theta-(k+1)\omega)\left( \prod_{j=0}^{k} f'(x_{M_0-j}) \right).
\end{align}
Using the bounds $\displaystyle \left|\prod_{j=0}^{k} f'(x_{M_0}-j)\right| \leq \alpha_c^{k}$ and $|g'(\theta-(k+1))|< \lambda^{-1/4 + 1/128}$, for all $k \in \{0,\dots ,M_0-1\}$, we obtain 
\begin{align*}
    \left| \sum_{k=0}^{M_0-1} g'(\theta-(k+1)\omega) \prod_{j=0}^{k} |f'(x_{M_0-j})| \right| & \leq \sum_{k=0}^{M_0-1}\left( \lambda^{-1/4+ 1/128} \alpha_c^{k+1}\right) < \lambda^{-1/4+ 1/64}.
\end{align*}
Further, using the estimates from proposition \ref{proposition_structure I0_I0 subsets} \ref{lemma item_structure I0_I0 components_b}, for each $i\in \{1,2\}$ we can derive
\begin{align*}
    - \lambda^{1 + 1/64} & < \partial_{\theta} \phi_{0,i}^{\pm}(\theta) < -\lambda^{1 - 1/64},  \text{ for all } \theta \in I_0^{(i,1)} \\
    \lambda^{1-1/64} & < \partial_{\theta} \phi_{0,i}^{\pm}(\theta) < \lambda^{1+1/64},   \text{ for all } \theta \in I_0^{(i,2)}.
\end{align*}
Similarly, we obtain 
\begin{align*}
 |\partial_{\omega} \phi_{0,i}^{\pm}(\theta)| & < \sum_{k=0}^{M_0-1} (k+1) \lambda^{-1/4} \alpha_c^{k} < \lambda^{-1/4+1/64}.
\end{align*}
\end{proof}
%------------------- PROP END -------------------

Note that the estimates in propositions \ref{proposition_structure I1_backward iterate curves} and \ref{proposition_structure I1_forward iterate curves} prove that the curves $\phi_{0,i}^{\pm}$ and $\psi_{0,i}^{\pm}$ depend robustly on the parameter $\omega$.

%------------------- LEMMA BEGIN -------------------
\begin{lemma}\label{lemma_structure I1_match Ai Bi}
    For all $i\in \{1,2\}$, we have
\begin{align*}
    F^{-M_0+1}(B_0^{(i)}) & \subset B_{-1}^{(i)} \\
    F^{M_0+1}(A_0^{(i)}) & \subset A_{-1}^{(i)}.
\end{align*}
\end{lemma}

\begin{proof}
Let us fix some $i\in \{1,2\}$. By proposition \ref{proposition_structure I1_forward iterate curves} \ref{lemma item_structure I1_forward iterate curves_a}, we have $F^{M_0}(A_0^{(i)}) \subset (I_0\times C_i)$. By definition \ref{definition_structure I0_iterating blocks} of the set $A_{-1}^{(i)}$ and the property , we have $F(I_0\times C_i) \subset A_{-1}^{(i)}$. Thus, $F^{M_0+1}(A_0^{(i)}) \subset  A_{-1}^{(i)}$. 
Similarly, by proposition \ref{proposition_structure I1_backward iterate curves} \ref{lemma item_structure I1_backward iterate curves_a} we have $F^{-M_0+1}(B_0^{(i)}) \subset  B_{-1}^{(i)}$. 
\end{proof}

%------------------- LEMMA END -------------------

The following proposition shows that, for each $\omega \in \mathcal{F}_0(M_0)$, the set $I_1(\omega)$ specified in definition \ref{definition_structure I1_I1 definition} satisfies the criteria of a critical set introduced in definition \ref{definition_block n_Critical sets In}. Furthermore, we characterize the geometric structure of $I_1(\omega)$.

%------------------- PROP BEGIN -------------------
\begin{proposition}\label{proposition_structure I1_critical set I_1}
   For each $\omega \in \mathcal{F}_0(M_0)$,
    \begin{enumerate}[label=\alph*)]
        \item \label{lemma item_structure I1_I1 estimates_a} The set $I_1(\omega)$, as in definition \ref{definition_structure I1_I1 definition}, satisfies 
        $$I_1(\omega) = \bigcup_{i,j =1}^{2} I_1^{(i,j)}(\omega),$$ where for each $i,j\in \{1,2\}$, the set $I_1^{(i,j)}$ is a non-empty closed interval.
        \item \label{lemma item_structure I1_I1 estimates_b} For each $i,j\in \{1,2\}$,
        $$\lambda^{-2}\alpha_c^{b(M_0-1)} < |I^{(i,j)}_1(\omega)| < \lambda^{-1+1/32} \alpha_{c}^{M_0-1}.$$
        \item \label{lemma item_structure I1_I1 estimates_c} For each $\omega \in \mathcal{F}_0(M_0)$, $|\partial_{\omega}I_1(\omega)| < \lambda^{-1}$.
    \end{enumerate}
\end{proposition}

\begin{proof}
%\textit{Proof of \ref{lemma item_structure I1_I1 estimates_a}:}
Let us fix $i,j \in \{1,2\}$. It follows from lemma \ref{lemma_structure I1_match Ai Bi} and proposition \ref{proposition_structure I0_I0 subsets} that, 
\begin{align*}
     \left(F^{M_0+1}(A_0^{(i)}) \cap F^{-M_0+1}(B_0^{(i)})\right) \subset \left(A_{-1}^{(i)} \cap B_{-1}^{(i)}\right),
\end{align*}
which gives $I_1  \subset \bigcup_{i,j=1}^{2} I_0^{(i,j)}$. Since $I_1^{(i,j)} = I_1 \cap I_0^{(i,j)}$ by definition \ref{definition_structure I1_I1 definition}, 
this proves that $\displaystyle I_1 = \bigcup_{i,j=1}^{2} I_1^{(i,j)}$. In proposition \ref{proposition_structure I0_I0 subsets}, we proved that $(I_0^{(i,j)}+\omega)$ is a non-empty closed interval, which is the projection of the transverse intersection of the sets $A_{-1}^{(i)}$ and $B_{-1}^{(i)}$. From the estimates in proposition \ref{proposition_structure I1_backward iterate curves} \ref{lemma item_structure I1_backward iterate curves_c} and proposition \ref{proposition_structure I1_forward iterate curves} \ref{lemma item_structure I1_forward iterate curves_c}, we can infer that the curves $\phi_{1,i}^{\pm}(\theta)$ and $\psi_{1,i}^{\pm}(\theta)$ are monotonic over each interval $I_0^{(i,j)}(\omega)$ and also intersect each other transversally. Since 
\begin{multline*}
    I_1^{(i,j)}(\omega) \\ = \pi_1 \left( \left\{(\theta,x): \theta \in I_0^{(i,j)}, x \in [\psi_{0,i}^{-}(\theta), \psi_{0,i}^{+}(\theta)] \right\} \bigcap \left\{ (\theta,x): \theta \in I_0^{(i,j)}, x \in [\phi_{0,i}^{-}(\theta), \phi_{0,i}^{+}(\theta)] \right\} \right),
\end{multline*}
this proves that $I_1^{(i,j)}$ must be a non-empty closed interval. 

%\textit{Proof of \ref{lemma item_structure I1_I1 estimates_b}:}
 For convenience, we use the notations $\phi_{0,i}^{\pm} = \phi_0^{\pm}$, $\psi_{0,i}^{\pm} = \psi_0^{\pm}$. Without loss of generality, assume that $j=2$ and thus $\partial_{\theta}\phi_0^{\pm}>0$. By proposition \ref{proposition_structure I1_critical set I_1} \ref{lemma item_structure I1_I1 estimates_a}, the set $I_1^{(i,j)}$is a non-empty closed interval. Let us denote the interval $I_1^{(i,j)}$ by $[a_1, b_1]$, where
\begin{align}\label{I_1 end points}
    \phi_0^{-}(a_1) < \phi_0^{+}(a_1) & = \psi_0^{-}(a_1) < \psi_0^{+}(a_1)\\ 
    \phi_0^{+}(b_1) > \phi_0^{-}(b_1) & = \psi_0^{+}(b_1) > \psi_0^{-}(b_1).
\end{align}
From the geometry of the curves and by the Mean Valued Theorem, we obtain the following inequalities:
\begin{align*}
   \frac{2\inf\{|\phi_0^{+}(\theta) - \phi_0^{-}(\theta)|: \theta\in I_0^{(i,j)}\}}{|\partial_{\theta}\phi_0^{\pm}|_{\sup} + |\partial_{\theta}\psi_0^{\pm}|_{\sup}} < b_1-a_1 < \frac{2\sup\{|\phi_0^{+}(\theta) - \phi_0^{-}(\theta)|: \theta\in I_0^{(i,j)}\}}{|\partial_{\theta}\phi_0^{\pm}|_{\inf} - |\partial_{\theta}\psi_0^{\pm}|_{\sup}}
\end{align*}
    From proposition \ref{proposition_structure I1_forward iterate curves} and proposition \ref{proposition_structure I1_backward iterate curves}, we have the bounds
    \begin{align*}
      \lambda^{1-1/64} < \partial_{\theta}\phi_0^{\pm} & < \lambda^{1+1/64} \\
       0 < |\partial_{\theta}\psi_0^{\pm}| & < \lambda^{-1/4+1/64},
\end{align*}
for any $\theta \in I_0^{(i,j)}$, and this gives us the estimate
\begin{align*}
    \lambda^{-1-1/32}\alpha_c^{b(M_0-1)}\min_{i=1,2}\{|e_i^{+}-e_{i}^{-}|, |c_i^{+}-c_{i}^{-}|\} < |I_1^{(i,j)}| <  \lambda^{-1+1/32} \alpha_{c}^{M_0-1}.
\end{align*}
For a sufficiently large $\lambda$ with respect to the regions $E_i$ and $C_i$, this proves that
\begin{align*}
    \lambda^{-2}\alpha_c^{b(M_0-1)} < |I_1^{(i,j)}| < \lambda^{-1+1/32} \alpha_{c}^{M_0-1}.
\end{align*}
The proof of the case $j=1$ follows similarly.
%\textit{Proof of \ref{lemma item_structure I1_I1 estimates_b}:} 

If $a_1(\omega)$ and $b_{1}(\omega)$ denotes the end points of the critical set $I^{(i,j)}_1(\omega)$, then by equations (\ref{I_1 end points}) the points satisfy the equations
\begin{align*}
   \phi_0^{+}(a_1(\omega)) & = \psi_0^{-}(a_1(\omega)) \\ 
    \phi_0^{-}(b_1(\omega)) & = \psi_0^{+}(b_1(\omega)).
\end{align*}
By the implicit function theorem for the function $a_1(\omega)$, we can derive 
\begin{align*}
a_1'(\omega) & = \frac{\partial_{\omega}\psi_{0}^{-}(a_1(\omega), \omega) - \partial_{\omega}\phi_{0}^{+}(a_1(\omega), \omega)}{\partial_{\theta}\phi_{0}^{+}(a_1(\omega), \omega) - \partial_{\theta}\psi_{0}^{-}(a_1(\omega). \omega)}.
\end{align*}
Then, by the estimates of partial derivatives in proposition \ref{proposition_structure I1_forward iterate curves} \ref{lemma item_structure I1_forward iterate curves_c} and proposition \ref{proposition_structure I1_backward iterate curves} \ref{lemma item_structure I1_backward iterate curves_c}, it follows that
\begin{align*}
   a_1'(\omega) & < \frac{|\partial_{\omega}\psi_{0}^{-} - \partial_{\omega}\phi_{0}^{+}|_{\sup}}{|\partial_{\theta}\phi_{0}^{+}|_{\inf} - |\partial_{\theta}\psi_{0}^{-}|_{\sup}} \\
    & < \frac{2\lambda^{-1/4+1/64}}{\lambda^{1-1/64}-\lambda^{-1/4+1/64}} < \lambda^{-1},
\end{align*}
and similarly, we can show that $b_1'(\omega) < \lambda^{-1}$. This proves the bound $|\partial_{\omega}I_1(\omega)|<\lambda^{-1}$.
\end{proof}
%------------------- PROP END -------------------

%%%%%%%%%%%%%%%%%%%%%%%%%%% CHARACTERISTICS OF G and S %%%%%%%%%%%%%%%%%%%

\subsection{Characteristics of subsets of ${G}$ and ${S}$}\label{subsection_characteristics G and S}

Recall the subsets $G$ and $S$ of the critical set $I_0$, given by 
\begin{align*}
    G & = (\theta^{(2,1)}, \theta^{(1,1)}) \cup (\theta^{(2,1)}, \theta^{(2,2)}) \\
    S & = (\theta^{(1,1)}, \theta^{(1,2)}) 
\end{align*}
which were introduced in definition \ref{definition_preliminaries F(1,1)_G and S sets}, and are characterised by the strength of forcing along the fibre. In this subsection, we define subsets of $G$ and $S$ that will be crucial for analysing the behaviour of orbits along the fibre direction, and particularly for determining the conditions for synchronisation. We now introduce subsets of $G$ and then, show them to be suitable choices for the sets $\hat{G}$ and $\bar{G}$, an in definitions \ref{definition_block n_sync condn parameter} and \ref{definition_block n_sync condn set_g1 g2}.

\begin{definition}[Synchronisation sets]\label{definition_structure I1_sync set Gs}
    We define the following subsets of $G$:
    \begin{align}\label{equation_structure I1_setting G_1, G_2}
    G_1 & := \left( \theta^{(2,1)}+ \lambda^{-1+1/64}, \theta^{(1,1)}-\lambda^{-1+1/64} \right) \cup \left( \theta^{(1,2)}+\lambda^{-1+1/64}, \theta^{(2,2)}-\lambda^{-1+1/64} \right)\\
    G_2 & := \left( \theta^{(2,1)}+3\lambda^{-1+1/64}, \theta^{(1,1)}-3\lambda^{-1+1/64} \right) \cup \left( \theta^{(1,2)}+3\lambda^{-1+1/64}, \theta^{(2,2)}-3\lambda^{-1+1/64} \right),
\end{align}
For each $i\in \{1,2\}$, we define $G_1^{(i)} := G_1 \cap G^{(i)}$ and $G_2^{(i)} := G_2 \cap G^{(i)}$.
\end{definition}

For each $\omega \in \mathcal{F}_0$ and $i,j\in \{1,2\}$, we recall the estimates $I^{(i,j)}_1(\omega) \subset I^{(i,j)}_0(\omega)$, $|I^{(i,j)}_1(\omega)| < \alpha_c^{M_0-1}$ and $|I^{(i,j)}_0(\omega)| < \lambda^{-1+1/64}$, from propositions \ref{proposition_structure I0_I0 subsets} and \ref{proposition_structure I1_critical set I_1}. 
Then, by the definition of the set $G$, it follows that 
\begin{align*}
    G_2 \subset B_{\varepsilon}(G_2) \subset G_1 \subset G,
\end{align*} 
for all $\varepsilon < |I_1(\omega)|$, and that 
 $G_1 = G_1^{(1)} \cup G_1^{(2)}$ and $G_2 = G_2^{(1)} \cup G_2^{(2)}$. 

\begin{definition}[Switching sets]\label{definition_structure I1_switch set Ss}
    We define the following subsets of $S$ and $I_0$:
    \begin{align}\label{equation_structure I1_setting S_1, S_2, S_0}
    S_2 & := \left( \theta^{(1,1)} + 3\lambda^{-1+1/64}, \theta^{(1,2)} - 3\lambda^{-1+1/64} \right)\\
    S_3 & := \left( \theta^{(1,1)} - 3\lambda^{-1+1/64}, \theta^{(1,2)} + 3\lambda^{-1+1/64} \right).\label{equation_structure I1_setting S_3}
\end{align}
For each $\omega \in \mathcal{F}_0$, let $S_1(\omega)$ denote the largest open interval in $I_0\setminus I_1(\omega)$ which contains the interval $S_2$.
\end{definition}
Note that $S_2 \subset S \subset S_3 \subset I_0$. The following proposition proves that the sets $G_1$ and $G_2$ are in fact suitable choices for the sets $\bar{G}$ and $\hat{G}$,  especially by verifying that $G_1$ satisfies the conditions $(g1-3)$ as in definition \ref{definition_block n_sync condn set_g1 g2}.

%------------------- PROP BEGIN -------------------
\begin{proposition}\label{proposition_structure I1_sync set Gs}
   The sets $G_1$, $G_2$ defined in definition \ref{definition_structure I1_sync set Gs} satisfy the following conditions:
\begin{enumerate}[label=\alph*)]
    \item\label{lemma item_structure I1_sync set Gs_a} For all $\omega \in \mathcal{F}_0$, we have 
    \begin{align*}
        G_1 & \subset I_0 \setminus \left( \bigcup_{i,j=1}^{2} I_0^{(i,j)} \right) \subset I_0\setminus I_1(\omega) \\
        G_2 & \subset I_0 \setminus \left( \bigcup_{i,j=1}^{2} 3I_0^{(i,j)} \right) \subset I_0\setminus 3 I_1(\omega).
    \end{align*}
    \item\label{lemma item_structure I1_sync set Gs_b} $|G_2| \geq \lambda^{-1/2-1/8}$.
    \item\label{lemma item_structure I1_sync set Gs_c} The set $G_1$ satisfies the conditions  
    \begin{align*}
\pi_2 \circ F^{2M_0+S_0}((G_1-M_0\omega)\times (C_1 \cup C_2)) & \subset C_1 \\
\pi_2 \circ F^{-2M_0-S_0}((G_1+M_0\omega)\times (E_1 \cup E_2)) & \subset E_1.
    \end{align*}
\end{enumerate}
\end{proposition}

\begin{proof}
%\textit{Proof of \ref{lemma item_structure I1_sync set Gs_a} and \ref{lemma item_structure I1_sync set Gs_b}:}
Let us fix $i,j\in \{1,2\}$. We have $\theta^{(i,j)}\in I_0^{(i,j)}$ and $|I_0^{(i,j)}|< \lambda^{-1+1/64}$, from proposition \ref{proposition_structure I0_I0 subsets}, and the estimate $d(\theta^{(2,1)}, \theta^{(1,1)}) = d(\theta^{1,2}, \theta^{(2,2)}) > 8\lambda^{-1/2-1/32}$, from lemma \ref{lemma_preliminary estimates_|G| and |S|}. We can infer from  these estimates that $G_1$ and $G_2$ are both non-empty sets that satisfy the inclusions: 
\begin{align*}
    G_1 & \subset I_0 \setminus \left( \bigcup_{i,j=1}^{2} I_0^{(i,j)} \right) \subset I_0 \setminus I_1 \\
    G_2 & \subset I_0 \setminus \left( \bigcup_{i,j=1}^{2} 3I_0^{(i,j)} \right) \subset I_0 \setminus 3I_1.
\end{align*}
Then, from the definition of $G_1$ and $G_2$, it follows that
\begin{align*}
    |G_1| &\geq 2(d(\theta^{(2,1)}, \theta^{(1,1)}) - 2|I_0^{(i,j)}|) > \lambda^{-1/2-1/16} \text{ and } \\
    |G_2| &\geq 2(d(\theta^{(2,1)}, \theta^{(1,1)}) - 6|I_0^{(i,j)}|) > \lambda^{-1/2-1/8}.
    \end{align*}
%\textit{Proof of \ref{lemma item_structure I1_sync set Gs_c}:} 
Let us fix a point $\theta \in G_1$. By proposition \ref{proposition_structure I1_sync set Gs} \ref{lemma item_structure I1_sync set Gs_a} we have $\theta-M_0\omega \in I_0-M_0\omega \subset \Theta_0$. Since $\theta \notin I_1$, it follows from $(C1-2)_0$ (in proposition \ref{proposition_block n contraction}) that 
\begin{align*}
\pi_2 \circ F^{2M_0+S_0}((\theta-M_0\omega)\times [c_1^{-}, c_1^{+}] ) & \subset (C_1\cup C_2)\\
    \pi_2 \circ F^{2M_0+S_0}((\theta-M_0\omega)\times [c_2^{-}, c_2^{+}] ) & \subset (C_1\cup C_2).
\end{align*}
Since the segment $\{\theta-k\omega\}_{k=1}^{M_0}$ does not intersect $I_0$, it follows from in equation (\ref{setup_I_0 contraction charac}) that
\begin{align*}
 F^{M_0+1}(\{\theta-M_0\omega\}\times [c_1^{-}, c_1^{+}]) 
    & \subset \left( \{\theta+\omega\} \times [c_1^{-} +g(\theta), c_1^{+}+g(\theta)] \right) \\
     F^{M_0+1}(\{\theta-M_0\omega\}\times [c_2^{-}, c_2^{+}]) 
    & \subset \left( \{\theta+\omega\} \times [c_2^{-} +g(\theta), c_2^{+}+g(\theta)] \right).
\end{align*}
Since $\theta \notin \bigcup_{i,j=1}^{2}I_0^{(i,j)}$, by our construction of the sets $I_0^{(i,j)}$ in proposition \ref{proposition_structure I0_I0 subsets}, we can infer that
\begin{align}\label{equation_structure I1_sync set Gs_1}
    \left\{ (\theta+\omega) \times [c_1^{-}+g(\theta), c_1^{+}+g(\theta)] \right\} \cap \left\{(\theta+\omega) \times [e_1^{-}, e_1^{+}]\right\} & = \emptyset.\\
    \left\{ (\theta+\omega) \times [c_2^{-}+g(\theta), c_2^{+}+g(\theta)] \right\} \cap \left\{(\theta+\omega) \times [e_2^{-}, e_2^{+}]\right\} & = \emptyset. \label{equation_structure I1_sync set Gs_2}
\end{align}
For each $i\in \{1,2\}$ and the points $c_i \in [c_i^-, c_i^+]$, by the definition of the sets $G_1$ and $G$ we have the bounds $c_1+g(\theta) > e_1$ and $c_2+g(\theta) < e_2$. Further, equations  \eqref{equation_structure I1_sync set Gs_1} and \eqref{equation_structure I1_sync set Gs_2} then gives us 
\begin{align*}
    \{(\theta+\omega) \times [c_1^{-}+g(\theta), c_1^{+}+g(\theta)] \} & \subset \{(\theta + \omega) \times [e_1^{+}, c_1^{+}]\} \\
    \{(\theta+\omega) \times [c_2^{-}+g(\theta), c_2^{+}+g(\theta)] \} & \subset \{(\theta + \omega) \times [c_1^{+}, e_2^{-}]\}.
\end{align*}
Due to the choice of $S_0$ (in equation \eqref{equation_preliminaries_f contraction S0 1}) and the property of contraction of $C_1$, as in definition \ref{definition_preliminaries_expanding contracting regions}, we obtain
\begin{align*}
    F^{M_0+S_0}\left((\theta+\omega) \times [c_1^{-}+g(\theta), c_1^{+}+g(\theta)] \right) & \subset \left((\theta + M_0\omega) \times C_1\right) \\
     F^{M_0+S_0}\left((\theta+\omega) \times [c_2^{-}+g(\theta), c_2^{+}+g(\theta)] \right) & \subset \left((\theta + M_0\omega) \times C_1\right).
\end{align*}
By the continuity of $F$ over the interval $G_1^{(i)}$ for each $i\in \{1,2\}$, this proves that
\begin{align*}
    \pi_2 \circ F^{2M_0+S_0}\left((G_1-M_0\omega) \times C_1 \right) & \subset C_1 \\
    \pi_2 \circ F^{2M_0+S_0}\left((G_1-M_0\omega) \times C_2 \right) & \subset C_1.
\end{align*}
The proof for the expanding regions under the backward iterations follows similarly.

\end{proof}
%------------------- PROP END -------------------

%------------------- LEMMA BEGIN -------------------
\begin{lemma}\label{lemma_structure I1_switch set S_2 size}
For each $\omega \in \mathbb{T}$ we have the estimate $|S_2| > \lambda^{-1/2-1/16}$ and $|S_3| > \lambda^{-1/2-1/32}$.
\end{lemma}
\begin{proof}
    By lemma \ref{lemma_preliminary estimates_|G| and |S|}, we have $|S| = d(\theta^{(1,1)}, \theta^{(1,2)}) > \lambda^{-1/2-1/32}$. Then, the proof follows directly from the definition of the sets $S_2$ and $S_3$.
\end{proof}
%------------------- LEMMA END -------------------

%------------------- PROP BEGIN -------------------
\begin{proposition}\label{proposition_structure I1_switch set Ss}
For each $\omega \in \mathcal{F}_0$, the set $S_1(\omega) \subset I_0\setminus I_1(\omega)$, as in definition \ref{definition_structure I1_switch set Ss}, satisfies the following conditions:

\begin{enumerate}[label=\alph*)]
\item \label{lemma item_structure I1_switch set Ss a} For each $\omega \in \mathcal{F}_0$,
\begin{align}\label{equation_structure I1_setting S_1}
   \left( \theta^{(1,1)}+ \lambda^{-1+1/64}, \theta^{(1,2)}-\lambda^{-1+1/64} \right) \subset S_1(\omega).
\end{align}

\item \label{lemma item_structure I1_switch set Ss b} For each $i \in \{1,2\}$, 
\begin{align*}
    F^{2M_0+S_0} \left( (S_1-M_0\omega) \times C_i \right) \subset (S_1+(M_0+S_0)\omega) \times C_{i+1}, \\
    F^{-2M_0-S_0} \left( (S_1+M_0\omega) \times E_i \right) \subset (S_1 -(M_0+S_0)\omega) \times E_{i+1}, 
\end{align*}
where $(i+1)$ denotes $(i+1) \pmod{2}$.

\item \label{lemma item_structure I1_switch set Ss c} For each $\theta \in I_0 \setminus(S_1\cup I_1)$ and $i\in\{1,2\}$, we have either
\begin{align*}
        \pi_2 \circ F^{2M_0+S_0}(\{\theta-M_0\omega\} \times C_i) & \subset C_{i} \\ \text{ or  } \ \ \
        \pi_2 \circ F^{2M_0+S_0}\{\theta-M_0\omega\} \times C_i) & \subset C_1.
\end{align*}

\item \label{lemma item_structure I1_switch set Ss d} For each $\theta \in I_0 \setminus(S_1\cup I_1)$ and $i\in\{1,2\}$, we have either
\begin{align*}
        \pi_2 \circ F^{-2M_0-S_0}(\{\theta+M_0\omega\} \times E_i) & \subset E_{i} \\ \text{ or  } \ \ \ 
        \pi_2 \circ F^{-2M_0-S_0}\{\theta+M_0\omega\} \times E_i) & \subset E_1.
\end{align*}
    
\end{enumerate}
\end{proposition}

\begin{proof}
Recall that for each $i,j \in \{1,2\}$, the critical sets are distributed such that $\theta^{(i,j)} \in I_0^{(i,j)}$, $|I_0^{(i,j)}| < \lambda^{-1+1/64}$ and  $I_1^{(i,j)}(\omega) \subset I_0^{(i,j)}$, due to propositions \ref{proposition_structure I0_I0 subsets} and \ref{proposition_structure I1_critical set I_1}.

\vspace{0.2cm}
\noindent\textit{Proof of \ref{lemma item_structure I1_switch set Ss a}:} 
By definition of $S_2$ (in definition \ref{definition_structure I1_switch set Ss}) and the estimates of $|I_1|$, it is clear that $S_2 \subset I_0\setminus I_1(\omega)$, for all $\omega \in \mathcal{F}_0$. Further, by the estimates of $|I_1|$, the largest open interval $S_1(\omega)$ in $I_0\setminus I_1(\omega)$ which contains $S_2$, must be non-empty and satisfy the condition
\begin{align*}
    \left( \theta^{(1,1)}+\lambda^{-1+1/64}, \theta^{(1,2)}-\lambda^{-1+1/64}\right) \subset S_1(\omega) \subset I_0\setminus I_1(\omega).
\end{align*}
\vspace{0.2cm}
\noindent\textit{Proof of \ref{lemma item_structure I1_switch set Ss b}:}
Let us fix a point $\theta \in S_2$. Since we have $(\theta-M_0\omega) \in \Theta_0$ and $(\theta-m\omega)_{m=0}^{2M_0+S_0} \cap I_1(\omega) = \emptyset$, it follows from $(C1)_0$ (in proposition \ref{proposition_block 0 sync' contraction}) that 
\begin{align*}
    F^{M_0}\left( \{\theta-M_0\omega\} \times C_i \right) & \subset \{\theta\} \times C_i \\ 
    F^{2M_0+S_0}\left( \{\theta-M_0\omega\} \times C_i \right) & \subset \{\theta+(M_0+S_0)\omega\} \times  (C_1\cup C_2),    
\end{align*}
where $i\in \{1,2\}$. Note that 
\begin{align}\label{equation_structure I1_subset S_eqn 1}
    F^{M_0+1}\left( \{\theta-M_0\omega\} \times [c_1^-, c_1^+] \right) & \subset \{\theta+\omega\} \times [c_1^- + g(\theta), c_1^+ + g(\theta)], \\
    F^{M_0+1}\left( \{\theta-M_0\omega\} \times [c_2^{-}, c_2^{+}]\right) & \subset \{\theta+\omega\} \times [c_2^- + g(\theta), c_2^+ + g(\theta)], 
\end{align}
and since $\displaystyle \theta \notin \bigcup_{i,j=1}^{2}I_0^{(i,j)}$, 
\begin{align}\label{equation_structure I1_subset S_eqn 3}
    [c_1^- + g(\theta), c_1^+ + g(\theta)] \cap (E_1\cup E_2) = \emptyset.\\
    [c_2^- + g(\theta), c_2^+ + g(\theta)] \cap (E_1 \cup E_2) = \emptyset. 
\end{align}
 Since $\theta \in S_2 \subset S$, by the definition of the set $S$ we can deduce that $c_1 + g(\theta) < e_1$ and $c_2 + g(\theta) < e_2$. Then, by equation \eqref{equation_structure I1_subset S_eqn 3} and the construction of the regions $\{C_i,E_i:i=1,2\}$ (in definition \ref{definition_preliminaries_expanding contracting regions}), we must have $c_1+g(\theta) \in [e_2^{+}, e_1^{-}]$ and $c_2+g(\theta) \in [e_1^{+}, e_2^{-}]$. This gives us
\begin{align*}
    [c_1^- + g(\theta), c_1^+ + g(\theta)] & \subset [e_2^+, e_1^-]\\
    [c_2^- + g(\theta), c_2^+ + g(\theta)] & \subset [e_1^+, e_2^-].
\end{align*}
Then, by the choice of $S_0$ in equation \eqref{equation_preliminaries_f contraction S0 2}, it follows that
\begin{align}\label{equation_structure I1_subset S_final 1}
    F^{2M_0+S_0}\left( \{\theta-M_0\omega\} \times C_1 \right) & \subset \{\theta+(M_0+S_0)\omega\} \times  C_2 \\
    F^{2M_0+S_0}\left( \{\theta-M_0\omega\} \times C_2 \right) & \subset \{\theta+(M_0+S_0)\omega\} \times  C_1. \label{equation_structure I1_subset S_final 2}
\end{align}
Since the equations \eqref{equation_structure I1_subset S_final 1} and \eqref{equation_structure I1_subset S_final 2} hold for all $\theta \in S_2$, and the set $S_1(\omega)$ is an interval that contains $S_2$, this proves that
\begin{align*}
    F^{2M_0+S_0}\left( (S_1-M_0\omega) \times C_1 \right) & \subset (S_1+(M_0+S_0)\omega) \times  C_2.     \\
    F^{2M_0+S_0}\left( (S_1-M_0\omega) \times C_1 \right) & \subset (S_1+(M_0+S_0)\omega) \times  C_2.  
\end{align*}
The proof for the backward iterations of the expanding regions follows similarly.

\vspace{0.2cm}
\noindent\textit{Proof of \ref{lemma item_structure I1_switch set Ss c} and \ref{lemma item_structure I1_switch set Ss d}:}
The set $I_1(\omega)$ consists of four intervals. Then, it follows from proposition \ref{proposition_structure I1_switch set Ss} \ref{lemma item_structure I1_switch set Ss a} that the set $I_0 \setminus (S_1\cup I_1)$ consists of exactly four intervals. For a sufficiently large $\lambda$, by the definition of the set $I_0^{(i,j)}$ and estimate $\displaystyle |g(1/2 \pm \lambda^{-1/4})| < \lambda^{-1/4+1/128}$, we can infer that $(1/2 \pm \lambda^{-1/4}) \notin \bigcup_{i,j=1}^{2} I_0^{(i,j)}$ and thus
\begin{align*}
    \{1/2-\lambda^{-1/4},1/2+\lambda^{-1/4}\} \cap (I_1\cup S_1) = \emptyset.
\end{align*}
 By the definition of the synchronisation sets $G_1, G_2$ and proposition \ref{proposition_structure I1_sync set Gs}, we can deduce that $I_0 \setminus (S_1\cup I_1)$ consists of an interval containing $G_1^{(1)}$ and another interval containing $G_2^{(2)}$. We now denote $I_0 \setminus (S_1\cup I_1)$ by the disjoint union
\begin{align}\label{equation_structure I1_switch set Ss_I0/I1 S1}
    I_0 \setminus (S_1\cup I_1) := \bigcup_{k=1}^{4} P_k,
\end{align}
where the interval $P_1$ contains $G^{(1)}_1$, $P_2$ contains $G^{(1)}_2$, $P_3$ contains a neighbourhood of $(1/2-\lambda^{-1/4})$ and $P_4$ contains a neighbourhood of $\displaystyle (1/2+\lambda^{-1/4})$. For any $\theta \in P_3\cup P_4$, it is clear by the disjointness condition $(\mathcal{F}1)_0$ that 
\begin{align*}
    F^{M_0+1}\left((\theta-M_0\omega) \times C_1 \right) & \subset \left(\{\theta+\omega\} \times [c_1^{-}+g(\theta), c_1^+ +g(\theta)]\right) \\
    F^{M_0+1}\left((\theta-M_0\omega) \times C_2 \right) & \subset \left(\{\theta+\omega\} \times [c_2^{-}+g(\theta), c_2^+ +g(\theta)]\right).
\end{align*}
We can choose a sufficiently large $\lambda$ and $\varepsilon >0$ such that $(1/2-\lambda^{-1/4}, 1/2-\lambda^{-1/4}+\varepsilon) \subset P_3$,  $(1/2+\lambda^{-1/4}-\varepsilon, 1/2+\lambda^{-1/4}) \subset P_4$, and for which we have the bound
\begin{align*}
    |g({\theta})| < \lambda^{-1/4+1/64} < \min\{|c_1^{-}-e_1^{+}|, |c_2^{-} - e_2^{+}|\}/2,
\end{align*}
for all ${\theta} \in (1/2-\lambda^{-1/4}, 1/2-\lambda^{-1/4}+\varepsilon) \cup (1/2+\lambda^{-1/4}-\varepsilon, 1/2+\lambda^{-1/4})$. Then, for any ${\theta} \in (1/2-\lambda^{-1/4}, 1/2-\lambda^{-1/4}+\varepsilon) \bigcup (1/2+\lambda^{-1/4}-\varepsilon, 1/2+\lambda^{-1/4})$, the rate of contraction in contracting regions and the choice of $S_0$ (in equations (\ref{setup_I_0 contraction charac} - \ref{equation_preliminaries_f contraction S0 1})) gives us 
\begin{align*}
    \pi_2\circ F^{2M_0+S_0}\left(({\theta}-M_0\omega) \times C_1 \right) & \subset C_1 \\
    \pi_2\circ F^{2M_0+S_0}\left(({\theta}-M_0\omega) \times C_2 \right) & \subset C_2.
\end{align*}
By $(C1)_0$ and the continuity of $F$ over the intervals $P_3$ and $P_4$, it then follows that 
\begin{align*}
    \pi_2\circ F^{2M_0+S_0}\left(({\theta}-M_0\omega) \times C_1 \right) & \subset C_1 \\
    \pi_2\circ F^{2M_0+S_0}\left(({\theta}-M_0\omega) \times C_2 \right) & \subset C_2,
\end{align*}
for all $\theta \in P_3\cup P_4$.  By proposition \ref{proposition_structure I1_sync set Gs} \ref{lemma item_structure I1_sync set Gs_c}, for each $\theta \in P_1\cup P_2$, 
\begin{align*}
    \pi_2 \circ F^{2M_0+S_0}(\{\theta-M_0\omega\} \times (C_1 \cup C_2)) \subset C_1.
\end{align*} 
This proves proposition \ref{proposition_structure I1_switch set Ss} \ref{lemma item_structure I1_switch set Ss c}. The proof for the iteration of the expanding regions follows similarly.

\end{proof}
%------------------- PROP END -------------------

%--------------------- Critical Set I_2 ---------------------------

\subsection{Structure of critical set $I_2$}\label{subsection_structure of critical set I2}

In this section, for each $\omega \in \tilde{\mathcal{F}}_1$, we construct a critical set $I_2(\omega)$ which satisfies definition \ref{definition_block n_Critical sets In}. Unlike the critical sets $I_0$ and $I_1$, there exists a dichotomy for the critical set $I_2(\omega)$, where  for different parameters $\omega \in \tilde{\mathcal{F}}_1$ the critical set $I_2(\omega)$ can either be empty or non-empty. In subsubsection \ref{subsubsection_empty vs nonempty I2}, this dichotomy will be derived as a characterisation of $\omega$. In subsubsection \ref{subsubsection_estimate nonempty I2}, we establish the structure of the critical set $I_2(\omega)$, when it is non-empty.  

\begin{definition}[Preliminary critical set]\label{definition_structure I2_ set I2 prime}
    Let $\omega \in \tilde{\mathcal{F}}_1$ and $M_1,M_0$ be the constants associated with $\omega$, as in definition \ref{definition_block n_disjointness condn}. Let $I_1(\omega)$ be the specific critical set associated with $\omega$ as defined in definition \ref{definition_structure I1_I1 definition}. We now define the following: 
    \begin{description}
    \item[$(\tilde{\mathcal{I}1)}_2$] We define the preliminary critical set $\Tilde{I}_{2}(\omega)$ as
    \begin{align*}
        \Tilde{I}_{2}(\omega) := \pi_1 \left( F^{M_{1}}(A_{1}) 
       \cap F^{-M_{1}}(B_{1}) \right),
    \end{align*}
    where 
    \begin{align*}
    A_1 & = \left( I_1-M_1\omega\right) \times (C_1\cup C_2) \\
    B_1 & = \left(I_1+M_1\omega\right)\times (E_1 \cup E_2).
    \end{align*}

    \item[$(\tilde{\mathcal{I}2)}_2$]  For each $i,j \in \{1,2\}$, let us define  
    $$\tilde{I}^{(i,j)}_{2}(\omega) := \Tilde{I}_{2}(\omega) \cap I^{(i,j)}_1(\omega).$$
\end{description}
\end{definition}
By the definition of the preliminary critical set $\tilde{I}_2(\omega)$ and the structure of the critical set $I_1(\omega)$ established in proposition \ref{proposition_structure I1_critical set I_1}, we infer that 
\begin{align*}
    \tilde{I_2}(\omega) = \bigcup_{i,j=1}^{2} \tilde{I}_2^{(i,j)}(\omega).
\end{align*}
We now define a specific critical set $I_2(\omega)$ using the preliminary critical set $\tilde{I}_2(\omega)$. Throughout the remainder of the article, the notation $I_2(\omega)$ is reserved exclusively for the sets introduced in the following definition.

\begin{definition}[Critical subsets]\label{definition_structure I2_critical set I2}
     For each $\omega \in \tilde{\mathcal{F}}_1$ and for each $i,j \in \{1,2\}$, let $I_2^{(i,j)}(\omega)$ denote the smallest closed interval in $\mathbb{T}$ such that $\tilde{I}_2^{(i,j)}(\omega) \subset I_2^{(i,j)}(\omega)$. Then, we fix the set $I_2$ to be
    $$I_2(\omega) := \bigcup_{i,j=1}^{2} I_2^{(i,j)}(\omega).$$
\end{definition}

It follows from definitions \ref{definition_structure I2_ set I2 prime} and \ref{definition_structure I2_critical set I2} that, the set $I_2(\omega)$ defined above is in fact a critical set as in definition \ref{definition_block n_Critical sets In}. 

%------------------- LEMMA BEGIN -------------------
\begin{lemma}\label{lemma_structure I2_match Ai Bi}
    For $\omega\in \tilde{\mathcal{F}}_1$, we have 
    \begin{align*}
        F^{M_1-M_{0}}(A_1) & \subset A_{0} \\
        F^{-M_1+M_{0}}(B_1) & \subset B_{0}.
    \end{align*}
\end{lemma}

\begin{proof}

By proposition \ref{proposition_structure I1_critical set I_1}, we have $I_1 \subset I_0$ and thus 
    \begin{align*}
        \pi_1 \circ F^{M_1-M_0}(A_1) = (I_1-M_0\omega) \subset (I_0-M_0\omega).
    \end{align*}
  Let us fix $(\theta, x)\in \left((I_1-M_1\omega)\times C_1\cup C_2\right)$. Using the disjointness conditions $(\mathcal{F}1)_1$ and $(\mathcal{F}2)_1$ in definition \ref{definition_block n_disjointness condn}, we can deduce that $\theta \in \Theta_0$ and $M_1$ is the smallest positive integer such that $\theta_{M_1} \in I_1$. Since $\theta_{M_1-M_0} \in I_0-M_0\omega$, it follows from $(C1)_1$ (in proposition \ref{proposition_block 0 contraction}) that $x_{M_1-M_0} \in C_1 \cup C_2$. By the continuity of the map $F$, this proves that $F^{M_1-M_0}(A_1) \subset A_0$. The proof for $F^{-M_1+M_0}(B_1) \subset B_0$ follows similarly using $(E1)_1$ (in proposition \ref{proposition_block 0 expansion}).

\end{proof}
%------------------- LEMMA END -------------------

For each $i\in \{1,2\}$, we denote the components of $A_1$ and $B_1$ by
\begin{align*}
    A_1^{(i)} & := \left(I_1-M_1\omega\right) \times C_i \\
    B_1^{(i)} & := \left(I_1+M_1\omega\right) \times E_i.
\end{align*}
For each $\omega \in \tilde{\mathcal{F}}_1$ and $i\in \{1,2\}$, we define the curves $\phi_{1,i}^{\pm}(\cdot,\omega), \psi_{1,i}^{\pm}(\cdot,\omega): I_1(\omega) \rightarrow \mathbb{T}$, by
\begin{align*}
    \phi_{1,i}^{\pm}(\theta, \omega) & = \pi_2 \circ F^{M_1+1}(\theta-M_1\omega, c_i^{\pm}) \\
    \psi_{1,i}^{\pm}(\theta, \omega) & = \pi_2 \circ F^{-M_1+1}(\theta + M_1\omega, e_i^{\pm}).
\end{align*}
Whenever convenient, for a fixed $\omega \in \tilde{\mathcal{F}}_1$, we use the notation $\phi_{1,i}^{\pm}(\theta) = \phi_{1,i}^{\pm}(\theta, \omega)$ and $\psi_{1,i}^{\pm}(\theta) = \psi_{1,i}^{\pm}(\theta, \omega)$.

\begin{lemma}
For any $\omega \in \tilde{\mathcal{F}}_1$ and $i \in \{1,2\}$, 
\begin{align*}
    F^{-M_1+1}(B_1^{(i)}) & = \left\{(\theta,x): \theta \in (I_1+\omega),\ x\in \left[\psi_{1,i}^{-}(\theta), \psi_{1,i}^{+}(\theta)\right] \right\} \\
    F^{M_1+1}(A_1^{(i)}) & = \left\{(\theta,x): \theta \in (I_1+\omega),\ x\in \left[\phi_{1,i}^{-}(\theta), \phi_{1,i}^{+}(\theta)\right] \right\}.
\end{align*}
\end{lemma}

\begin{proof}
For each $\theta \in I_1$, by the continuity of $F$ and since $f_{\theta}(x)$ is an orientation preserving diffeomorphism on the torus, we obtain
   \begin{align*}
       F^{-M_1+1}\left( \{\theta + M_1\omega \} \times [e_i^{-}, e_i^{+}]\right) & \subset \{\theta + \omega \} \times [\psi_{1,i}^{-}(\theta), \psi_{1,i}^{+}(\theta)] \\
       F^{M_1+1}\left( \{\theta - M_1\omega \} \times [c_i^{-}, c_i^{+}]\right) & \subset \{\theta + \omega \} \times [\phi_{1,i}^{-}(\theta), \phi_{1,i}^{+}(\theta)],
   \end{align*}
for each $i \in \{1,2\}$.
    
\end{proof}
%------------------- LEMMA END -------------------

The preceding lemma shows that the geometry of the critical set $I_2(\omega)$ is closely connected to the structure of the curves $\phi_{1,i}^{\pm}$ and $\psi_{1,i}^{\pm}$, where $i \in \{1,2\}$. In the following propositions we establish the necessary estimates for these curves.

%------------------- PROP BEGIN -------------------
\begin{proposition}\label{proposition_structure I2_backward iterate curves}
   For all sufficiently large $\lambda$, $\omega \in \tilde{\mathcal{F}}_1$ and for each $i \in \{1,2\}$, the maps $\psi_{1,i}^{\pm}(\theta)$ satisfy the following conditions:  
\begin{enumerate}[label=\alph*)]
\item \label{lemma item_structure I2_backward iterate curves_a} For all $\theta \in I_1$, we have $[\psi_{1,i}^{-}(\theta), \psi_{1,i}^{+}(\theta)] \subset E_1 \cup E_2.$
\item \label{lemma item_structure I2_backward iterate curves_b} For all $\theta \in I_1$,  
\begin{align*}
   \alpha_c^{b(M_1-1)}|e_i^+ - e_i^-| & \leq |\psi_{1,i}^{+}(\theta) - \psi_{1,i}^{-}(\theta)| \leq \left(\alpha_c^{1-(b+1)\beta_{0}}\right)^{M_1-1}|e_i^{+} -e_i^{-}|.
\end{align*}
\item \label{lemma item_structure I2_backward iterate curves_c}  For all $\theta\in I_1$,
$$|\partial_{\theta}\psi_{1,i}^{\pm}(\theta)|< 1.$$
\item \label{lemma item_structure I2_backward iterate curves_d}  For any $\omega \in \tilde{\mathcal{F}}_1$ and $\theta \in I_1$,  
$$ |\partial_{\omega}\psi_{1,i}^{\pm}(\theta, \omega)| < 1. $$
\end{enumerate}
\end{proposition}

\begin{proof}

%\textit{Proof of \ref{lemma item_structure I2_backward iterate curves_a}:} 
Let us fix a point $\theta \in I_1(\omega)$. By the disjointness condition $(\mathcal{F}1-2)_1$ (in definition \ref{definition_block n_disjointness condn}), it follows that $(\theta+M_1\omega) \in \Theta_0$ and 
\begin{align*}
    (\theta+k\omega)_{k=1}^{M_1} \cap \left(\bigcup_{l=0}^{1} \bigcup_{m=1}^{M_1+S_0}(I_l-m\omega) \right) = \emptyset.
\end{align*} 
Then, by $(E1)_0$ (in proposition \ref{proposition_block 0 expansion}), we obtain 
\begin{align*}
    F^{-M_1+1}\left( (\theta+M_1\omega) \times (E_1\cup E_2) \right) \subset (\theta+\omega) \times (E_1\cup E_2).
\end{align*}
By the continuity of $F^{-1}$ and since $f^{-1}_{\theta}$ is orientation preserving for each $\theta \in \mathbb{T}$, this proves that $$\bigcup_{i=1}^{2}[\psi_{1,i}^{-}(\theta), \psi_{1,i}^{-}(\theta)] \subset (E_1\cup E_2).$$
%\textit{Proof of \ref{lemma item_structure I2_backward iterate curves_b}:}  
For each $x\in (E_1\cup E_2)$, by $(E3)_0$ (in proposition \ref{proposition_block 0 expansion}), we have the estimate
    \begin{align*}
        \alpha_c^{b(M_1-1)} & \leq |\partial_x \left(\pi_2 \circ F^{-M_1+1}(\theta+M_1\omega, x)\right)| \leq (\alpha_c^{1-(1+b)\beta_0})^{M_1-1}. 
    \end{align*}
This proves that, 
\begin{align*}
    \alpha_c^{b(M_1-1)}|e_i^+ - e_i^-| & \leq |\psi_{1,i}^{+}(\theta)-\psi_{1,i}^{-}(\theta)| \leq (\alpha_c^{1-(1+b)\beta_0})^{M_1-1}|e_i^+- e_i^-|.
\end{align*}
%\textit{Proof of \ref{lemma item_structure I2_backward iterate curves_c} and \ref{lemma item_structure I2_backward iterate curves_d}:}  
 From the definition of the curves $\psi_{1,i}^{\pm}$, we can derive the expressions: 
 \begin{align*}
     \partial_{\theta} \psi_{1,i}^{\pm}(\theta, \omega)
    & = \sum_{k=1}^{M_1-1}\left((-g'(\theta+k\omega)) \prod_{j=1}^{k}(f^{-1})'(f(x_{-M_1+j})) \right) \text{ and } \\
    \partial_{\omega}\psi_{1,i}^{\pm}(\theta, \omega) & =  \sum_{k=1}^{M_1-1}\left((-kg'(\theta+ k\omega)) \prod_{j=1}^{k}(f^{-1})'(f(x_{-M_1+j})) \right).
\end{align*}
For every $x \in E_1 \cup E_2$, by disjointness condition $(\mathcal{F}2)_1$ and $(E2)_0$ (in proposition \ref{proposition_block 0 expansion}), we obtain the following estimates for any $k \in \{1,\dots ,(M_1-1)\}$:
\begin{align*}
    \left|\prod_{j=1}^{k}(f^{-1})'(f(x_{-M_1+j})) \right| < (\alpha_c^{1-(1+b)\beta_0})^{k}.
\end{align*}
Since $I_1 \subset I_0$, by the disjointness condition $(\mathcal{F}1)_0$, we have the bound $|g'(\theta+k\omega)|< \lambda^{-1/4+1/128}$ for all $k \in \{1,\dots, M_0\}$, and the bound $|g'(\theta+k\omega )|< \lambda^{1+1/128}$ for all $k \in \{M_0+1, \dots ,M_1\}$. Using these bounds, we derive the estimates
\begin{align*}
    |\partial_{\theta}\psi_{1,i}^{\pm}(\theta)| 
    & \leq \sum_{k=1}^{M_0} \lambda^{-1/4+1/128}(\alpha_c^{1-(1+b)\beta_0})^k + \sum_{k=M_0+1}^{M_1} \lambda^{1+1/128}(\alpha_c^{1-(1+b)\beta_0})^k.
\end{align*}
The constant $M_0 = \lambda^{1/16}$ and $\beta_0$ is small enough that $\beta_0 < 1/(4(b+1))$, which gives us $\alpha_c^{1-(1+b)\beta_0} < \alpha_c^{1/2}$. Then, for a big enough $\lambda$ with respect to constant $\alpha_c$, we obtain
\begin{align*}
    |\partial_{\theta}\psi_{1,i}^{\pm}(\theta)| 
    & < \lambda^{-1/4+1/128} \sum_{k=1}^{M_0} \alpha_c^{k/2} + \lambda^{1+1/128}\alpha_c^{M_0/2} \sum_{k=1}^{M_1-M_0} \alpha_c^{k/2} \\
    & < \frac{\alpha_c^{1/2}}{1-\alpha_c^{1/2}} \left(\lambda^{-1/4+1/128} + \lambda^{1+1/128}\alpha_c^{M_0/2}\right) < 1. 
    \stepcounter{equation}\tag{\theequation} \label{equation_backward curve I2_theta slope}
\end{align*}
Similarly, we can choose $\lambda$ big enough such that   
\begin{align*}
 |\partial_{\omega}\psi_{1,i}^{\pm}(\theta, \omega)| 
    & < \sum_{k=1}^{M_0} k\lambda^{-1/4+1/128} \alpha_c^{k/2} + \sum_{k=M_0+1}^{M_1} k\lambda^{1+1/128} \alpha_c^{k/2} \\
    & < \sum_{k=1}^{M_0} k\lambda^{-1/4+1/128} \alpha_c^{k/2} + \alpha_c^{M_0/2}\sum_{k=1}^{\infty} (M_0+k)\lambda^{1+1/128} \alpha_c^{k/2} < 1.
    \stepcounter{equation}\tag{\theequation} 
    \label{equation_backward curve I2_theta w slope}
\end{align*}
    
\end{proof}
%------------------- PROP END -------------------

%------------------- PROP BEGIN -------------------
\begin{proposition}\label{proposition_structure I2_forward iterate curves}
 For all sufficiently large $\lambda$, $\omega \in \tilde{\mathcal{F}}_1$ and for each $i \in \{1,2\}$, the maps $\phi_{1,i}^{\pm}(\theta)$ satisfy the following conditions:
\begin{enumerate}[label=\alph*)]
\item \label{lemma item_structure I2_forward iterate curves_a} For all $\theta \in I_1$,
\begin{align*}
[\pi_2 \circ F^{M_1} (\theta- M_1\omega, c_i^{-}),\pi_2 \circ F^{M_1}(\theta-M_1\omega, c_i^{+})] & \subset C_1 \cup C_2. 
\end{align*}
\item \label{lemma item_structure I2_forward iterate curves_b} For all $\theta \in I_1$, 
\begin{align*} 
\alpha_c^{b(M_1+1)}|c_i^+ - c_i^-| & \leq |\phi_{1,i}^+(\theta) - \phi_{1,i}^-(\theta)| \leq (\alpha_c^{1-(1+b)\beta_0})^{M_1+1}|c_i^{+}-c_i^{-}|.
\end{align*}

\item  \label{lemma item_structure I2_forward iterate curves_c} For all $\theta \in \bigcup_{i'=1}^{2} I^{(i',1)}_1$,
\begin{align*}
     -\lambda^{1+1/64} & < \partial_{\theta} \phi_{1,i}^{\pm}(\theta) < -\lambda^{1-1/64}, 
\end{align*}
and for all $\theta \in \bigcup_{i'=1}^{2} I^{(i',2)}_1$,
\begin{align*}
    \lambda^{1-1/64} & < \partial_{\theta} \phi_{1,i}^{\pm}(\theta) < \lambda^{1+1/64}.
\end{align*}

\item  \label{lemma item_structure I2_forward iterate curves_d} For any $\omega \in \tilde{\mathcal{F}}_1$ and $\theta \in I_1(\omega)$, we have 
\begin{align*}
    |\partial_{\omega}\phi_{1,i}^{\pm}(\theta, \omega)| < 1.
\end{align*}
\end{enumerate}
\end{proposition}

\begin{proof}
    %\textit{Proof of \ref{lemma item_structure I2_forward iterate curves_a}:} 
   Let us fix a point $\theta \in I_1(\omega)$. By the disjointness condition $(\mathcal{F}1-2)_1$ (in definition \ref{definition_block n_disjointness condn}), it follows that $(\theta - M_1\omega) \in \Theta_0$ and  
\begin{align*}
    (\theta-k\omega)_{k=1}^{M_1} \cap \left(\bigcup_{l=0}^{1} \bigcup_{m=1}^{M_1+S_0-1}(I_l+m\omega) \right) = \emptyset.
\end{align*} 
Then, by $(C1)_0$ (in proposition \ref{proposition_block 0 contraction}), we obtain
\begin{align*}
    F^{M_1-1}\left( (\theta-M_1\omega) \times (C_1\cup C_2) \right) \subset \{\theta\} \times (C_1\cup C_2).
\end{align*}
By the continuity of $F$ and since $f_{\theta}$ is orientation preserving for each $\theta \in \mathbb{T}$, this proves that
$$\pi_2\circ F^{M_1-1}\left( (I_1-M_1\omega) \times (C_1\cup C_2) \right) \subset (C_1\cup C_2).$$
For each $x\in C_1\cup C_2$, by $(C3)_0$ (in proposition \ref{proposition_block 0 contraction}), we have the estimate
    \begin{align*}
        \alpha_c^{b(M_1+1)} & \leq |\partial_x \left(\pi_2 \circ F^{M_1+1}(\theta-M_1\omega, x)\right)| \leq (\alpha_c^{1-(1+b)\beta_0})^{M_1+1}. 
    \end{align*}
This proves that, 
\begin{align*}
    \alpha_c^{b(M_1+1)}|c_i^+ - c_i^-| & \leq |\phi_{1,i}^{+}(\theta)-\phi_{1,i}^{-}(\theta)| \leq (\alpha_c^{1-(1+b)\beta_0})^{M_1+1}|c_i^+- c_i^-|.
\end{align*}
 From the definition of the curves $\phi_{1,i}^{\pm}$, we can derive the expression 
\begin{align*}
      \partial_{\theta} \phi_{1,i}^{\pm}(\theta) & = g'(\theta) + \sum_{k=1}^{M_1}\left( g'(\theta-k\omega)\prod_{j=0}^{k-1}f'(x_{M_1-j}) \right)\\
      \partial_{\omega}\phi_{1,i}^{\pm}(\theta, \omega) & =  \sum_{k=1}^{M_1}\left( (-k)g'(\theta-k\omega) \prod_{j=0}^{k-1} f'(x_{M_1-j}) \right).
\end{align*}
For every $x \in C_1 \cup C_2$, by disjointness condition $(\mathcal{F}2)_1$ and $(C2)_0$ (in proposition \ref{proposition_block 0 contraction}), we obtain the following estimates for any $k \in \{1,\dots ,M_1\}$:
\begin{align*}
    \left| \prod_{j=0}^{k-1} f'(x_{M_1-j}) \right| < (\alpha_c^{1-(1+b)\beta_0})^k.
\end{align*}
Since $I_1 \subset I_0$, by the disjointness condition $(\mathcal{F}1)_0$, we have the bound $|g'(\theta-k\omega)|< \lambda^{-1/4+1/128}$ for all $k \in \{1,\dots, M_0\}$, and the bound $|g'(\theta-k\omega )|< \lambda^{1+1/128}$ for all $k \in \{M_0+1, \dots ,M_1\}$. Using these bounds, we derive the estimates
\begin{align*}
    \left|\sum_{k=1}^{M_0} g'(\theta-k\omega) \prod_{j=0}^{k-1}f'(x_{M_1-j}) \right| & < \lambda^{-1/4+1/128} \frac{(\alpha_c^{1-(1+b)\beta_0})}{1-(\alpha_c^{1-(1+b)\beta_0})}\\
    \left|\sum_{k=M_0+1}^{M_1} g'(\theta-k\omega) \prod_{j=0}^{k-1}f'(x_{M_1-j}) \right| 
    & < \lambda^{1+1/128} (\alpha_c^{1-(1+b)\beta_0})^{M_0} \frac{\alpha_c^{1-(1+b)\beta_0}}{1-\alpha_c^{1-(1+b)\beta_0}}. 
\end{align*}
The constant $M_0 = \lambda^{1/16}$ and $\beta_0$ is small enough that $\beta_0 < 1/(4(b+1))$, which gives us $\alpha_c^{1-(1+b)\beta_0} < \alpha_c^{1/2}$. Then, for a big enough $\lambda$ with respect to constant $\alpha_c$, it is clear that
\begin{align*}
    \left|\sum_{k=1}^{M_1} g'(\theta-k\omega) \prod_{j=0}^{k-1}f'(x_{M_1-j}) \right| < \lambda^{1/2}.
\end{align*}
For every $\theta \in \bigcup_{i=1}^{2}I_1^{(i,2)}$, since $I_1^{(i,2)} \subset I_0^{(i,2)}$, by proposition \ref{proposition_structure I0_I0 subsets} \ref{lemma item_structure I0_I0 components_b} we have $\lambda^{1-1/128} < g'(\theta) < \lambda^{1+1/128},$ and thus
\begin{align*}
    \lambda^{1-1/128} - \lambda^{1/2} & < \partial_{\theta} \phi_{1,i}^{\pm}(\theta) < \lambda^{1+1/128} + \lambda^{1/2} \\
     \implies \lambda^{1-1/64} & < \partial_{\theta} \phi_{1,i}^{\pm}(\theta) < \lambda^{1+1/64}.
    \stepcounter{equation}\tag{\theequation} \label{equation_forward curve I2_theta slope}
\end{align*}
Similarly, for every $\theta \in \bigcup_{i=1}^{2}I_1^{(i,1)}$, since $I_1^{(i,2)} \subset I_0^{(i,2)}$, it follows that
\begin{align*}
   -\lambda^{1+1/64} & < \partial_{\theta} \phi_{1,i}^{\pm}(\theta) < -\lambda^{1-1/64}.
    \stepcounter{equation}\tag{\theequation} \label{equation_forward curve I2_theta slope}
\end{align*}
Similarly, we can choose $\lambda$ big enough such that for each $\theta \in I_1$ and $\omega \in \tilde{\mathcal{F}}_1$,  
\begin{align*}
    |\partial_{\omega}\phi_{1,i}^{\pm}(\theta, \omega)| 
    & < \left|\sum_{k=1}^{M_0} (-k\lambda^{-1/4+1/128})\alpha_c^{k/2}\right| + \left|\sum_{k=M_0+1}^{M_1} (-k\lambda^{1+1/128})\alpha_c^{k/2}\right| < 1.
\end{align*}

\end{proof}
%------------------- PROP END -------------------

In the following propositions, we focus on parameters $\omega \in \mathcal{F}_1$ and establish stronger results regarding the structure of the iterated curves $\psi_{1,i}^{\pm}$ and $\phi_{1,i}^{\pm}$. Recall that $\mathcal{F}_1 \subset \tilde{\mathcal{F}}_1$. 

%------------------- PROP BEGIN -------------------
\begin{proposition}\label{proposition_structure I2_backward iterate curves sync}
Assume there exists a non-empty set $\bar{G} \subset I_0 \setminus 3I_1(\omega)$ and an associated set $\hat{G}$ as in definition \ref{definition_block n_sync condn set_g1 g2}.
     Then, for all sufficiently large $\lambda$ and $\omega \in \mathcal{F}_1(f,\lambda,{\bar{G}})$, we have the following conditions: 
\begin{enumerate}[label=\alph*)]
\item \label{lemma item_structure I2_backward iterate curves sync_a} For each $i,j \in \{1,2\}$, there exists $k \in \{1,2\}$ such that
$$\bigcup_{l=1}^{2}[\psi_{1,l}^{-}(\theta), \psi_{1,l}^{+}(\theta)] \subset E_k, \ \ \ \text{for all } \theta \in I_1^{(i,j)}.$$

\item \label{lemma item_structure I2_backward iterate curves sync_b}  For all $\theta \in I_1$, 
\begin{align*}
\alpha_c^{b(M_1-1)}|e_1^+ - e_1^-| \leq \max_{i,j \in \{1,2\}} |\psi_{1,i}^{+}(\theta) - \psi_{1,j}^{-}(\theta)| \leq \left(\alpha_c^{1-(b+1)\beta_{0}} \right)^{M_1-M_0^{60}-1}.
\end{align*}

\end{enumerate}
\end{proposition}

\begin{proof}
%\textit{Proof of \ref{lemma item_structure I2_backward iterate curves sync_a}:} 
 Let us fix $i,j \in \{1,2\}$ and $\theta \in I_1^{(i,j)}$. By the disjointness condition $(\mathcal{F}1)_1$, we have $(\theta+M_1\omega) \in \Theta_1 \subset \Theta_0$. Since $\omega \in \mathcal{F}_1(f,\lambda, \bar{G})$, by the condition $(\mathcal{G})$ (in definition \ref{definition_block n_sync condn parameter}), there exists a constant $N_s$ such that $N_s < M_0^{60}<M_1$ and $(\theta + (M_1-N_s+M_0+S_0)\omega) \in \bar{G}$. 
 The sets $\hat{G}$ satisfies the condition $(g1)$ (in definition \ref{definition_block n_sync condn set_g1 g2}) which is that $B_{|I_1(\omega)|}(\bar{G}) \subset \hat{G}$ and this gives us $(I_1^{(i,j)}+(M_1-N_s+M_0+S_0)\omega) \subset \hat{G}$. Then, it follows from $(\mathcal{F}1)_0$ and $(E1)_0$ (in proposition \ref{proposition_block 0 expansion}) that 
 \begin{align*}
     F^{-N_s+2M_0+S_0}((I_1^{(i,j)}+M_1\omega) \times (E_1 \cup E_2)) \subset (\hat{G}+M_0\omega) \times (E_1 \cup E_2).
 \end{align*}
Since $\hat{G}$ satisfies condition $(g3)$ (in definition \ref{definition_block n_sync condn set_g1 g2}), we obtain 
 \begin{align}\label{equation_structure I2_backward iterate curves sync_1}
     F^{-N_s}((I_1^{(i,j)}+M_1\omega) \times (E_1 \cup E_2)) & \subset (I_1+(M_1-N_s)\omega) \times E_1 \\
     & \subset (I_0-(M_0+S_0)\omega) \times E_1.\label{equation_structure I2_backward iterate curves sync_1 0}
 \end{align}
From the above, we can also observe that $(I_1^{(i,j)}+(M_1-N_s)\omega) \subset \Theta_0$. Then, since
$$(I_1^{(i,j)}+\omega) \cap \left(\bigcup_{m=0}^{M_0+S_0-1}(I_0-m\omega)\right) = \emptyset,$$ 
it follows from $(Se1)_0$ (in proposition \ref{proposition_Block 0 sync expansion}) and continuity of $F$ that, there exists $k\in \{1,2\}$ such that
\begin{align*}
    \pi_2 \circ F^{-M_1+1}\left( (I_1^{(i,j)}+M_1\omega) \times (E_1\cup E_2)\right) & \subset E_k. 
\end{align*}
This proves that $\bigcup_{l=1}^{2} [\psi_{1,l}^{-}(\theta), \psi_{1,l}^{+}(\theta)] \subset E_k,$ for all $\theta \in I_1^{(i,j)}$.

%\textit{Proof of \ref{lemma item_structure I2_backward iterate curves sync_b}:} 
For fixed $i',j' \in \{1,2\}$, by $(E2)_0$ we obtain the estimate
\begin{align*}
    \left|\partial_x \left(\pi_2 \circ F^{-M_1+N_s+1}(\theta+(M_1-N_s)\omega, x)\right)\right| < (\alpha_c^{1-(1+b)\beta_0})^{M_1-N_s-1},
\end{align*}
for any $\theta \in I_1^{(i',j')}$ and $x\in E_1\cup E_2$. Since $F$ is a continuous and orientation preserving map, it follows from equation \eqref{equation_structure I2_backward iterate curves sync_1} that
\begin{align*}
   \alpha_c^{b(M_1-1)}|e_1^{+}-e_1^{-}| < \max_{i,j \in \{1,2\}} |\psi_{1,i}^{-}(\theta) - \psi_{1,j}^{+}(\theta)| & < (\alpha_c^{1-(1+b)\beta_0})^{M_1-N_s-1}, 
\end{align*}
for every $\theta \in I_1^{(i',j')}$ and since $N_s< M_0^{60}$, this gives us 
\begin{align*}
   \alpha_c^{b(M_1-1)}|e_1^{+}-e_1^{-}| < \max_{i,j \in \{1,2\}} |\psi_{1,i}^{-}(\theta) - \psi_{1,j}^{+}(\theta)| & < (\alpha_c^{1-(1+b)\beta_0})^{M_1-M_0^{60}-1}. 
\end{align*}
for every $\theta \in I_1$.
\end{proof}
%------------------- PROP END -------------------

%------------------- PROP BEGIN -------------------
\begin{proposition}\label{proposition_structure I2_forward iterate curves sync}
        Assume there exists a non-empty set $\bar{G} \subset I_0 \setminus 3I_1(\omega)$ and an associated set $\hat{G}$ as in definition \ref{definition_block n_sync condn set_g1 g2}. Then, for all sufficiently large $\lambda$ and $\omega \in \mathcal{F}_1(f,\lambda,{\bar{G}})$, we have the following conditions:    
\begin{enumerate}[label=\alph*)]
\item \label{lemma item_structure I2_forward iterate curves sync_a} For each $i,j \in \{1,2\}$, there exists $k \in \{1,2\}$ such that 
\begin{align*}
\bigcup _{l=1}^{2} [\pi_2 \circ F^{M_1} (\theta- M_1\omega, c_l^{-}),\pi_2 \circ F^{M_1}(\theta-M_1\omega, c_l^{+})] & \subset C_k, \ \text{ for all } \ \theta \in I_1^{(i,j)}. 
\end{align*}
\item \label{lemma item_structure I2_forward iterate curves sync_b} For each $\theta\in I_1$, 
\begin{align*}
    \alpha_c^{b(M_1+1)}|c_1^{+}-c_1^{-}| \leq \max_{i,j \in \{1,2\} }|\phi_{1,i}^{+}(\theta)-\phi_{1,j}^{-}(\theta)| \leq (\alpha_c^{1-(1+b)\beta_0})^{M_1-M_0^{60}+1}.
\end{align*}
\end{enumerate}
\end{proposition}

\begin{proof}
%\textit{Proof of \ref{lemma item_structure I2_forward iterate curves sync_a}:} 
 Let us fix $i,j \in \{1,2\}$ and $\theta \in I_1^{(i,j)}$. By the disjointness condition $(\mathcal{F}1)_1$, we have $(\theta-M_1\omega) \in \Theta_1 \subset \Theta_0$. Since $\omega \in \mathcal{F}_1(f,\lambda, \bar{G})$, by the condition $(\mathcal{G})$ (in definition \ref{definition_block n_sync condn parameter}), there exists a constant $N_s$ such that $N_s < M_0^{60}<M_1$ and $(\theta - (M_1-N_s+M_0+S_0)\omega) \in \bar{G}$. 
 The sets $\hat{G}$ satisfies the condition $(g1)$ (in definition \ref{definition_block n_sync condn set_g1 g2}) which is that $B_{|I_1(\omega)|}(\bar{G}) \subset \hat{G}$ and this gives us $(I_1^{(i,j)}-(M_1-N_s+M_0+S_0)\omega) \subset \hat{G}$. Then, by $(\mathcal{F}1)_0$ and $(C1)_0$ (in proposition \ref{proposition_block 0 contraction}), it follows that 
 \begin{align*}
     F^{N_s-2M_0-S_0}((I_1^{(i,j)}-M_1\omega) \times (C_1 \cup C_2)) \subset (\hat{G}-M_0\omega) \times (C_1 \cup C_2).
 \end{align*}
Since $\hat{G}$ satisfies condition $(g3)$ (in definition \ref{definition_block n_sync condn set_g1 g2}), we obtain 
 \begin{align}
     F^{N_s}((I_1^{(i,j)}-M_1\omega) \times (C_1 \cup C_2)) 
     & \subset (I_1-(M_1-N_s)\omega) \times C_1 \\
     & \subset (I_0+(M_0+S_0)\omega) \times C_1.\label{equation_structure I2_forward iterate curves sync_1}
 \end{align}
From the above, we can also observe that $(I_1^{(i,j)}-(M_1-N_s)\omega) \subset \Theta_0$. Then, since
$$I_1^{(i,j)} \cap \left(\bigcup_{m=1}^{M_0+S_0-1}(I_0+m\omega)\right) = \emptyset,$$ 
it follows from $(Sc1)_0$ (in proposition \ref{proposition_Block 0 sync contraction}) and continuity of $F$ that, there exists $k\in \{1,2\}$ such that
\begin{align*}
    \pi_2 \circ F^{M_1}\left( (I_1^{(i,j)}-M_1\omega) \times (C_1\cup C_2)\right) & \subset C_k.
\end{align*}
Since $N_s< M_0^{60}$, this proves that $\bigcup _{l=1}^{2} [\pi_2 \circ F^{M_1} (\theta- M_1\omega, c_l^{-}),\pi_2 \circ F^{M_1}(\theta-M_1\omega, c_l^{+})] \subset C_k,$ for all $\theta \in I_1^{(i,j)}$.

%\textit{Proof of \ref{lemma item_structure I2_forward iterate curves sync_b}:} 
For fixed $i',j' \in \{1,2\}$, by $(C2)_0$ we obtain the estimate
\begin{align*}
    \left|\partial_x \left(\pi_2 \circ F^{M_1-N_s+1}(\theta-(M_1-N_s)\omega, x)\right)\right| < (\alpha_c^{1-(1+b)\beta_0})^{M_1-N_s+1},
\end{align*}
for any $\theta \in I_1^{(i',j')}$ and $x\in C_1\cup C_2$. Since $F$ is continuous and orientation preserving, it follows from equation \eqref{equation_structure I2_forward iterate curves sync_1} that 
\begin{align*}
   \alpha_c^{b(M_1+1)}|c_1^{+}-c_1^{-}| < \max_{i,j \in \{1,2\}} |\phi_{1,i}^{-}(\theta) - \phi_{1,j}^{+}(\theta)| & < (\alpha_c^{1-(1+b)\beta_0})^{M_1-N_s+1}, 
\end{align*}
for every $\theta \in I_1^{(i',j')}$, and since $N_s< M_0^{60}$, this proves that  
\begin{align*}
   \alpha_c^{b(M_1-1)}|c_1^{+}-c_1^{-}| < \max_{i,j \in \{1,2\}} |\phi_{1,i}^{-}(\theta) - \phi_{1,j}^{+}(\theta)| & < (\alpha_c^{1-(1+b)\beta_0})^{M_1-M_0^{60}+1}, 
\end{align*}
for every $\theta \in I_1$.

\end{proof}
%------------------- PROP END -------------------

%%%%%%%%%%%%%%%%% EMPTY vs NON_EMPTY DICHOTOMY for I2 %%%%%%%%%%%%%%%%%%%%%%%%%
\subsubsection{The empty and non-empty cases for $I_2$}\label{subsubsection_empty vs nonempty I2}

We now characterize the conditions on $\omega \in \tilde{\mathcal{F}}_1$ under which the critical set $I_2(\omega)$ is either empty or non-empty. It is also worth noting that the propositions established below will also play a crucial role in the construction of the parameter sets $\Omega_1$ and $\Omega_2$ in Section \ref{section_parameter space}.

\begin{definition}[Non-emptiness condition]\label{definition_structure I2_SNA character of parameters}
For each $i,j \in \{1,2\}$, let $\tilde{\mathcal{F}}_1^{(i,j)}(M_1,M_0)$ denote the set of all $\omega \in \tilde{\mathcal{F}}_1(M_1, M_0)$ for which 
\begin{align*}
    \pi_2 \circ F^{M_1-M_0}\left( \left( I^{(i,j)}_1-M_1\omega\right) \times C_k \right) & \subset C_i, \text{ for some } k\in \{1,2\},  \\
        \text{and  } \pi_2 \circ F^{-M_1+M_0}\left( \left( I^{(i,j)}_1+M_1\omega\right) \times E_k \right) & \subset E_i, \text{ for some } k\in \{1,2\}. 
\end{align*}
Further, we define 
$$ \tilde{\mathcal{F}}^{(i,j)}_1 := \left\{ \omega\in \mathbb{T}: \omega \in \tilde{\mathcal{F}}^{(i,j)}_1(M_1, M_0) \text{ for some } M_1 \in [\tilde{M}_1, 2\tilde{M}_1]\right\}.$$
\end{definition}

The following proposition establishes the dichotomy in the structure of the critical set $I_2(\omega)$.

%------------------- PROP BEGIN -------------------
\begin{proposition}\label{proposition_structure I2_I2 empty nonempty}
    For each $\omega \in \tilde{\mathcal{F}}_1$ and $i,j \in \{1,2\}$, we have $I_2^{(i,j)}(\omega) \neq \emptyset$ if and only if $\omega \in \tilde{\mathcal{F}}_1^{(i,j)}$. Consequently, $I_2(\omega) \neq \emptyset$ if and only if $\displaystyle \omega \in \bigcup_{i,j=1}^{2} \tilde{\mathcal{F}}_1^{(i,j)}.$
\end{proposition}

\begin{proof}
Let us fix $\omega \in \tilde{\mathcal{F}}_1^{(i,j)}$, for some fixed $i,j\in \{1,2\}$. By the definition of $\tilde{\mathcal{F}}_1^{(i,j)}$ (in definition \ref{definition_structure I2_SNA character of parameters}) and proposition \ref{proposition_structure I1_critical set I_1}, we have a non-empty interval $I_1^{(i,j)}$ that satisfies the conditions
    \begin{align*}
        F^{M_1-M_0} \left( (I_1^{(i,j)}-M_1\omega) \times C_{k_1} \right) & \subset (I_0^{(i,j)}-M_0\omega) \times C_i \\
        \text{ and } F^{-M_1+M_0} \left((I_1^{(i,j)}-M_1\omega) \times E_{k_2} \right) & \subset (I_0^{(i,j)}+M_0\omega) \times E_i,  
    \end{align*}
for some $k_1,k_2 \in \{1,2\}$, and this gives us
\begin{align*}
        F^{M_1+1} \left( (I_1^{(i,j)}-M_1\omega) \times C_{k_1} \right) & \subset F^{M_0+1}\left( (I_0^{(i,j)}-M_0\omega) \times C_i \right) \\
        \text{ and } F^{-M_1+1} \left( (I_1^{(i,j)}-M_1\omega) \times E_{k_2} \right) & \subset F^{-M_0+1} \left( (I_0^{(i,j)}+M_0\omega) \times E_i \right).  
    \end{align*}
Recall the definition of set $I_1^{(i,j)}$, which is
    \begin{align*}
         I_1^{(i,j)}+\omega = \pi_1 \circ \left( F^{M_0+1} \left( (I_0^{(i,j)}-M_0\omega)\times C_i \right) \bigcap F^{-M_0+1} \left( (I_0^{(i,j)}+M_0\omega)\times E_i \right) \right). 
    \end{align*}
By the estimates in proposition \ref{proposition_structure I1_backward iterate curves} \ref{lemma item_structure I1_backward iterate curves_c} and proposition \ref{proposition_structure I1_forward iterate curves} \ref{lemma item_structure I1_forward iterate curves_c}, we know that the sets $F^{M_0+1}\left( (I_0^{(i,j)}-M_0\omega)\times C_i \right)$ and $F^{-M_0+1} \left( (I_0^{(i,j)}+M_0\omega)\times E_i \right)$ intersect transversally over the interval $(I_1^{(i,j)}+\omega)$. We can also deduce from  propositions \ref{proposition_structure I2_backward iterate curves} \ref{lemma item_structure I2_backward iterate curves_c} and \ref{proposition_structure I2_forward iterate curves} \ref{lemma item_structure I1_forward iterate curves_c} that the curves $\phi_{1,k_1}^{\pm}$ and $\psi_{1,k_2}^{\pm}$ are monotonic over the interval $I_1^{(i,j)}$, with slopes $|\partial_{\theta}\phi_{1,k_1}^{\pm}|> \lambda^{1/2}$ and $|\partial_{\theta}\psi_{1,k_2}^{\pm}| < 1$. Hence, the curves $\phi_{1,k_1}^{\pm}$ and $\psi_{1,k_2}^{\pm}$ must intersect transversally over $I_1^{(i,j)}$, which proves that the set $\tilde{I}_2^{(i,j)}$ must be non-empty. Since $I_2^{(i,j)}(\omega)$ contains $\tilde{I}_2^{(i,j)}$, this proves that $I_2^{(i,j)}(\omega)$ is non-empty. 

For a fixed $\omega \in \tilde{\mathcal{F}}_1$ and $i,j\in \{1,2\}$, say $\omega \notin \tilde{\mathcal{F}}_1^{(i,j)}$, that is, either of the following conditions hold true: 
\begin{align*}
    \pi_2 \circ F^{M_1-M_0}\left( \left( I^{(i,j)}_1-M_1\omega\right) \times (C_1 \cup C_2) \right) & \subset C_{i+1} \\
    \pi_2 \circ F^{-M_1+M_0}\left( \left( I^{(i,j)}_1+M_1\omega\right) \times (E_1\cup E_2) \right) & \subset E_{i+1},
\end{align*}
where $(i+1)$ denotes $(i+1) \pmod{2}$. Without loss of generality, let us assume that 
\begin{align*}
    \pi_2 \circ F^{M_1-M_0}\left( \left( I^{(i,j)}_1-M_1\omega\right) \times (C_1 \cup C_2) \right) & \subset C_{i+1}.
\end{align*}
Then, it follows from the above condition that
\begin{align}\label{equation_structure I2_Fij_1}
    F^{M_1+1} \left( (I_1^{(i,j)}-M_1\omega) \times (C_1\cup C_2) \right) & \subset F^{M_0+1}\left( (I_0^{(i,j)}-M_0\omega) \times C_{i+1} \right) \subset A_{-1}^{(i+1)}.
\end{align}
Since the sets $B_1$, $B_0$ and $B_{-1}$ are nested, as shown in lemmas \ref{lemma_structure I1_match Ai Bi} and lemma \ref{lemma_structure I0_match Ai Bi}, we have
\begin{align}\label{equation_structure I2_Fij_2}
    F^{-M_1+1} \left( (I_1^{(i,j)}+M_1\omega) \times (E_1\cup E_2) \right) & \subset F^{-M_0+1}\left( (I_0^{(i,j)}+M_0\omega) \times (E_1\cup E_2) \right) \\
    & \subset B_{-1},
\end{align}
and also that $\left( A_{-1}^{(i+1)} \cap B_{-1}^{(i)} \right) = \emptyset$. Then, by the definition of set $I_2^{(i,j)}(\omega)$ and equations (\ref{equation_structure I2_Fij_1} - \ref{equation_structure I2_Fij_2}), we obtain
\begin{align}\label{equation_structure I2_Fij_2}
    I_2^{(i,j)}(\omega) \subset \pi_1 \left( A_{-1}^{(i+1)} \cap B_{-1}^{(i+1)} \right) = \bigcup_{j=1}^{2} I_0^{(i+1,j)}.
\end{align}
However, by construction $I_2^{(i,j)} \subset I_0^{(i,j)}$ and, by proposition \ref{proposition_structure I0_I0 subsets}, the subsets in $\{I_0^{(i',j')}\}_{i',j'=1}^{2}$ are disjoint subsets of $I_0$. Thus, we can conclude that $I_2^{(i,j)} = \emptyset$.

\end{proof}
%------------------- PROP END -------------------

In the following proposition, the condition on the  parameter $\omega$ that leads to a non-empty critical set $I_2(\omega)$, as in definition \ref{definition_structure I2_SNA character of parameters}, is realised as a characteristic of the orbit of $\theta^{(1,1)}$ with respect to $\omega$. This will be essential for determining parameters in the set $\tilde{\mathcal{F}}^{(1,1)}_1$ in section \ref{section_parameter space}.

%------------------- PROP BEGIN -------------------
\begin{proposition}\label{proposition_structure I2_I2  nonempty F^{(1,1)}}
    For all sufficiently large $\lambda$ and any $\omega \in \tilde{\mathcal{F}}_1$, if there exists $1 \leq q  \leq M_0^9$ such that $(\theta^{(1,1)} \pm q\omega)  \in G_2$ and
    \begin{align*}
        \left\{(\theta^{(1,1)} + m\omega): |m|\leq (q-1) \right\} \cap S_3 & = \emptyset,
    \end{align*}
    then $\omega \in \tilde{\mathcal{F}}_1^{(1,1)}$.
\end{proposition}

\begin{proof}

For a sufficiently large $\lambda$, we have $ M_0^{9} < M_1$. Recall that $I_1^{(1,1)}(\omega) \subset I_0^{(1,1)}$ (by definition \ref{definition_structure I1_I1 definition}). Then, by the definition of the sets $G_1$ and $G_2$ (in definition \ref{definition_structure I1_sync set Gs}) and the structure of the critical set $I_0^{(1,1)}$ (in proposition \ref{proposition_structure I0_I0 subsets}), it is clear that $B_{|I_1(\omega)|}(G_2) \subset G_1$ and thus
\begin{align*}
    (\theta^{(1,1)}\pm q\omega) \in G_2 \implies (I_1^{(1,1)} \pm q\omega) \subset G_1 \subset I_0.
\end{align*}
By disjointness condition $(\mathcal{F}1)_0$, we have $(I_1^{(1,1)} - (q+M_0)\omega) \subset (G_1-M_0\omega) \subset (I_0-M_0\omega)$. Then, by property $(C1)_0$ (in proposition \ref{proposition_block 0 contraction}) we infer that 
\begin{align*}
F^{M_1 - q - M_0} \left((I_1^{(1,1)} - M_1\omega) \times C_1\cup C_2\right)  \subset \left((G_1-M_0\omega) \times (C_1\cup C_2)\right), 
\end{align*}
and thus $\pi_2 \circ F^{M_1 - q +M_0+S_0} \left((I_1^{(1,1)} - M_1\omega) \times C_1\cup C_2\right)  \subset (C_1\cup C_2)$.
Moreover, by the property of $G_1$ which was derived in proposition \ref{proposition_structure I1_sync set Gs} \ref{lemma item_structure I1_sync set Gs_c}, 
\begin{align}\label{equation_F(1,1) sync 0}
    \pi_2 \circ F^{M_1 - q +M_0+S_0} \left((I_1^{(1,1)} - M_1\omega) \times C_1\cup C_2\right) \subset \pi_2 \circ F^{2M_0+S_0}\left((G_1-M_0\omega) \times (C_1\cup C_2)\right) \subset C_1.
\end{align}
Recall that, by proposition \ref{proposition_structure I1_switch set Ss}, the sets $S_3 = (\theta^{(1,1)} - 3\lambda^{-1+1/64}, \theta^{(1,2)} + 3\lambda^{-1+1/64})$ and $S_1(\omega)$ satisfy the condition
\begin{align*}
    S_1(\omega) \subset (\theta^{(1,1)}-\lambda^{-1+1/64}, \theta^{(1,2)} + \lambda^{-1+1/64}) \subset S_3.
\end{align*}
Then, it follows that $ (I_1^{(1,1)} + m\omega) \cap S_1(\omega)  = \emptyset$, for all $|m| \leq q-1$. By equation \eqref{equation_F(1,1) sync 0}, condition $(C1)_0$ and the property of $S_1(\omega)$ established in proposition \ref{proposition_structure I1_switch set Ss}, the above condition gives us
\begin{align}\label{equation_F(1,1) sync}
    \pi_2 \circ F^{M_1 - M_0} \left((I_1^{(1,1)} - M_1\omega) \times C_1\cup C_2\right) \subset C_1.
\end{align}
Similarly, we can also derive 
\begin{align}\label{equation_F(1,1) sync}
    \pi_2 \circ F^{-M_1 + M_0} \left((I_1^{(1,1)} + M_1\omega) \times E_1\cup E_2\right) \subset E_1,
\end{align}
and this proves $\omega \in \tilde{\mathcal{F}}_1^{(1,1)}$.
\end{proof}
%------------------- PROP END -------------------

In the following proposition, we establish the conditions for which the critical set $I_2(\omega)$ will be empty. This will be essential for the construction of the parameter set $\Omega_2$ in section \ref{section_parameter space}.

%------------------- PROP BEGIN -------------------
\begin{proposition}\label{proposition_structure I2_I2  empty}
    For any $\omega \in \tilde{\mathcal{F}}_1$, let $N_1$ and $N_2$ be the smallest positive integers for which $N_1 < N_2 < M_1$, 
    \begin{align*}
        \{(I_1 - N_1\omega), (I_1+N_1\omega)\} \cap I_0 \neq \emptyset,
    \end{align*}
    and 
    \begin{align*}
        \{(I_1 - N_2\omega), (I_1+N_2\omega)\} \cap I_0 \neq \emptyset.
    \end{align*}
Then, the following conditions are satisfied:
\begin{enumerate}[label=\alph*)]
    \item \label{lemma item_structure I2_empty set a} If $(\theta^{(1,1)} - N_2\omega) \in G_2^{(2)}$ and $(\theta^{(1,1)} - N_1\omega) \in S_2$, then $I_2^{(1,1)}(\omega) = \emptyset$.
    
\item \label{lemma item_structure I2_empty set b} If  $(\theta^{(1,2)} + N_2\omega) \in G_2^{(1)}$ and $(\theta^{(1,2)} + N_1\omega) \in S_2$, then $I_2^{(1,2)}(\omega) = \emptyset$.

    \item  \label{lemma item_structure I2_empty set c} If $ (\theta^{(2,1)} - {N}_1) \in G_2^{(1)}$, then $I_2^{(2,1)}(\omega) = \emptyset$.

\item  \label{lemma item_structure I2_empty set d} If  $(\theta^{(2,2)} + {N}_1) \in G_2^{(2)}$, then $I_2^{(2,2)}(\omega) = \emptyset$.
\end{enumerate}
\end{proposition}

\begin{proof}
\textit{Proof of \ref{lemma item_structure I2_empty set a}:} 
Recall that $I_1^{(1,1)} \subset I_0^{(1,1)}$ and by proposition \ref{proposition_structure I0_I0 subsets}, we have the estimate $|I_0^{(i,j)}| < \lambda^{-1+1/64}$. Since $(\theta^{(1,1)} - N_2\omega) \in G_2^{(2)}$, by the definition of the sets $G_1$ and $G_2$ (in definition \ref{definition_structure I1_sync set Gs}), we obtain $(I_1^{(1,1)} - N_2\omega) \subset G_1$. Since $G_1 \subset (I_0\setminus 3 I_1)$ and $N_2 < M_1$, it follows from $(C1)_0$ that
\begin{align*}
    F^{M_1-N_2-M_0} \left( (I_1^{(1,1)} -M_1\omega) \times (C_1\cup C_2) \right) \subset (G_1-M_0\omega) \times (C_1 \cup C_2).
\end{align*}
Then, by the main property of synchronisation which was shown for $G_1$ in proposition \ref{proposition_structure I1_sync set Gs},
\begin{align*}
    F^{M_1-N_2+M_0+S_0} \left( (I_1^{(1,1)} -M_1\omega) \times (C_1\cup C_2) \right) \subset (I_1^{(1,1)}-(N_2+M_0+S_0)\omega) \times C_1.
\end{align*}
Since $(\theta^{(1,1)} - N_1\omega) \in S_2$, by the definition of $S_1$, $S_2$ and proposition \ref{proposition_structure I1_switch set Ss}, we deduce that $(I_1^{(1,1)} - N_1\omega) \subset S_1(\omega)$. By the hypothesis of the lemma, we have $(I_1^{(1,1)} -m\omega)_{m=N_2+M_0+S_0}^{N_1-1} \cap I_0 = \emptyset$ and this gives us
\begin{align*}
    F^{M_1-N_1-M_0} \left( (I_1^{(1,1)} -M_1\omega) \times (C_1\cup C_2) \right) \subset (I_0 -M_0\omega) \times C_1. 
\end{align*}
Then, by the property of $S_1$ established in proposition \ref{proposition_structure I1_switch set Ss}, we obtain
\begin{align*}
    F^{M_1-N_1+M_0+S_0} \left( (I_1^{(1,1)} -M_1\omega) \times (C_1\cup C_2) \right) \subset (I_1^{(1,1)} - (N_1+M_0+S_0)\omega) \times C_2.
\end{align*}
Since $(I_1^{(1,1)} -m\omega)_{m=1}^{N_1+M_0+S_0} \cap I_0 = \emptyset$,  
\begin{align*}
    F^{M_1-M_0} \left( (I_1^{(1,1)} -M_1\omega) \times (C_1\cup C_2) \right) \subset (I_1^{(1,1)} - M_0\omega) \times C_2.
\end{align*}
This proves that $\omega \notin \tilde{\mathcal{F}}_1^{(1,1)}$ (as defined in \ref{definition_structure I2_SNA character of parameters}), and by proposition \ref{proposition_structure I2_I2 empty nonempty}, it then follows that $I_2^{(1,1)}(\omega) = \emptyset$.

\vspace{0.2cm}
\noindent\textit{Proof of \ref{lemma item_structure I2_empty set b}:} Similarly, since $(\theta^{(1,2)} + N_2\omega) \in G_2^{(1)}$ this gives us $ (I_1^{(1,2)} + N_2\omega) \subset G_1^{(1)}$. Then, using the properties of $G_1$ shown in proposition \ref{proposition_structure I1_sync set Gs}, and $(E1)_0$, we obtain 
\begin{align*}
    F^{-M_1+N_2-M_0-S_0} \left( (I_1^{(1,2)} +M_1\omega) \times (E_1\cup E_2) \right) \subset (I_1^{(1,2)}+(N_2-M_0-S_0)\omega) \times E_1.
\end{align*}
Since $(I_1^{(1,2)} +m\omega)_{m=N_2-M_0-S_0}^{N_1-1} \cap I_0 = \emptyset$ due to the hypothesis, we also have
\begin{align*}
    F^{-M_1+N_1+M_0} \left( (I_1^{(1,2)} +M_1\omega) \times (E_1\cup E_2) \right) \subset (I_0 -M_0\omega) \times E_1.
\end{align*}
Since $(\theta^{(1,2)} - N_1\omega) \in S_2$, by the definition of sets $S_1(\omega)$ and $S_2$, we infer $(I_1^{(1,2)} - N_1\omega) \subset S_1$. Then, using the properties of $S_1$ established in proposition \ref{proposition_structure I1_switch set Ss}, and $(E1)_0$, we obtain
\begin{align*}
    F^{-M_1+N_1-M_0-S_0} \left( (I_1^{(1,2)} +M_1\omega) \times (E_1\cup E_2) \right) \subset (I_1^{(1,2)} + (N_1-M_0-S_0)\omega) \times E_2.
\end{align*}
Since $(I_1^{(1,2)} +m\omega)_{m=1}^{N_1-M_0-S_0} \cap I_0 = \emptyset$, this gives us 
\begin{align*}
    F^{-M_1+M_0} \left( (I_1^{(1,2)} + M_1\omega) \times (E_1\cup E_2) \right) \subset (I_1^{(1,2)} + M_0\omega) \times E_2.
\end{align*}
This proves that $\omega \notin \tilde{\mathcal{F}}_1^{(1,2)}$ (as defined in \ref{definition_structure I2_SNA character of parameters}), and by proposition \ref{proposition_structure I2_I2 empty nonempty}, it then follows that $I_2^{(1,2)}(\omega) = \emptyset$.

\vspace{0.2cm}
\noindent\textit{Proof of \ref{lemma item_structure I2_empty set c}:} Similar to the preceding parts of the proof, by the definition of the sets $G_1^{(1)}$ and $G_2^{(1)}$ (in definition \ref{definition_structure I1_sync set Gs}), we can deduce that $(I_1^{(2,1)} - {N}_1) \subset G_1^{(1)}$, which then gives us
\begin{align*}
    F^{M_1-{N}_1+M_0+S_0} \left( (I_1^{(2,1)} - M_1\omega) \times (C_1\cup C_2) \right) \subset (I_1^{(2,1)}-({N}_1+M_0+S_0)\omega) \times C_1,
\end{align*}
by the properties of $G_1$ established in proposition \ref{proposition_structure I1_sync set Gs}, and $(C1)_0$. Since we have $(I_1^{(2,1)} - m\omega)_{m=1}^{{N}_1-1} \cap I_0 = \emptyset$, we obtain 
\begin{align*}
    F^{M_1-M_0} \left( (I_1^{(2,1)} - M_1\omega) \times (C_1\cup C_2) \right) \subset (I_1^{(2,1)}-M_0\omega) \times C_1.
\end{align*}
This proves that $\omega \notin \tilde{\mathcal{F}}_1^{(2,1)}$ (as defined in \ref{definition_structure I2_SNA character of parameters}), and by proposition \ref{proposition_structure I2_I2 empty nonempty}, it follows that $I_2^{(2,1)}(\omega) = \emptyset$.

\vspace{0.2cm}
\noindent\textit{Proof of \ref{lemma item_structure I2_empty set d}:} Similarly, by the definition of the sets $G_1^{(2)}$ and $G_2^{(2)}$ (in definition \ref{definition_structure I1_sync set Gs}), we can deduce that $(I_1^{(2,2)} + {N}_1) \subset G_1^{(2)}$
which then gives us
\begin{align*}
     F^{-M_1+{N}_1-M_0-S_0} \left( (I_1^{(2,2)} + M_1\omega) \times (E_1\cup E_2) \right) \subset (I_1^{(2,2)}+({N}_1+M_0+S_0)\omega) \times E_1.
\end{align*}
by the properties of $G_1$ established in proposition \ref{proposition_structure I1_sync set Gs}, and $(E1)_0$. Since we have $(I_1^{(2,2)} + m\omega)_{m=1}^{{N}_1-1} \cap I_0 = \emptyset$, we obtain 
\begin{align*}
    F^{-M_1+M_0} \left( (I_1^{(2,2)} + M_1\omega) \times (E_1\cup E_2) \right) \subset (I_1^{(2,2)}-M_0\omega) \times E_1.
\end{align*}
This proves that $\omega \notin \tilde{\mathcal{F}}_1^{(2,2)}$ (as defined in \ref{definition_structure I2_SNA character of parameters}) and by proposition \ref{proposition_structure I2_I2 empty nonempty}, we then have $I_2^{(2,2)}(\omega) = \emptyset$.

\end{proof}
%------------------- PROP END -------------------

%%%%%%%%%%%%%%%%%%%%%%%%%% ESTIMATES FOR NON_EMPTY I2 %%%%%%%%%%%%%%
\subsubsection{Estimates for a non-empty $I_2(\omega)$}\label{subsubsection_estimate nonempty I2}

For any $\omega \in \mathcal{F}_1$, the following proposition describes the structure of the critical set $I_2(\omega)$ when it is non-empty. As a consequence of the following estimate, one can infer that the disjointness condition $(\bar{\mathcal{F}}1)_1$, which was introduced in definition \ref{definition_block n_disjointness condn}, implies the condition $(\mathcal{F}1)_1$.

%------------------- PROP BEGIN -------------------
\begin{proposition}\label{proposition_structure I2_critical set I2}
   Let $\lambda$ be sufficiently large and let $i,j \in \{1,2\}$. For any $\omega \in \tilde{\mathcal{F}}^{(i,j)}_1 \cap \mathcal{F}_1$, we have the following estimates for the set $I_2^{(i,j)}(\omega)$:
    \begin{enumerate}[label=\alph*)]
        \item \label{lemma item_structure I2_critical set I2 a}
        $ \lambda^{-1-1/32} \alpha_{c}^{3b\tilde{M}_1} < |I^{(i,j)}_2(\omega)| < \alpha_c^{\tilde{M}_1/3}\lambda^{-1+1/32}.$
        \item \label{lemma item_structure I2_critical set I2 b} $|\partial_{\omega}I^{(i,j)}_2(\omega)| < 1/4$.
    \end{enumerate}
\end{proposition}

\begin{proof}
    %\textit{Proof of \ref{lemma item_structure I2_critical set I2 a}:} 
    Let us fix $i,j \in \{1,2\}$ and $\omega \in \mathcal{F}_1 \cap \tilde{\mathcal{F}}_1^{(i,j)}$. By proposition \ref{proposition_structure I2_I2 empty nonempty}, the set $I_2^{(i,j)}(\omega)$ is then a non-empty closed interval. Recall that, by definition \ref{definition_structure I2_ set I2 prime} and \ref{definition_structure I2_critical set I2}, the critical set $I_2^{(i,j)}$ satisfies  
    \begin{align*}
         \pi_2 \circ \left( F^{M_1}((I^{(i,j)}_1-M_1\omega) \times (C_1\cup C_2)) \cap F^{-M_1}((I^{(i,j)}_1+M_1\omega) \times (E_1\cup E_2)) \right) \subset I_2^{(i,j)},
    \end{align*}
    where $I_2^{(i,j)}$ is the smallest such interval. Since $f_{\theta}$ is orientation preserving for each $\theta \in \mathbb{T}$ and $I_1^{(i,j)}$ is a closed interval, there exists distinct $i_1,i_2 \in \{1,2\}$ such that 
    \begin{align}\label{equation_structure I2_backward curves}
        \psi_{1,i_1}^{-}(\theta) < \psi_{1,i_1}^{+}(\theta) < \psi_{1,i_2}^{-}(\theta) < \psi_{1,i_2}^{+}(\theta),
    \end{align}
    for all $\theta \in I_1^{(i,j)}$. Similarly, there exists distinct $i'_1, i'_2 \in \{1,2\}$ such that
\begin{align}\label{equation_structure I2_forward curves}
    \phi_{1,i'_1}^{-}(\theta) < \phi_{1,i'_1}^{+}(\theta) < \phi_{1,i'_2}^{-}(\theta) < \phi_{1,i'_2}^{+}(\theta),
\end{align}
for all $\theta \in I_1^{(i,j)}$. Let us define new sets
\begin{align*}
    \Tilde{B}_1 & := \left\{(\theta, [\psi_{1,i_1}^{-}(\theta), \psi_{1,i_2}^{+}(\theta)]): \theta \in I_1^{(i,j)} \right\} \\
    \Tilde{A}_1 & := \left\{(\theta, [\phi_{1,i'_1}^{-}(\theta), \phi_{1,i'_2}^{+}(\theta)]): \theta \in I_1^{(i,j)} \right\}.
\end{align*}
Note that, by equations \eqref{equation_structure I2_backward curves} and \eqref{equation_structure I2_forward curves}, we have $F^{M_1}(A_1) \subset \Tilde{A}_1$ and $F^{-M_1}(B_1) \subset \Tilde{B}_1$. This gives us 
$$\left( I_1^{(i,j)} \cap (\pi_1 (F^{M_1}(A_1) \cap F^{-M_1}(B_1))) \right) \subset \pi_1 (\tilde{A}_1 \cap \tilde{B}_1).$$
and thus $\tilde{I}_2^{(i,j)} \subset \pi_1 (\tilde{A}_1 \cap \tilde{B}_1)$. 

%------------------- FIGURE TRANSVERSE INTERSECTION -----------

% trim = left bottom right top
\begin{figure}[t]
    \centering
    \includegraphics[
        width=0.9\linewidth,
        trim=1cm 0.1cm 0.8cm 1cm,
        clip
    ]{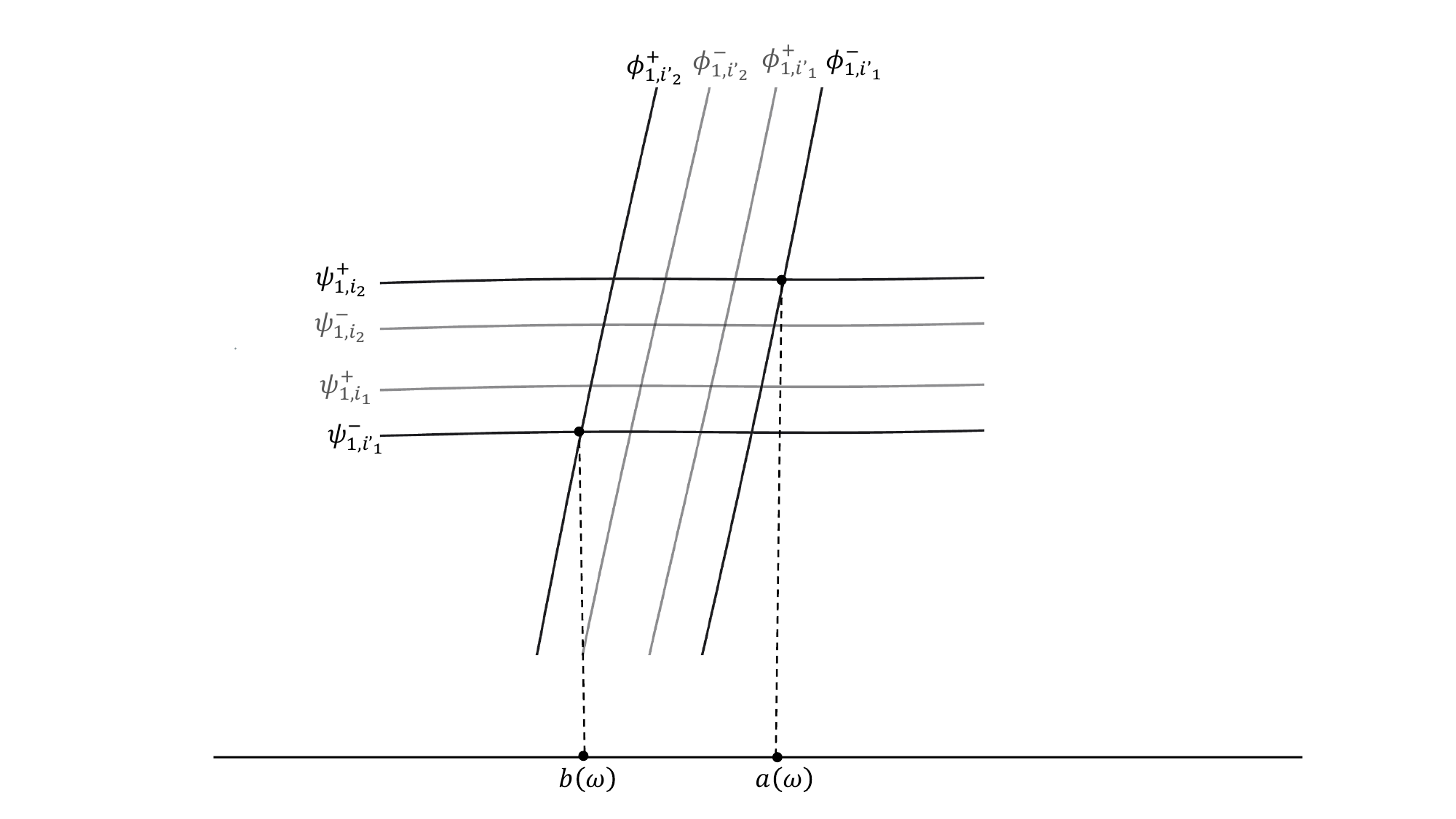}
    \caption{The above figure illustrates the transverse intersection of the curves $\{\phi_{1,i}^{\pm}\}_{i=1}^{2}$ and $\{\psi_{1,i}^{\pm}\}_{i=1}^{2}$ over a fixed interval $I_1^{(i,j)}$, as discussed in equations \eqref{equation_structure I2_backward curves} and \eqref{equation_structure I2_forward curves}.}
    \label{figure_transverse intersection}
\end{figure}

%--------------------------------------

By proposition \ref{proposition_structure I2_backward iterate curves} \ref{lemma item_structure I2_backward iterate curves_c} and \ref{proposition_structure I2_forward iterate curves} \ref{lemma item_structure I2_forward iterate curves_c}, we have the estimate $|\partial_{\theta}\phi_{1,l}^{\pm}|> \lambda^{1-1/64}$ and $|\partial_{\theta}\psi_{1,l}^{\pm}| < 1$, for all $\theta \in I_1^{(i,j)}$ and $l\in \{1,2\}$. This shows that the curves in $\{\phi_{1,l}^{\pm}\}_{l=1}^{2}$ and $\{\psi_{1,l}^{\pm}\}_{l=1}^{2}$ intersect transversally over the interval $I_1^{(i,j)}$. Then, it follows that the set $\pi_1 (\tilde{A}_1 \cap \tilde{B}_1)$ is a closed interval, for which the end points $a$ and $b$ satisfy the equations 
\begin{align}\label{equation_structure I2_transverse intersection}
    \phi_{1,i'_1}^{-}(a) & = \psi_{1,i_2}^{+}(a) \\
    \phi_{1,i'_2}^{+}(b) & = \psi_{1,i_1}^{-}(b). \label{equation_structure I2_transverse intersection-2}
\end{align}
By construction, there exists no closed interval which can be smaller than $\pi_1 (\tilde{A}_1 \cap \tilde{B}_1)$ and also contain $\tilde{I}^{(i,j)}_2$. This proves that 
$$I_2^{(i,j)} = \pi_1 (\tilde{A}_1 \cap \tilde{B}_1). $$
For the transverse intersection of the sets $\tilde{A}_1$ and $\tilde{B}_1$, we derive the following estimates using the Mean Value Theorem:
\begin{align*}
    \frac{2\inf\{|\phi_{1,i'_1}^{-}(\theta)-\phi_{1,i'_2}^{+}(\theta)|: \theta \in I_1^{(i,j)}\}}{|\partial_{\theta}\phi_{1}^{\pm}|_{sup}+|\partial_{\theta}\psi^{\pm}_{1}|_{sup}} 
    & < |I_{2}^{(i,j)}| < \frac{2\sup\{|\psi_{1,i_1}^{-}(\theta)-\psi_{1,i_2}^{+}(\theta)|: \theta \in I_1^{(i,j)}\}}{|\partial_{\theta}\phi^{\pm}_{1}|_{inf}-|\partial_{\theta} \psi_{1}^{\pm}|_{sup}},
\end{align*}
where $|\partial_{\theta}\phi^{\pm}_{1}|_{sup}$ denotes the quantity $\sup \{ |\partial_{\theta}\phi^{\pm}_{1,i}|_{sup}: i=1,2\}$, and the analogous notations for the curves $\{\psi_{1,i}^{\pm}\}$ and for infimum are defined similarly. Using the estimates in propositions \ref{proposition_structure I2_backward iterate curves}, \ref{proposition_structure I2_forward iterate curves}, \ref{proposition_structure I2_backward iterate curves sync} and \ref{proposition_structure I2_forward iterate curves sync}, we obtain 
\begin{align*}
    \frac{2 \alpha_{c}^{b(M_1-1)} |e_1^{+}-e_1^{-}|}{\lambda^{1+1/64}+1} & < |I_{2}^{(i,j)}| < \frac{2 \left(\alpha_c^{(1-(1+b)\beta_0)}\right)^{M_1-M_0^{60}+1}}{\lambda^{1-1/64}-1}.
\end{align*}
Since $M_1 < 2\alpha_c^{M_0/8}$, for a sufficiently large $\lambda$ with respect to $\alpha_c$, this proves the estimate 
$$\lambda^{-1-1/32} \alpha_{c}^{b(M_1+1)} < |I^{(i,j)}_2(\omega)| < \alpha_c^{M_1/3}\lambda^{-1+1/32}.$$

%\textit{Proof of \ref{lemma item_structure I2_critical set I2 b}:} 
For each $\omega \in \mathcal{F}_1\cap\tilde{\mathcal{F}}_1^{(i,j)}$, by equations (\ref{equation_structure I2_transverse intersection} - \ref{equation_structure I2_transverse intersection-2}), the end points of the set $I_2^{(i,j)}(\omega)$, denoted by $a(\omega), b(\omega) \in I_1^{(i,j)}$, satisfy the conditions 
\begin{align*}
     \phi_{1,i'_1}^{-}(a(\omega), \omega) & = \psi_{1,i_2}^{+}(a(\omega), \omega) \\
    \phi_{1,i'_2}^{+}(b(\omega), \omega) & = \psi_{1,i_1}^{-}(b(\omega), \omega).
\end{align*}
By the implicit function theorem, the functions $a(\omega)$ and $b(\omega)$ are continuously differentiable and we have the estimate
\begin{align*}
|a'(\omega)| & < \frac{|\partial_{\omega}\phi_{1,i'_1}^{-}(\theta, \omega)|_{sup} + |\partial_{\omega}\psi_{1,i_2}^{+}(\theta, \omega)|_{sup}}{|\partial_{\theta}\psi_{1,i_2}^{+}(\theta, \omega) - \partial_{\theta}\phi_{1,i'_1}^{-}(\theta, \omega)|_{inf}}.
\end{align*}
By the propositions \ref{proposition_structure I2_backward iterate curves} and \ref{proposition_structure I2_forward iterate curves}, we have the estimates $|\partial_{\theta} \phi_{1,i}^{\pm}| > \lambda^{1-1/64}$, $|\partial_{\theta} \psi_{1,i}^{\pm}| <1$, $|\partial_{\omega}\phi^{\pm}_{1,i}|<1$ and $|\partial_{\omega}\psi^{\pm}_{1,i}|$, for all $i \in \{1,2\}$. This shows that the end point depend robustly on the parameter $\omega$ and we obtain $|a'(\omega)| < \lambda^{-1+1/32}$.
By symmetry, we similarly have $|b'(\omega)| < \lambda^{-1+1/32}$. This proves that $|\partial_{\omega}I_2(\omega)| < \lambda^{-1+1/32}$.
\end{proof}

%------------------- PROP END -------------------

%--------------------- Critical Set I_n ---------------------------

\subsection{Structure of critical set $I_{n+1}$}\label{subsection_structure of critical set I(n+1)}

In this section, we describe an inductive construction of a nested sequence of non-empty sets $(I_n(\omega))_{n>2}$, satisfying the conditions $(\mathcal{I}1-3)_n$ for critical sets (in definition \ref{definition_block n_Critical sets In}). For each $n>2$, the dichotomy for the critical set $I_n(\omega)$ being empty or non-empty will be shown to follow directly from the results in subsection \ref{subsection_structure of critical set I2} for critical set $I_2(\omega)$. 

For the rest of the subsection, we assume the existence of a non-empty set $\bar{G} \subset I_0 \setminus 3I_1(\omega)$ and an associated set $\hat{G}$, as in definition \ref{definition_block n_sync condn set_g1 g2}. Recall that, by definition \ref{definition_block n_sync condn parameter}, a parameter $\omega$ is said to satisfy condition $(\mathcal{G})$ if
\begin{align*}
     \left\{\theta+k\omega \right\}_{k=0}^{M_0^{56}} \cap \bar{G} \neq \emptyset,
\end{align*}
for every $\theta \in \mathbb{T}$. Further, for each $n\geq 1$ we also have the associated parameter set $\mathcal{F}_n$ given by
\begin{align}
    \mathcal{F}_n = \mathcal{F}_n(f,\lambda,{\bar{G}}) := \tilde{\mathcal{F}}_n(f,\lambda) \cap \left\{\omega \in \mathbb{T}: \omega \text{ satisfies } (\mathcal{G}) \text{ for the set } \bar{G}\right\}. 
\end{align}
In this section, we focus exclusively on parameters in $\mathcal{F}_n$, where $n\geq 1$.

\begin{definition}[Preliminary Critical Set]\label{definition_structure In_ set In prime}
    For any $l \geq 2$ and $\omega \in \mathcal{F}_{l-1}(M_{l-1},\dots , M_0, \bar{G})$, given a sequence of non-empty critical sets $\{I_i(\omega)\}_{i=0}^{l-1}$, that satisfy the conditions in definition \ref{definition_block n_Critical sets In}, and for which $I_1(\omega)$ is as in definition \ref{definition_structure I1_I1 definition}, we define the following: 
    \begin{description}
    \item[$(\tilde{\mathcal{I}1)}_{l}$] We define a preliminary critical set $\tilde{I}_{l}(\omega)$ by $$\Tilde{I}_{l}(\omega) := \pi_1 \left( F^{M_{l-1}}(A_{l-1}) 
       \cap F^{-M_{l-1}}(B_{l-1}) \right),$$ where 
    \begin{align*}
    A_{l-1} & = \left( I_{l-1}-M_{l-1}\omega\right) \times (C_1\cup C_2) \\
    B_{l-1} & = \left(I_{l-1}+M_{l-1}\omega\right)\times (E_1 \cup E_2).
\end{align*}

    \item[$(\tilde{\mathcal{I}2)}_{l}$]  For each $i,j \in \{1,2\}$, let us define  $$\tilde{I}^{(i,j)}_{l}(\omega) := \Tilde{I}_{l}(\omega) \cap I^{(i,j)}_1(\omega).$$
\end{description}
\end{definition} 

For each $l\geq 1$, we have $\Tilde{I}_{l} \subset {I}_{1}$ by construction, and thus $\displaystyle \Tilde{I}_{l} = \bigcup_{i,j=1}^{2} \Tilde{I}_{l}^{(i,j)}$. Throughout the remainder of the article, the notation $I_l(\omega)$ for $l\geq 3$, is reserved exclusively for the sets introduced in the following definition.

\begin{definition}[Critical subsets]\label{definition_structure In_critical set In}
    For any $l\geq 2$ and $i,j \in \{1,2\}$, under the assumptions of definition \ref{definition_structure In_ set In prime}, let $I_{l}^{(i,j)}(\omega)$ denote the smallest closed interval in $\mathbb{T}$ containing the set $\tilde{I}_{l}^{(i,j)}(\omega)$. Then, we fix the set $I_{l}(\omega)$ to be
    $$I_{l}(\omega) := \bigcup_{i,j=1}^{2} I_{l}^{(i,j)}(\omega).$$
\end{definition}

Observe that, for the case $l=2$ the above definition coincides with the definition of critical set $I_2(\omega)$ in definition \ref{definition_structure I2_critical set I2}. Then, it follows from proposition \ref{proposition_structure I2_critical set I2} that, for each $i,j\in \{1,2\}$ and $\omega \in (\mathcal{F}_1 \cap \tilde{\mathcal{F}}_1^{(i,j)})$ we have a non-empty critical set $I_2(\omega)$.

For each $n > 2$ and $\displaystyle \omega \in \mathcal{F}_n$, the set $I_n(\omega)$ introduced in definition \ref{definition_structure In_critical set In} also satisfies the conditions of critical set in definition \ref{definition_block n_Critical sets In}. This ensures the validity of the contraction rates obtained in section \ref{section_Mechanism of slow-returns}.  For each $n\geq 2$, we prove that such a construction will yield a non-empty critical set $I_{n+1}(\omega)$, for which we describe its geometric structure.

%------------------- LEMMA BEGIN -------------------
\begin{lemma}\label{lemma_structure In_match Ai Bi}
    Let $n\geq 2$, $\omega\in \mathcal{F}_n$ and $(I_l(\omega))_{l=0}^{n}$ be a sequence of non-empty critical sets, as defined in definition \ref{definition_structure In_critical set In}. Then, 
    \begin{align*}
        F^{M_n-M_{n-1}}(A_n) & \subset A_{n-1} \\
        F^{-M_n+M_{n-1}}(B_n) & \subset B_{n-1}.
    \end{align*}
\end{lemma}
\begin{proof}

 By the definition of critical sets in definition \ref{definition_structure In_ set In prime} and \ref{definition_structure In_critical set In}, we have $\pi_1 \circ F^{M_n-M_{n-1}}(A_n) = (I_n-M_{n-1}\omega)$. By the disjointness condition $(\mathcal{F}2)_n$ (in definition \ref{definition_block n_disjointness condn}), it follows that
    \begin{align*}
        (I_n-M_n\omega) \cap \left( \bigcup_{l=0}^{n-1} \bigcup_{|m| \leq M_l} (I_l-m\omega) \right) & = \emptyset.
    \end{align*}
Thus, for any $\theta \in I_n-M_n\omega$, we also have $\theta \in \Theta_{n-1}$. By the disjointness condition $(\mathcal{F}1)_n$, $M_n$ is the smallest positive integer such that $(\theta + M_n \omega) \in I_n$. Since $(\theta+(M_n-M_{n-1})\omega) \in (I_{n-1}-M_{n-1}\omega)$, by the disjointness condition $(\mathcal{F}1-2)_{n-1}$ (in definition \ref{definition_block n_disjointness condn}),    
     \begin{align*}
         (\theta + (M_n-M_{n-1}) \omega) \notin \bigcup_{l=0}^{n-1} \bigcup_{m=1}^{M_l+S_0-1} (I_l +m\omega).
     \end{align*}
     From the above and by $(C1)_{n-1}$ (in proposition \ref{proposition_block n contraction}), for any $x \in C_1 \cup C_2$ this gives us $ \pi_2 \circ F^{M_n-M_{n-1}}(\theta, x) \in (C_1 \cup C_2)$. Hence, 
     $$ F^{M_n-M_{n-1}}((I_{n}-M_n\omega) \times C_1 \cup C_2) \subset (I_{n-1}-M_{n-1}\omega) \times (C_1 \cup C_2) = A_{n-1}. $$
It follows similarly that $F^{-M_n-M_{n-1}}(B_n) \subset B_{n-1}$.

\end{proof}
%------------------- LEMMA END -------------------

For each $i\in \{1,2\}$, we denote the components of $A_n$ and $B_n$ by
\begin{align*}
    A_n^{(i)} & := \left(I_n-M_n\omega\right) \times C_k \\
 B_n^{(i)} & := \left(I_n+M_n\omega\right) \times E_k.
\end{align*}
For each $\omega \in \mathcal{F}_n$ and $i\in \{1,2\}$, let us define the curves $\phi_{n,i}^{\pm}(\cdot,\omega), \psi_{n,i}^{\pm}(\cdot,\omega): I_n \rightarrow \mathbb{T}$ by
\begin{align*}
    \phi_{n,i}^{\pm}(\theta, \omega) & = \pi_2 \circ F^{M_n+1}(\theta-M_n\omega, c_i^{\pm}) \\
    \psi_{n,i}^{\pm}(\theta, \omega) & = \pi_2 \circ F^{-M_n+1}(\theta + M_n\omega, e_i^{\pm}). 
\end{align*}
Whenever convenient, we denote the curves as $\phi_{n,i}^{\pm}(\theta, \omega) = \phi_{n,i}^{\pm}(\theta)$ and $\psi_{n,i}^{\pm}(\theta, \omega) = \psi_{n,i}^{\pm}(\theta)$.

\begin{lemma}\label{lemma_structure In_iterate curves}
Let $n\geq 2$, $\omega\in \mathcal{F}_n$ and $(I_l(\omega))_{l=0}^{n}$ be the sequence of non-empty critical sets, as in definition \ref{definition_structure In_critical set In}. For any $i\in \{1,2\}$,  
\begin{align*}
    F^{-M_n+1}(B_n^{(i)}) = \left\{(\theta,x): \theta \in (I_n+\omega),\ x\in \left[\psi_{n,i}^{-}(\theta), \psi_{n,i}^{+}(\theta)\right] \right\} \\
    F^{M_n+1}(A_n^{(i)}) = \left\{(\theta,x): \theta \in (I_n+\omega),\ x\in \left[\phi_{n,i}^{-}(\theta), \phi_{n,i}^{+}(\theta)\right] \right\}.
\end{align*}
\end{lemma}

\begin{proof}
For each $\theta \in I_n$, by the continuity of $F$ and since $f_{\theta}(x)$ is an orientation preserving diffeomorphism on $\mathbb{T}$, we obtain 
\begin{align*}
    F^{M_n+1}\left((I_n-M_n\omega) \times [c_i^{-}, c_i^{+}]\right) & = \left\{ (\theta, x):\theta \in (I_n+\omega),\ x\in [\phi_{n,i}^{-}(\theta, \omega), \phi_{n,i}^{+}(\theta, \omega)]\right\},\\
    F^{-M_n+1}\left((I_n+M_n\omega) \times [e_i^{-}, e_i^{+}]\right) & = \left\{(\theta,x): \theta \in (I_n+\omega), \ x\in [\psi_{n,i}^{-}(\theta, \omega), \psi_{n,i}^{+}(\theta, \omega)] \right\}.
\end{align*}
    
\end{proof}
%------------------- LEMMA END -------------------

The preceding lemma shows that, for each $n\geq 2$, the geometry of the critical set $I_{n+1}(\omega)$ is closely connected to the structure of the curves $\phi_{n,i}^{\pm}$ and $\psi_{n,i}^{\pm}$, where $i \in \{1,2\}$. In the following propositions we establish the necessary estimates for these curves.

%------------------- PROP BEGIN -------------------
\begin{proposition}\label{proposition_structure In_backward iterate curves}
    Let $n\geq 2$, $\omega\in \mathcal{F}_n$ and $(I_l(\omega))_{l=0}^{n}$ be a sequence of non-empty critical sets, as in definition \ref{definition_structure In_critical set In}. For each $i \in \{1,2\}$, the maps $\psi_{n,i}^{\pm}(\theta)$ satisfy the following conditions: 
\begin{enumerate}[label=\alph*)]
\item \label{lemma item_structure In_backward iterate curves a} For each ${i'},{j'} \in \{1,2\}$, there exists $k \in \{1,2\}$ such that
$$\bigcup_{i=1}^{2}[\psi_{n,i}^{-}(\theta), \psi_{n,i}^{+}(\theta)] \subset E_k, \ \ \ \text{for all } \theta \in I_n^{({i'},{j'})}.$$

\item \label{lemma item_structure In_backward iterate curves b} For all $\theta \in I_n$ and $i \in \{1,2\}$, 
\begin{align*}
   \alpha_c^{b(M_n-1)}|e_i^+ - e_i^-| & \leq |\psi_{n,i}^{+}(\theta) - \psi_{n,i}^{-}(\theta)| \leq \left(\alpha_c^{1-(b+1)\beta_{n-1}}\right)^{M_n-1}|e_i^{+} -e_i^{-}|.
\end{align*}

\item \label{lemma item_structure In_backward iterate curves c} For all $\theta \in I_n$,  
\begin{align*}
\alpha_c^{b(M_n-1)}|e_1^+ - e_1^-| \leq \max_{i,j \in \{1,2\}} |\psi_{n,i}^{+}(\theta) - \psi_{n,j}^{-}(\theta)| \leq \left(\alpha_c^{1-(b+1)\beta_{n-1}} \right)^{M_n-M_0^{60}-1}.
\end{align*}

\item \label{lemma item_structure In_backward iterate curves d} For all $\theta\in I_n$ and $i \in \{1,2\}$, 
$$|\partial_{\theta}\psi_{n,i}^{\pm}(\theta)|< 1.$$

\item \label{lemma item_structure In_backward iterate curves e} For any $\omega \in \mathcal{F}_n$, $\theta \in I_n$ and $i\in \{1,2\}$, 
$$ |\partial_{\omega}\psi_{n,i}^{\pm}(\theta, \omega)| < 1. $$

\end{enumerate}
\end{proposition}

\begin{proof}
\textit{Proof of \ref{lemma item_structure In_backward iterate curves a}:} 
 Let us fix $i',j' \in \{1,2\}$. Since $\omega \in \mathcal{F}_n$, by the disjointness condition $(\mathcal{F}1)_n$, we have $(I_n^{(i',j')}+M_n \omega) \subset \Theta_n$ and by condition $(\mathcal{F}2)_n$, we also have $(I_n^{(i',j')}+M_n \omega) \subset \tilde{\Theta}_{n-1}$. Since $\omega$ satisfies the condition $(\mathcal{G})$, for a fixed $\theta \in I_n^{(i',j')}$ there exists a constant $N_s < M_0^{60} < M_1$ such that $(\theta +(M_n-N_s+M_0+S_0)\omega) \in \bar{G}$. Then, it follows by the property $(g1)$ (in definition \ref{definition_block n_sync condn set_g1 g2}) that $(I_n^{(i',j')} +(M_n-N_s+M_0+S_0)\omega) \subset \hat{G}$. 
 
By the property $(g3)$ (in definition \ref{definition_block n_sync condn set_g1 g2}), this gives us 
\begin{align}\label{equation_structure In_backward iterate_1}
\pi_2 \circ F^{-N_s} \left((I^{(i',j')}_n+M_n\omega)\times (E_1 \cup E_2) \right) & \subset E_1.
\end{align} 
Since $(I_n^{(i',j')}+M_n \omega) \subset \tilde{\Theta}_{n-1}$ and $N_s < M_1$, we infer that $(I_n^{(i',j')} + (M_n-N_s)\omega) \subset \Theta_{n-1}$. Note that the points in the set $\left((I_n^{(i',j')} + (M_n-N_s)\omega) \times E_1\right)$ satisfy the hypothesis of proposition \ref{proposition_block n sync expansion}, and by $(\mathcal{F}1)_n$, it follows that 
\begin{align*}
    (I_n^{(i',j')}+\omega) \cap \left( \bigcup_{l=0}^{n-1}\bigcup_{m=0}^{M_l+S_0-1}(I_l-m\omega) \right) = \emptyset.
\end{align*}
Since $I_n^{(i',j')}$ is a closed interval, by the continuity of the map $F$ and $(Se1)_{n-1}$ (in proposition \ref{proposition_block n sync expansion}), we conclude that
\begin{align*}
\pi_2 \circ F^{-M_n+1}\left((I_n^{(i',j')}+M_n\omega) \times (E_1 \cup E_2) \right) & \subset E_k, \text{ for some } k\in \{1,2\}.
\end{align*}
This proves the property $\bigcup_{i=1}^{2} [\psi_{n,i}^{-}(\theta), \psi_{n,i}^{+}(\theta)] \subset E_k$, for some $k\in \{1,2\}$.

\vspace{0.2cm}
\noindent\textit{Proof of \ref{lemma item_structure In_backward iterate curves b}:}  Let us fix $\theta \in I_n$ and $i\in \{1,2\}$. Recall that
\begin{align*}
    [\psi_{n,i}^{-}(\theta), \psi_{n,i}^{+}(\theta)] = \pi_2 \circ F^{-M_n+1} \left( \left\{(\theta+ M_n \omega, x): x\in [e_i^-, e_i^+] \right\} \right).
\end{align*}
By $(\mathcal{F}2)_n$, we have $(\theta+M_n\omega)\in \Theta_{n-1}$ and $M_n-1$ is the smallest positive integer $k$ for which $((\theta+M_n\omega)-k\omega) \in (I_n + \omega)$. Then, by $(E2)_{n-1}$ (in proposition \ref{proposition_block n expansion}), it follows that 
    \begin{align*}
    \alpha_c^{b(M_n-1)}|e_i^+ - e_i^-| & \leq |\psi_{n,i}^{+}(\theta)-\psi_{n,i}^{-}(\theta)| \leq (\alpha_c^{1-(1+b)\beta_{n-1}})^{M_n-1}|e_i^+- e_i^-|.
    \end{align*}

\vspace{0.2cm}
\noindent\textit{Proof of \ref{lemma item_structure In_backward iterate curves c}:} Let us fix $i',j' \in \{1,2\}$. By equation \eqref{equation_structure In_backward iterate_1}, there exists a constant $N_s < M_0^{60} < M_1$ such that $(I_n^{(i',j')} +(M_n-N_s)\omega) \subset \Theta_{n-1}$ and
\begin{align*}
F^{-N_s}\left( (I_1^{(i',j')}+M_n\omega) \times (E_1 \cup E_2 ) \right) & \subset (I_1^{(i',j')} +(M_1-N_s)\omega) \times E_1.
\end{align*}
That is,
\begin{align*}
     \bigcup_{i=1}^{2} \left\{ (\theta+\omega,x): \theta \in I_n^{(i',j')}, x\in [\psi_{n,i}^{-}(\theta), \psi_{n,i}^{+}(\theta)]  \right\} \subset F^{-M_n+N_s+1} \left( (I_n^{(i',j')} +(M_n-N_s)\omega) \times E_1 \right). 
\end{align*}
     Since $F$ is an orientation preserving map, by $(E2)_{n-1}$ (in proposition \ref{proposition_block n expansion}), for any $\theta \in I_n^{(i',j')}$ we have the estimate
\begin{align*}
    \max_{i,j \in \{1,2\}} |\psi_{n,i}^{-}(\theta) - \psi_{n,j}^{+}(\theta)| 
    & < \max_{x\in E_1} |\partial_x \left(\pi_2 \circ F^{-M_n+N_s+1}(\theta+(M_n-N_s)\omega, x)\right) | \\
    & < (\alpha_c^{1-(1+b)\beta_{n-1}})^{M_n-N_s-1}.
\end{align*}
Similarly, for any $\theta \in I_n^{(i',j')}$, we obtain
\begin{align*}
     \alpha_c^{b(M_n-1)}|e_1^{+}-e_1^{-}| < \max_{i,j\in \{1,2\}} |\psi_{n,i}^{+}(\theta) - \psi_{n,j}^{-}(\theta)|.
\end{align*}
This gives us
\begin{align*}
    \alpha_c^{b(M_n-1)}|e_1^{+}-e_1^{-}| \leq \max_{i,j \in \{1,2\}} |\psi_{n,i}^{-}(\theta)-\psi_{n,j}^{+}(\theta)| \leq (\alpha_c^{1-(1+b)\beta_{n-1}})^{M_n-M_0^{60}-1},
\end{align*}
for each $\theta \in I_n$.

\vspace{0.2cm}
\noindent\textit{Proof of \ref{lemma item_structure In_backward iterate curves d}:} Let us fix some $\theta \in I_n(\omega)$. For the curves $\{\psi_{n,i}^{\pm}\}_{i=1}^{2}$, we have the expression 
 \begin{align*}
     \partial_{\theta} \psi_{n,i}^{\pm}(\theta, \omega)
    & = \sum_{k=1}^{M_n-1}\left((-g'(\theta+k\omega)) \prod_{j=1}^{k}(f^{-1})'(f(x_{-M_n+j})) \right).
\end{align*}
Since $(\theta+M_n\omega) \in \Theta_{n-1}$, for any $x \in E_1 \cup E_2$, we have the following estimate by $(E2)_{n-1}$ (in proposition \ref{proposition_block 0 expansion}), 
\begin{align*}
    \left|\prod_{j=1}^{k}(f^{-1})'(f(x_{-M_n+j}))\right| < (\alpha_c^{1-(1+b)\beta_{n-1}})^{k},
\end{align*}
for any $k \in \{1,\dots, M_n-1\}$. Since $I_n \subset I_0$, the bound $|g'(\theta+k\omega)|< \lambda^{-1/4+1/128}$, for all $k=1,\dots, M_0$, follows from $(\mathcal{F}1)_0$. For $k=M_0+1, \dots ,M_n$, we consider the maximum upper bound $|g'(\theta+k\omega )|< \lambda^{1+1/128}$. This gives us 
\begin{align*}
    |\partial_{\theta}\psi_{n,i}^{\pm}(\theta)| 
    & \leq \lambda^{-1/4+1/128} \sum_{k=1}^{M_0} (\alpha_c^{1-(1+b)\beta_{n-1}})^k + \lambda^{1+1/128}(\alpha_c^{1-(1+b)\beta_{n-1}})^{M_0} \sum_{k=1}^{M_n-M_0} (\alpha_c^{1-(1+b)\beta_{n-1}})^k.
\end{align*}
Recall that $\alpha_c^{1-(1+b)\beta_{n-1}} < \alpha_c^{1/2}$ and from equation \eqref{equation_backward curve I2_theta slope} (in proposition \ref{proposition_structure I2_backward iterate curves})  the constant $\lambda$ satisfies the bound
\begin{align*}
    \frac{\alpha_c^{1/2}}{1-\alpha_c^{1/2}} \left(\lambda^{-1/4+1/128} + \lambda^{1+1/128}\alpha_c^{M_0/2}\right) < \frac{1}{2}.
\end{align*}
This proves that
\begin{align*}
    |\partial_{\theta}\psi_{n,i}^{\pm}(\theta)| 
    & < \frac{\alpha_c^{1/2}}{1-\alpha_c^{1/2}} \left(\lambda^{-1/4+1/128} + \lambda^{1+1/128}\alpha_c^{M_0/2}\right) < 1.
    \stepcounter{equation}\tag{\theequation} \label{equation_backward curve In_theta slope}
\end{align*}

\vspace{0.2cm}
\noindent\textit{Proof of \ref{lemma item_structure In_backward iterate curves e}:} 
For the curves $\{\psi_{n,i}^{\pm}\}_{i=1}^{2}$, we have the expression 
\begin{align*}
    \partial_{\omega}\psi_{n,i}^{\pm}(\theta, \omega) & =  \sum_{k=1}^{M_n-1}\left((-kg'(\theta+ k\omega)) \prod_{j=1}^{k}(f^{-1})'(f(x_{-M_n+j})) \right)
\end{align*}
For any $\theta \in I_n$ and $\omega \in \mathcal{F}_n$, by equation \eqref{equation_backward curve I2_theta w slope} (in proposition \ref{proposition_structure I2_backward iterate curves}) it follows that the constant $\lambda$ satisfies the bound
\begin{align*}
    |\partial_{\omega}\psi_{n,i}^{\pm}(\theta, \omega)| 
    & < \left|\sum_{k=1}^{M_0} (-k\lambda^{-1/4+1/128})(\alpha_c^{1-(1+b)\beta_{n-1}})^{k}\right| + \left|\sum_{k=M_0+1}^{M_n} (-k\lambda^{1+1/128})(\alpha_c^{1-(1+b)\beta_{n-1}})^{k}\right| \\
     & < \sum_{k=1}^{M_0} k\lambda^{-1/4+1/128} \alpha_c^{k/2} + \alpha_c^{M_0/2}\sum_{k=1}^{\infty} (M_0+k)\lambda^{1+1/128} \alpha_c^{k/2} < 1,
\end{align*}
which proves the lemma.
\end{proof}
%------------------- PROP END ----------------------

%------------------- PROP BEGIN --------------------
\begin{proposition}\label{proposition_structure In_forward iterate curves}
         Let $n\geq 2$, $\omega\in \mathcal{F}_n$ and $(I_l(\omega))_{l=0}^{n}$ be the sequence of non-empty critical sets defined in \ref{definition_structure In_critical set In}. For each $i\in \{1,2\}$, the curves $\phi_{n,i}^{\pm}(\theta)$ satisfy the following conditions:
    
\begin{enumerate}[label=\alph*)]
\item \label{lemma item_structure In_forward iterate curves a} For each ${i'},{j'}\in \{1,2\}$, there exists $k \in \{1,2\}$ such that
\begin{align*}
\bigcup _{i=1}^{2} [\pi_2 \circ F^{M_n} (\theta- M_n\omega, c_i^{-}),\pi_2 \circ F^{M_n}(\theta+M_n\omega, c_i^{+})] & \subset C_k,
\end{align*}
for all $\theta \in I_n^{({i'},{j'})}$.

\item \label{lemma item_structure In_forward iterate curves b} For all $\theta \in I_n$ and $i \in \{1,2\}$, 
\begin{align*} 
\alpha_c^{b(M_n+1)}|c_i^+ - c_i^-| & \leq |\phi_{n,i}^+(\theta) - \phi_{n,i}^-(\theta)| \leq \left(\alpha_c^{1-(b+1)\beta_{n-1}} \right)^{M_n+1}|c_i^{+}-c_i^{-}|.
\end{align*}

\item \label{lemma item_structure In_forward iterate curves c} For all $\theta\in I_n$, 
\begin{align*}
    \alpha_c^{b(M_n+1)}|c_1^{+}-c_1^{-}| \leq \max_{i,j \in \{1,2\} }|\phi_{n,i}^{+}(\theta)-\phi_{n,j}^{-}(\theta)| \leq \left(\alpha_c^{1-(b+1)\beta_{n-1}} \right)^{M_n-M_0^{60}+1}.
\end{align*}

\item \label{lemma item_structure In_forward iterate curves d} For any $i\in \{1,2\}$ and $\theta \in I_n$, if $\theta \in \bigcup_{{i'}=1}^{2} I^{({i'},1)}_n$ then 
\begin{align*}
     -\lambda^{1+1/64} & < \partial_{\theta} \phi_{n,i}^{\pm}(\theta) < -\lambda^{1-1/64},
\end{align*}
and if $\theta \in \bigcup_{{i'}=1}^{2} I^{({i'},2)}_n$ then 
\begin{align*}
    \lambda^{1-1/64} & < \partial_{\theta} \phi_{n,i}^{\pm}(\theta) < \lambda^{1+1/64}.
\end{align*}

\item \label{lemma item_structure In_forward iterate curves e} For any $\theta \in I_n(\omega)$ and $i\in \{1,2\}$, 
\begin{align*}
    |\partial_{\omega}\phi_{n,i}^{\pm}(\theta, \omega)| < 1.
\end{align*}
\end{enumerate}

\end{proposition}

\begin{proof}
\textit{Proof of \ref{lemma item_structure In_forward iterate curves a}:} 
Let us fix $i',j'\in \{1,2\}$ and $\theta \in I_n^{(i',j')}(\omega)$. By the disjointness condition $(\mathcal{F}2)_n$, we have $(\theta-M_n \omega) \in \tilde{\Theta}_{n-1}$ and by $(\mathcal{F}1)_n$, it follows that 
$$\{\theta - k \omega\}_{k=1}^{M_n} \cap I_n = \emptyset.$$
For the fixed $\omega \in \mathcal{F}_n$, by definition \ref{definition_block n_sync condn parameter} there exists a constant $N_s < M_0^{60} < M_1$ such that $(\theta - (M_n-N_s+M_0+S_0)\omega) \in \bar{G}$. By the definition of $I_n^{(i',j')} \subset I_0^{(i',j')}$ and the measure of the sets $I_0^{(i',j')}$ (in proposition \ref{proposition_structure I0_I0 subsets}), we can also deduce that $(I_n^{(i',j')}-(M_n-N_s+M_0+S_0)\omega) \subset \hat{G}$. Since $N_s<M_1$, we have $(I_n^{(i',j')}-(M_n-N_s)\omega) \subset \Theta_{n-1}$. Then, by the condition $(g2)$ in definition \ref{definition_block n_sync condn set_g1 g2} and $(C1)_{n-1}$ (in proposition \ref{proposition_block n contraction}), it follows that 
    \begin{align}\label{equation_forward curves In_sync N_s}
     F^{N_s} \left( (I_n^{(i',j')}- M_n\omega) \times (C_1\cup C_2) \right) & \subset (I_n^{(i',j')}-(M_n-N_s)\omega) \times C_1.  
     \end{align}
By $(\mathcal{F}1)_n$, it is clear that
\begin{align*}
    I_n^{(i',j')} \cap \left( \bigcup_{l=0}^{n}\bigcup_{m=1}^{M_l+S_0-1} (I_l+m\omega) \right) = \emptyset,
\end{align*}
and thus, by $(S1)_{n-1}$ (in proposition \ref{proposition_block n sync' contraction}), there exists $k\in \{1,2\}$ such that 
 \begin{align*}
     \pi_2 \circ F^{M_1-N_s} \left( (I_1^{(i',j')}- (M_1-N_s)\omega) \times C_1 \right) & \subset \pi_2 \circ F^{M_1} \left( (I_1^{(i',j')}-M_1\omega) \times (C_1 \cup C_2) \right) \subset C_k. 
 \end{align*} 

\vspace{0.2cm}
\noindent\textit{Proof of \ref{lemma item_structure In_forward iterate curves b}:}  
Let $\theta \in I_n$ and $i\in \{1,2\}$. By disjointness condition $(\mathcal{F}2)_n$, we have
    $$(I_n-M_n\omega) \cap \left(\bigcup_{l=1}^{n-1}\bigcup_{|m|\leq 2M_{n-1}} (I_l+m\omega) \right) = \emptyset $$ 
    and thus, $\theta \in \Theta_{n-1}$. By condition $(\mathcal{F}1)_n$, it follows that $(\theta-(M_n-k)\omega) \notin I_n$, for all $0\leq k < M_n$. By the estimate in $(C3)_n$ (in proposition \ref{proposition_block n contraction}), for any $x_0 \in C_i$ we derive the estimate 
    \begin{align*}
    \alpha_c^{b(M_n+1)} & \leq |\partial_x(\pi_2 \circ F^{M_n+1}((\theta-M_n\omega),x_0))| \leq (\alpha_c^{1-(b+1)\beta_{n-1}})^{M_n+1}.
\end{align*}
Thus, for any $\theta \in I_n$ and $i \in \{1,2\}$, 
\begin{align*}
    \alpha_c^{b(M_n+1)} |c_i^{+}-c_i^{-}| \leq & |\phi_{n,i}^{+}(\theta)-\phi_{n,i}^{-}(\theta)| \leq (\alpha_c^{1-(b+1)\beta_{n-1}})^{M_n+1}|c_i^{+}-c_i^{-}|.
\end{align*}

\vspace{0.2cm}
\noindent\textit{Proof of \ref{lemma item_structure In_forward iterate curves c}:} For each fixed $i',j' \in \{1,2\}$, by equation \eqref{equation_forward curves In_sync N_s}, there exists $N_s < M_0^{60}<M_1$ such that 
    \begin{align*}
        \pi_1 \circ F^{N_s}(I_n^{(i',j')}-M_n\omega \times C_1\cup C_2) & \subset \Theta_{n-1} \text{ and }\\
    \pi_2 \circ F^{N_s}(I_n^{(i',j')}-M_n\omega \times C_1\cup C_2) & \subset C_1,
    \end{align*}
where $k=(M_n-N_s)$ is the smallest integer for which $(\theta-(M_n-N_s)\omega+k\omega) \in I_n$. By $(C2)_{n-1}$ (in proposition \ref{proposition_block n contraction}), this gives the estimate:
\begin{align*}
    \alpha_c^{b(M_n-N_s+1)} & \leq |\partial_x(\pi_2 \circ F^{M_n-N_s+1}((\theta-(M_n-N_s)\omega),x))| \leq (\alpha_c^{1-(1+b)\beta_{n-1}})^{M_n-N_s+1},
\end{align*}
for all $\theta \in I_{n}^{(i',j')}$ and $x\in C_1$. Then, by $(Sc1)_{n-1}$, it follows that
\begin{align*}
     \alpha_u^{-M_n} |c_1^{+}-c_1^{-}| & \leq \max_{i,j \in \{1,2\}} |\phi_{n,i}^{+}(\theta)- \phi_{n,j}^{-}(\theta)| 
     \leq (\alpha_c^{1-(1+b)\beta_{n-1}})^{M_n-N_s+1}. 
\end{align*}
Since $N_s<M_0^{60}$, for all $i',j' \in \{1,2\}$, we derive the following estimate for all $\theta \in I_n$:
\begin{align*}
    \alpha_c^{bM_n} |c_1^{+}-c_1^{-}| & \leq \max_{i,j \in \{1,2\}} |\phi_{n,i}^{+}(\theta)- \phi_{n,j}^{-}(\theta)| 
     \leq (\alpha_c^{1-(1+b)\beta_{n-1}})^{M_n-M_0^{60}+1}.
\end{align*}

\vspace{0.2cm}
\noindent\textit{Proof of \ref{lemma item_structure In_forward iterate curves d}:} For any $\theta \in I_n$ and $x_0 = c_i^{\pm}$, we have the expression
\begin{align}\label{equation_forward curves n_theta slope}
      \partial_{\theta} \phi_{n,i}^{\pm}(\theta) = g'(\theta) + \sum_{k=1}^{M_n}\left( g'(\theta-k\omega)\prod_{j=0}^{k-1}f'(x_{M_n-j}) \right).
\end{align}
For any $(\theta-M_n\omega) \in (I_n-M_n\omega) \subset \Theta_{n-1}$ and $x_0 = c_i^{\pm}$, by $(C1-2)_{n-1}$ (in proposition \ref{proposition_block n contraction}) we derive the estimate
\begin{align*}
    \left| \prod_{j=0}^{k-1} f'(x_{M_n-j}) \right| < (\alpha_c^{1-(1+b)\beta_{n-1}})^k,
\end{align*}
for $k=1,\dots ,M_n$. By the definition of critical set $I_n$, we have $(I_n-k\omega) \subset (I_0-k\omega)$, for $k=1,\dots M_0$. Then, for any $\theta \in I_n$ and $k=1,\dots , M_0$, we consider the estimate $|g'(\theta-k\omega)|< \lambda^{-1/4+1/128}$. This gives us
\begin{align*}
    \left|\sum_{k=1}^{M_0} g'(\theta-k\omega) \prod_{j=0}^{k-1}f'(x_{M_n-j}) \right| < \lambda^{-1/4+1/128} \alpha_c^{1/2} \frac{1-\alpha_c^{M_0/2}}{1-\alpha_c^{1/2}}.
\end{align*}
and 
\begin{align*}
    \left|\sum_{k=M_0+1}^{M_n} g'(\theta-k\omega) \prod_{j=0}^{k-1}f'(x_{M_n-j}) \right| & < \lambda^{1+1/128} \alpha_c^{M_0/2} \frac{\alpha_c^{1/2}(1-\alpha_c^{(M_n-M_0)/2})}{1-\alpha_c^{1/2}}.
\end{align*}
Similar to our estimates for $\psi^{\pm}_{n,i}$ in equation \eqref{equation_backward curve In_theta slope}, we conclude that 
\begin{align*}
    \left|\sum_{k=1}^{M_n} g'(\theta-k\omega) \prod_{j=0}^{k-1}f'(x_{M_n-j}) \right| < 1.
\end{align*}
By definition of the set $I_n$, we have $I_n^{(i,j)} \subset I_0^{(i,j)}$, for all $i,j\in \{1,2\}$. Then, by estimates for $g'(\theta)$ for $\theta \in  \bigcup_{i,j=1}^{2}I_0^{(i,j)}$ in proposition \ref{proposition_structure I0_I0 subsets} \ref{lemma item_structure I0_I0 components_c}, we obtain 
\begin{align*}
    \lambda^{1-1/64} & < \partial_{\theta} \phi_{n,i}^{\pm}(\theta) < \lambda^{1+1/64},
\end{align*}
for any $\theta \in \bigcup_{i=1}^{2}I_n^{(i,2)}$, and 
\begin{align*}
    -\lambda^{1+1/64} & < \partial_{\theta} \phi_{n,i}^{\pm}(\theta) < -\lambda^{1-1/64},
\end{align*}
for any $\theta \in \bigcup_{i=1}^{2}I_n^{(i,1)}$.

\vspace{0.2cm}
\noindent\textit{Proof of \ref{lemma item_structure In_forward iterate curves e}:} For the curves $\phi_{n,i}^{\pm}$, we have the expression  
\begin{align*}
    \partial_{\omega}\phi_{n,i}^{\pm}(\theta, \omega) & =  \sum_{k=1}^{M_n}\left( (-k)g'(\theta-k\omega) \prod_{j=0}^{k-1} f'(x_{M_n-j}) \right).
\end{align*}
Similar to the estimate in the proof of proposition \ref{proposition_structure In_forward iterate curves} \ref{lemma item_structure In_forward iterate curves d}, we prove that 
\begin{align*}
    |\partial_{\omega}\phi_{n,i}^{\pm}(\theta, \omega)| 
    & < \left|\sum_{k=1}^{M_0} (-k\lambda^{-1/4+1/128})\alpha_c^{k/2} + \sum_{k=M_0+1}^{M_n} (-k\lambda^{1+1/128})\alpha_c^{k/2}\right| \\
    & < \sum_{k=1}^{M_0} k\lambda^{-1/4+1/128} \alpha_c^{k/2} + \alpha_c^{M_0/2}\sum_{k=1}^{\infty} (M_0+k)\lambda^{1+1/128} \alpha_c^{k/2} < 1.
\end{align*}

\end{proof}
%------------------- PROP END -------------------

%%%%%%%%%%%% EMPTY vs NONEMPTY DICHOTOMY FOR In+1 %%%%%%%%%%%%%%%%%

\subsubsection{The empty and non-empty cases for $I_{n+1}$}\label{subsubsection_empty vs nonempty In+1}

The following proposition proves that, if the critical set $I_2(\omega)$ is non-empty for some $\omega \in \mathcal{F}_n$, then the critical set $I_n(\omega)$ must also be non-empty.

%------------------- PROP BEGIN -------------------
\begin{proposition}\label{proposition_structure In_In empty nonempty}
       For each $n\geq 2$, $\omega \in \mathcal{F}_n$ and ${i'},{j'} \in \{1,2\}$, if $\displaystyle I_n^{(i',j')}(\omega) \neq \emptyset$ then $\displaystyle I_{n+1}^{(i',j')} \neq \emptyset$.\\
       Equivalently, $\displaystyle I_{n+1}^{({i'},{j'})}(\omega) \neq \emptyset$ if and only if $\displaystyle \omega \in \mathcal{F}_n \cap \tilde{\mathcal{F}}_1^{({i'},{j'})}$.
\end{proposition}

\begin{proof}
Say $I_n^{(i', j')}(\omega) \neq \emptyset,$ for some $i',j' \in \{1,2\}$. By definition $\ref{definition_structure In_ set In prime}$ and \ref{definition_structure In_critical set In}, we must then have $\tilde{I}_n^{(i',j')}(\omega) \neq \emptyset$. By $(\tilde{\mathcal{I}}1)_{n}$, recall the definitions 
    \begin{align*}
        \tilde{I}_{n}^{(i',j')}(\omega) & = I_1^{(i',j')}(\omega) \cap \pi_1  \left( F^{M_{n-1}}(A_{n-1}) \cap F^{-M_{n-1}}(B_{n-1}) \right) \\
        \tilde{I}_{n+1}^{(i',j')}(\omega) & = I_1^{(i',j')}(\omega) \cap \pi_1  \left( F^{M_{n}}(A_{n}) \cap F^{-M_{n}}(B_{n}) \right).
    \end{align*}
By proposition \ref{proposition_structure I0_I0 subsets} \ref{lemma item_structure I0_I0 components_d}, since the sets $\{I_0^{(i,j)}\}_{i,j=1}^{2}$ are disjoint, $\tilde{I}_{n+1}^{(i',j')}(\omega)$ is equivalently given by
\begin{align*}
    \tilde{I}_{n+1}^{(i',j')}(\omega) & = \pi_1 \circ \left( F^{M_{n}}((I_n^{(i', j')}-M_n\omega) \times (C_1 \cup C_2))  \cap F^{-M_{n}}((I_n^{(i', j')}+M_n\omega) \times (E_1\cup E_2)) \right). 
\end{align*}
By lemma \ref{lemma_structure In_match Ai Bi}, 
\begin{align}\label{equation_structure In_nested transverse curves_1}
    F^{M_n+1}(A_n) \subset F^{M_{n-1}+1}(A_{n-1}) \\
    F^{M_n-1}(B_n) \subset F^{M_{n-1}-1}(B_{n-1}),
    \stepcounter{equation}\tag{\theequation} \label{equation_structure In_nested transverse curves_2}
\end{align}
where the sets $F^{M_{n-1}}(A_{n-1})$ and $F^{-M_{n-1}}(B_{n-1})$ are bounded by the curves $\psi_{n-1,i}^{\pm}, \phi_{n-1,i}^{\pm}: I_{n-1} \rightarrow \mathbb{T}$. Since $F$ is an orientation preserving diffeomorphism, the curves in the set $\{\phi_{n-1,i}^{\pm}:i=1,2\}$ are parallel to each other and similarly the curves in the set $\{\psi_{n-1,i}^{\pm}: i=1,2\}$ are parallel to each other. For any $\theta \in I_{n-1}^{(i',j')}$, the estimates of $\partial_{\theta}\phi_{n-1,i}^{\pm}(\theta)$ and $\partial_{\theta}\psi_{n-1,i}^{\pm}(\theta)$ established in proposition \ref{proposition_structure In_forward iterate curves} \ref{lemma item_structure In_forward iterate curves d} and proposition \ref{proposition_structure In_backward iterate curves} \ref{lemma item_structure In_backward iterate curves d} show that the curves $\phi_{n-1,i}^{\pm}$ and $\psi_{n-1,i}^{\pm}$ intersect transversally over the interval $I_{n-1}^{(i',j')}$. By definition \ref{definition_structure In_critical set In}, since $I_{n}^{(i',j')}$ is the smallest interval containing $\tilde{I}_{n}^{(i',j')}$, the end points of $I_{n}^{(i',j')}$ must exactly be the points of intersection of the curves $\phi_{n-1,i_1}^{\pm}$ and $\psi_{n-1,i_2}^{\pm}$, for some $i_1, i_2 \in \{1,2\}$. Since $\omega \in \mathcal{F}_n$, by the estimates of $\partial_{\theta}\phi_{n,i}^{\pm}$ and $\partial_{\theta}\psi_{n,i}^{\pm}$ and equations (\ref{equation_structure In_nested transverse curves_1} - \ref{equation_structure In_nested transverse curves_2}), the curves $\phi_{n,i}^{\pm}$ and $\psi_{n,i}^{\pm}$ must also intersect transversally over the interval $I_{n}^{(i', j')}$. This proves that $\tilde{I}_{n+1}^{(i', j')}(\omega) \neq \emptyset$ and consequently the set $I_{n+1}^{(i',j')}(\omega)$ are non-empty.   

Say $I_n^{(i', j')}(\omega)$ is non-empty for some $i',j' \in \{1,2\}$. Since $\displaystyle \left(I_l^{(i',j')}(\omega)\right)_{l=1}^{n}$ are a sequence of nested intervals, the set $I_2^{(i',j')}(\omega)$ is also non-empty. Then, by proposition \ref{proposition_structure I2_I2 empty nonempty}, it follows that $\omega \in \tilde{\mathcal{F}}_1^{(i', j')}$.

\end{proof}
%------------------- PROP END -------------------

%%%%%%%%%%%%%%%%%%%%% ESTIMATES FOR NON-EMPTY %%%%%%%%%%%%%%%%%%%%%%%%%%%%%%%%%%%%%%%%%%%%%%%

\subsubsection{Estimates for a non-empty $I_{n+1}(\omega)$}\label{subsubsection_estimate nonempty In+1}

For any $n\geq 2$ and $\omega \in \mathcal{F}_n$ for which the critical set $I_{n+1}(\omega)$ is non-empty, the following proposition describes its geometric structure. As a consequence of the following estimate, one can infer that the disjointness condition $(\bar{\mathcal{F}}1)_n$, which was introduced in definition \ref{definition_block n_disjointness condn}, implies the condition $(\mathcal{F}1)_n$.

%------------------- PROP BEGIN -------------------
\begin{proposition}\label{proposition_structure In+1_critical set I_n+1}
    Let $\lambda$ be sufficiently large and let ${i'},{j'} \in \{1,2\}$. For each $n\geq 2$ and $\omega \in \tilde{\mathcal{F}}^{(i',j')}_1 \cap \mathcal{F}_n$, we have the following estimates for the set $I_{n}^{({i'},{j'})}(\omega)$:
    \begin{enumerate}[label=\alph*)]
        \item \label{lemma item_structure In+1_critical set I_n+1 a} 
        $ \lambda^{-1-1/32} \alpha_{c}^{3b\tilde{M}_n} < |I^{({i'},{j'})}_{n+1}(\omega)| < \alpha_c^{\tilde{M}_n/3}\lambda^{-1+1/32}.$
        \item \label{lemma item_structure In+1_critical set I_n+1 b} $|\partial_{\omega}I^{({i'},{j'})}_{n+1}(\omega)| < 1/4 $.
    \end{enumerate}
\end{proposition}

\begin{proof}
%\textit{Proof of \ref{lemma item_structure In+1_critical set I_n+1 a}:} 
For $\omega \in \mathcal{F}_1^{(i',j')}\cap \mathcal{F}_n$, the critical $I_{n+1}^{(i',j')}(\omega)$ is non-empty due to proposition \ref{proposition_structure In_In empty nonempty}.
From the definition of the set $I_{n+1}^{(i',j')}$, in definition \ref{definition_structure In_ set In prime} and \ref{definition_structure In_critical set In}, recall that
    \begin{align*}
        \pi_2 \circ \left( F^{M_n}((I^{(i',j')}_n-M_n\omega) \times C_1\cup C_2) \cap F^{-M_n}((I^{(i',j')}_n+M_n\omega) \times E_1\cup E_2) \right) \subset I_{n+1}^{(i',j')}
    \end{align*}
    and $I_{n+1}^{(i',j')}$ is the smallest such interval. Since $f_{\theta}$ is an orientation preserving diffeomorphism, for each $\theta \in \mathbb{T}$, and $I_{n}^{(i',j')}$ is an interval, there exists distinct $i_1,i_2 \in \{1,2\}$ such that 
    \begin{align}\label{equation_structure In_backward curves}
        \psi_{n,i_1}^{-}(\theta) < \psi_{n,i_1}^{+}(\theta) < \psi_{n,i_2}^{-}(\theta) < \psi_{n,i_2}^{+}(\theta),
    \end{align}
    for all $\theta \in I_n^{(i',j')}$. Similarly, there exists distinct $i'_1, i'_2 \in \{1,2\}$ such that
\begin{align}\label{equation_structure In_forward curves}
    \phi_{n,i'_1}^{-}(\theta) < \phi_{n,i'_1}^{+}(\theta) < \phi_{n,i'_2}^{-}(\theta) < \phi_{n,i'_2}^{+}(\theta),
\end{align}
for all $\theta \in I_n^{(i',j')}$. Let us define sets
\begin{align*}
    \Tilde{B}_n & := \left\{(\theta, x): \theta \in I_n^{(i',j')}, x \in [\psi_{n,i_1}^{-}(\theta), \psi_{n,i_2}^{+}(\theta)] \right\} \\
    \Tilde{A}_n & := \left\{(\theta, x): \theta \in I_n^{(i',j')}, x \in [\phi_{n,i'_1}^{-}(\theta), \phi_{n,i'_2}^{+}(\theta)] \right\}.
\end{align*}
By equations \eqref{equation_structure In_backward curves} and \eqref{equation_structure In_forward curves}, we observe that $F^{M_n}(A_n) \subset \Tilde{A}_n$ and $F^{-M_n}(B_n) \subset \Tilde{B}_n$ and hence, 
$$\left( I_n^{(i',j')} \cap (\pi_1 (F^{M_n}(A_n) \cap F^{-M_n}(B_n))) \right) \subset \pi_1 (\tilde{A}_n \cap \tilde{B}_n).$$
This gives us $\tilde{I}_{n+1}^{(i',j')} \subset \pi_1 (\tilde{A}_n \cap \tilde{B}_n)$. In propositions \ref{proposition_structure In_backward iterate curves} \ref{lemma item_structure In_backward iterate curves d} and \ref{proposition_structure In_forward iterate curves} \ref{lemma item_structure In_forward iterate curves d}, we derived the estimates $|\partial_{\theta}\phi_{n,i}^{\pm}|> \lambda^{1-1/128}$ and $|\partial_{\theta}\psi_{n,i}^{\pm}| < 1$, for all $\theta \in I_n^{(i',j')}$. This shows that the curves $\phi_{n,i}^{\pm}$ and $\psi_{n,i}^{\pm}$ intersect transversally  over the interval $I_n^{(i',j')}$, which yields a closed interval $\pi_1 (\tilde{A}_n \cap \tilde{B}_n)$ for which the end points $a$ and $b$ satisfy the conditions:
\begin{align}\label{equation_structure I2_transverse intersection 1}
    \phi_{n,i'_1}^{-}(a) & = \psi_{n,i_2}^{+}(a) \\
    \phi_{n,i'_2}^{+}(b) & = \psi_{n,i_1}^{-}(b).
    \stepcounter{equation}\tag{\theequation} \label{equation_structure I2_transverse intersection 2}
\end{align}
This also shows that there exists no smaller interval than $\pi_1 (\tilde{A}_n \cap \tilde{B}_n)$ which contains $\tilde{I}^{(i,j)}_{n+1}$, and this give us 
$$I_{n+1}^{(i,j)} = \pi_1 (\tilde{A}_n \cap \tilde{B}_n). $$
For the transverse intersection, we derive the estimates
\begin{align*}
    \frac{2\inf\{|\phi_{n,i'_1}^{-}(\theta)-\phi_{n,i'_2}^{+}(\theta)|: \theta \in I_n^{(i',j')}\}|}{|\partial_{\theta}\phi_{n}^{\pm}|_{sup}+|\partial_{\theta}\psi^{\pm}_{n}|_{sup}} 
    & < |I_{n+1}^{(i',j')}| < \frac{2\sup\{|\psi_{n,i_1}^{-}(\theta)-\psi_{n,i_2}^{+}(\theta)|: \theta \in I_n^{(i',j')}\}}{|\partial_{\theta}\phi^{\pm}_{n}|_{inf}-|\partial_{\theta} \psi_{n}^{\pm}|_{sup}},
\end{align*}
where $|\partial_{\theta}\phi^{\pm}_{n}|_{sup}$ denotes the quantity $\sup \{ |\partial_{\theta}\phi^{\pm}_{n,i}|_{sup}: i=1,2\}$, and the analogous notations for the curves $\{\psi_{n,i}^{\pm}\}$ and for infimum are defined similarly. By the estimates in proposition \ref{proposition_structure In_backward iterate curves} and \ref{proposition_structure In_forward iterate curves}, we then obtain
\begin{align*}
    \lambda^{-1-1/64} \alpha_{c}^{b(M_n+1)} & < |I^{(i',j')}_{n+1}(\omega)| < \alpha_c^{M_n/3}\lambda^{-1+1/64}.
\end{align*}

For $\omega \in \mathcal{F}_n \cap \mathcal{F}_1^{(i',j')}$, by equation (\ref{equation_structure I2_transverse intersection 1} - \ref{equation_structure I2_transverse intersection 2}), the end points $a(\omega), b(\omega)$ of the set $I_{n+1}^{(i',j')}$ satisfy the conditions
\begin{align*}
     \phi_{n,i'_1}^{-}(a(\omega), \omega) & = \psi_{n,i_2}^{+}(a(\omega), \omega) \\
    \phi_{n,i'_2}^{+}(b(\omega), \omega) & = \psi_{n,i_1}^{-}(b(\omega), \omega).
\end{align*}
Since $|\partial_{\theta} \phi_{n,i}^{\pm}| > \lambda^{1-1/128}$, for all $i \in \{1,2\}$, it follows by the implicit function theorem that the functions $a(\omega)$ and $b(\omega)$ are continuously differentiable and we have the expression
\begin{align*}
    a'(\omega) & = \frac{\partial_{\omega}\phi_{n,i_1'}^{-}(a(\omega),\omega) - \partial_{\omega}\psi_{n,i_2}^{+}(a(\omega),\omega)}{\partial_{\theta}\psi_{n,i_2}^{+}(a(\omega),\omega) - \partial_{\theta}\phi_{n,i_1'}^{-}(a(\omega),\omega)}.
\end{align*}
Then, by the estimates in proposition \ref{proposition_structure In_backward iterate curves} and \ref{proposition_structure In_forward iterate curves}, we obtain $|a'(\omega)| < \lambda^{-1+1/32}$. By symmetry, we also have $ |b'(\omega)| < \lambda^{-1+1/32}$and this proves the bound $|\partial_{\omega}I_{n+1}^{(i',j')}(\omega)| < \lambda^{-1+1/32}$.

\end{proof}
%------------------- PROP END -------------------

%%%%%%%%%%%%%%%%%%%%%%%%%%%%%%%%%%%%%%%%%%%%%%%%%%%%%%%%%%%

%%%%%%%%%%%%%%%%%%%%% PARAMETER SPACE %%%%%%%%%%%%%%%%%%%%%%%%%%%%%%%%%%%%%%
\section{Parameter Space}\label{section_parameter space}
The goal of this section is to construct the parameter sets $\Omega_1$ and $\Omega_2$ of positive measure, as introduced in Theorem \ref{Theorem SNA}. That is, for each $\omega \in \Omega_1$ the map $F_{\omega}$ will be proved to admit measurable discontinuous invariant graphs $\Phi$ and $\Psi$, and for each $\omega \in \Omega_2$ the map $F_{\omega}$ will be proved to admit two continuous invariant graphs. The exact formulation and estimation of the sets $\Omega_1$ and $\Omega_2$  is contained in subsection \ref{subsection_parameter set Omega1 and Omega2} . 

The construction of $\Omega_1$ and $\Omega_2$ builds on the sets of parameters $\mathcal{F}_0$, $\mathcal{F}^{s}(G_2)$, $\tilde{\mathcal{F}}_1^{(1,1)}$, $\tilde{\mathcal{F}}_1$ and ${\mathcal{F}}_n$ for $n\geq 2$, which were defined in sections \ref{section_Mechanism of slow-returns} and \ref{section_structure of critical sets}. The constructions and estimates related to these parameter sets, contained in subsections (\ref{subsection_Parameters satisfying synchronisation condition} - \ref{subsection_Parameter exclusion for the construction of Fn}), are thus central to the construction of $\Omega_1$ and $\Omega_2$. Before we proceed, let us state the following lemma which contains measure theoretic estimates for a set of parameters satisfying an abstract slow-return condition. This will be used widely in the inductive parameter exclusions in subsections (\ref{subsection_Parameter exclusion for the construction of F0} - \ref{subsection_Parameter exclusion for the construction of Fn}).

%%%%%%%%%%%%%%%%%%%%% PARAMETER SPACE MAIN  %%%%%%%%%%%%%%%%%%%%%%%%%%%%%%

%%%%%%%%%%%%%%%%%%%%% ABSTRACT LEMMA - PARAMETER SPACE %%%%%%%%%%%%
%------------------- LEMMA BEGIN -------------------
\begin{lemma}\label{lemma_parameters_prelim_disjoint condn}
Assume that, for each $\omega \in \mathbb{T}$, there exists corresponding intervals $I^1(\omega), \dots , I^{\mathcal{N}}(\omega)$ in $\mathbb{T}$, each of length at most $\delta>0$ and satisfying $|\partial_{\omega}I^{i}(\omega)| < 1/4$, for $i \in \{1,\dots,\mathcal{N}\}$. Then, for $M>0$ and $\varepsilon > 0$, the set 
\begin{align*}
    \left\{ \omega \in \mathbb{T}: dist \left( \bigcup_{i=1}^{\mathcal{N}} I^i(\omega), \bigcup_{m=1}^{M} \bigcup_{i=1}^{\mathcal{N}} (I^i(\omega) + m\omega) \right) < \varepsilon \right\}
\end{align*}
has measure $\leq 4 \mathcal{N}^2 M (\delta + \varepsilon)$ and consists of at most $\mathcal{N}^2 M^2$ components.
\end{lemma}
\begin{proof}
    For the proof, we refer to Lemma 7.1 in \cite{BJERKLOV_2005}.
\end{proof}
%------------------- LEMMA END -------------------

%%%%%%%%%%%%%%%%%%%%% F0 PARAMETER SPACE %%%%%%%%%%%%

\subsection{Parameter exclusion for the construction of $\mathcal{F}_0$}\label{subsection_Parameter exclusion for the construction of F0}

This subsection consists of the construction of the set of parameters $\mathcal{F}_0$ given by 
\begin{align*}
   \mathcal{F}_0(M_0,f,\lambda) = \left\{ \omega:  d\left(3I_0, \bigcup_{m=1}^{2K_0M_0}3(I_0+m\omega)\right) > 0 \right\},
\end{align*}
as in definition \ref{definition_block 0_disjointness condn}, where the constants are $\lambda>0$, $M_0 = \lambda^{1/16}$ and $K_0 = 2^8(b+1)$. 

%------------------- PROP BEGIN -------------------
\begin{proposition}\label{proposition_parameters_block I0_lemma}
The set
\begin{equation}
 \mathbb{T} \setminus \mathcal{F}_0(M_0,f,\lambda) = \left\{\omega\in \mathbb{T}: dist\left( 3I_0, \bigcup_{m=1}^{2K_0M_0}3(I_0+m\omega) \right) \leq 0 \right\}    
\end{equation} 
has measure at most $16K_0M_0|I_0|$ and consists of at most $(2K_0M_0)^2$ components.
\end{proposition}
\begin{proof}
    Note that the set $I_0$ is independent of $\omega$. For any two intervals $J_1$ and $J_2$ each of length $\delta$ and $m\neq 0$, the set
    $$ \{\omega\in \mathbb{T}: dist(J_1,J_2+m\omega)<\varepsilon\} $$
consists of at most $m$ intervals, each of length at most $\displaystyle \frac{2\delta+2\varepsilon}{|m|}.$ Thus, for the set 
\begin{equation*}
    \bigcup_{m=1}^{M_0}\{\omega \in \mathbb{T}: dist(3I_0, 3(I_0+m\omega))< 0\}
\end{equation*}
we have the estimates 
\begin{align*}
    \sum_{m=1}^{2K_0M_0}|\{\omega\in \mathbb{T}: dist(3I_0, 3(I_0+m\omega))< 0\}| & \leq \sum_{m=1}^{2K_0M_0} \frac{6m|I_0|}{m} < 16K_0M_0|I_0|,
\end{align*}
and 
\begin{align*}
   \sum_{m=1}^{2K_0M_0}|m| < (2K_0M_0)^2 = 4K_0^2M_0^2.
\end{align*}
This proves the lemma.
\end{proof}

%------------------- PROP END -------------------

%%%%%%%%%%%%%%%%%%%%% tilde Fs(G2) PARAMETER SPACE %%%%%%%%%%%%
\subsection{Parameters satisfying synchronisation condition}\label{subsection_Parameters satisfying synchronisation condition}

In this subsection, we construct a measurable subset of $\mathcal{F}^s(G_2)$, which has a positive measure. Recall the set $G_2 \subset I_0$ in definition \ref{definition_structure I1_sync set Gs}, given by
\begin{align}\label{equation_parameters_sync_G set}
    {G_2} & = \left( \theta^{(2,1)}+3\lambda^{-1+1/64}, \theta^{(1,1)}-3\lambda^{-1+1/64} \right) \bigcup \left( \theta^{(1,2)}+3\lambda^{-1+1/64}, \theta^{(2,2)}-3\lambda^{-1+1/64} \right).
\end{align}
By proposition \ref{proposition_structure I1_sync set Gs},  we have the estimate $|{G_2}| > \lambda^{-1/2-1/8} > \lambda^{-1}$ and the set satisfies $G_2 \subset I_0\setminus I_1(\omega)$, for every $\omega \in \mathcal{F}_0$. Then, by definition \ref{definition_block n_sync condn parameter}, there exists a well-defined set $\mathcal{F}^s(G_2)$ given by
\begin{align}\label{equation_parameters_sync condn parameter}
    \mathcal{F}^s(G_2) = \left\{ \omega\in \mathbb{T}: \left\{\theta+k\omega \right\}_{k=0}^{M_0^{56}} \bigcap G_2 \neq \emptyset, \text{ for all } \theta \in \mathbb{T} \right\}.
\end{align}

\begin{definition}[Synchronisation parameters]\label{definition_parameters_sync_FsG}
    For a fixed $\displaystyle M' := \left[\frac{1}{8|I_0|(|{G_2}|/3)^{3}}\right]$, let us define 
    \begin{align}\label{equation_parameters_sync condn parameter_new}
        \bar{\mathcal{F}}^s(G_2)  :=  \left\{ \omega\in \mathbb{T}: |m\omega| \geq \frac{8|I_0|}{m^2},\ \text{ for all } \ 0\leq m\leq \frac{1}{8|I_0|(|{G_2}|/3)^{3}}\right\}.
    \end{align}
\end{definition}

The following lemma proves that the number theoretic condition on $\omega$ that is specified in equation \eqref{equation_parameters_sync condn parameter}, corresponds to the condition of distribution of orbits in equation \eqref{equation_parameters_sync condn parameter_new}. 

%------------------- LEMMA BEGIN -------------------
\begin{lemma}\label{lemma_parameters_sync_lemma 0}
    For any fixed $0<K < 1$, $0<\delta < 1$ and $\theta\in \mathbb{T}$, if 
    \begin{align}\label{equation_parameters_sync_lemma FsG}
        \omega\in \left\{ \omega\in \mathbb{T}: d(0,m\omega) \geq \frac{K}{m^2},\ \text{ for all } 0 < m\leq \frac{1}{K\delta^3}\right\},
    \end{align}
    then the set 
    \begin{align*}
        \left\{ \theta+m\omega: 0 \leq m\leq \frac{1}{K\delta^3}\right\}
    \end{align*}
    is $\delta$ dense in $\mathbb{T}$.
\end{lemma}
\begin{proof}
Let us fix $w\in \mathbb{T}$ that satisfies the condition in equation \eqref{equation_parameters_sync_lemma FsG}. By Dirichlet's approximation theorem, there exists an integer $ k\leq {1}/{\delta}$ such that 
        \begin{equation*}
            d(0,k\omega) < \delta. 
        \end{equation*}
        For one such fixed $k$, let us take 
        $$ \alpha:= d(0,k\omega).$$
         From the conditions on $\omega$ and since $k\leq {1}/{\delta} < 1/(K\delta^3)$, 
        \begin{align*}
        d(0,k\omega) & \geq \frac{K}{k^2} \geq K\delta^2 >0
        \implies 0 < K\delta^2 < \alpha < \delta.        
        \end{align*}
Since $\alpha > 0$, the points in the set $\displaystyle \{nk\omega\}_{n=1}^{\left[ {1}/{\alpha}\right]+1}$ are all distinct and between two consecutive points in the set we have the distance to be at most $\alpha$. Then, since $K\delta^2 < \alpha < \delta$, it follows that for fixed $\displaystyle N := \left[ \frac{1}{K\delta^2}\right]$, the set $\displaystyle \{nk\omega\}_{n=1}^{N}$ is $\delta$ dense in $\mathbb{T}$. This proves that the set $\left\{ m\omega \right\}_{m=1}^{M'},$ for $\displaystyle M' = \frac{1}{K\delta^3}$, is $\delta$ dense in $\mathbb{T}$. 
        
\end{proof}
%------------------- LEMMA END -------------------

%------------------- PROP BEGIN -------------------
\begin{proposition}\label{proposition_parameters_sync_lemma 1}
For the set $\bar{\mathcal{F}}^s({G_2})$ defined in definition \ref{definition_parameters_sync_FsG}, we have $\bar{\mathcal{F}}^s({G_2}) \subset \mathcal{F}^s({G_2})$. The measure of $\bar{\mathcal{F}}^s({G_2})$ is at least $1 - 32|I_0| $ and the set consists of at most $M_0^{114}$ intervals.
\end{proposition}

\begin{proof}    
For the fixed set ${G_2}$, by proposition \ref{proposition_structure I1_sync set Gs} \ref{lemma item_structure I1_sync set Gs_b} we have estimate $|{G_2}|/3 = \lambda^{-1/2-1/8}/3 > M_0^{-16}$. By the estimate $|I_0| = 2\lambda^{-1/4} = 2M_0^{-4}$, we then have $8|I_0|(|{G_2}|/3)^3 > M_0^{-56}$.

By the definition of the sets $\mathcal{F}^s({G_2})$ and $\bar{\mathcal{F}}^s(G_2)$, for the choice of constants $K = 8|I_0|$, $\displaystyle M' = \left[\frac{1}{8|I_0|(|{G_2}|/3)^{3}}\right] < M_0^{56}$ and $\delta = |{G_2}|/3$ we have 
\begin{align*}
    \mathcal{F}^s(G_2) & \supset \{\omega \in \mathbb{T}: \left(\theta + m\omega\right)_{m=0}^{M'} \text{ is } |{G_2}|/3\text{ dense in } \mathbb{T}, \text{ for all } \theta \in \mathbb{T}\}.
\end{align*}
Then, by lemma \ref{lemma_parameters_sync_lemma 0}, it follows that
\begin{align*}
    \mathcal{F}^s(G_2)  & \supset  \left\{ \omega\in \mathbb{T}: |m\omega| \geq \frac{8|I_0|}{m^2},\ \text{ for all } \ 1 \leq m\leq \frac{1}{8|I_0|(|{G_2}|/3)^{3}}\right\} = \bar{\mathcal{F}}^s(G_2).
\end{align*}
This gives us the bound
\begin{align*}
      |{\mathcal{F}}^s(G_2)| & \geq 1 - \left| \left\{ \omega\in \mathbb{T}: |m\omega| < \frac{8|I_0|}{m^2},\ \text{ for some } 1 \leq m\leq \frac{1}{8|I_0|(|{G_2}|/3)^{3}}\right\} \right|.
\end{align*}
Since $\displaystyle \frac{1}{8|I_0|(|{G_2}|/3)^{3}} \leq M_0^{56}$, 
\begin{align*}
    \left| \left\{ \omega\in \mathbb{T}: |m\omega| < \frac{8|I_0|}{m^2},\ \text{ for some } 1 \leq m\leq M_0^{56}\right\} \right| & \leq  \sum_{m=1}^{M_0^{56}} \sum_{k=0}^{m-1} \frac{16|I_0|}{m^3} < 32|I_0|.
\end{align*}
This gives the estimate $ |\bar{\mathcal{F}}^s(G_2)| > 1 - 32|I_0|$. Similarly, by lemma \ref{lemma_parameters_sync_lemma 0} we have
\begin{align*}
 \# \bar{\mathcal{F}}^s(G_2) & = \# \left\{ \omega\in \mathbb{T}: |m\omega| < \frac{8|I_0|}{m^2},\ \text{ for some } 0\leq m\leq \frac{1}{8|I_0|(|{G_2}|/3)^{3}}\right\} \\
 & \leq \sum_{m=1}^{M_0^{56}} m < M_0^{114}.
\end{align*}

\end{proof}
%------------------- PROP END -------------------

%%--------------------- \tilde{F}_1 parameters -------------------------
\subsection{Parameter exclusion for the construction of $\tilde{\mathcal{F}}_1$}\label{subsection_Parameter exclusion for the construction of F1}

%%%%%%%%%%%%%%%%%%%%% F1 PARAMETER SPACE %%%%%%%%%%%%%%%%%%%%
This subsection contains the construction and estimation of the parameter set $\tilde{\mathcal{F}}_1$, as a subset of the parameter set ${\mathcal{F}}_0$. For the fixed constants $\tilde{M}_1 = \alpha_c^{-{M}_{0}/8}$ and $K_1 = (b+1)2^{9}$, the set $\tilde{\mathcal{F}}_1$ is as in definition \ref{definition_block n_disjointness condn} and is given by
\begin{align*}
    \tilde{\mathcal{F}}_1 = \left\{\omega\in \mathbb{T}: \omega \in  \tilde{\mathcal{F}}_1(M_1,M_0), \text{ for some } M_1\in [\tilde{M}_1, 2\tilde{M}_1]\right\},
\end{align*}
where $\tilde{\mathcal{F}}_1(M_1, M_0)$ is the set of parameters for which the following disjointness conditions are satisfied:
\begin{description} 
\item[$(\mathcal{F}1)_1$] For each $j \in \{0,1\}$,
\begin{align*}
    d\left(3I_j, \bigcup_{m=1}^{2K_jM_j}3(I_j+m\omega)\right) > 0.
\end{align*}

\item[$(\bar{\mathcal{F}}1)_1$] 
         \begin{align*}
    d\left(I_1, \bigcup_{m=1}^{2K_1M_1}(I_1+m\omega)\right) > \varepsilon_1,
\end{align*}
where $\varepsilon_1 := 3\alpha_c^{M_0}$.

\item[$(\mathcal{F}2)_1$] 
         \begin{align*}
             d\left(I_1 \pm (M_1+S_0)\omega, \bigcup_{|m|\leq 2M_0} (I_0 + m\omega)\right) > 0.
         \end{align*}
\end{description}

The construction of the set $\tilde{\mathcal{F}}_1$ comprises two  steps. The first step is to prove that for each $\omega \in \mathcal{F}_0$ there exists a constant $M_1 \in [\tilde{M}_1, 2\tilde{M}_1]$ for which conditions $(\bar{\mathcal{F}}1)_1$ and $(\mathcal{F}2)_1$ are satisfied, and this is settled in proposition \ref{proposition_parameters_block I1_lemma 1}. The second step is to exclude all the parameters for which the iterates of critical sets $I_1(\omega)$ exhibit resonances and thus do not satisfy the condition $(\mathcal{F}1)_1$, and this is derived in proposition \ref{proposition_parameters_block I1_lemma 2}.

%-------------------------- PROP BEGIN -------------------------
\begin{proposition}\label{proposition_parameters_block I1_lemma 1}
    For each $\omega \in \mathcal{F}_0$, there exists $M_1 \in [\tilde{M}_1, 2\tilde{M}_1]$ such that 
    \begin{align*}
        d \left( I_1(\omega) \pm (M_1+ S_0)\omega, \bigcup_{|m|\leq 2M_0}(I_0+m\omega) \right) >  |I_0|/2.
    \end{align*}
\end{proposition}

\begin{proof}
Let us fix $\omega \in \mathcal{F}_0$. By $(\mathcal{F}1)_0$, we have
        \begin{align*}
            d\left(I_0, \bigcup_{m=1}^{2K_0M_0} (I_0+m\omega)\right) > 2|I_0|.
        \end{align*}
Thus, for any $L \in \mathbb{Z}^{+}$, there exists at most one $p\in [L, L+2K_0M_0]$ such that 
        $$ d((I_0+p\omega), I_0) < |I_0|/2 $$
and there exists at most one $q \in [L, L+2K_0M_0]$ such that 
$$ d((I_0-q\omega), I_0) < |I_0|/2. $$
Then, there exists at most $2(4M_0+1)$ $p$'s in $[L, L+2K_0M_0]$ such that 
        $$ d\left( I_0\pm p\omega, \bigcup_{|m|\leq 2M_0} (I_0+m\omega) \right) < |I_0|/2.$$
       For the fixed constants, we have the bounds $2(4M_0+1) \leq 10M_0$, $2K_0M_0 >16M_0$ and $\tilde{M}_1 = \alpha_c^{-M_0/8}$. Thus, the number of iterates of $I_0$ that return close to itself is significantly smaller than $2K_0M_0$.   
        Since $2K_0M_0 < \tilde{M}_1$, we have $[\tilde{M}_1, \tilde{M}_1+2K_0M_0] \subset [\tilde{M}_1, 2\tilde{M}_1]$. Thus, there exists $M_1 \in [\tilde{M}_1, \tilde{M}_1+2K_0M_0]$ such that 
        \begin{align*}
        d\left( I_0\pm (M_1+S_0) \omega, \bigcup_{|m|\leq 2M_0} (I_0+m\omega)\right) & > |I_0|/2.
        \end{align*}
Since $I_1 \subset I_0$, this proves that there exists a constant $M_1\in [\tilde{M}_1, 2\tilde{M}_1]$ such that
\begin{align*}
    d\left(I_1(\omega)\pm (M_1+S_0)\omega, \bigcup_{|m| \leq 2M_0}(I_0+m\omega)\right) > |I_0|/2.
\end{align*}
\end{proof}
%-------------------------- PROP END -------------------------

%-------------------------- PROP BEGIN -------------------------
\begin{proposition}\label{proposition_parameters_block I1_lemma 2}
For each interval $\Omega \subset \mathcal{F}_0$, there exists a set $\Omega^{+} \subset \Omega$ such that 
     \begin{align*}
         |\Omega^{+}| & \geq |\Omega| - \frac{1}{\tilde{M}_1^{2}} \\
         \# \Omega^{+} & \leq \tilde{M}_1^{4},
     \end{align*}
     and for each interval $\Omega_{i}^{+} \subset \Omega^{+}$ there is a corresponding $M_1 \in [\tilde{M}_1, 2\tilde{M}_1]$ such that the following conditions $(\mathcal{F}1)_1$, $(\mathcal{F}2)_1$ are satisfied for every $\omega \in \Omega^{+}_i$:
     \begin{enumerate}
      \item[$(\mathcal{F}1)_1$] 
        \begin{align*}
            d\left(3I_1(\omega) ,\bigcup_{m=1}^{2K_1M_1} 3(I_1(\omega) + m\omega)\right) > 0.
        \end{align*}

\item[$(\bar{\mathcal{F}}1)_1$] 
         \begin{align*}
    d\left(I_1, \bigcup_{m=1}^{2K_1M_1}(I_1+m\omega)\right) > \varepsilon_1,
\end{align*}
where $\varepsilon_1 := 3\alpha_c^{M_0}$.
        
         \item[$(\mathcal{F}2)_1$] 
         \begin{align*}
             d\left(I_1(\omega) \pm (M_1+S_0)\omega, \bigcup_{|m|\leq 2M_0} (I_0+m\omega)\right) > 0.
         \end{align*}
     \end{enumerate}
Equivalently, we have $\Omega_i^+ \subset \tilde{\mathcal{F}}_1(M_0, M_1)$.
\end{proposition}

\begin{proof}
Let us fix a constant $\displaystyle \delta:= \frac{|I_0|}{2^4\tilde{M}_1}$. Let us subdivide the fixed interval $\Omega$ such that we have $[1/\delta]$  intervals each of size $\delta$ and possibly one interval extra interval of smaller size. Then, we have 
    \begin{align*}
        \Omega = \bigcup_{i=1}^{T} \Omega_i, 
    \end{align*} 
    where $|\Omega_i| = \delta$ for each $i \in \{1,\dots , T-1\}$, and $|\Omega_T| < \delta$. For each $\Omega_i$, let $\omega_i$ denote its midpoint and $M_1^i \in [\tilde{M}_1, 2\tilde{M}_1]$ be the constant corresponding to $\omega_i$, as defined in proposition \ref{proposition_parameters_block I1_lemma 1}, for which 
    \begin{align*}
        d\left(I_1(\omega_i)\pm (M_1^{i}+S_0)\omega_i, \bigcup_{|m|\leq 2M_0} (I_0 + m\omega_i) \right) > |I_0|/2. 
    \end{align*}
 Recall that, for each $i \in \{1,\dots , T-1\}$ we have $\displaystyle |\Omega_i|\leq \delta = \frac{|I_0|} {2^4\tilde{M}_1}$, $M_1^i \in [\tilde{M}_1, 2\tilde{M}_1]$ and $|\partial_{\omega}I_1| < 1/4$. Then, for each $i\in \{1,\dots , T-1\}$ and any $\omega, \tilde{\omega} \in \Omega_i$, 
\begin{align*}
    |(I_1(\omega) \pm (M_1^{i}+S_0)\omega) - (I_1(\tilde{\omega}) \pm (M_1^{i}+S_0)\tilde{\omega})| & \leq |I_1(\omega) - I_1(\tilde{\omega})| + (M_1^{i} + S_0)|\omega - \tilde{\omega}| \\
    & \leq \frac{\delta}{4} + 2\tilde{M}_1\delta \ < \ \frac{|I_0|}{4}.
\end{align*}
Thus, for every $\omega \in \Omega_{i}$, 
\begin{align*}
   d\left( I_1(\omega) \pm M_1^i\omega,  \bigcup_{|m|\leq 2M_0} (I_0 + m\omega) \right) > 0. 
\end{align*}
By proposition \ref{proposition_structure I1_critical set I_1}, for each $\omega \in \Omega_i$ we have $\lambda^{-2}\alpha_c^{b(M_0-1)} < |I^{(i',j')}_1(\omega)| < \lambda^{-1+1/32} \alpha_{c}^{M_0-1}$, for each $i',j'\in \{1,2\}$, and $|\partial_{\omega}I_1(\omega)| < 1/4$. This shows that, for each $\omega \in \Omega_i$ the condition $(\bar{\mathcal{F}}1)_1$ implies  $(\mathcal{F}1)_1$. Hence, it is in fact sufficient consider the condition $(\bar{\mathcal{F}}1)_1$. 

For each $i \in \{1,\dots ,T\}$, let $\Omega_i^{+}$ denote the set of all parameters in $\Omega_i$ that satisfy the condition $(\bar{\mathcal{F}}1)_1$ with respect to $M_1^i$. By lemma \ref{lemma_parameters_prelim_disjoint condn}, the set 
\begin{align*}
        \left\{ \omega\in \Omega_i: d \left( \bigcup_{i',j'=1}^{2} I_1^{(i',j')}(\omega), \bigcup_{m=1}^{2K_1M_1^i} \bigcup_{i',j'=1}^{2} (I_1^{(i',j')}(\omega) + m\omega)\right) < \varepsilon_1  \right\}
\end{align*}
has measure at most $(2^8K_1\tilde{M}_1)(\lambda^{-1+1/32}\alpha_{c}^{M_0-1} + 3\alpha_c^{M_0})$ and consists of at most $2^{8}(K_1\tilde{M}_1)^2$. Thus, for each $i \in \{1,\dots , T\}$, 
    \begin{align*}
        |\Omega_i^+| & \geq |\Omega_i| - 2^{10}K_1\tilde{M}_1\alpha_{c}^{M_0-1}. \\
        \# \Omega_i^+ & \leq  2^8 (K_1\tilde{M}_1)^2
    \end{align*}
Since $\Omega = \bigcup_{i=1}^{T}\Omega_i$ and $\displaystyle T \leq \frac{1}{\delta} + 1  = \frac{2^4\tilde{M}_1}{|I_0|}+1$, we 
    have 
\begin{align*}
    |\Omega^+| & \geq |\Omega| - \frac{2^{25} (b+1)}{\lambda^{-1/4}} \tilde{M}^2_1 \alpha_{c}^{M_0-1} \\
    & > |\Omega| - (\alpha_c^{M_0/4}) = |\Omega| - \frac{1}{\tilde{M}_1^2},
\end{align*}
for a sufficiently large $\lambda$ with respect to the constant $b$. Similarly, the maximum number of intervals in $\Omega^+$ is given by
\begin{align*}
    \# \Omega^+ & \leq \frac{2^5\tilde{M}_1}{ \lambda^{-1/4}} 2^8 (K_1\tilde{M}_1)^2 < \tilde{M}_1^4.  
\end{align*}

\end{proof}
%-------------------------- PROP END -------------------------

%%%%%%%%%%%%%%%%%%%%% tilde F(1,1) PARAMETER SPACE - nonempty I2 %%%%%%%%%%%%
\subsection{Parameters corresponding to non-empty critical set $I_2$}\label{subsection_Parameters corresponding to non-empty critical sets}

In this subsection, we describe a specific measurable subset of $\tilde{\mathcal{F}}_1^{(1,1)}$.
Recall that the set $\tilde{\mathcal{F}}_1^{(1,1)}$ was introduced in definition \ref{definition_structure I2_SNA character of parameters} and by proposition \ref{proposition_structure I2_I2 empty nonempty}, it is equivalently given by
\begin{align*}
    \tilde{\mathcal{F}}_1^{(1,1)} = \left\{ \omega \in \mathbb{T}: I_2^{(1,1)}(\omega) \neq \emptyset \right\}.
\end{align*}
Let us recall the following sets as in definition \ref{definition_structure I1_sync set Gs} and definition \ref{definition_structure I1_switch set Ss}, which are central to our construction in this subsection,
\begin{align*}
G_2^{(1)} = ( \theta^{(2,1)}+3\lambda^{-1+1/64}, \theta^{(1,2)} - 3\lambda^{-1+1/64})\\
G_2^{(2)} = ( \theta^{(1,2)}+3\lambda^{-1+1/64}, \theta^{(1,2)} -3\lambda^{-1+1/64})\\
    {S_3} = ( \theta^{(1,1)}-3\lambda^{-1+1/64}, \theta^{(1,2)} +3\lambda^{-1+1/64}),
\end{align*}
and $G_2 = G_2^{(1)} \cup G_2^{(2)}$.

\begin{definition}\label{definition_parameter space F(1,1)_Rk and R sets}
    Let $k \geq 1$ be an integer such that $\displaystyle \frac{1}{10|{G_2^{(i)}}|} < k < \frac{1}{40 |{S_3}|}$. Then, we define the following sets:
\begin{align*}
    \tilde{R}_k & = \left\{ \frac{p}{q} : 1\leq q \leq k, 1 \leq p \leq (q-1), (p,q)=1 \right\} \\
    R_k & = \left\{ \frac{p}{q} : \frac{k}{2} < q \leq k, 1 \leq p \leq (q-1), (p,q)=1 \right\}, 
\end{align*}
where $(p,q) =1$ indicates that the integers $p$ and $q$ are coprime.
\end{definition}

Note that, by the above definition, for each $k$ we have $R_k \subset \tilde{R}_k$. 

\begin{definition}\label{definition_parameter space F(1,1)_Omega sets}
    For each $k\geq 1$ such that $\displaystyle \frac{1}{10|{G_2^{(i)}}|} < k < \frac{1}{40 |{S_3}|}$, and $p/q \in R_k$, let 
    \begin{align*}
        \Omega(p/q) = \left( \frac{p}{q} -\frac{1}{10k^2}, \frac{p}{q} - \frac{1}{20 k^2} \right).
    \end{align*}
    Let $k'>1$ be fixed such that $\displaystyle \frac{1}{10|{G_2^{(i)}}|} < k' < \frac{1}{40 |{S_3}|}$, then we define 
    \begin{align*}
        \Omega = \Omega_{k'} := \bigcup_{p/q \in R_{k'}} \Omega(p/q).  
    \end{align*}
\end{definition}

We show that there exists integers $k\geq 1$ that satisfy the condition $\displaystyle \frac{1}{10|{G_2^{(i)}}|} < k < \frac{1}{40 |{S_3}|}$ and thus, has corresponding non-trivial sets $R_k$ and $\tilde{R}_k$ through the following lemma.

%------------------- LEMMA BEGIN -------------------
\begin{lemma}\label{lemma_parameters_non-empty I2_welldefinedness of k}
    For a sufficiently large $\lambda>0$, the measure of the sets $G_2$ and $S_3$ satisfy $\displaystyle \frac{1}{40|S_3|} - \frac{1}{10|G_2^{(i)}|} > \frac{1}{10|G_2^{(i)}|}$ and 
    \begin{align} \label{equation_parameters_non-empty critical sets}
    \left[\frac{1}{10|G_2^{(i)}|}, \frac{1}{40|{S_3}|} \right] \subset [M_0^3, M_0^9].
\end{align}
\end{lemma}

\begin{proof}
By lemma \ref{lemma_preliminary estimates_|G| and |S|}, we have $|G^{(1)}| = |G^{(2)}| > 9|S| > 9\lambda^{-1/2-1/32}$. Then, the measure of the sets $G_2^{(i)}$ for $i\in \{1,2\}$ and $S_3$ satisfy the relation $|G_2^{(i)}| > 8|S_3|$. This gives us
\begin{align*}
 \frac{1}{40|S_3|} - \frac{1}{10|G_2^{(i)}|}  = \frac{|G_2^{(i)}| - 4|S_3|}{40|G_2^{(i)}||S_3|} > \frac{1}{10|G_2^{(i)}|}, \text{ for } i\in \{1,2\}.
\end{align*}
Since $|G_2^{(i)}| > 8|S_3| > 8\lambda^{-1/2-1/32}$ and $M_0 = \lambda^{-1/16}$, for a sufficiently large $\lambda$ we have 
\begin{align}
    \left[\frac{1}{10|G_2^{(i)}|}, \frac{1}{40|{S_3}|} \right] \subset [M_0^3, M_0^9].
\end{align}
Thus, there exists many positive integers $k \in [M_0^3, M_0^9]$ for which $\displaystyle \frac{1}{10|{G_2^{(i)}}|} < k < \frac{1}{40 |{S_3}|}$.
\end{proof}
%------------------- LEMMA END -------------------

%------------------- PROP BEGIN -------------------
\begin{proposition}\label{proposition_parameters_non-empty I2_characterise}
Let us fix an integer $k\geq 1$ such that $\displaystyle \frac{1}{10|{G_2^{(i)}}|} < k < \frac{1}{40 |{S_3}|}$, and a rational $p/q \in R_k$. For any $\omega \in \Omega(p/q)$, the orbits $(\theta^{(1,1)}+m\omega)_{m=1}^{q}$ and $(\theta^{(1,1)}-m\omega)_{m=1}^{q}$ satisfy the following conditions: 
    \begin{enumerate}[label=\alph*)]
         \item\label{lemma item_parameters_non-empty I2_characterise a} For $m \in \{1,\dots , q-1\}$, 
         $$d(\theta^{(1,1)}+m\omega, \theta^{(1,1)}) > |{S_3}|$$
         and thus $(\theta^{(1,1)}+m\omega) \notin {S_3}$. 
        \item\label{lemma item_parameters_non-empty I2_characterise b} We have $$|{S_3}|< d(\theta^{(1,1)}+q\omega, \theta^{(1,1)}) < |G_2^{(1)}|$$ and thus $(\theta^{(1,1)} + q\omega) \in G_2^{(1)}$. 
        \item\label{lemma item_parameters_non-empty I2_characterise c} For $m\in \{1,\dots , q-1\}$,  $$d(\theta^{(1,1)}-m\omega, \theta^{(1,1)}) > |{S_3}|$$ and thus $(\theta^{(1,1)}-m\omega) \notin {S_3}$. 
        \item\label{lemma item_parameters_non-empty I2_characterise d} We have $$|{S_3}| < d(\theta^{(1,1)}-q\omega, \theta^{(1,1)}) < |G_2^{(2)}|$$ and thus $(\theta^{(1,1)} - q\omega) \in G_2^{(2)}$. 
    \end{enumerate}
\end{proposition}
\begin{proof}
\textit{Distribution of the forward iterates:} For each $p/q \in R_k$, the points in $\displaystyle \left( \theta^{(1,1)} + {mp}/{q} \right)_{m=1}^{q}$ are equidistributed on $\mathbb{T}$ by a distance of $1/q$. Then, for any $\displaystyle \omega_1, \omega_2 \in \Omega(p/q)$ and $m\in \{1,\dots ,q\}$, we have 
    \begin{align*}
        d\left((\theta^{(1,1)}+m\omega_1), (\theta^{(1,1)}+m\omega_2)\right) & \leq \frac{m}{20k^2} < \frac{1}{10k} < \frac{1}{q}. 
    \end{align*}
For any $\displaystyle \omega \in \Omega(p/q)$ and $m \in \{1,\dots , q\}$, 
\begin{align*}
    \frac{m}{20k^2} \leq d \left(\left(\theta^{(1,1)}+ \frac{mp}{q}\right), (\theta^{(1,1)}+m\omega)\right) \leq \frac{m}{10k^2}.
\end{align*}
Thus, if we consider $0 \leq m_1, m_2 \leq q$ such that $m_1 \neq m_2$ and $\displaystyle d\left(\frac{m_1p}{q}, \frac{m_2p}{q}\right) = \frac{1}{q}$, then 
\begin{align*}
    d(m_1\omega, m_2\omega) \leq \frac{1}{q} + \frac{\max\{m_1, m_2\}}{10k^2} \leq \frac{1}{q} + \frac{1}{10k}. 
\end{align*}
For any $0 \leq m_1, m_2 \leq q$ such that $m_1 \neq m_2$, from the above estimates we get 
\begin{align*}
    d(m_1\omega, m_2\omega) & \geq \frac{1}{q} - \frac{1}{10k} = \frac{9}{10k},
\end{align*}
and thus, 
\begin{align}\label{equation_parameters_non-empty I2_characterise_1}
    \frac{9}{10k} \leq d(m_1\omega, m_2\omega) \leq \frac{1}{q}+\frac{1}{10k}.
\end{align}

\vspace{0.2cm}
\noindent\textit{Proof of \ref{lemma item_parameters_non-empty I2_characterise a} and \ref{lemma item_parameters_non-empty I2_characterise b}:} For the fixed $k$, and for each $m \in \{1,\dots , q\}$, by equation \eqref{equation_parameters_non-empty I2_characterise_1}, we have 
\begin{align*}
    d(0,m\omega) & > \frac{9}{10k} > \frac{9}{10} 40 |S_3| > |S_3|,
\end{align*}
which gives us the estimate $d(\theta^{(1,1)}+m\omega, \theta^{(1,1)}) > |S_3|$. Since $\displaystyle {S_3} = ( \theta^{(1,1)}-3\lambda^{-1+1/64}, \theta^{(1,2)} +3\lambda^{-1+1/64})$ and $|S_3|> \lambda^{-1/2-1/16}$, by lemma \ref{lemma_structure I1_switch set S_2 size}, this proves that $(\theta^{(1,1)}+m\omega) \notin S_3$, for all $m \in \{1,\dots ,q\}$.

Since $\displaystyle \omega \in \Omega_{p/q} = \left( \frac{p}{q} - \frac{1}{10k^2}, \frac{p}{q}- \frac{1}{20k^2} \right)$, note that 
\begin{align*}
\theta^{(1,1)}+q\omega < \theta^{(1,1)},
\end{align*}
where the order corresponds to that on $[0,1)$. By the choice of $k$, 
\begin{align*}
    d(\theta^{(1,1)}+q\omega, \theta^{(1,1)}) < \frac{q}{10k^2} < \frac{1}{10k} < |G_2^{(1)}|.
\end{align*}
Then, since $G_2^{(1)} = \left( \theta^{(2,1)}+3\lambda^{-1+1/64}, \theta^{(1,1)}-3\lambda^{-1+1/64} \right)$ and $d(\theta^{(1,1)}+q\omega, \theta^{(1,1)}) > |S_3|$, this gives us $\theta^{(1,1)} \in G_2^{(1)}$.

%--------------------- Case m<0 ---------------------------

\textit{Distribution of backward iterations:}
For each $p/q \in R_k$, the points in $\displaystyle \left( \theta^{(1,1)} - {mp}/{q} \right)_{m=1}^{q}$ are equidistributed on $\mathbb{T}$ by a distance of $1/q$. Then, for any $\displaystyle \omega_1, \omega_2 \in \Omega(p/q)$ and $m\in \{1,\dots ,q\}$, we have 
    \begin{align*}
        d\left((\theta^{(1,1)} - m\omega_1), (\theta^{(1,1)} - m\omega_2)\right) & \leq \frac{m}{20k^2} < \frac{1}{10k} < \frac{1}{q}. 
    \end{align*}
For any $\displaystyle \omega \in \Omega(p/q)$ and $m \in \{1,\dots , q\}$, 
\begin{align*}
    \frac{m}{20k^2} \leq d \left(\left(\theta^{(1,1)} -  \frac{mp}{q}\right), (\theta^{(1,1)} - m\omega)\right) \leq \frac{m}{10k^2}.
\end{align*}
Similar to equation \eqref{equation_parameters_non-empty I2_characterise_1}, for any $0 \leq m_1, m_2 \leq q$ such that $m_1 \neq m_2$, we have
\begin{align}\label{equation_parameters_non-empty I2_characterise_2}
    \frac{9}{10k} \leq d(-m_1\omega, -m_2\omega) \leq \frac{1}{q}+\frac{1}{10k}.
\end{align}

\vspace{0.2cm}
\noindent\textit{Proof of \ref{lemma item_parameters_non-empty I2_characterise c} and \ref{lemma item_parameters_non-empty I2_characterise d}:} For the fixed $k$, and for each $m \in \{1,\dots , q\}$, by equation \eqref{equation_parameters_non-empty I2_characterise_2}, we have 
\begin{align*}
    d(0,-m\omega) & > |S_3|.
\end{align*}
This gives us the estimate $d(\theta^{(1,1)}-m\omega, \theta^{(1,1)}) > |S_3|$. By the definition of the set ${S_3}$ and the estimate  $|S_3|> \lambda^{-1/2-1/16}$ by lemma \ref{lemma_structure I1_switch set S_2 size}, it follows that $(\theta^{(1,1)}-m\omega) \notin S_3$, for all $m\in \{1,\dots ,q\}$.

Since $\displaystyle \omega \in \Omega_{p/q} = \left( \frac{p}{q} - \frac{1}{10k^2}, \frac{p}{q}- \frac{1}{20k^2} \right)$, note that
\begin{align*}
\theta^{(1,1)} < \theta^{(1,1)}-q\omega,
\end{align*}
where the order corresponds to that on $[0,1)$. By the choice of $k$, 
\begin{align*}
    d(\theta^{(1,1)}-q\omega, \theta^{(1,1)}) & = d(q\omega,0) < \frac{q}{10k^2} < \frac{1}{10k} < |G_2^{(2)}|.
\end{align*}
Then, since $G_2^{(2)} = \left( \theta^{(1,2)}+3\lambda^{-1+1/64}, \theta^{(2,2)}-3\lambda^{-1+1/64} \right)$ and $d(\theta^{(1,1)}-q\omega, \theta^{(1,1)}) > |S_3|$, this proves that $\theta^{(1,1)} \in G_2^{(2)}$.
\end{proof}
%------------------- PROP END -------------------

%------------------- PROP BEGIN -------------------
\begin{proposition}\label{lemma_parameters_non-empty I2_estimates}
    Let the set $\Omega$ be as defined in definition \ref{definition_parameter space F(1,1)_Omega sets}. Then, 
\begin{align*}
    \left(\Omega \cap \mathcal{F}_0 \cap \tilde{\mathcal{F}}_1 \cap \mathcal{F}^s(G_2)\right) \subset \mathcal{F}^{(1,1)}_1,
\end{align*}
and we have the estimates
    \begin{align*}
        |\Omega| & \geq \frac{1}{640} \\
      \# \Omega & < M_0^{20}.
    \end{align*}
\end{proposition}

\begin{proof} 
Let us fix $\omega \in \Omega\cap \mathcal{F}_0 \cap \tilde{\mathcal{F}}_1$. By definition of the set $\Omega$ in definition \ref{definition_parameter space F(1,1)_Omega sets} and proposition \ref{proposition_parameters_non-empty I2_characterise}, there exists a constant $1 \leq q  \leq M_0^9$ such that 
\begin{align*}
    (\theta^{(1,1)} \pm q\omega) \in G_2, \\
    (\theta^{(1,1)} \pm k\omega)_{|m|< q-1} \cap S_3 = \emptyset.
\end{align*}
By proposition \ref{proposition_structure I2_I2  nonempty F^{(1,1)}}, this proves that $\omega \in \tilde{\mathcal{F}}_1^{(1,1)}$ and thus, $\left(\Omega\cap \mathcal{F}_0 \cap \tilde{\mathcal{F}}_1 \cap \mathcal{F}^s(G_2)\right) \subset \mathcal{F}_1^{(1,1)}$.

For any two distinct rationals $\displaystyle \frac{p_1}{q_1}, \frac{p_2}{q_2} \in R_k$, 
        \begin{align*}
            d \left( \frac{p_1}{q_1}, \frac{p_2}{q_2} \right) > \frac{1}{k^2},
        \end{align*}
which shows that $\Omega_{p_1/q_1} \cap \Omega_{p_2/q_2} = \emptyset$. Thus, the intervals $\Omega_{p/q}$, for $p/q \in R_k$, are all disjoint and the number of intervals in $\Omega$ is at most the cardinality of $R_k$. Since $\displaystyle |\Omega_{p/q}| = \frac{1}{20k^2}$, it follows that
\begin{align*}
    |\Omega| = (\# R_k) \frac{1}{20k^2}.
\end{align*}
For the fixed $k$, we have 
\begin{align*}
        \# \tilde{R}_k & = \# \left\{ \frac{p}{q} : 1\leq q \leq k, 1\leq p \leq (q-1), (p,q) =1 \right\} \\
        & \geq \frac{k^2}{2} \left(1- \frac{1}{4} - \frac{1}{2}\right) \geq \frac{k^2}{2^3}.
\end{align*}        
and $\displaystyle \# \tilde{R}_k \leq \frac{3k^2}{2^3}$. Similarly, we have 
\begin{align*}
        \frac{k^2}{2^3.2^2} & \leq \# \tilde{R}_{k/2} \leq \frac{3k^2}{2^3 2^2}.
\end{align*}
By the definition of sets $R_k$ and $\tilde{R}_k$, we have $\# R_k = \# \tilde{R}_k - \# \tilde{R}_{k/2}$ and thus,
    \begin{align*}
        \frac{k^2}{2^3} \left( 1-\frac{3}{2^2} \right) \leq & \# R_k \leq \frac{3k^2}{2^3}
    \end{align*}
 By the above bounds we obtain that 
    \begin{align*}
        |\Omega| & = (\# R_k) \left( 1-\frac{1}{2} \right) \frac{1}{10k^2} \geq \frac{1}{640}. 
    \end{align*}
and $\# \Omega < {3k^2}/{8} < M_0^{20}$.
\end{proof}
%------------------- PROP END -------------------

%%%%%%%%%%%%%%%%%%%%% tilde Omega PARAMETER SPACE - empty I2 %%%%%%%%%%%%
\subsection{Parameters corresponding to empty critical sets}\label{subsection_Parameters corresponding to empty critical sets}

In this subsection, we construct an open interval of parameters for which the critical set $I_2(\omega)$, when it is well-defined, is empty. More specifically, we construct a open interval contained in the set $\mathcal{F}_0$ which satisfies the conditions in proposition \ref{proposition_structure I2_I2  empty}. For the purpose of this subsection, we recall the following sets as in definition \ref{definition_structure I1_sync set Gs} and definition \ref{definition_structure I1_switch set Ss}, 
\begin{align*}
G_2^{(1)} & = ( \theta^{(2,1)}+3\lambda^{-1+1/64}, \theta^{(1,2)} - 3\lambda^{-1+1/64})\\
G_2^{(2)} & = ( \theta^{(1,2)}+3\lambda^{-1+1/64}, \theta^{(1,2)} -3\lambda^{-1+1/64})\\
S_2 & = \left( \theta^{(1,1)} + 3\lambda^{-1+1/64}, \theta^{(1,2)} - 3\lambda^{-1+1/64} \right).
\end{align*}
By the estimates in lemma \ref{lemma_structure I1_switch set S_2 size}, we have $|S_2| > \lambda^{-1/2-1/16} = M_0^{-9}$.

\begin{definition}\label{definition_parameter space F empty_Omega sets}
    Let us define the interval
    \begin{align*}
        \tilde{\Omega} := \left( \frac{1}{5K_0M_0} -\frac{|S_2|}{5K_0M_0}, \frac{1}{5K_0M_0} -\frac{3|S_2|}{20K_0M_0} \right).
    \end{align*}
\end{definition}

For each $\omega \in \tilde{\Omega}$, the following result describes the behaviour of the orbits of the points $\{\theta^{(i,j)}\}_{ij=1}^{2}$ on the base, particularly at their close returns to the critical set $I_0$.

%------------------- PROP BEGIN -------------------
\begin{proposition}\label{proposition_parameters_empty I2_lemma 00}
    For each $\omega \in \tilde{\Omega}$, we have $\omega \in \mathcal{F}_0$ and the following conditions are satisfied:
    \begin{enumerate}[label=\alph*)]
         \item \label{lemma item_parameters_empty I2_lemma 0 a} For each $i \in \{1,2\}$, 
\begin{align}\label{equation_parameters_empty I2_returns}
\left\{(\theta^{(i,1)} - m\omega)\right\}_{m=1}^{10K_0M_0} \cap 2I_0 & = \left\{ (\theta^{(i,1)} - 5K_0M_0\omega), (\theta^{(i,1)} - 10K_0M_0\omega) \right\} \\
\left\{(\theta^{(i,2)} + m\omega)\right\}_{m=1}^{10K_0M_0} \cap 2I_0 & = \left\{ (\theta^{(i,2)} + 5K_0M_0\omega), (\theta^{(i,2)} + 10K_0M_0\omega) \right\}. 
\end{align}    

        \item\label{lemma item_parameters_empty I2_lemma 0 b} The behaviour of the forward and backward orbits at the first return to $I_0$ is:
        \begin{align}\label{equation_parameters_empty I2_first returns}
            (\theta^{(1,1)} - 5K_0M_0\omega) \in S_2, \ \ (\theta^{(2,1)} - 5K_0M_0\omega) \in G_2^{(1)} \\
            (\theta^{(1,2)} + 5K_0M_0\omega) \in S_2, \ \ (\theta^{(2,2)} + 5K_0M_0\omega) \in G_2^{(2)}.
        \end{align}
        
        \item\label{lemma item_parameters_empty I2_lemma 0 c} The behaviour of the forward and backward orbits at the second return to $I_0$ is:
        \begin{align}\label{equation_parameters_empty I2_second returns}
            (\theta^{(1,1)} - 10K_0M_0\omega) \in G_2^{(2)}, \\
            (\theta^{(1,2)} + 10K_0M_0\omega) \in G_2^{(1)}.
        \end{align}

        \item\label{lemma item_parameters_empty I2_lemma 0 d} For each $\omega \in \tilde{\Omega} \cap \tilde{\mathcal{F}}_1$, we have $I_2(\omega) = \emptyset$.
    \end{enumerate}
\end{proposition}

\begin{proof}

 Let us fix $\omega \in \tilde{\Omega}$. Note that the points in the sets $\displaystyle\left\{\frac{m}{5K_0M_0}\right\}_{m=1}^{5K_0M_0}$ and $\displaystyle\left\{\frac{m}{5K_0M_0}\right\}_{m=6K_0M_0}^{10K_0M_0}$ coincide and the points are equidistributed on $\mathbb{T}$ by a distance of $\displaystyle \frac{1}{5K_0M_0}$. A similar observation can be made for $10K_0M_0 \leq m < 0$. For each  $1 \leq m \leq 10K_0M_0$,  
\begin{align}\label{equation_parameter_I2 empty equidistribution of iterates}
  \frac{3|S_2|}{20K_0M_0} \leq \frac{m3|S_2|}{20K_0M_0} < d\left( \frac{m}{5K_0M_0}, m\omega \right) < \frac{m|S_2|}{5K_0M_0} \leq 2|S_2|.
\end{align} 
Similarly, for each  $1 \leq m \leq 10K_0M_0$,  
\begin{align}\label{equation_parameter_I2 empty equidistribution of iterates negative}
  \frac{3|S_2|}{20K_0M_0} \leq \frac{m3|S_2|}{20K_0M_0} < d\left( \frac{-m}{5K_0M_0}, -m\omega \right) < \frac{m|S_2|}{5K_0M_0} \leq 2|S_2|.
\end{align} 
By the distribution of the points $\displaystyle \left\{\frac{m}{5K_0M_0}\right\}_{m=1}^{2K_0M_0}$ and equation \eqref{equation_parameter_I2 empty equidistribution of iterates}, we obtain 
\begin{align}
    d\left(3I_0, 3(I_0+m\omega)\right) > \frac{1}{5K_0M_0} - 2|S_2| - 6|I_0| > \frac{1}{5K_0M_0} - \frac{16}{M_0^{4}},
\end{align}
for all $m \in \{1,\dots ,2K_0M_0\}$. This proves that 
\begin{align*}
    d\left(3I_0, \bigcup_{m=1}^{2K_0M_0}3(I_0+m\omega)\right) > 0,
\end{align*}
and thus, $\omega \in \mathcal{F}_0$.

\vspace{0.2cm}
\noindent\textit{Proof of \ref{lemma item_parameters_empty I2_lemma 0 a}:}  Recall that $\displaystyle\{\theta^{(i,j)}\}_{i,j=1}^{2} \subset I_0$ and $|I_0| = 2M_0^{-4}$. Then, by equation \eqref{equation_parameter_I2 empty equidistribution of iterates} and the equidistribution of the points $\displaystyle \left\{\frac{m}{5K_0M_0}\right\}_{m=1}^{5K_0M_0}$, we have 
\begin{align*}
    \left\{(\theta^{(i,2)} + m\omega)\right\}_{m=1}^{10K_0M_0} \cap 2I_0 & = \left\{ (\theta^{(i,2)} + 5K_0M_0\omega), (\theta^{(i,2)} + 10K_0M_0\omega) \right\},
\end{align*}
for each $i\in \{1,2\}$. Similarly, we obtain 
\begin{align*}
    \left\{(\theta^{(i,1)} - m\omega)\right\}_{m=1}^{10K_0M_0} \cap 2I_0 & = \left\{ (\theta^{(i,1)} - 5K_0M_0\omega), (\theta^{(i,1)} - 10K_0M_0\omega) \right\},
\end{align*}
for each $i\in \{1,2\}$. 

\vspace{0.2cm}
\noindent\textit{Proof of \ref{lemma item_parameters_empty I2_lemma 0 b}:}
By the above estimates and equation \eqref{equation_parameter_I2 empty equidistribution of iterates negative}, we obtain $ \theta^{(1,1)} <( \theta^{(1,1)}-5K_0M_0\omega)$ and
\begin{align}
     |S_2|/2 < d\left(\theta^{(1,1)}, (\theta^{(1,1)}-5K_0M_0\omega) \right) < |S_2|.
\end{align}
Similarly, for the iterates of $\theta^{(2,1)}$ we have $\theta^{(2,1)} <( \theta^{(2,1)}-5K_0M_0\omega)$ and
\begin{align*}
   |S_2|/2 < d\left(\theta^{(2,1)}, (\theta^{(2,1)}-5K_0M_0\omega) \right) < |S_2| < |G_2^{(1)}|. 
\end{align*}
By the definition of the sets $S_2$ and $G_2$, and the estimate $|S_2| > \lambda^{-1/2-1/16}$, it follows that $(\theta^{(1,1)} - 5K_0M_0\omega) \in S_2$ and $(\theta^{(2,1)} - 5K_0M_0\omega) \in G_2^{(1)}$. By equation \eqref{equation_parameter_I2 empty equidistribution of iterates}, we obtain $(\theta^{(1,2)}+5K_0M_0\omega) < \theta^{(1,2)}$ and 
\begin{align}
    |S_2|/2 < d\left(\theta^{(1,2)}, (\theta^{(1,2)} + 5K_0M_0\omega) \right) < |S_2|.
\end{align}
Similarly, for the iterates of $\theta^{(2,2)}$ we have $(\theta^{(2,1)} + 5K_0M_0\omega) < \theta^{(2,2)}$ and
\begin{align*}
  |S_2|/2 < d\left(\theta^{(2,2)}, (\theta^{(2,2)}-5K_0M_0\omega) \right) < |S_2| < |G_2^{(2)}|. 
\end{align*}
Then, by the definition of the sets $S_2$ and $G_2$, and the estimate $|S_2| > \lambda^{-1/2-1/16}$, it follows that $(\theta^{(1,2)} + 5K_0M_0\omega) \in S_2$ and $(\theta^{(2,2)} + 5K_0M_0\omega) \in G_2^{(2)}$.

\vspace{0.2cm}
\noindent\textit{Proof of \ref{lemma item_parameters_empty I2_lemma 0 c}:} By the distribution of the iterates $\displaystyle \left\{\frac{m}{5K_0M_0}\right\}_{m=6K_0M_0}^{10K_0M_0}$ and the estimates in equation \eqref{equation_parameter_I2 empty equidistribution of iterates negative}, we obtain $\theta^{(1,1)} <( \theta^{(1,1)} - 10K_0M_0\omega)$ and
\begin{align}
   3|S_2|/2 < d\left(\theta^{(1,1)}, (\theta^{(1,1)}-10K_0M_0\omega) \right) < 2|S_2| < |G_2^{(2)}|.
\end{align}
Similarly, we have $ (\theta^{(1,2)} + 10K_0M_0\omega) < \theta^{(1,2)}$ and 
\begin{align*}
   3|S_2|/2 < d\left(\theta^{(1,2)}, (\theta^{(1,2)} + 10K_0M_0\omega) \right) < 2|S_2| < |G_2^{(1)}|. 
\end{align*}
Then, by the definition of the sets $S_2$ and $G_2$, and the estimate $|S_2| > \lambda^{-1/2-1/16}$, it follows that $(\theta^{(1,1)} - 10K_0M_0\omega) \in G_2^{(2)}$ and $(\theta^{(1,2)} + 10K_0M_0\omega) \in G_2^{(1)}$.

\vspace{0.2cm}
\noindent\textit{Proof of \ref{lemma item_parameters_empty I2_lemma 0 d}:} For each $\omega \in \tilde{\Omega}\cap \tilde{\mathcal{F}}_1$, it follows from the proof of proposition \ref{proposition_parameters_empty I2_lemma 0} (\ref{lemma item_parameters_empty I2_lemma 0 a}-\ref{lemma item_parameters_empty I2_lemma 0 d}) and proposition \ref{proposition_structure I2_I2  empty} that $I_2(\omega) = \emptyset$. 

\end{proof}
%------------------- PROP END -------------------

%------------------- PROP BEGIN -------------------
\begin{proposition}\label{proposition_parameters_empty I2_lemma 0}
    Let the set $\tilde{\Omega}$ be as defined in definition \ref{definition_parameter space F empty_Omega sets}. Then, we have 
    \begin{align*}
        |\tilde{\Omega}| & > \frac{1}{K_0M_0^{11}} \\
      \# \tilde{\Omega} & = 1
    \end{align*}
\end{proposition}

\begin{proof}
Since $\displaystyle|\tilde{\Omega}| = \frac{|S_2|}{20K_0M_0}$, the proof follows from definition \ref{definition_parameter space F empty_Omega sets}, and the bound $|S_2|>M_0^{-9}$, due to lemma \ref{lemma_structure I1_switch set S_2 size}.
\end{proof}
%------------------- PROP END -------------------

%%%%%%%%%%%%%%%%%%%%% Fn PARAMETER SPACE %%%%%%%%%%%%%%%%%%%%
\subsection{Parameter exclusion for the construction of ${\mathcal{F}}_n$}\label{subsection_Parameter exclusion for the construction of Fn}

%%%%%%%%%%%%%%%%%%%%% Fn PARAMETER SPACE %%%%%%%%%%%%%%%%%%%%
In this subsection, we establish an inductive construction corresponding to the parameter sets $({\mathcal{F}}_n)_{n\geq 2}$. Recall that, for each $n\geq 1$, we have the constants $\tilde{M}_n = \alpha_c^{-\tilde{M}_{n-1}/8}$, $K_n = 2^{n+8}(b+1)$, and by definition \ref{definition_block n_disjointness condn}, 
\begin{align*}
    {\mathcal{F}}_n = \left\{\omega\in \mathbb{T}: \omega \in  \tilde{\mathcal{F}}_n(M_n,\dots,M_0), \text{ for some } M_n\in [\tilde{M}_n, 2\tilde{M}_n]\right\} \cap \mathcal{F}^s(G_2),
\end{align*}
where $\tilde{\mathcal{F}}_n(M_n,\dots , M_0)$ is the set of all parameters $\omega \in \mathbb{T}$ for which the following disjointness conditions are satisfied:
\begin{description} 
\item[$(\mathcal{F}1)_n$]\label{(F1)_n} For every $j\in \{0, 1,\dots ,n\}$,
\begin{align*}
    d\left(3I_j, \bigcup_{m=1}^{2K_jM_j}3(I_j+m\omega)\right) > 0.
\end{align*}

\item[$(\bar{\mathcal{F}}1)_n$] For each $j \in \{1,\dots,n\}$,
         \begin{align*}
    d\left(I_j, \bigcup_{m=1}^{2K_jM_j}(I_j+m\omega)\right) > \varepsilon_j,
\end{align*}
where $\varepsilon_1 := 3\alpha_c^{M_0}$ and $\varepsilon_j := 3\alpha_c^{\tilde{M}_{j-1}/3}$, for $j\in \{2,\dots ,n\}$.

\item[$(\mathcal{F}2)_n$] For every $j \in \{1,\dots,n\}$,
         \begin{align*}
             d\left(I_j \pm (M_j+S_0)\omega, \bigcup_{l=0}^{j-1} \bigcup_{|m|\leq 2M_l} (I_l + m\omega)\right) > 0.
         \end{align*}

\end{description}

Similar to the construction of $\tilde{\mathcal{F}}_1$ in subsection \ref{subsection_Parameter exclusion for the construction of F1}, for each $n\geq 2$ the construction of ${\mathcal{F}}_n$ comprises two  steps. The first step is to prove that for each $\omega \in {\mathcal{F}}_{n-1}$ there exists a constant $M_n\in [\tilde{M}_n, 2\tilde{M}_n]$ for which condition $(\mathcal{F}2)_n$ is satisfied, and this is settled in proposition \ref{proposition_parameters_Block In_lemma 0}. The second step is to exclude all the parameters for which the iterates of critical sets $I_n$ exhibit resonances and thus do not satisfy the conditions $(\bar{\mathcal{F}}1)_n$ and $(\mathcal{F}1)_n$, and this is derived in proposition \ref{proposition_parameters_Block In_lemma 1}.

%-------------------------- PROP BEGIN -------------------------
\begin{proposition}\label{proposition_parameters_Block In_lemma 0}
Let $n\geq 2$ and $\omega \in  {\mathcal{F}}_{n-1}$. Assume that the critical set $I_n(\omega)$ is non-empty. Then, there exists $\displaystyle M_n \in [\tilde{M}_n, 2\tilde{M}_n]$ such that 
    \begin{align*}
        d \left( I_n(\omega) \pm (M_n + S_0)\omega, \bigcup_{l=0}^{n-1} \bigcup_{|m|\leq 2M_l} (I_l + m\omega) \right) > \varepsilon_{n-1}/3,
    \end{align*}
    where $\displaystyle \varepsilon_{n-1} = 3 \alpha_c^{\tilde{M}_{n-2}/3}$.
\end{proposition}

\begin{proof}
  For the fixed $\omega \in {\mathcal{F}}_{n-1}$, let $(M_l(\omega))_{l=0}^{n-1}$ be the corresponding constants as in definition \ref{definition_block n_disjointness condn}. For each $l\in \{1,\dots ,n\}$ we have the estimate
  \begin{align*}
      \max\left\{|I_{l}^{(i,j)}(\omega)|: i,j=1,2\right\} < \lambda^{-1+1/32}\alpha_c^{\tilde{M}_{l-1}/3} < \varepsilon_l/3,
  \end{align*}
 due to propositions \ref{proposition_structure I1_critical set I_1}, \ref{proposition_structure I2_critical set I2} and \ref{proposition_structure In+1_critical set I_n+1}. Since $\omega$ satisfies the disjointness condition $(\bar{\mathcal{F}}1)_{n-1}$ (in definition \ref{definition_block n_disjointness condn}), we have
        \begin{align*}
            d\left(I_l, \bigcup_{m=1}^{2K_lM_l} (I_l + m\omega)\right) > {\varepsilon}_l,
        \end{align*}
for each $l \in \{0,\dots , n-1\}$. Then, it follows from the above conditions that, for each $l \in \{0,\dots , n-1\}$, $\tau, \kappa \in \{(i,j)\}_{i,j=1}^{2}$ and  $L \in \mathbb{Z}^{+}$, there exists at most one $p\in [L, L+2K_lM_l]$ such that 
        $$ d((I^{\kappa}_l + p\omega), I^{\tau}_l) < \varepsilon_l/3. $$
Similarly, for each $l \in \{0,\dots , n-1\}$, there exists at most one $q \in [L, L+2K_lM_l]$ such that 
$$ d((I^{\tau}_l - q\omega), I^{\kappa}_l) < \varepsilon_l/3, $$
and thus, there exists at most $2^5(4M_l+1)$ $p$'s in $[L, L+2K_lM_l]$ such that 
        \begin{align}\label{parameter exclusion F1_equation_returns of block l}
            d\left( I_l\pm p\omega, \bigcup_{|m|\leq 2M_l} (I_l+m\omega) \right) < \varepsilon_l/3.
        \end{align}
        
Consider the interval $[\tilde{M}_n, \tilde{M}_n + 2K_{n-1}M_{n-1}] \subset [\tilde{M}_n, 2 \tilde{M}_n]$. For a fixed $l \in \{0,\dots ,n-1\}$, we partition the interval $[\tilde{M}_n, \tilde{M}_n + 2K_{n-1}M_{n-1}]$ such that each subinterval contains $2K_lM_l$ integers. In this case, we have at most $\displaystyle  \frac{2^{n-l}M_{n-1}}{M_{l}} + 1$ such subintervals. Since $I_n \subset I_l$ and $\varepsilon_{n-1} \leq \varepsilon_l$, in each of the intervals $[L, L + 2K_lM_l]$, there are at most $2^5(4M_l + 1)$ integers for which 
\begin{align*}
     d\left( I_n\pm p\omega, \bigcup_{|m|\leq 2M_l} (I_l+m\omega) \right) & < \varepsilon_{n-1}/3.
\end{align*}
Then, the number of integers $p\in [\tilde{M}_n, \tilde{M}_n + 2K_{n-1}M_{n-1}]$ for which equation \eqref{parameter exclusion F1_equation_returns of block l} is satisfied for the fixed $l$, is at most $2^5(4M_l+1)(2^{n-l}M_{n-1}/M_l + 1) < 2^8 2^{n-l} M_{n-1}$. Thus, the maximum number of integers $p\in [\tilde{M}_n, \tilde{M}_n + 2K_{n-1}M_{n-1}]$ for which equation \eqref{parameter exclusion F1_equation_returns of block l} is satisfied for some $l \in \{0,\dots ,n-1\}$, is given by 
\begin{align*}
    \sum_{l=0}^{n-1}  2^8 2^{n-l} M_{n-1} & \leq 2^{n+9} M_{n-1} < K_{n-1}M_{n-1}. 
\end{align*}
Further, due to the super exponential growth of the sequence $(\tilde{M}_{n})_{n=0}^{\infty}$, we have $2K_{n-1}M_{n-1} < \tilde{M_{n}}$. Thus, there exists $M_n \in [\tilde{M}_n, \tilde{M}_n + 2K_{n-1}M_{n-1}]$ such that 
        \begin{align*}
        d\left( I_n \pm (M_n+S_0) \omega,  \bigcup_{l=0}^{n-1} \bigcup_{|m|\leq 2M_l} (I_l + m\omega)\right) & > \varepsilon_{n-1}/3.
        \end{align*}   
\end{proof}
%-------------------------- PROP END -------------------------

%-------------------------- PROP BEGIN -------------------------
\begin{proposition}\label{proposition_parameters_Block In_lemma 1}
Let $n\geq 2$. For each interval $\Omega \subset {\mathcal{F}}_{n-1}\cap \tilde{\mathcal{F}}^{(1,1)}$, there exists a set $\Omega^{+} \subset \Omega$ such that 
     \begin{align*}
         |\Omega^{+}| & \geq |\Omega| - \frac{1}{\tilde{M}_{n}^{1/3}} \\
         \# \Omega^{+} & \leq \tilde{M}_{n}^{4}.
     \end{align*}
     and for each interval $\Omega_{i}^{+} \subset \Omega^{+}$ there is a corresponding $M_n \in [\tilde{M}_n, 2\tilde{M}_n]$ such that the conditions $(\mathcal{F}2)_n$ and $(\mathcal{F}1)_n$ are satisfied for every $\omega \in \Omega^{+}_i$:
     \begin{enumerate}
     
        \item[$(\mathcal{F}1)_n$] 
        \begin{align*}
            d\left(3I_n(\omega) ,\bigcup_{m=1}^{2K_n M_n} 3(I_n(\omega) + m\omega)\right) > 0.
        \end{align*}
        
        \item[$(\bar{\mathcal{F}}1)_n$] 
         \begin{align*}
    d\left(I_n, \bigcup_{m=1}^{2K_nM_n}(I_n+m\omega)\right) > \varepsilon_n,
\end{align*}
where $\varepsilon_n := 3\alpha_c^{\tilde{M}_{n-1}/3}$.

         \item[$(\mathcal{F}2)_n$] 
         \begin{align*}
             d\left(I_n(\omega) \pm (M_n+S_0)\omega, \bigcup_{l=0}^{n-1} \bigcup_{|m|\leq 2M_l} (I_l + m\omega)\right) > 0.
         \end{align*}         
     \end{enumerate}
Equivalently, we have $\Omega_i^+ \subset \tilde{\mathcal{F}}_n(M_0,\dots ,M_n)$.
\end{proposition}

\begin{proof}
For each $\omega \in \Omega$, since $\omega \in \mathcal{F}_{n-1} \cap \tilde{\mathcal{F}}^{(1,1)}_1$, by proposition \ref{proposition_structure I2_I2  nonempty F^{(1,1)}} and proposition \ref{proposition_structure In_In empty nonempty} we have a non-empty critical set $I_n(\omega)$. Thus, for each $\omega \in \Omega$ we have non-empty critical sets $(I_l(\omega))_{l=0}^{n}$. 

Let us fix $\displaystyle \delta := \frac{\alpha_c^{\tilde{M}_{n-2}}}{2^4\tilde{M}_n}$. Recall that $\varepsilon_{n-1} = 3\alpha_c^{\tilde{M}_{n-2}/3} > \max\left\{|I_{n-1}^{(i,j)}(\omega)|: i,j=1,2\right\}$. Let us subdivide the fixed interval $\Omega$ such that we have $[1/\delta]$ intervals each of size $\delta$ and possibly one interval extra interval of smaller size. Then, we have 
    \begin{align*}
        \Omega = \bigcup_{i=1}^{T} \Omega_i, 
    \end{align*} 
    where $|\Omega_i| = \delta$ for every $i \in \{1,\dots , T-1\}$, and $|\Omega_T| < \delta$. For each $\Omega_i$, let $\omega_i$ denote its midpoint and $M_n^i \in [\tilde{M}_n, 2\tilde{M}_n]$ be the constant corresponding to $\omega_i$, as defined in proposition \ref{proposition_parameters_Block In_lemma 0}, for which 
    \begin{align*}
        d\left(I_n(\omega_i)\pm (M_n^{i}+S_0)\omega_i, \bigcup_{l=0}^{n-1} \bigcup_{|m|\leq 2M_l} (I_l + m\omega_i) \right) > \varepsilon_{n-1}/3. 
    \end{align*}
 Recall that, for each $i \in \{1,\dots , T\}$ we have $\displaystyle |\Omega_i|\leq \delta$, $M_n^i \in [\tilde{M}_n, 2\tilde{M}_n]$ and $|\partial_{\omega}I_n| < 1/4$. Then, for each $i$ and any $\omega, \tilde{\omega} \in \Omega_i$, 
\begin{align*}
    |(I_n(\omega) \pm (M_n^{i}+S_0)\omega) - (I_n(\tilde{\omega}) \pm (M_n^{i}+S_0)\tilde{\omega})| & \leq |I_n(\omega) - I_n(\tilde{\omega})| + (M_n^{i} + S_0)|\omega - \tilde{\omega}| \\
    & \leq \frac{\delta}{4} + 2\tilde{M}_n\delta \ < \ {\varepsilon_{n-1}}/{4}.
\end{align*}
Thus, for each $\omega \in \Omega_{i}$, we have 
\begin{align*}
   d\left( I_n(\omega) \pm (M_n^i+S_0)\omega, \bigcup_{l=0}^{n-1} \bigcup_{|m|\leq 2M_l} (I_l + m\omega) \right) > 0. 
\end{align*}

By proposition \ref{proposition_structure In+1_critical set I_n+1}, for each $\omega \in \Omega_i$ we have the estimates $\lambda^{-1-1/32} \alpha_{c}^{b(2\tilde{M}_{n-1}+1)} < |I^{({i'},{j'})}_{n}(\omega)| < \alpha_c^{\tilde{M}_{n-1}/3}\lambda^{-1+1/32}$, for each $i',j'=1,2$, and $|\partial_{\omega}I_n(\omega)| < 1/4$. This shows that, for each $\omega \in \Omega_i$ the condition $(\bar{\mathcal{F}}1)_n$ implies  $(\mathcal{F}1)_n$. Hence, it is in fact sufficient consider the condition $(\bar{\mathcal{F}}1)_n$.

For each $i \in \{1,\dots ,T\}$, let $\Omega_i^{+}$ denote the set of all parameters in $\Omega_i$ that satisfy the condition $(\bar{\mathcal{F}}1)_n$ with respect to $M_n^i$. By lemma \ref{lemma_parameters_prelim_disjoint condn}, the set
    \begin{align*}
        \left\{ \omega\in \Omega_i: d \left( \bigcup_{i,j=1}^{2} I_n^{(i,j)}(\omega), \bigcup_{m=1}^{2K_nM_n^i} \bigcup_{i,j=1}^{2} (I_n^{(i,j)}(\omega) + m\omega)\right) < \varepsilon_n  \right\}
    \end{align*}
    has measure at most $(2^8K_n\tilde{M}_n)(\alpha_c^{\tilde{M}_{n-1}/3}\lambda^{-1+1/32} + 3\alpha_c^{\tilde{M}_{n-1}/3})$ and consists of at most $2^8 (K_n\tilde{M}_n)^2$ components. Thus, for each $i \in \{1,\dots , T\}$,
    \begin{align*}
        |\Omega_i^+| & \geq |\Omega_i| - (2^8K_n\tilde{M}_n)(\alpha_c^{\tilde{M}_{n-1}/4}), \\
       \# \Omega_i^+ & \leq  2^8 (K_n\tilde{M}_n)^2.
    \end{align*}
 Since $\Omega = \bigcup_{i=1}^{T}\Omega_i$ and  $\displaystyle T \leq 1/\delta  = {2^4 \tilde{M}_n}{\alpha_c^{-\tilde{M}_{n-2}}}$, we have 
\begin{align*}
    |\Omega^+| & \geq |\Omega| - {2^4\tilde{M}_n} (2^8K_n\tilde{M}_n)(\alpha_c^{\tilde{M}_{n-1}/4}{\alpha_c^{-\tilde{M}_{n-2}}}),\\
   & > |\Omega| - 2^{11} K_n \tilde{M}_{n}^2 \alpha_{c}^{\tilde{M}_{n-1}/4} \alpha_{c}^{-\tilde{M}_{n-2}}.
\end{align*}
Recall that the sequence $(\tilde{M}_n)_{n\geq 1}$ increases super exponentially, where $\Tilde{M}_n = \alpha_c^{-\tilde{M}_{n-1}/8}$. This gives us the estimate 
\begin{align*}
    |\Omega^+| > |\Omega| - \alpha_c^{\tilde{M}_{n-1}/24} = |\Omega| -\frac{1}{\tilde{M}_{n}^{1/3}}.
\end{align*}
Similarly, the maximum number of intervals in $\Omega^+$ can be estimated to be
\begin{align*}
    \# \Omega^+ & \leq ({2^4 \tilde{M}_n}{\alpha_c^{-\tilde{M}_{n-2}}}) 2^8 (K_n\tilde{M}_n)^2  < \tilde{M}_n^4.
\end{align*}

\end{proof}

%-------------------------- PROP END -------------------------

%%%%%%%%%%%%%%%%%%%%% FINAL PARAMETER SPACE %%%%%%%%%%%%%%%%%%%%
\subsection{Final set of parameters}\label{subsection_parameter set Omega1 and Omega2}

In this subsection, we finally define and estimate the parameter sets $\Omega_1$ and $\Omega_2$, for which we prove the theorem \ref{Theorem SNA}. The construction in this subsection is built upon the parameter sets $\bar{\mathcal{F}}^s(G_2)$ (in definition \ref{definition_parameters_sync_FsG}), $\tilde{\Omega}$ (in definition \ref{definition_parameter space F(1,1)_Omega sets}) and $\Omega$ (in definition \ref{definition_parameter space F empty_Omega sets}), together with the inductive exclusion of parameters established for the sets $\mathcal{F}_0$, $\tilde{\mathcal{F}}_1$ and $\mathcal{F}_n$, for $n\geq 2$. 

\begin{definition}[Parameter set $\Omega_1$]\label{definition_parameters_final_omega1}
    For $\lambda > 1$, we define the parameter set $\Omega_1 \subset \mathbb{T}\setminus \mathbb{Q}$ as
    \begin{align*}
        \Omega_1(\lambda) := \Omega \cap \bar{\mathcal{F}}^s(G_2) \cap \tilde{\mathcal{F}}_1 \cap \left( \bigcap_{l=2}^{\infty} \mathcal{F}_l \right),
    \end{align*}
    where $\bar{\mathcal{F}}^s(G_2)$ is as in definition \ref{definition_parameters_sync_FsG}, and $\Omega$ is as in definition \ref{definition_parameter space F empty_Omega sets}.  
\end{definition}
For each rational $p/q \in \mathbb{T}$, by the definition of the sets $\mathcal{F}_n$, there exists an integer $N\geq 1$ such that $p/q \notin \mathcal{F}_N$ and this shows that $\Omega_1 \subset \mathbb{T} \setminus \mathbb{Q}$.

\begin{definition}[Parameter set $\Omega_2$]\label{definition_parameters_final_omega2}
    For $\lambda > 1$, let us define the set 
    \begin{align*}
        \Omega_2 := \tilde{\Omega} \cap \tilde{\mathcal{F}}_1 \cap (\mathbb{T}\setminus \mathbb{Q}),
    \end{align*}
     where $\tilde{\Omega}$ is the open interval as defined in definition \ref{definition_parameter space F empty_Omega sets}.
\end{definition}

%-------------------- PROP BEGIN ----------------------
\begin{proposition}\label{proposition_parameters_final_omega1}
For all sufficiently large $\lambda$, the set $\Omega_1(\lambda) \subset \mathbb{T}\setminus \mathbb{Q}$ is measurable and satisfies the following:
\begin{enumerate} 
\item $\displaystyle |\Omega_1(\lambda)| > \frac{1}{800}$.
    \item For each $\omega \in \Omega_1$, there exists a sequence of integers $(M_n)_{n=0}^{\infty}$, such that we have $M_n \in [\tilde{M}_n, 2\tilde{M}_n]$ for each $n \geq 1$, and
\begin{align*}
    \omega \in \mathcal{F}^s(G_2) \cap {\mathcal{F}}_1^{(1,1)} \cap \left(\bigcap_{n\geq 0} {\mathcal{F}}_n(M_0,\dots,M_n)\right).
\end{align*}
\end{enumerate}
\end{proposition}

\begin{proof}
By lemma \ref{lemma_parameters_non-empty I2_estimates}, we have the estimates
\begin{align*}
        |\Omega| & \geq \frac{1}{640} \\
      \# \Omega & < M_0^{20}.
\end{align*}
By proposition \ref{proposition_parameters_block I0_lemma} and proposition \ref{proposition_parameters_sync_lemma 1}, we have 
\begin{align*}
    \left|\Omega \cap \bar{\mathcal{F}}^s({G_2}) \cap \mathcal{F}_0\right| & \geq \frac{1}{640} - 32|I_0| - 16K_0M_0|I_0|,\\
    \# \left( \Omega \cap \bar{\mathcal{F}}^s({G_2}) \cap \mathcal{F}_0 \right) & < 4K_0^2M_0^{136}.
\end{align*}
Since $K_0=2^8(b+1)$, $M_0 = \lambda^{1/16}$, $\tilde{M}_1 = \alpha_c^{-M_0/8}$ and $|I_0| = 2\lambda^{-1/4}$, for a large enough $\lambda$, the above estimates can be reduced to
\begin{align}\label{equation_paramters omega1_partial exclusion}
    \left|\Omega \cap \bar{\mathcal{F}}^s({G_2}) \cap \mathcal{F}_0\right| & \geq \frac{1}{650},\\
    \# \left( \Omega \cap \bar{\mathcal{F}}^s({G_2}) \cap \mathcal{F}_0 \right) & < \tilde{M}_1. \label{equation_paramters omega1_partial exclusion 1}
\end{align}
For each interval in the set $\left(\Omega \cap \bar{\mathcal{F}}^s({G_2}) \cap \mathcal{F}_0\right)$, we have the estimates in proposition \ref{proposition_parameters_block I1_lemma 2} for the parameter exclusion with respect to the condition $(\mathcal{F}1-2)_1$ and $(\bar{\mathcal{F}1})_1$. Thus, it follows from the estimates in equations (\ref{equation_paramters omega1_partial exclusion} - \ref{equation_paramters omega1_partial exclusion 1}) and proposition \ref{proposition_parameters_block I1_lemma 2} that
\begin{align}\label{equation_paramters omega1_partial exclusion 2}
    \left| \Omega \cap \bar{\mathcal{F}}^s({G_2}) \cap \mathcal{F}_0 \cap \tilde{\mathcal{F}}_1 \right| & \geq \frac{1}{650} - \frac{\tilde{M}_1}{\tilde{M}_1^{2}} = \frac{1}{650} - \frac{1}{\tilde{M}_1}, \\
    \# \left( \Omega \cap \bar{\mathcal{F}}^s({G_2}) \cap  \mathcal{F}_0 \cap \tilde{\mathcal{F}}_1 \right) & < \tilde{M}_1^{5}.\label{equation_paramters omega1_partial exclusion 3}
\end{align}
By lemma \ref{lemma_parameters_non-empty I2_estimates}, we have $\left(\Omega \cap \bar{\mathcal{F}}^s({G_2}) \cap \tilde{\mathcal{F}}_1\right) \subset \tilde{\mathcal{F}}^{(1,1)}_1(M_0, M_1)$. Then, for every $\omega \in \left(\Omega \cap \bar{\mathcal{F}}^s({G_2}) \cap \tilde{\mathcal{F}}_1\right)$, by proposition \ref{proposition_structure I2_I2  nonempty F^{(1,1)}} and proposition \ref{proposition_structure In_In empty nonempty}, the critical set $I_2^{(1,1)}(\omega)$ is non-empty. Similarly, for each $n\geq 2$ and $\omega \in \Omega \cap \bar{\mathcal{F}}^s({G_2}) \cap \tilde{\mathcal{F}}_1 \cap \left(\cap_{l=2}^{n-1}\mathcal{F}_l\right)$, by proposition \ref{proposition_structure In+1_critical set I_n+1}, the critical set $I_n^{(1,1)}(\omega)$ must be non-empty. This enables us to build the following series of parameter exclusions.

For each $n\geq 2$, and each interval in the set $\left| \Omega \cap \bar{\mathcal{F}}^s({G_2}) \cap \tilde{\mathcal{F}}_1 \cap \left(\cap_{l=2}^{n-1}\mathcal{F}_l\right)  \right|$, we have the estimates in proposition \ref{proposition_parameters_Block In_lemma 1} for the parameter exclusion with respect to the condition $(\mathcal{F}1-2)_n$ and $(\bar{\mathcal{F}1})_n$. Since the sequence $(\tilde{M}_n)_{n\geq 0}$ increases super exponentially for every $\lambda>1$, it follows from the estimates in equations (\ref{equation_paramters omega1_partial exclusion 2} - \ref{equation_paramters omega1_partial exclusion 3}) and proposition \ref{proposition_parameters_Block In_lemma 1} that
\begin{align*}
    \left| \Omega \cap \bar{\mathcal{F}}^s({G_2}) \cap \tilde{\mathcal{F}}_1 \cap \left(\cap_{l=2}^{n}\mathcal{F}_l\right)  \right| & \geq \frac{1}{650} - \frac{1}{\tilde{M}_1} - \sum_{l=2}^{n} \frac{\tilde{M}_{n-1}^5}{\tilde{M}_n^{1/3}}\\
    & > \frac{1}{650} -  \sum_{l=1}^{n} \frac{1}{\tilde{M}_n^{1/6}}. 
\end{align*}
and 
\begin{align*}
    \# \left( \Omega \cap \bar{\mathcal{F}}^s({G_2}) \cap \tilde{\mathcal{F}}_1 \cap \left(\cap_{l=2}^{n}\mathcal{F}_l\right) \right) & < \tilde{M}_1^{5} \prod_{l=2}^{n} \tilde{M}_{l}^4 < \tilde{M}_n^5.
\end{align*}
For a sufficiently enough $\lambda$ with respect to $\alpha_c$ and $b$, the sum $\displaystyle \sum_{n=1}^{\infty} \frac{1}{\tilde{M}_n^{1/3}}$ is bounded and we have the estimate
\begin{align*}
     \left| \Omega \cap \bar{\mathcal{F}}^s({G_2}) \cap \tilde{\mathcal{F}}_1 \cap \left(\cap_{l=2}^{\infty}\mathcal{F}_l\right) \right| \geq \frac{1}{650} - \sum_{n=1}^{\infty}\frac{1}{\tilde{M}_n^{1/6}} > \frac{1}{800}.
\end{align*}
It is worth noting that the lower bound conditions imposed on $\lambda$ so far, remains a finite condition. By construction the set $ \displaystyle \Omega \cap \bar{\mathcal{F}}^s({G_2}) \cap \tilde{\mathcal{F}}_1 \cap \left(\cap_{l=2}^{\infty}\mathcal{F}_l\right)$ is measurable. By the definition of the set $\Omega_1$, 
\begin{align*}
    \Omega_1 = \Omega \cap \bar{\mathcal{F}}^s(G_2) \cap \tilde{\mathcal{F}}_1 \cap \left( \bigcap_{l=2}^{\infty} \mathcal{F}_l \right) & \subset  \mathcal{F}^s({G_2}) \cap \mathcal{F}^{(1,1)}_1 \cap \left(\bigcap_{n=0}^{\infty} \mathcal{F}_n\right),
\end{align*}
and this proves that $\displaystyle |\Omega_1(\lambda)| > {1}/{800}$, for all large $\lambda$.

\end{proof}
%-------------------- PROP END ----------------------

%-------------------- PROP BEGIN ----------------------
\begin{proposition}\label{proposition_parameters_final_omega2}
For all sufficiently large $\lambda$, 
\begin{enumerate}
  \item The set $(\Omega_2(\lambda) \cup \mathbb{Q})$ is open and 
  \begin{align*}
      |\Omega_2| & > \frac{1}{M_0^{12}}.
  \end{align*}
    \item For each $\omega \in \Omega_2$, we have $\omega \in \mathcal{F}_0(M_0)\cap \tilde{\mathcal{F}}_1$ and $I_2(\omega) = \emptyset$.
\end{enumerate}
\end{proposition}

\begin{proof}
By propositions \ref{proposition_parameters_empty I2_lemma 00}, we have $\tilde{\Omega} \subset \mathcal{F}_0$. Since $\tilde{\Omega}$ is an open interval, by proposition \ref{proposition_parameters_block I1_lemma 2}, we obtain the estimates
\begin{align*}
   |\Omega_2| = |\tilde{\Omega} \cap \tilde{\mathcal{F}}_1| &
    > \frac{1}{K_0M_0^{11}} - \frac{1}{\tilde{M}_1^2} > \frac{1}{M_0^{12}}. \\
    \# \Omega_2 & \leq \tilde{M}_1^4.
\end{align*}
Thus, the set $\Omega_2$ is measurable and $(\Omega_2 \cup \mathbb{Q})$ is open. By proposition \ref{proposition_parameters_empty I2_lemma 00}, for each  $\omega \in \tilde{\Omega} \cap \tilde{\mathcal{F}}_1$ the critical set $I_2(\omega)$ is empty. 
\end{proof}

%-------------------- PROP END ----------------------

It is worth observing that, for the purpose of theorem \ref{Theorem SNA}, we do not seek an optimal set of parameters $\Omega_2$ for which the system has continuous invariant curves. Instead, we obtain a specific non-empty open set for which the conditions in theorem \ref{Theorem SNA} \ref{Theorem SNA b} are satisfied.

%%%%%%%%%%%%%%%%%%%%%%%%%%%%%%%%%%%%%%%%%%%%%%%%%%%%%%%%%%%

%%%%%%%%%%%%%%%%%%%%% PROOF OF THE MAIN THEOREM %%%%%%%%%%%%%%%%%%%%%%%%%%%%%
\section{Proof of the Theorem \ref{Theorem SNA}}\label{section_proof of thm}

This section contains the proof of theorem \ref{Theorem SNA}. For a fixed $a>0$ satisfying the hypothesis of the theorem, recall the parameter set $\Omega_1(\lambda)$ as in definition \ref{definition_parameters_final_omega1}, and $\Omega_2(\lambda)$ as in definition \ref{definition_parameters_final_omega2}. In propositions \ref{proposition_parameters_final_omega1} and \ref{proposition_parameters_final_omega2}, for sufficiently large $\lambda$, we proved that  $\Omega_1$ has positive measure, and that $(\Omega_2 \cup \mathbb{Q})$ is an open set. For parameters in these sets, the geometric structure of critical sets derived in section \ref{section_structure of critical sets} and the mechanism of slow-returns introduced in section \ref{section_Mechanism of slow-returns} collectively form a cohesive induction. This enables the construction and analysis of the invariant graphs below.

%%%%%%%%%%%%%%%%%%%%%%%%%%%%%% (a) SNA CONSTRUCTION %%%%%%%%%%%%%%%%%%%%%%%%%

\vspace{0.2cm}
\noindent\textbf{Proof of the Theorem \ref{Theorem SNA} (a):} Let us fix a parameter $\omega \in \Omega_1$, which by definition \ref{definition_parameters_final_omega1} must be irrational. By proposition \ref{proposition_parameters_final_omega1}, there exists a corresponding sequence $\displaystyle (M_n)_{n=0}^{\infty} = (M_n(\omega))_{n=0}^{\infty}$ such that $\displaystyle \omega \in {\mathcal{F}^s(G_2)} \cap \tilde{\mathcal{F}}^{(1,1)} \cap \left( \cap_{n\geq 0} \mathcal{F}_n(M_0,\dots ,M_n) \right)$ and $\displaystyle M_n \in [\tilde{M}_n, 2\tilde{M}_n]$ for $\tilde{M}_n := \alpha_c^{-\tilde{M}_{n-1}/8}$, for all $n\geq 1$. Then, with respect to the constants $(K_n)_{n\geq 1}$ and $(M_n)_{n\geq 0}$, we can define a sequence of critical sets $\displaystyle (I_n(\omega))_{n\geq 0}$, where $I_0(\omega) = I_0$ and $I_n(\omega)$ for each $n\geq 1$ is as in definition \ref{definition_structure In_critical set In}. For $n=0$, we have $|I_0| = 2M_0^{-4} = 2\lambda^{-1/4}$. By propositions \ref{proposition_structure I1_critical set I_1}, \ref{proposition_structure I2_critical set I2} and \ref{proposition_structure In+1_critical set I_n+1}, all the critical sets in the sequence are in fact non-empty and satisfy the conditions 
\begin{align}
    \text{comp}(I_{n}) \leq 4, \ \ I_{n}(\omega) \subset I_{n-1}(\omega) \ \ \mbox{   and   } \ \ 0 < |I_n(\omega)| < \lambda^{-1/8} \alpha_c^{M_{n-1}/3},
\end{align}
for each $n\geq 1$, where $\comp(\cdot)$ denote the number of components. From definitions \ref{definition_Block 0_Big Theta} and \ref{definition_block n_Big theta}, recall the sequence of nested sets $\left( \tilde{\Theta}_n \right)_{n\geq 1}$ in $\mathbb{T}$ given by 
\begin{align*}
    \tilde{\Theta}_{n} & = \Theta_0 \setminus \bigcup_{l=1}^{n}\bigcup_{|m|\leq 2M_l}(I_l + m\omega),
\end{align*}
where $\Theta_{0} = \mathbb{T}\setminus \bigcup_{|m|\leq M_0 - 1}(I_0 + m\omega)$. Using the sets in $\left( \tilde{\Theta}_n \right)_{n\geq 1}$, we now define the following measurable subsets in $\mathbb{T}$:
\begin{align}
        {\Theta}_{\infty} & := \mathbb{T} \setminus \bigcup_{l=0}^{\infty} \bigcup_{|m| \leq 2M_l} (I_l + m\omega) 
\stepcounter{equation}\tag{\theequation}\label{equation_final proof_Theta infty}\\
        \Theta & := \bigcap_{m= -S_0}^{S_0} ({\Theta}_{\infty} + m\omega) \stepcounter{equation}\tag{\theequation}\label{equation_final proof_Theta}\\
        B & := \bigcup_{l=0}^{\infty} \bigcup_{m =1 }^{M_l +S_0} (I_l + m\omega) \stepcounter{equation}\tag{\theequation}\label{equation_final proof_B} \\
        \tilde{B} & := \bigcup_{l=0}^{\infty} \bigcup_{m = 0}^{M_l + S_0} (I_l - m\omega).
    \stepcounter{equation}\tag{\theequation}\label{equation_final proof_tilde B}
\end{align}
By construction, the sets defined above satisfy the relation $\Theta \subset {\Theta}_{\infty} \subset  \bigcap_{n\geq 1} \tilde{\Theta}_n$, and the bounds
\begin{align*}
     |{\Theta}_{\infty}| & > 1 -  \lambda^{-1/8} - \sum_{k=0}^{\infty}  \alpha_c^{\tilde{M}_k / 8} > 0 \\
    |{\Theta}| & > 1 - (3S_0 + 1) \left( \lambda^{-1/8} + \sum_{k=0}^{\infty} \alpha_c^{\tilde{M}_k/8} \right) > 0 \\
   0 < |B| & < 4\lambda^{-1/8} + \sum_{k=0}^{\infty} 8\alpha_c^{\tilde{M}_{k}/8} \\
   0 < |\tilde{B}| & < 4\lambda^{-1/8} + \sum_{k=0}^{\infty} 8\alpha_c^{\tilde{M}_{k}/8} 
\end{align*}

From the above estimates, we observe that the measure of $\Theta$ and $\Theta_{\infty}$ both converge to 1 as $\lambda$ goes to infinity. Whereas the measure of $B$ and $\tilde{B}$ both converge to $0$ as $\lambda$ goes to $\infty$. Then, for a sufficiently large $\lambda$, we can ensure that $|B|, |\tilde{B}| < {1}/{6b}$, where $b$ is the smallest positive integer for which $\alpha_u\alpha_c^{b}<1$. By ergodicity of the base dynamics, which is an irrational rotation, and the Birkhoff ergodic theorem, there exists a set $\tilde{\Theta}_{\infty} \subset {\Theta}_{\infty}$ such that $|\tilde{\Theta}_{\infty}| = |{\Theta}_{\infty}|$ and for each $\theta \in \tilde{\Theta}_{\infty}$, the following holds:
\begin{align}\label{equation_final proof_ergodicity and 1/6b visits 1}
    \lim_{n\rightarrow \infty} \frac{1}{n} \sum_{k=0}^{n} 1_{B}(\theta - k\omega) = |B| & < \frac{1}{6b} \\
    \lim_{n\rightarrow \infty} \frac{1}{n} \sum_{k=0}^{n} 1_{\tilde{B}}(\theta + k\omega) = |\tilde{B}| & < \frac{1}{6b}.
    \stepcounter{equation}\tag{\theequation}\label{equation_final proof_ergodicity and 1/6b visits 2}
\end{align}
Note that, for each $\theta \in \mathbb{T}$ both the orbits $(\theta - m\omega)_{m\geq 1}$ and $(\theta + m\omega)_{m\geq 1}$ intersect the set $\tilde{\Theta}_{\infty}$ infinitely often. Similarly, let us define a set $\tilde{\Theta} \subset \Theta$ as
\begin{align}\label{equation_final proof_tilde Theta 1}
    \tilde{\Theta} := \bigcap_{m=-S_0}^{S_0} (\tilde{\Theta}_{\infty} + m\omega) \subset \bigcap_{m=-S_0}^{S_0} (\Theta_{\infty} + m\omega),
\end{align}
for which by construction $|\tilde{\Theta}| = |\Theta|$.

\vspace{0.2cm}
\noindent\textit{Construction of the graphs $\Phi(\theta)$ and $\Psi(\theta)$:} 
\vspace{0.1cm}

Let us fix a point $\theta \in \tilde{\Theta}_{\infty}$. Then, by equations \eqref{equation_final proof_ergodicity and 1/6b visits 1} and \eqref{equation_final proof_ergodicity and 1/6b visits 2}, there exists a constant $N_1 > 2\tilde{M_1}$ such that $(\theta - N_1\omega) \in \tilde{\Theta}_{\infty}$ and
\begin{align}\label{equation_final proof_bad iterates}
    0 \leq \# \left\{ 0\leq m \leq n : (\theta - m\omega) \in B \right\} < \frac{n}{3b}, \ \ \text{ for all } \ n\geq N_1.
\end{align}
Further, there exists an increasing sequence $(N_k)_{k \geq 1}$ corresponding to $\theta$ such that for each $k\geq 2$ we have $(\theta - N_k\omega) \in \tilde{\Theta}_{\infty}$ and $N_{k-1} - N_{k} > 2\tilde{M}_1$. 

For the fixed $\theta$ and sequence $(N_k)_{k\geq 1}$, by the definition of $\tilde{\Theta}_{\infty}$, there exists an increasing sequence of positive integers $(n_k)_{k\geq 1}$ for which 
\begin{align*}
    \{\theta - m\omega: m=0, \dots ,N_k\} \cap I_{n_k} = \emptyset.
\end{align*}
For each $k\geq 2$, since $\tilde{\Theta}_{\infty} \subset \tilde{\Theta}_{n_k}$, it follows that $(\theta - N_k\omega), (\theta - N_{k-1}\omega) \in \tilde{\Theta}_{n_k}$. Since $N_k - N_{k-1} \geq 2\tilde{M}_1$, by $(S1)_{n_k}$ (in proposition \ref{proposition_block n sync' contraction}) and the continuity of the map $F$, we obtain
\begin{align}\label{equation_final proof_ N_k nested sets on fiber}
    F^{N_k - N_{k-1}} \left( \{\theta - N_k\omega\} \times C_1 \cup C_2 \right) & \subset \{\theta - N_{k-1}\omega \} \times C_{i_k}, \text{ for some } i_k \in \{1,2\},
\end{align}
for each $k\geq 2$. Since $\theta \notin B$, by $(Sc1)_{n_{k-1}}$ (in proposition \ref{proposition_block n sync contraction}) and $(C1)_{n_{k-1}}$ (in proposition \ref{proposition_block n contraction}), there exists $i\in \{1,2\}$ such that
\begin{align}
F^{N_{k-1}} \left( \{\theta - N_{k-1}\omega\} \times C_{i_k} \right) \subset \{\theta\} \times C_i,
\end{align}
for each $k\geq 2$. Then, by equation \eqref{equation_final proof_ N_k nested sets on fiber}, this gives us
\begin{align}
F^{N_k} \left( \{\theta - N_k\omega\} \times C_1 \cup C_2 \right) & \subset F^{N_{k-1}} \left( \{\theta - N_{k-1}\omega\} \times C_{i_k} \right) \subset \{\theta\} \times C_i,
\end{align}
for each $k\geq 2$, which proves that $\left(F^{N_k} (\{\theta - N_k \omega\} \times (C_1 \cup C_2))\right)_{k\geq 2}$ is a sequence of nested intervals in $\{\theta\} \times C_i$, for some $i\in \{1,2\}$.

For any $k\geq 1$ and $x\in C_1\cup C_2$, by the definition of set $B$ (in equation \eqref{equation_final proof_B}) and property $(C1)_{n_{k-1}}$ (in proposition \ref{proposition_block n  contraction}), for every $m \in \{0, \dots ,N_{k-1}\}$ we have the following condition:
\begin{align*}
    \pi_2 \circ F^m(\theta - N_{k-1}\omega, x) \notin C_1\cup C_2 \implies (\theta - (N_{k-1}-m)\omega) \in B. 
\end{align*}
Recall that an orbit experiences a uniform contraction along the fibre when it intersects the region $\mathbb{T}\times (C_1\cup C_2)$, as defined in definition \ref{definition_preliminaries_expanding contracting regions}. Then, by equation \eqref{equation_final proof_bad iterates} and property $(C1)_{n_{k-1}}$(in proposition \ref{proposition_block n  contraction}), for each $k\geq 2$ and $x\in C_1\cup C_2$, we have the following rate of contraction under forward iteration:
\begin{align}\label{equation_final proof_contraction rate N_k-1}
    \left|\partial_x \left(\pi_2 \circ F^{N_{k-1}} \left((\theta - N_{k-1}\omega), x \right) \right)\right| & < \left(\alpha_u^{1/3b} \alpha_c^{1 - 1/3b}\right)^{N_{k-1}} < \alpha_c^{(-1/3 + 1 - 1/3b)N_{k-1}} < \alpha_c^{N_{k-1}/3}. 
\end{align}
By equation \eqref{equation_final proof_ N_k nested sets on fiber}, the estimates in equation \eqref{equation_final proof_contraction rate N_k-1} and continuity of $F$, it follows that the one-dimensional Lebesgue measure of each set in the sequence $\left(\pi_2 \circ F^{N_k} \left( \{\theta - N_k\omega\} \times (C_1 \cup C_2) \right)\right)_{k\geq 1}$ satisfies the following bound:
\begin{align*}
 \left| \pi_2 \circ F^{N_k} \left( \{\theta - N_k\omega\} \times (C_1 \cup C_2) \right) \right| & < \alpha_c^{N_{k-1}/3}, \text{ for all } k\geq 2,
\end{align*}
which gives us
\begin{align*}
    \lim_{k\rightarrow \infty} \left| \pi_2 \circ F^{N_k} \left( \{\theta - N_k\omega\} \times (C_1 \cup C_2) \right) \right| = 0.
\end{align*}
This proves that the sequence of nested intervals on the fibre $\left(\pi_2 \circ F^{N_k} \left( \{\theta - N_k\omega\} \times (C_1 \cup C_2) \right)\right)_{k\geq 1}$ converges to a single point. 

For each $k\geq 1$ and any $N_k < n < N_{k+1}$ for which $(\theta - n\omega) \in \tilde{\Theta}_{\infty}$, we have 
\begin{align*}
F^{n - N_k} \left( \{\theta - n\omega\} \times (C_1 \cup C_2) \right) & \subset \left( \{\theta - N_{k}\omega\} \times (C_1 \cup C_2) \right),
\end{align*}
by $(C1)_{n_{k+1}}$ (in proposition \ref{proposition_block n contraction}). This proves that the point of convergence of the sets $$\left(F^{N_k} \left( \{\theta - N_k\omega\} \times C_1 \cup C_2 \right)\right)_{k\geq 1}$$
is independent of the choice of sequence $(N_k)_{k\geq 1}$.

Thus, for each $\theta \in \tilde{\Theta}_{\infty}$ there exists a unique $\Phi(\theta) \in \mathbb{T}$ such that
\begin{align}\label{equation_final proof_partial invariant graph on Theta_infty}
    \Phi(\theta) := \lim_{k\rightarrow \infty}  \pi_2 \circ F^{N_k} \left( (\theta - N_k\omega), x \right), 
\end{align}
for all $x\in C_1 \cup C_2$. By the construction of the function $\Phi$ and $(C1)_n$ (in proposition \ref{proposition_block n contraction}), for every $\theta \in \tilde{\Theta}_{\infty}$ we have $\Phi(\theta) \in C_1 \cup C_2$. Further, for each $\theta \in \tilde{\Theta}_{\infty}$ the map $\Phi$ satisfies the relation
\begin{align}\label{equation_final proof_partial invariance on Theta_infty}
   F^n(\theta, \Phi(\theta)) =  \left( \theta + n\omega, \Phi(\theta + n\omega) \right), 
\end{align}
whenever $(\theta + n\omega) \in \tilde{\Theta}_{\infty}$. For any $\theta \in \tilde{\Theta}_{\infty}$ and $x\in C_1\cup C_2$, by equation \eqref{equation_final proof_ergodicity and 1/6b visits 1}, we have the following rate of contraction under forward iteration:
\begin{align}\label{equation_final proof_contraction rate}
    \left|\partial_x \left(\pi_2 \circ F^{n} \left(\theta , x \right) \right)\right| & < \left(\alpha_u^{1/3b} \alpha_c^{1 - 1/3b}\right)^{n} < \alpha_c^{(-1/3 + 1 - 1/3b)n} < \alpha_c^{n/3}, 
\end{align}
for all sufficiently large $n\geq 1$. This gives us the following bound for the forward Lyapunov exponent:
\begin{align}\label{equation_final proof_Lyapunov in contracting region}
    \overline{L}^{+}(\theta, x) < \frac{\log \alpha_c}{3}  < 0.
\end{align}
Since $|\tilde{\Theta}_{\infty}| > 0$ and the function $\Phi: \tilde{\Theta}_{\infty} \rightarrow \mathbb{T}$ satisfies the equation \eqref{equation_final proof_partial invariance on Theta_infty}, $\Phi$ can be naturally extended to a measurable function $\Phi:\mathbb{T} \rightarrow \mathbb{T}$ which satisfies
\begin{align}\label{equation_final proof_full invariance of graph phi}
    F(\theta, \Phi(\theta)) = (\theta+\omega, \Phi(\theta + \omega)),
\end{align}
for almost every $\theta \in \mathbb{T}$. By the construction of $\Phi$ (in equation \eqref{equation_final proof_partial invariant graph on Theta_infty} and equation \eqref{equation_final proof_full invariance of graph phi}), for every $\theta \in \tilde{\Theta}_{\infty}$ and $x \in C_1 \cup C_2$ the point $(\theta,x)$ converges to the invariant graph $\Phi$ under forward iteration, that is, 
\begin{align}\label{equation_final proof_partial attraction from C}
    \lim_{n \rightarrow \infty} |\pi_2 \circ F^n(\theta, x) - \Phi(\theta + n\omega)| = 0.
\end{align}

Let us now fix a point $\theta \in \tilde{\Theta}$, where the set $\tilde{\Theta}$ given by 
\begin{align*}
    \tilde{\Theta} = \bigcap_{m=-S_0}^{S_0}(\tilde{\Theta}_{\infty}+m\omega),
\end{align*}
as in equation \eqref{equation_final proof_tilde Theta 1}. Then, we have $\{\theta + n\omega\}_{n=0}^{S_0} \cap I_0 = \emptyset$ and $(\theta + S_0\omega) \in \tilde{\Theta}_{\infty}$. By the choice of constant $S_0$ and equations (\ref{preliminaries_f contraction and expansion S01} - \ref{preliminaries_f contraction and expansion S04}), it follows that
\begin{align*}
    x_{S_0} := \pi_2 \circ F^{S_0}(\theta, x) \in C_1 \cup C_2, \ \ \text{ for all } x\in E_1^c \cap E_2^c.
\end{align*}
For every $x\in (E_1^{c} \cap E_2^{c})$, since $(\theta +S_0\omega) \in \tilde{\Theta}_{\infty}$ and $x_{S_0} \in (C_1 \cup C_2)$, the orbit $(F^n(\theta+{S_0}\omega,x_{S_0}))_{n\geq 1}$ must converge to the invariant graph $\Phi$ under the forward iteration, that is,  
\begin{align*}
    \lim_{n \rightarrow \infty} |\pi_2 \circ F^n(\theta + S_0, x_{S_0}) - \Phi(\theta + (S_0+ n)\omega)| = 0,
\end{align*}
which follows from equation \eqref{equation_final proof_partial attraction from C}. This proves that, for all $\theta \in \tilde{\Theta}$ and $x\in E_1^{c} \cap E_2^{c}$, 
\begin{align}\label{equation_final proof_partial attraction from E^c}
    \lim_{n \rightarrow \infty} |\pi_2 \circ F^n(\theta, x) - \Phi(\theta+n\omega)| = 0.
\end{align}

Similarly, the estimate $|\tilde{B}| < 1/6b$, the properties  $(E1-2)_n$ (in proposition \ref{proposition_block n expansion}), $(\bar{S}1)_n$ (in proposition \ref{proposition_block n sync' expansion}) and $(Se1)_n$ (in proposition \ref{proposition_block n sync expansion}), enable the construction of a well-defined invariant measurable function $\Psi: \mathbb{T} \rightarrow \mathbb{T}$ for which
\begin{align}\label{equation_final proof_partial repulsion E}
    \lim_{n \rightarrow \infty} |\pi_2 \circ F^{-n}(\theta, x) - \Psi(\theta - n\omega)| = 0,
\end{align} 
for every $\theta \in \tilde{\Theta}_{\infty}$ and $x \in E_1 \cup E_2$. For such an invariant graph, for every $\theta \in \tilde{\Theta}$ and any $x \in C^c_1 \cap C^c_2$, the orbit of $(\theta,x)$ under backward iteration converges to the graph $\Psi$, that is,
\begin{align}\label{equation_final proof_partial repulsion C^c}
    \lim_{n \rightarrow \infty} |\pi_2 \circ F^{-n}(\theta, x) - \Psi(\theta - n\omega)| = 0. 
\end{align} 
By $(E1)_n$, for every $\theta \in \tilde{\Theta}_{\infty}$ we have $\Psi(\theta) \in E_1 \cup E_2$. For every $\theta \in \tilde{\Theta}_{\infty}$ and $\Psi(\theta) \in E_1 \cup E_2$, we also have the following bound for the backward Lyapunov exponent:
\begin{align}\label{equation_final proof_Lyapunov in expanding region}
    \overline{L}^{-}(\theta, x) = \overline{L}^{-}(\theta, \Psi(\theta)) < \frac{\log \alpha_c}{3}.
\end{align}

%------------------------- Uniqueness of invariant graph -------------------- 
\vspace{0.2cm}
\noindent\textit{Uniqueness of the attracting(repelling) invariant graph:} 
\vspace{0.1cm}

Let us fix a point $\theta \in \tilde{\Theta}$. Then, by equation \eqref{equation_final proof_ergodicity and 1/6b visits 2}, there exists $N > 0$ such that  
\begin{align}\label{equation_final proof_bad region visits 0}
    \# \{ 0 \leq j \leq n: (\theta + j\omega) \in \tilde{B}\} < \frac{n}{3b}, \text{ for all } n\geq N.
\end{align}
Let us fix $y \neq \Psi(\theta)$. Then, there exists $\varepsilon > 0$ such that $|\psi(\theta) - y| > \varepsilon$. Let us choose a constant $\tilde{k} > N + 2\tilde{M}_1$ that satisfies $(\theta + \tilde{k}\omega) \in \tilde{\Theta}$ and $\alpha_c^{(\tilde{k}-2\tilde{M}_1)/3} < \varepsilon$.

%Then, by equation \eqref{equation_final proof_bad region visits 0}, for any $x\in  E_1 \cup E_2$, we have the following estimates for the rate of contraction for the backward iteration 
%\begin{align*}
%    \left|\partial_x \left(\pi_2 \circ F^{-\tilde{k}} (\theta + (\tilde{k}-2\tilde{M})\omega, x)\right)\right| =|(f^{-1})'(f(x_{-(\tilde{k}-2\tilde{M})})) \dots (f^{-1})'(f(x_{-1}))| &  < \alpha_c^{(\tilde{k}-2\tilde{M})/3}.
%\end{align*}
Since $(\theta +\tilde{k}\omega) \in \tilde{\Theta}$, by the construction of the invariant graph $\Psi$, it follows that $\Psi(\theta) \in E_1\cup E_2$ and $\Psi(\theta+\tilde{k}\omega) \in E_1\cup E_2$. For any $x \in E_1 \cup E_2$, by $(E1)_n$ (in proposition \ref{proposition_block n expansion}),  
\begin{align*}
    \pi_2 \circ F^{-j} (\theta + \tilde{k}\omega, x) \notin E_1 \cup E_2 \implies (\theta - j\omega) \in \tilde{B}, \text{ for all } j = 0,\dots ,\tilde{k}.
\end{align*}

Note that, for $\theta \in \tilde{\Theta}$ and fixed $\tilde{k}$, there exists $n\geq 2$ such that $(\theta+\tilde{k}\omega) \in \tilde{\Theta}_n$. Then, by $(\bar{S}1)_n$ (in proposition \ref{proposition_block n sync' expansion}) and $(Se1)_n$ (in proposition \ref{proposition_block n sync expansion}), there exists $N_s<2\tilde{M}_1$ such that $(\theta+(\tilde{k}-N_s)\omega)\in \Theta_n$ and 
\begin{align*}
    \pi_2\circ F^{-N_s}\left(\{\theta+\tilde{k}\omega\}\times (E_1\cup E_2)\right) \subset E_i, 
\end{align*}
for some $i\in \{1,2\}$. Then, by equation \eqref{equation_final proof_bad region visits 0} and $(E1)_n$ (in proposition \ref{proposition_block n expansion}), for every $x\in E_1\cup E_2$ we obtain the estimate 
\begin{align*}
   \left|\pi_2\circ F^{-\tilde{k}}\left((\theta+\tilde{k}\omega), \Psi(\theta+\tilde{k}\omega)\right) - \pi_2 \circ F^{-\tilde{k}} (\theta + \tilde{k}\omega, x)\right| & = |\Psi(\theta) - \pi_2 \circ F^{-\tilde{k}} (\theta + \tilde{k}\omega, x)| \\
   & < (\alpha_c^{1/3})^{\tilde{k} - 2\tilde{M}_1} < \varepsilon. 
\end{align*}
 Since $\varepsilon$ was chosen such that $|\Psi(\theta) - y| > \varepsilon$, the above estimate implies that 
\begin{align*}
    F^{\tilde{k}}(\theta, y) & \in \tilde{\Theta} \times (E_1^c\cap E_2^c).
\end{align*} 
Then, by the attracting property of the invariant graph $\Phi$ as in equation \eqref{equation_final proof_partial attraction from E^c},
\begin{align}\label{equation_final proof_synceqn1a}
    \lim_{n\rightarrow \infty} |\pi_2 \circ F^n(\theta, y) - \Phi(\theta + n\omega)| =  \lim_{n\rightarrow \infty} |\pi_2 \circ F^{\tilde{k}+n}(\theta, y) - \Phi(\theta + (\tilde{k}+n)\omega)| = 0.
\end{align}
This proves that for every $\theta \in \tilde{\Theta}$ and $y \neq \Psi(\theta)$,
\begin{align*}
    \lim_{n\rightarrow \infty} |\pi_2 \circ F^n(\theta, y) - \Phi(\theta + n\omega)| = 0.
\end{align*}
Since the set $\tilde{\Theta}$ has positive measure and $\Psi$ is an invariant graph, this gives us 
\begin{align*}
    \lim_{n\rightarrow \infty} |\pi_2 \circ F^n(\theta, y) - \Phi(\theta + n\omega)| = 0,
\end{align*}
for almost every $\theta \in \mathbb{T}$ and for every $y\neq \Psi(\theta)$, which proves Theorem \ref{Theorem SNA} \ref{Theorem SNA_a}. The proof of Theorem \ref{Theorem SNA} \ref{Theorem SNA b} follows similarly. Note that, by equations \eqref{equation_final proof_Lyapunov in contracting region} and \eqref{equation_final proof_Lyapunov in expanding region}, for almost every $(\theta,x)\in \mathbb{T}^2$,
\begin{align}\label{equation_final proof_Lyapunov final bound}
    \overline{L}^{+}(\theta, x) & = \overline{L}^{+}(\theta, \Phi(\theta))  < \frac{\log \alpha_c}{3}  < 0 \\
    \overline{L}^{-}(\theta, x) & = \overline{L}^{-}(\theta, \Psi(\theta)) < \frac{\log \alpha_c}{3}  < 0. 
\end{align}

%%%%%%%%%%%%%%%%%%%%%%%% SINK-SOURCE ORBIT %%%%%%%%%%%%%%%%%%%%%%%%%%%%%
\vspace{0.2cm}
\noindent\textit{Construction of sink-source orbit:} 
\vspace{0.1cm}

Recall that an orbit $\displaystyle (F_{\omega}^n(\theta,x))_{n\in \mathbb{Z}}$ is called a sink-source orbit if its corresponding Lyapunov exponents satisfy the bounds $\underline{L}^{+}(\theta, x) > 0$ and $\underline{L}^{-}(\theta, x) > 0$, where the lower Lyapunov exponents are as defined in equations (\ref{equation_intro_lower_forward Lyapunov} - \ref{equation_intro_lower_backward Lyapunov}).

For the fixed parameter $\omega \in \Omega_1$ and for each $n\geq 1$, we consider the set
\begin{align}\label{equation_final proof_critical regions}
   V_n := \left( F^{M_{n}}(A_{n}) \cap F^{-M_{n}}(B_{n}) \right),
\end{align}
 where 
\begin{align*}
    A_{n} & = \left( I_{n}-M_{n}\omega\right) \times (C_1\cup C_2) \\
    B_{n} & = \left(I_{n}+M_{n}\omega\right)\times (E_1 \cup E_2).
\end{align*}
By lemma \ref{lemma_structure I2_match Ai Bi} and lemma \ref{lemma_structure In_match Ai Bi}, for each $n\geq 1$ we have $V_{n+1} \subset V_n$. From definition \ref{definition_block n_Critical sets In}, recall the  preliminary critical sets given by 
\begin{align*}
    \tilde{I}_{n+1}(\omega) = \pi_1 \left( F^{M_{n}}(A_{n}) \cap F^{-M_{n}}(B_{n}) \right) =  \pi_1 (V_n).
\end{align*}
Since the critical set $I_n(\omega)$ is non-empty for the fixed $\omega \in \Omega_1$, by the definition of critical set $I_n(\omega)$ in definition \ref{definition_block n_Critical sets In}, the set $\tilde{I}_{n+1}(\omega)$ must also be non-empty. This implies that the nested sequence of closed sets $(V_n)_{n\geq 1}$ are all non-empty. Then, there exists a point $\displaystyle (\theta, x) \in \cap_{n\geq 0} V_n$, where $\theta \in \cap_{n\geq 0} I_n$. For an arbitrary $n\geq 1$, by the disjointness conditions $(\mathcal{F}1-2)_n$ it follows that $(\theta - M_n\omega) \in \Theta_n$ and the iterates in the set $\{\theta+k\omega\}_{|k| = 1}^{M_n}$ do not intersect $I_n$. Since $(\theta, x) \in V_n$, on the fibre we have $\pi_2 \circ F^{-M_n}(\theta, x) \in (C_1 \cup C_2)$. Then, by $(C2)_n$ (in proposition \ref{proposition_block n contraction}), we obtain the following bound for the expansion under backward iteration:
\begin{align*}
\left|\partial_{x}\left(\pi_2\circ F^{-k}(\theta,x)\right)\right| = \left|\prod_{i=1}^{k}(f'(x_{-i}))^{-1}\right| \geq \alpha_c^{-k/4}, \text{ for every } 0< k \leq M_n.
\end{align*}
Since we have the above estimate for every $n\geq 1$, this proves that 
\begin{align*}
    \underline{L}^{-}(\theta, x) > -{(\log \alpha_c)}/{4} > 0.
\end{align*}
 Thus the orbit of $(\theta, x)$ is a sink. Similarly, for the fixed $(\theta, x) \in V_n$, we have $\theta_{M_n} \in (I_n+M_n\omega)$ and $x_{M_n} \in (E_1\cup E_2)$, and hence by $(E2)_n$ (in proposition \ref{proposition_block n expansion}) we obtain the following bound for the contraction under forward iteration: 
\begin{align*}
\left|\partial_{x}\left(\pi_2\circ F^{k}(\theta,x)\right)\right| = \left|\prod_{i=0}^{k-1}\left((f^{-1})'(f(x_{i}))\right)^{-1}\right| \geq \alpha_c^{-k/4}, \text{ for every } 1<k\leq M_n.
\end{align*}
Since we have the bound for every $n\geq 1$, this proves that $\displaystyle \underline{L}^{+}(\theta, x) > -{(\log \alpha_c)}/{4} > 0$. Thus, the orbit of $(\theta,x)$ is a source. 

This proves that, for every $\omega \in \Omega_1$, the map $F_{\omega}$ has a sink-source orbit $(F_{\omega}^n(\theta,x))_{n \in \mathbb{Z}}$. For such a point $(\theta,x)$, since $\underline{L}^{+}(\theta, x) > 0$ and  $\displaystyle \underline{L}^{-}(\theta, x) > 0$, it must also be a point of discontinuity of the graphs $\Phi$ and $\Psi$. By the minimality of the irrational translation, the existence of the sink-source orbit  proves that the invariant graphs are discontinuous on a dense set on the base $\mathbb{T}$. Hence, the synchronisation and convergence of orbits along the fibre is a strictly non-uniform behaviour.

%%%%%%%%%%%%%%%%%%%%%%%%%%%%%% (b) CURVES CONSTRUCTION %%%%%%%%%%%%%%%%%%%%%%
\vspace{0.2cm}
\noindent\textbf{Proof of the Theorem \ref{Theorem SNA} (b):} Let us fix $\omega \in \Omega_2$, which is irrational by definition \ref{definition_parameters_final_omega2}. Then, there exist constants $M_0, M_1(\omega) \in [\tilde{M}_1, 2\tilde{M}_1]$ for $\tilde{M}_1 := \alpha_c^{-{M}_{0}/8}$, such that $\omega \in \tilde{\mathcal{F}}_1(M_0,M_1)$. By proposition \ref{proposition_structure I1_critical set I_1}, there exists a non-empty critical set $I_1(\omega)$ as in definition \ref{definition_structure In_critical set In}, for which
\begin{align}
    \text{comp}(I_{1}) = 4, \ \ I_{1}(\omega) \subset I_{0}(\omega) \ \ \mbox{   and   } \ \ 0 < |I_1(\omega)| <  \alpha_c^{M_{0}/2}.
\end{align}
Recall the sets $G_2 \subset G_1 \subset I_0$ as in definition \ref{definition_structure I1_sync set Gs}. By proposition \ref{proposition_structure I1_sync set Gs}, the sets satisfy the condition
\begin{align}\label{equation_final proof_invariant curve_G2 sync prop 1}
    \pi_2 \circ F^{2M_0+S_0}((G_1-M_0\omega)\times (C_1 \cup C_2)) & \subset C_1 \\
    \pi_2 \circ F^{-2M_0-S_0}((G_1+M_0\omega)\times (E_1 \cup E_2)) & \subset E_1,
    \stepcounter{equation}\tag{\theequation} \label{equation_final proof_invariant curve_G2 sync prop 2}
\end{align}
and $|G_2| > \lambda^{-3/4}$. Note that, by $(\mathcal{F}1)_0$ (in definition \ref{definition_block 0_disjointness condn}), we have $(G_2-3M_0\omega) \subset \Theta_0$. By the estimates $|G_2| > \lambda^{-3/4}$,  $|I_1(\omega)| < \alpha_c^{M_0/2}$, and 
\begin{align*}
    \left|\bigcup_{|m| \leq 2M_1} 3(I_1 + m\omega) \right| < 60\tilde{M}_1\alpha_c^{M_0/2} < \alpha_c^{M_0/4}
\end{align*}
for a sufficiently large $\lambda$, this gives us
\begin{align*}
   |G_2| - \left|\bigcup_{|m| \leq 2M_1} 3(I_1 + m\omega) \right| > \lambda^{-3/4}/2.
\end{align*}
By the above estimates, there exists a non-empty closed interval $\Delta \subset \mathbb{T}$ such that $3\Delta \subset (G_2-3M_0\omega)$,
\begin{align}\label{equation_final proof_invariant curve_Delta prop}
    \Delta \cap \left(\bigcup_{|m| \leq 2M_1} 3(I_1 + m\omega) \right) = \emptyset,
\end{align}
and $|\Delta| = \delta < \min\{|I_1^{(i,j)}|:i,j=1,2\}$. By $(\mathcal{F}1)_0$, for such a fixed interval $\Delta \subset (G_2-3M_0\omega)$, we then also have 
\begin{align}\label{equation_final proof_invariant curve_Delta prop 2}
    \Delta \cap \left(\bigcup_{|m| \leq 2M_0} 3(I_0 + m\omega) \right) = \emptyset.
\end{align}
Then, it follows from the definition of $\tilde{\Theta}_1$ that $\Delta \subset \tilde{\Theta}_1$. Let us now define the following new measurable sets:
\begin{align}
 \bar{\Theta}_{1} & := \mathbb{T} \setminus \bigcup_{l=0}^{1} \bigcup_{|m| \leq 2M_l} (I_l + m\omega)  \label{equation_final proof_invariant curve_bar Theta 1 set} \\
 \bar{\Theta} & := \bigcap_{m= -3S_0}^{3S_0} (\bar{\Theta}_{1} + m\omega) \\
        D & := \bigcup_{l=0}^{1} \bigcup_{m =1 }^{M_l +S_0} (I_l + m\omega) \\
        \tilde{D} & := \mathbb{T} \setminus 3D,\\
        D' & := \bigcup_{l=0}^{1} \bigcup_{m =0 }^{M_l +S_0-1} (I_l - m\omega) 
        \stepcounter{equation}\tag{\theequation} \label{equation_final proof_invariant curve_D' set}
\end{align}
Note that $\bar{\Theta}_1 \subset \tilde{\Theta}_1$. For any fixed $\lambda$, each of the sets $D$ and $\tilde{D}$ are a finite union of open intervals, and the fixed $\Delta$ satisfies $\Delta \cap D = \emptyset$. Using the estimates of the critical sets in $\{I_n\}_{n=0,1}$, we obtain the following estimates:
\begin{align*}
|\bar{\Theta}_{1}| & > 1 -  \lambda^{-1/8} -  \alpha_c^{\tilde{M}_0 / 8} > 0 \\
    |\bar{\Theta}| & > 1 - (7S_0 + 1) \left( \lambda^{-1/8} + \alpha_c^{\tilde{M}_0/8} \right) > 0 \\
   0 < |D| & < \lambda^{-1/8} + \alpha_c^{\tilde{M}_{0}/8} \\
   0 < |D'| & < \lambda^{-1/8} + \alpha_c^{\tilde{M}_{0}/8}  \\
    |\tilde{D}| & > 1 - 3\lambda^{-1/8} - 3\alpha_c^{\tilde{M}_{0}/8}.
\end{align*}
We observe that the measure of sets $\bar{\Theta}_1$, $\bar{\Theta}$ and $\tilde{D}$ all converge to 1 as $\lambda$ goes to infinity. Whereas the measure of $D$ and $D'$ converge to $1$ as $\lambda$ goes to infinity.

By the Birkhoff ergodic theorem, there exists a measurable set $\tilde{\Delta} \subset \Delta$ such that $|\tilde{\Delta}| = |\Delta|$ and for each $\theta \in \tilde{\Delta}$ there exists an integer $N>0$ large enough such that 
\begin{align}\label{equation_final proof_invariant curve_rare visits to tilde D}
   \# \left\{ 0 \leq k \leq n: ({\theta} -k\omega) \in \tilde{D} \right\} > n\left(1 - \frac{1}{3b}\right),
\end{align}
for all $n\geq N$. By the definition of $\tilde{D}$ and since $|\Delta| = \delta < \min\{|I_1^{(i,j)}|:i,j=1,2\}$, for any $\theta \in \tilde{\Delta}$ and $k\geq 1$,
\begin{align}\label{equation_final proof_invariant curve_rare visits to D}
    ({\theta} - k\omega) \in \tilde{D} \implies \left(\Delta-k\omega\right) \cap D = \emptyset,
\end{align}
which gives us
\begin{align*}
    \# \left\{ 0 \leq k \leq n: (\Delta-k\omega) \cap D = \emptyset \right\} > n\left(1 - \frac{1}{3b}\right), \text{ for all } n\geq N.
\end{align*}

\vspace{0.2cm}
\noindent\textit{Construction of invariant curves $\Phi(\theta)$ and $\Psi(\theta)$:} 
\vspace{0.1cm}

Let us fix a point $\bar{\theta} \in \tilde{\Delta}$. Then, corresponding to the fixed $\bar{\theta}$, there exists an increasing sequence $(N_k)_{k\geq 0}$, such that $N_0 := 0$, $N_1 \geq N$, $N_k - N_{k-1} > 2\tilde{M}_1$ for all $k\geq 1$, and $(\bar{\theta}-N_k\omega) \in \tilde{\Delta}$.
For each $k\geq 1$, since $N_k - N_{k-1} > 2\tilde{M}_1$ and $(\bar{\theta}-N_k\omega) \in (G_2-3M_0\omega)$, it follows from equation \eqref{equation_final proof_invariant curve_G2 sync prop 1} that
\begin{align*}
    F^{4M_0+S_0}\left( \{\bar{\theta}-N_k\omega\} \times (C_1 \cup C_2) \right) \subset \left( \{\bar{\theta}-(N_k-4M_0-S_0)\omega\} \times C_1 \right).
\end{align*}
For each $k\geq 1$, since $(\bar{\theta}-(N_k-2M_0-S_0)\omega) \in \Theta_1$, by $(Sc1)_1$ (in proposition \ref{proposition_block n sync contraction}) and $(C3)_1$ (in proposition \ref{proposition_block n contraction}) we obtain
\begin{align}\label{equation_final proof_invariant curve_nestedness on point thetabar}
    F^{N_k-N_{k-1}}\left( \{\bar{\theta}-N_k\omega\} \times (C_1 \cup C_2) \right) \subset \left( \{\bar{\theta}-N_{k-1}\omega\} \times C_{i_{k-1}} \right),
\end{align}
for some $i_{k-1} \in \{1,2\}$. Similarly, there exists $i_0\in \{1,2\}$ such that 
\begin{align*}
     F^{N_k}\left( \{\bar{\theta}-N_k\omega\} \times (C_1 \cup C_2) \right) \subset F^{N_{k-1}}\left( \{\bar{\theta}-N_{k-1}\omega\} \times C_{i_{k-1}} \right) \subset \left( \{\bar{\theta}\} \times C_{i_0}\right).
\end{align*}
For each $k\geq 0$, by the conditions imposed on $\Delta$ in the equations \eqref{equation_final proof_invariant curve_Delta prop} and \eqref{equation_final proof_invariant curve_Delta prop 2}, 
\begin{align*}
    (\bar{\theta}-N_{k}\omega) \in \mathbb{T} \setminus \left(\bigcup_{l=0}^{1}\bigcup_{|m|\leq 2M_l}3(I_l+m\omega) \right) \implies (\Delta - N_{k}\omega) \subset \tilde{\Theta}_1.
\end{align*}
Then, for every $k\geq 1$, by equation \eqref{equation_final proof_invariant curve_nestedness on point thetabar},  the continuity of $F$, $(C1)_0$ (in proposition \ref{proposition_block 0 contraction}) and $(C3)_1$ (in proposition \ref{proposition_block n contraction}), 
\begin{align}\label{equation_final proof_curves_nested blocks}
    F^{N_k - N_{k-1}}\left( (\Delta - N_{k}\omega) \times (C_1 \cup C_2) \right) \subset \{ (\Delta - N_{k-1}\omega) \} \times C_{i_{k-1}}, \text{ for some } i_{k-1} \in \{1,2\}.
\end{align}
By $(C1)_0$ (in proposition \ref{proposition_block 0 contraction}) and $(C1)_1$ (in proposition \ref{proposition_block n contraction}), we know that for any $\theta \in {\Theta}_1$ and $x\in C_1 \cup C_2$, 
\begin{align*}
    \theta_m \notin D & \implies x_m \in C_1 \cup C_2  \implies |f'(x_m)| < \alpha_c. 
\end{align*}
Since $\bar{\theta} \in \tilde{\Delta}$ and the sequence of integers $(N_k)_{k\geq 0}$ satisfy the equation \eqref{equation_final proof_invariant curve_rare visits to tilde D}, it follows from equation \eqref{equation_final proof_invariant curve_rare visits to D} that for any $\theta \in \Delta$ and $x\in C_1 \cup C_2$ we have the following bound for the rate of contraction along the fibre:
\begin{align}\label{equation_final proof_curves_shrinking rate}
    |\partial_x \left( \pi_2 \circ F^{N_k}(\theta- N_k\omega, x) \right)| \leq \left(\alpha_c^{1-1/3b} \alpha_c^{b/3b} \right)^{N_k} < (\alpha_c^{1/3})^{N_k}.
\end{align}
For each $k\geq 1$ and $\theta \in \Delta$, by equation \eqref{equation_final proof_curves_nested blocks} and \eqref{equation_final proof_curves_shrinking rate}, 
\begin{align}\label{equation_final proof_curves_uniform shrinking rate}
    \left|\pi_2 \circ F^{N_k} \left( \{\theta - N_k\omega\} \times (C_1 \cup C_2) \right) \right| & < (\alpha_c^{1/3})^{N_{k-1}}.
\end{align}
Thus, for each $k\geq 1$ and $\theta \in \Delta$, the measure of the nested sets $\left(F^{N_k} \left( (\theta - N_{k}\omega) \times (C_1 \cup C_2) \right)\right)_{k\geq 1}$ satisfy the condition
\begin{align*}
    \lim_{k\rightarrow \infty}|\pi_2 \circ F^{N_k} \left( \{\theta - N_k\omega\} \times (C_1 \cup C_2) \right)| & = 0.
\end{align*}
This proves that the nested sets $\left(F^{N_k} \left( (\theta - N_{k}\omega) \times (C_1 \cup C_2) \right)\right)_{k\geq 1}$ converge to a single point which we denote by $\Phi(\theta) \in (C_1\cup C_2)$. For any $\theta \in \Delta$ and $m\geq 1$ such that $(\theta-m\omega) \in \Delta$ and $m > N_k$ for some $k\geq 1$, by $(C1)_0$ (in proposition \ref{proposition_block 0 contraction}) and $(C3)_1$ (in proposition \ref{proposition_block n contraction}), we have 
\begin{align*}
    F^m(\{\theta-m\omega\}\times(C_1\cup C_2)) \subset F^{N_k}(\{\theta-N_k\omega\}\times(C_1\cup C_2)). 
\end{align*}
This proves that, for each $\theta \in \Delta$ the point $\Phi(\theta)$ is independent of the specific choice of sequence $(N_k)_{k\geq 0}$.

The exponential rate of convergence of the points in equation \eqref{equation_final proof_curves_uniform shrinking rate} proves that the function $\Phi$ is in fact continuous on $\Delta$. Note that the function $\Phi$ satisfies the condition 
\begin{align*}
    F^n \left( \theta, \Phi(\theta) \right) = (\theta + n\omega, \Phi(\theta + n\omega)),
\end{align*}
whenever $\theta, (\theta + n\omega) \in \Delta$. Since $\Delta$ is a non-empty closed interval, the function $\Phi$ naturally extends to a continuous invariant graph on $\mathbb{T}$, which satisfies
\begin{align}
    F \left( \theta, \Phi(\theta) \right) = (\theta + \omega, \Phi(\theta + \omega)), \ \text{ for all } \theta \in \mathbb{T}.
\end{align}

Recall the set $\bar{\Theta}_1 \subset \tilde{\Theta}_1$ which was defined in equation \eqref{equation_final proof_invariant curve_bar Theta 1 set}. For every $\theta \in \bar{\Theta}_1$ and $x\in C_1 \cup C_2$, if $(m_k)_{k\geq 1}$ is the sequence of all integers for which $(\theta+m_k\omega) \in \Delta$ then, by $(C3)_1$ (in proposition \ref{proposition_block n contraction}), we have $ \pi_2 \circ F^{m_k}(\theta, x) \in (C_1\cup C_2)$, and thus
\begin{align*}
    \lim_{k \rightarrow \infty} |\pi_2 \circ F^{m_k} \left(\theta, x \right) - \Phi(\theta + m_k\omega)| = 0.
\end{align*}
Thus, for every $\theta \in \bar{\Theta}_1$ and $x\in C_1 \cup C_2$,
\begin{align}\label{equation_final proof_curves_partial attraction C}
    \lim_{n \rightarrow \infty} |\pi_2 \circ F^{n} \left(\theta, x \right) - \Phi(\theta + n\omega)| = 0.
\end{align}
Note that this also shows that $\Phi(\theta) \in (C_1\cup C_2)$, for every $\theta \in \bar{\Theta}_1$. For each $\theta \in \bar{\Theta}$ and $x\in E_1^c \cap E_2^c$, since $(\theta + m\omega)_{|m|\leq S_0} \cap I_0 = \emptyset$, it follows that
\begin{align*}
    F^{S_0}(\theta, x) & \in \left(\bar{\Theta}_1 \times (C_1 \cup C_2)\right).
\end{align*}
Thus, by equation \eqref{equation_final proof_curves_partial attraction C},
\begin{align}\label{equation_final proof_curves_partial attraction E^c}
     \lim_{n \rightarrow \infty} |\pi_2 \circ F^{n} \left(\theta, x \right) - \Phi(\theta + n\omega)| & = 0,
\end{align}
for all $\theta \in \bar{\Theta}$ and $x\in E_1^c \cap E_2^c$.

Similarly, we can construct a continuous invariant graph $\Psi: \mathbb{T} \rightarrow \mathbb{T}$ that satisfies 
\begin{align}\label{equation_final proof_curves_partial repulsion E}
    \lim_{n \rightarrow \infty} |\pi_2 \circ F^{-n} \left(\theta, x \right) - \Psi(\theta - n\omega)|,
\end{align}
for all $\theta \in \bar{\Theta}_1$ and $x\in C_1 \cup C_2$,  
\begin{align}\label{equation_final proof_curves_partial attraction C^c}
     \lim_{n \rightarrow \infty} |\pi_2 \circ F^{-n} \left(\theta, x \right) - \Psi(\theta - n\omega)| & = 0,
\end{align}
for all $\theta \in \bar{\Theta}$ and $x\in C_1^c \cap C_2^c$, and  $\Psi(\theta) \in (E_1\cup E_2)$ for every  $\theta \in \tilde{\Theta}_1$.

\vspace{0.2cm}
\noindent\textit{Uniqueness of the attracting(repelling) invariant curve:}
\vspace{0.1cm}

By the definition of the set $D'$ in equation \eqref{equation_final proof_invariant curve_D' set}, for a sufficiently large $\lambda$ with respect to $b$, we have $|D'|<1/(6b)$. By the ergodicity of the base dynamics, we can fix a point $\theta \in \bar{\Theta}$ and $N > 0$ such that  
\begin{align}\label{equation_final proof_bad region visits}
    \# \{ 0 \leq j \leq n: (\theta + j\omega) \in D'\} < \frac{n}{3b}, \text{ for all } n\geq N.
\end{align}
Since $\theta \in \bar{\Theta}$, by $(E3)_1$ (in proposition \ref{proposition_block n expansion}), we have $\Psi(\theta) \in E_1\cup E_2$. Let us fix $y \neq \Psi(\theta)$ and a constant $\varepsilon > 0$ such that $|\Psi(\theta) - y| > \varepsilon$. Let $1 <k < \tilde{k}$ be integers that satisfy conditions $(\theta + \tilde{k}\omega) \in \bar{\Theta}$, $(\theta+k\omega) \in \Delta$, $2N + 2\tilde{M}_1 < N+ 2\tilde{M}_1 +k < \tilde{k}$ and $\alpha_c^{k/3} < \varepsilon$. 

For any $x \in E_1 \cup E_2$, by $(E3)_1$ (in proposition \ref{proposition_block n expansion}), it follows that 
\begin{align*}
   (\theta + (\tilde{k}-j)\omega) \notin D'\implies \pi_2 \circ F^{-j} (\theta + \tilde{k}\omega, x) \in E_1 \cup E_2 , 
\end{align*}
for all $j = 0,\dots ,\tilde{k}$. Since $\tilde{k}>N$, by the bound in equation \eqref{equation_final proof_bad region visits}, for every $x\in E_1\cup E_2$, we obtain
\begin{align*}
    \left|\partial_x \left(\pi_2 \circ F^{-\tilde{k}} (\theta + \tilde{k}\omega, x)\right)\right| = \left|(f^{-1})'(f(x_{-\tilde{k}})) \dots (f^{-1})'(f(x_{-1}))\right| & < \alpha_c^{\tilde{k} - \frac{\tilde{k}}{3b}} \alpha_c^{\frac{-\tilde{k}b}{3b}} < (\alpha_c^{1/3})^{\tilde{k}}.  
\end{align*}
and similarly $\left|\partial_x \left(\pi_2 \circ F^{-{k}} (\theta + {k}\omega, x)\right)\right| <  (\alpha_c^{1/3})^{{k}}$. Since $(\theta+\tilde{k}\omega), (\theta+{k}\omega) \in \tilde{\Theta}_1$, by the construction of the invariant graph $\Psi$ and equation \eqref{equation_final proof_curves_partial repulsion E}, this gives us $\Psi(\theta+\tilde{k}\omega) \in E_1 \cup E_2$ and $\Psi(\theta+{k}\omega) \in E_i$, for some $i\in \{1,2\}$. Since $\Delta \subset (G_2-3M_0\omega) \cap \bar{\Theta}_1$ and $(\tilde{k}-k)>2\tilde{M}_1$, by $(E3)_1$ (in proposition \ref{proposition_block n expansion}) and equation \eqref{equation_final proof_invariant curve_G2 sync prop 2} it follows that
\begin{align*}
    F^{-\tilde{k}+k}(\{\theta+\tilde{k}\omega\}\times (E_1\cup E_2)) \subset (\{\theta+k\omega\}\times E_i),
\end{align*}
where $i \in \{1,2\}$ for which $\Psi(\theta+{k}\omega) \in E_i$. Then, for every $x\in E_1\cup E_2$, we have
\begin{align*}
    |\Psi(\theta) - \pi_2 \circ F^{-\tilde{k}} (\theta + \tilde{k}\omega, x)| < (\alpha_c^{1/3})^{k} < \varepsilon.
\end{align*}
By the choice of $\varepsilon$, since $|\Psi(\theta) - y| > \varepsilon$, we must have 
\begin{align*}
    F^{\tilde{k}}(\theta, y) & \in \bar{\Theta} \times (E_1^c\cap E_2^c).
\end{align*} 
Then, by equation \eqref{equation_final proof_curves_partial attraction E^c}, we obtain that
\begin{align}\label{equation_final proof_synceqn1a}
    \lim_{n\rightarrow \infty} |\pi_2 \circ F^n(\theta, y) - \Phi(\theta + n\omega)| = 0.
\end{align}
This proves the equation \eqref{statementofthem_thm_synceqn2a} in Theorem \ref{Theorem SNA}. The proof of equation \eqref{statementofthem_thm_synceqn2b} in Theorem \ref{Theorem SNA} follows similarly.

%--------------------------------------------- ACKNOWLEDGEMENTS ------------------------------------

\subsection*{Acknowledgements}
The present paper is a part of my thesis \cite{Rajasekar:2026}.
I want to express my sincere gratitude to my PhD advisor Kristian Bjerkl\"{o}v for his continued support and guidance. I also thank \`{A}lex Haro and Jordi-Lluís Figueras for the insightful conversations and your valuable perspectives.

% -------------------------------------------- Disclosure ----------------------------------

%%%%%%%%%%%%%%%%%%%%%%%%%%%%%%%%%%%%%%%%%%%%%%%%%%%%%%%%%%%
    
%%%%%%%%%%%%%%%%%%%%% THE END %%%%%%%%%%%%%%%%%%%%%%%%%%%%%

\bibliographystyle{abbrv}
\bibliography{bibliography_SNA} 

\end{document}